\documentclass[
  a4paper,
  twoside,
  reqno
]{amsart}
\usepackage[english]{babel}
\usepackage{amssymb}
\usepackage{amsmath}
\usepackage{amsthm}
\usepackage{subfiles}
\usepackage{mathtools}
\usepackage{nicefrac}
\input xy
\xyoption{all}
\usepackage{enumitem}
\usepackage{multicol}
\usepackage{tikz}

\newcommand{\R}{\mathbb{R}}
\newcommand{\Prm}{\mathbb{P}}
\newcommand{\Z}{\mathbb{Z}}
\newcommand{\Q}{\mathbb{Q}}
\newcommand{\E}{{[-\infty,+\infty]}}
\newcommand{\EX}{\smash{\E^X}}
\newcommand{\N}{\mathbb{N}}
\newcommand{\Lex}{\mathbb{L}}

\newcommand{\ra}{\rightarrow}
\newcommand{\ol}{\overline}
\newcommand{\ul}{\underline}
\newcommand{\bv}{\textstyle{\bigvee}}
\newcommand{\bw}{\textstyle{\bigwedge}}
\newcommand{\eqdf}{:=}

\newcommand{\qvL}[1]{#1/{\approx}}
\newcommand{\qvphi}[1]{#1/{\approx}}

\newcommand{\vs}[4]{#1 \supseteq #2\,\smash{\xrightarrow{#3}\,#4}}
\newcommand{\Lebesguesign}{\mathcal{L}}
\newcommand{\LA}{\mathcal{A}_{\Lebesguesign}}
\newcommand{\Lmu}{\mu_{\Lebesguesign}}
\newcommand{\LF}{F_{\Lebesguesign}}
\newcommand{\Lphi}{\varphi_{\Lebesguesign}}
\newcommand{\vsLA}{\vs{\wp{\R}}{\LA}{\Lmu}{\R}}
\newcommand{\vsLF}{\vs{\E^\R}{\LF}{\Lphi}{\R}}

\newcommand{\Borelsign}{\mathcal{B}}
\newcommand{\BA}{\mathcal{A}_{\Borelsign}}
\newcommand{\Bmu}{\mu_{\Borelsign}}
\newcommand{\vsBA}{\vs{\wp{\R}}{\BA}{\Bmu}{\R}}

\newcommand{\Simplesign}{\mathrm{S}}
\newcommand{\SA}{\mathcal{A}_{\Simplesign}}
\newcommand{\Smu}{\mu_{\Simplesign}}
\newcommand{\SF}{F_{\Simplesign}}
\newcommand{\Sphi}{\varphi_{\Simplesign}}
\newcommand{\vsSA}{\vs{\wp{\R}}{\SA}{\Smu}{\R}}
\newcommand{\vsSF}{\vs{\E^\R}{\SF}{\Sphi}{\R}}

\newcommand{\ulim}[1]{\smash{\overline{\lim}_{#1}}}
\newcommand{\llim}[1]{\smash{\underline{\lim}_{#1}}}
\newcommand{\pulim}[2]{{#1}\text{-}\!\ulim{#2}}
\newcommand{\pllim}[2]{{#1}\text{-}\!\llim{#2}}
\newcommand{\plim}[2]{{#1}\text{-}\!{\lim}_{#2}}

\newcommand{\ve}{\varepsilon}
\newcommand{\se}{\,\ve\,}
\newcommand{\dtn}[2]{\smash{\nicefrac{\textstyle{#1}}{#2}}}
\newcommand{\dt}[1]{\dtn{#1}{2}}

\newcommand{\ld}[1]{d_{\smash{#1}}}

\newcommand\lbb{\mathchoice{[\kern-1.42pt[}%
                    {[\kern-1.42pt[}%
                    {[\kern-1.2pt[}%
                    {[\kern-1.15pt[}}
\newcommand\rbb{\mathchoice{]\kern-1.42pt]}%
                    {]\kern-1.42pt]}%
                    {]\kern-1.2pt]}%
                    {]\kern-1.15pt]}}
\newcommand{\decode}[1]{\lbb#1\rbb}
\newcommand{\Stump}{\mathbf{Stp}}

\newcommand{\dv}{\preccurlyeq}
\newcommand{\lcm}{\operatorname{lcm}}
\newcommand{\keyword}[1]{\textbf{#1}}
\newcommand{\Zmod}[1]{\Z_{#1}}

\newcommand{\rsub}[2]{\phantom{#2}\mathllap{#1}}

\setlist[1]{label=(\roman*)}

\theoremstyle{definition}
\newtheorem{numbering}{Use to get sequential numbering}
\newtheorem{dfn}[numbering]{Definition}
\newtheorem{ex}[numbering]{Example}
\newtheorem{exs}[numbering]{Examples}

\theoremstyle{remark}
\newtheorem{rem}[numbering]{Remark}

\theoremstyle{notation}
\newtheorem{nt}[numbering]{Notation}
\newtheorem{cnv}[numbering]{Convention}
\newtheorem{sit}[numbering]{Situation}

\theoremstyle{theorem}
\newtheorem{lem}[numbering]{Lemma}
\newtheorem{cor}[numbering]{Corollary}
\newtheorem{thm}[numbering]{Theorem}
\newtheorem{prop}[numbering]{Proposition}

\newcommand{\shortversion}{ 2013-06-03 (bab2968)}

\usepackage[final]{hyperref}

\begin{document}
{ 
\title[A Generalisation of Measure and Integral]{Lattice Valuations,\\
A Generalisation of Measure and Integral}
\email{bram@westerbaan.name}

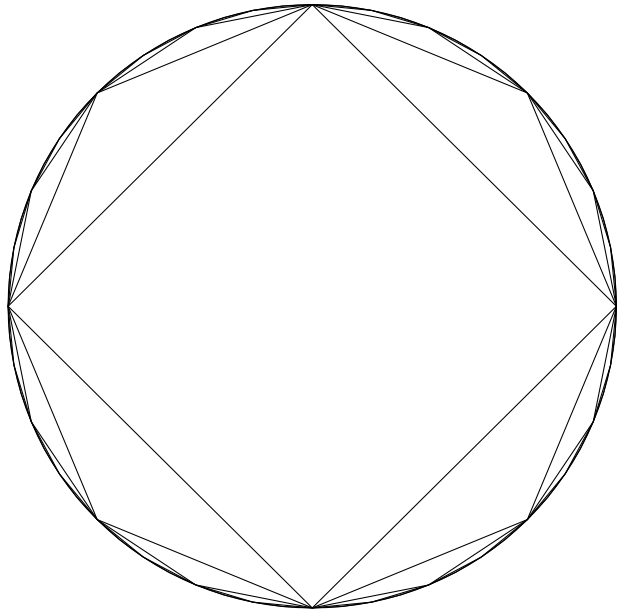
\author[A.A.~Westerbaan]{Bram Westerbaan}
%\date{\today \quad{\tiny \version}}
\maketitle

\vspace{1cm}
\begin{center}
\begin{tikzpicture}
\draw[thin,black!100] (3.99518182482068962003cm,0.19627069730967205974cm)
	 -- (3.98073890668878771493cm,0.39206856131824241452cm)
	 -- (3.95670603985912405776cm,0.58692189782144699173cm)
	 -- (3.92314112161292172232cm,0.78036128806451299234cm)
	 -- (3.88012501277817589695cm,0.97192071961305548378cm)
	 -- (3.82776134292883529753cm,1.16113870901784932421cm)
	 -- (3.76617626073208322524cm,1.34755941356888020444cm)
	 -- (3.69551813004514739802cm,1.53073372946035912712cm)
	 -- (3.61595717249377379687cm,1.71022037372112856168cm)
	 -- (3.52768505739342064231cm,1.88558694730399079020cm)
	 -- (3.43091444000108891643cm,2.05641097677288708923cm)
	 -- (3.32587844921018138677cm,2.22228093207840871059cm)
	 -- (3.21283012592258021556cm,2.38279721796973342762cm)
	 -- (3.09204181345094841760cm,2.53757313665458195118cm)
	 -- (2.96380450141983686763cm,2.68623581938807332037cm)
	 -- (2.82842712474619073504cm,2.82842712474618984686cm)
	 -- (2.68623581938807420855cm,2.96380450141983642354cm)
	 -- (2.53757313665458283936cm,3.09204181345094797351cm)
	 -- (2.38279721796973431580cm,3.21283012592257977147cm)
	 -- (2.22228093207840959877cm,3.32587844921018094269cm)
	 -- (2.05641097677288753331cm,3.43091444000108847234cm)
	 -- (1.88558694730399101225cm,3.52768505739342019822cm)
	 -- (1.71022037372112878373cm,3.61595717249377379687cm)
	 -- (1.53073372946035957121cm,3.69551813004514784211cm)
	 -- (1.34755941356888064853cm,3.76617626073208411341cm)
	 -- (1.16113870901784976830cm,3.82776134292883662980cm)
	 -- (0.97192071961305592787cm,3.88012501277817722922cm)
	 -- (0.78036128806451332540cm,3.92314112161292305458cm)
	 -- (0.58692189782144721377cm,3.95670603985912539002cm)
	 -- (0.39206856131824252554cm,3.98073890668878904719cm)
	 -- (0.19627069730967208749cm,3.99518182482069139638cm)
	 -- (-0.00000000000000005551cm,4.00000000000000177636cm)
	 -- (-0.19627069730967219852cm,3.99518182482069139638cm)
	 -- (-0.39206856131824263656cm,3.98073890668878949128cm)
	 -- (-0.58692189782144732479cm,3.95670603985912583411cm)
	 -- (-0.78036128806451343642cm,3.92314112161292349867cm)
	 -- (-0.97192071961305603889cm,3.88012501277817767331cm)
	 -- (-1.16113870901784999035cm,3.82776134292883707388cm)
	 -- (-1.34755941356888087057cm,3.76617626073208500159cm)
	 -- (-1.53073372946035979325cm,3.69551813004514917438cm)
	 -- (-1.71022037372112922782cm,3.61595717249377557323cm)
	 -- (-1.88558694730399167838cm,3.52768505739342241867cm)
	 -- (-2.05641097677288797740cm,3.43091444000109069279cm)
	 -- (-2.22228093207841004286cm,3.32587844921018316313cm)
	 -- (-2.38279721796973475989cm,3.21283012592258154783cm)
	 -- (-2.53757313665458328344cm,3.09204181345094974986cm)
	 -- (-2.68623581938807509673cm,2.96380450141983819989cm)
	 -- (-2.82842712474619162322cm,2.82842712474619206731cm)
	 -- (-2.96380450141983819989cm,2.68623581938807554081cm)
	 -- (-3.09204181345094974986cm,2.53757313665458417162cm)
	 -- (-3.21283012592258154783cm,2.38279721796973520398cm)
	 -- (-3.32587844921018271904cm,2.22228093207841048695cm)
	 -- (-3.43091444000109024870cm,2.05641097677288842149cm)
	 -- (-3.52768505739342197458cm,1.88558694730399190043cm)
	 -- (-3.61595717249377557323cm,1.71022037372112967191cm)
	 -- (-3.69551813004514961847cm,1.53073372946036023734cm)
	 -- (-3.76617626073208588977cm,1.34755941356888131466cm)
	 -- (-3.82776134292883840615cm,1.16113870901785043444cm)
	 -- (-3.88012501277817900558cm,0.97192071961305648298cm)
	 -- (-3.92314112161292483094cm,0.78036128806451376949cm)
	 -- (-3.95670603985912716638cm,0.58692189782144754684cm)
	 -- (-3.98073890668879082355cm,0.39206856131824274758cm)
	 -- (-3.99518182482069317274cm,0.19627069730967222627cm)
	 -- (-4.00000000000000355271cm,0.00000000000000000000cm)
	 -- (-3.99518182482069317274cm,-0.19627069730967222627cm)
	 -- (-3.98073890668879126764cm,-0.39206856131824274758cm)
	 -- (-3.95670603985912761047cm,-0.58692189782144754684cm)
	 -- (-3.92314112161292527503cm,-0.78036128806451376949cm)
	 -- (-3.88012501277817944967cm,-0.97192071961305648298cm)
	 -- (-3.82776134292883885024cm,-1.16113870901785043444cm)
	 -- (-3.76617626073208677795cm,-1.34755941356888153670cm)
	 -- (-3.69551813004515095074cm,-1.53073372946036045938cm)
	 -- (-3.61595717249377734959cm,-1.71022037372113011600cm)
	 -- (-3.52768505739342419503cm,-1.88558694730399256656cm)
	 -- (-3.43091444000109246915cm,-2.05641097677288930967cm)
	 -- (-3.32587844921018493949cm,-2.22228093207841137513cm)
	 -- (-3.21283012592258332418cm,-2.38279721796973609216cm)
	 -- (-3.09204181345095152622cm,-2.53757313665458461571cm)
	 -- (-2.96380450141983997625cm,-2.68623581938807642899cm)
	 -- (-2.82842712474619384366cm,-2.82842712474619339957cm)
	 -- (-2.68623581938807731717cm,-2.96380450141983997625cm)
	 -- (-2.53757313665458550389cm,-3.09204181345095197031cm)
	 -- (-2.38279721796973653625cm,-3.21283012592258376827cm)
	 -- (-2.22228093207841181922cm,-3.32587844921018493949cm)
	 -- (-2.05641097677288975376cm,-3.43091444000109246915cm)
	 -- (-1.88558694730399301065cm,-3.52768505739342419503cm)
	 -- (-1.71022037372113078213cm,-3.61595717249377779368cm)
	 -- (-1.53073372946036134756cm,-3.69551813004515183891cm)
	 -- (-1.34755941356888220284cm,-3.76617626073208811022cm)
	 -- (-1.16113870901785110057cm,-3.82776134292884062660cm)
	 -- (-0.97192071961305703809cm,-3.88012501277818122603cm)
	 -- (-0.78036128806451421358cm,-3.92314112161292705139cm)
	 -- (-0.58692189782144787991cm,-3.95670603985912938683cm)
	 -- (-0.39206856131824296963cm,-3.98073890668879304400cm)
	 -- (-0.19627069730967233729cm,-3.99518182482069539319cm)
	 -- (0.00000000000000000000cm,-4.00000000000000621725cm)
	 -- (0.19627069730967236505cm,-3.99518182482069583727cm)
	 -- (0.39206856131824302514cm,-3.98073890668879393218cm)
	 -- (0.58692189782144787991cm,-3.95670603985913027500cm)
	 -- (0.78036128806451421358cm,-3.92314112161292793957cm)
	 -- (0.97192071961305703809cm,-3.88012501277818211420cm)
	 -- (1.16113870901785132261cm,-3.82776134292884151478cm)
	 -- (1.34755941356888242488cm,-3.76617626073208944248cm)
	 -- (1.53073372946036156961cm,-3.69551813004515361527cm)
	 -- (1.71022037372113122622cm,-3.61595717249378001412cm)
	 -- (1.88558694730399389883cm,-3.52768505739342685956cm)
	 -- (2.05641097677289064194cm,-3.43091444000109513368cm)
	 -- (2.22228093207841270740cm,-3.32587844921018760402cm)
	 -- (2.38279721796973742443cm,-3.21283012592258598872cm)
	 -- (2.53757313665458639207cm,-3.09204181345095374667cm)
	 -- (2.68623581938807820535cm,-2.96380450141984219670cm)
	 -- (2.82842712474619517593cm,-2.82842712474619562002cm)
	 -- (2.96380450141984175261cm,-2.68623581938807909353cm)
	 -- (3.09204181345095374667cm,-2.53757313665458728025cm)
	 -- (3.21283012592258554463cm,-2.38279721796973831260cm)
	 -- (3.32587844921018715993cm,-2.22228093207841315149cm)
	 -- (3.43091444000109468959cm,-2.05641097677289108603cm)
	 -- (3.52768505739342641547cm,-1.88558694730399434292cm)
	 -- (3.61595717249378001412cm,-1.71022037372113189235cm)
	 -- (3.69551813004515405936cm,-1.53073372946036223574cm)
	 -- (3.76617626073209033066cm,-1.34755941356888309102cm)
	 -- (3.82776134292884284704cm,-1.16113870901785198875cm)
	 -- (3.88012501277818344647cm,-0.97192071961305781524cm)
	 -- (3.92314112161292927183cm,-0.78036128806451487971cm)
	 -- (3.95670603985913160727cm,-0.58692189782144854604cm)
	 -- (3.98073890668879526444cm,-0.39206856131824352474cm)
	 -- (3.99518182482069761363cm,-0.19627069730967278138cm)
	 -- (4.00000000000000799361cm,-0.00000000000000033307cm)
	 -- (3.99518182482068962003cm,0.19627069730967205974cm);

\draw[thin,black!100] (3.98073890668878771493cm,0.39206856131824241452cm)
	 -- (3.92314112161292216641cm,0.78036128806451310336cm)
	 -- (3.82776134292883574162cm,1.16113870901784954626cm)
	 -- (3.69551813004514739802cm,1.53073372946035934916cm)
	 -- (3.52768505739342064231cm,1.88558694730399101225cm)
	 -- (3.32587844921018138677cm,2.22228093207840959877cm)
	 -- (3.09204181345094841760cm,2.53757313665458283936cm)
	 -- (2.82842712474619073504cm,2.82842712474619117913cm)
	 -- (2.53757313665458283936cm,3.09204181345094930577cm)
	 -- (2.22228093207841004286cm,3.32587844921018227495cm)
	 -- (1.88558694730399190043cm,3.52768505739342153049cm)
	 -- (1.53073372946036023734cm,3.69551813004514873029cm)
	 -- (1.16113870901785043444cm,3.82776134292883751797cm)
	 -- (0.78036128806451388051cm,3.92314112161292394276cm)
	 -- (0.39206856131824302514cm,3.98073890668878993537cm)
	 -- (0.00000000000000038858cm,4.00000000000000266454cm)
	 -- (-0.39206856131824230349cm,3.98073890668879037946cm)
	 -- (-0.78036128806451321438cm,3.92314112161292483094cm)
	 -- (-1.16113870901784999035cm,3.82776134292883840615cm)
	 -- (-1.53073372946036001530cm,3.69551813004515006256cm)
	 -- (-1.88558694730399190043cm,3.52768505739342330685cm)
	 -- (-2.22228093207841048695cm,3.32587844921018405131cm)
	 -- (-2.53757313665458417162cm,3.09204181345095108213cm)
	 -- (-2.82842712474619295548cm,2.82842712474619295548cm)
	 -- (-3.09204181345095108213cm,2.53757313665458461571cm)
	 -- (-3.32587844921018449540cm,2.22228093207841137513cm)
	 -- (-3.52768505739342419503cm,1.88558694730399301065cm)
	 -- (-3.69551813004515139482cm,1.53073372946036112552cm)
	 -- (-3.82776134292884018251cm,1.16113870901785110057cm)
	 -- (-3.92314112161292660730cm,0.78036128806451432460cm)
	 -- (-3.98073890668879259991cm,0.39206856131824319167cm)
	 -- (-4.00000000000000532907cm,0.00000000000000027756cm)
	 -- (-3.98073890668879304400cm,-0.39206856131824263656cm)
	 -- (-3.92314112161292749548cm,-0.78036128806451388051cm)
	 -- (-3.82776134292884107069cm,-1.16113870901785087852cm)
	 -- (-3.69551813004515272709cm,-1.53073372946036112552cm)
	 -- (-3.52768505739342597138cm,-1.88558694730399323269cm)
	 -- (-3.32587844921018671585cm,-2.22228093207841226331cm)
	 -- (-3.09204181345095330258cm,-2.53757313665458639207cm)
	 -- (-2.82842712474619517593cm,-2.82842712474619517593cm)
	 -- (-2.53757313665458639207cm,-3.09204181345095374667cm)
	 -- (-2.22228093207841315149cm,-3.32587844921018715993cm)
	 -- (-1.88558694730399456496cm,-3.52768505739342685956cm)
	 -- (-1.53073372946036245779cm,-3.69551813004515405936cm)
	 -- (-1.16113870901785221079cm,-3.82776134292884329113cm)
	 -- (-0.78036128806451521278cm,-3.92314112161293016001cm)
	 -- (-0.39206856131824374678cm,-3.98073890668879615262cm)
	 -- (-0.00000000000000049960cm,-4.00000000000000888178cm)
	 -- (0.39206856131824280309cm,-3.98073890668879659671cm)
	 -- (0.78036128806451432460cm,-3.92314112161293104819cm)
	 -- (1.16113870901785176670cm,-3.82776134292884462340cm)
	 -- (1.53073372946036245779cm,-3.69551813004515583572cm)
	 -- (1.88558694730399500905cm,-3.52768505739342863592cm)
	 -- (2.22228093207841403967cm,-3.32587844921018893629cm)
	 -- (2.53757313665458816843cm,-3.09204181345095552302cm)
	 -- (2.82842712474619739638cm,-2.82842712474619695229cm)
	 -- (3.09204181345095596711cm,-2.53757313665458816843cm)
	 -- (3.32587844921018982447cm,-2.22228093207841448375cm)
	 -- (3.52768505739342952410cm,-1.88558694730399567518cm)
	 -- (3.69551813004515716798cm,-1.53073372946036334596cm)
	 -- (3.82776134292884639976cm,-1.16113870901785287693cm)
	 -- (3.92314112161293326864cm,-0.78036128806451554585cm)
	 -- (3.98073890668879926125cm,-0.39206856131824374678cm)
	 -- (4.00000000000001243450cm,-0.00000000000000016653cm)
	 -- (3.98073890668878771493cm,0.39206856131824241452cm);

\draw[thin,black!100] (3.92314112161292172232cm,0.78036128806451299234cm)
	 -- (3.69551813004514695393cm,1.53073372946035890507cm)
	 -- (3.32587844921018094269cm,2.22228093207840871059cm)
	 -- (2.82842712474618984686cm,2.82842712474618984686cm)
	 -- (2.22228093207840871059cm,3.32587844921018049860cm)
	 -- (1.53073372946035890507cm,3.69551813004514650984cm)
	 -- (0.78036128806451299234cm,3.92314112161292083414cm)
	 -- (0.00000000000000022204cm,3.99999999999999911182cm)
	 -- (-0.78036128806451254825cm,3.92314112161292083414cm)
	 -- (-1.53073372946035823894cm,3.69551813004514606575cm)
	 -- (-2.22228093207840782242cm,3.32587844921018005451cm)
	 -- (-2.82842712474618895868cm,2.82842712474618940277cm)
	 -- (-3.32587844921017961042cm,2.22228093207840826651cm)
	 -- (-3.69551813004514517758cm,1.53073372946035846098cm)
	 -- (-3.92314112161291950187cm,0.78036128806451277029cm)
	 -- (-3.99999999999999777955cm,0.00000000000000022204cm)
	 -- (-3.92314112161291950187cm,-0.78036128806451232620cm)
	 -- (-3.69551813004514473349cm,-1.53073372946035779485cm)
	 -- (-3.32587844921017872224cm,-2.22228093207840693424cm)
	 -- (-2.82842712474618807050cm,-2.82842712474618807050cm)
	 -- (-2.22228093207840737833cm,-3.32587844921017827815cm)
	 -- (-1.53073372946035846098cm,-3.69551813004514384531cm)
	 -- (-0.78036128806451310336cm,-3.92314112161291816960cm)
	 -- (-0.00000000000000077716cm,-3.99999999999999644729cm)
	 -- (0.78036128806451154905cm,-3.92314112161291816960cm)
	 -- (1.53073372946035668463cm,-3.69551813004514384531cm)
	 -- (2.22228093207840560197cm,-3.32587844921017783406cm)
	 -- (2.82842712474618629415cm,-2.82842712474618762641cm)
	 -- (3.32587844921017650179cm,-2.22228093207840693424cm)
	 -- (3.69551813004514206895cm,-1.53073372946035823894cm)
	 -- (3.92314112161291639325cm,-0.78036128806451321438cm)
	 -- (3.99999999999999467093cm,-0.00000000000000122125cm)
	 -- (3.92314112161292172232cm,0.78036128806451299234cm);

\draw[thin,black!100] (3.69551813004514695393cm,1.53073372946035912712cm)
	 -- (2.82842712474618984686cm,2.82842712474619029095cm)
	 -- (1.53073372946035868303cm,3.69551813004514695393cm)
	 -- (-0.00000000000000044409cm,4.00000000000000000000cm)
	 -- (-1.53073372946035957121cm,3.69551813004514695393cm)
	 -- (-2.82842712474619073504cm,2.82842712474618984686cm)
	 -- (-3.69551813004514739802cm,1.53073372946035846098cm)
	 -- (-4.00000000000000000000cm,-0.00000000000000066613cm)
	 -- (-3.69551813004514650984cm,-1.53073372946035979325cm)
	 -- (-2.82842712474618940277cm,-2.82842712474619073504cm)
	 -- (-1.53073372946035801689cm,-3.69551813004514739802cm)
	 -- (0.00000000000000111022cm,-4.00000000000000000000cm)
	 -- (1.53073372946036023734cm,-3.69551813004514650984cm)
	 -- (2.82842712474619073504cm,-2.82842712474618895868cm)
	 -- (3.69551813004514695393cm,-1.53073372946035757280cm)
	 -- (3.99999999999999911182cm,0.00000000000000155431cm)
	 -- (3.69551813004514695393cm,1.53073372946035912712cm);

\draw[thin,black!100] (2.82842712474619029095cm,2.82842712474618984686cm)
	 -- (0.00000000000000088818cm,4.00000000000000000000cm)
	 -- (-2.82842712474618940277cm,2.82842712474619073504cm)
	 -- (-4.00000000000000000000cm,0.00000000000000111022cm)
	 -- (-2.82842712474619117913cm,-2.82842712474618895868cm)
	 -- (-0.00000000000000177636cm,-4.00000000000000000000cm)
	 -- (2.82842712474618851459cm,-2.82842712474619162322cm)
	 -- (4.00000000000000000000cm,-0.00000000000000266454cm)
	 -- (2.82842712474619029095cm,2.82842712474618984686cm);

\draw[thin,black!100] (0.00000000000000024493cm,4.00000000000000000000cm)
	 -- (-4.00000000000000000000cm,0.00000000000000048986cm)
	 -- (-0.00000000000000073479cm,-4.00000000000000000000cm)
	 -- (4.00000000000000000000cm,-0.00000000000000097972cm)
	 -- (0.00000000000000024493cm,4.00000000000000000000cm);

%\draw (0cm,0cm) circle(1cm);

\end{tikzpicture}
\end{center}

\vfill
\noindent
{\tiny
Thesis for the Master's Examination
Mathematics at the Radboud University Nijmegen,\\
supervised by prof.~dr.~A.C.M.~van Rooij
with second reader dr.~O.W.~van Gaans,\\
written by Abraham A.~Westerbaan, 
student number 0613622,
on November 16th 2012.\\
\\
This is a revised version written on \shortversion.\\
For the latest version,
see \url{http://bram.westerbaan.name/master.pdf}.

% The whitespace above is needed.
% If it is removed, the text above will be tiny,
% but the spacing between lines will be normal.
}
\thispagestyle{empty}
\clearpage
\thispagestyle{empty}
$\,$
\clearpage
\begin{abstract}
Measure and integral are two closely related,
but distinct objects of study.
Nonetheless,
they are both real-valued \emph{lattice valuations}:
order preserving real-valued functions~$\varphi$
on a lattice~$L$
which are \emph{modular}, i.e.,
\begin{equation*}
\varphi(x) + \varphi(y) 
\,=\, 
\varphi(x\wedge y) + \varphi(x\vee y)\qquad(x,y\in L).
\end{equation*}
We unify measure and integral
by developing a theory for lattice valuations.
We allow these lattice valuations
to take their values from the reals,
or any suitable ordered Abelian group.
\end{abstract}
\clearpage
\thispagestyle{empty}
$\,$
\newpage
\thispagestyle{empty}
$\,$
\vfill
\begin{center}
{\em dedicated to my father}\\
Henk Westerbaan\\
 $\dagger$ june, 2012
\end{center}
\vfill
\vfill
\clearpage
\thispagestyle{empty}
%
%                  FOREWORD
%
\section*{Preface}
\noindent
In the summer of 2009  Bas Westerbaan
and I worked out an 
overly general approach
to the introduction of the Lebesgue measure and the Lebesgue integral
with the help of dr.~A.C.M.~van~Rooij.
The theory that is presented in this thesis
is based on the work done in that summer.

Since I was fortunate enough to be offered
a Ph.D.-position,
this thesis was written under time constraints.
Hence the text is not nearly as polished as 
I would like it to be,
and the proofs of some statements
have been left to the reader.
I hope the reader will be able to ignore the rough edges
and enjoy this fresh view on the old subject
of measure and integration.

I would like to thank  all my teachers
for showing me the beauty of mathematics.
In particular, I  thank
dr.~Mai Gehrke for showing me its elegance,
dr.~Wim Veldman for showing me its content, and
dr.~Henk Barendregt for showing me how it is written.
Furthermore, I am most grateful
to dr.~A.C.M.~van~Rooij
for his never relenting willingness 
to answer my questions
and note my errors.
\clearpage
\thispagestyle{empty}
\tableofcontents
\clearpage
 }

% Sections
{ \thispagestyle{empty}
\begin{flushright}
\begin{minipage}{.7\columnwidth}
\begin{flushright}
The theory of integration,
because of its 
central r\^ole 
in mathematical analysis and geometry,
continues to afford opportunities 
for serious investigation.\\
--- \textsc{M.H. Stone}, 1948, \cite{Stone48}
\end{flushright}
\end{minipage}
\end{flushright}
\clearpage
\section{Introduction}%
\label{S:intro}%
\noindent
There are many ways (some more popular than
others) to introduce the Lebesgue measure
and the Lebesgue integral.
For the purposes of this introduction,
we define the Lebesgue measure and integral
in such a way that the similarity between them is obvious.
This similarity is the basis of this thesis.
We leave it to the reader to 
compare the  definitions below
to those that are familiar to him/her.
\begin{dfn}
\label{D:lebesgue-measure}
The \keyword{Lebesgue measure}
$\Lmu \colon \LA \ra \R$
is the smallest\footnote{%
``Smallest'' with respect to the following order.
We say that $\mu_1$ \emph{is extended by}~$\mu_2$
where $\mu_i\colon \mathcal{A}_i \ra \R$
and $\mathcal{A}_i \subseteq \wp(\R)$
provided that
$\mathcal{A}_1 \subseteq \mathcal{A}_2$,
and $\mu_1(A) = \mu_2(A)$
for all~$A\in\mathcal{A}_1$.}
$\mu\colon \mathcal{A} \ra \R$
where $\mathcal{A}$ is a subset of~$\wp(\R)$
that has the following properties.
\begin{enumerate}
\item
\label{prop:measure-1}
Let $a,b\in\R$ with $a\leq b$. Then $[a,b]\in\mathcal{A}$
and $(a,b)\in\mathcal{A}$, 
and
\begin{equation*}
\mu(\,[a,b]\,)
\ = \ \mu(\,(a,b)\,)
\ =\ b-a.
\end{equation*}

\item
\label{prop:measure-2}
\textit{(Monotonicity)}\ 
Let $A,B\in \mathcal{A}$.
Then $\mu(A)\leq \mu(B)$
when  $A\subseteq B$.

\item
\label{prop:measure-3}
\textit{(Modularity)}\ 
Let $A,B\in \mathcal{A}$.
Then $A\cap B\in\mathcal{A}$ and $A\cup B \in \mathcal{A}$, and
\begin{equation*}
\mu(\,A\cap B\,)\,+\,\mu(\,A\cup B\,)\ =\ \mu(A)\,+\,\mu(B).
\end{equation*}

\item
\label{prop:measure-4}
\textit{($\Pi$-Completeness)}\ 
Let $A_1 \supseteq A_2 \supseteq \dotsb$
from~$\mathcal{A}$ be given.\\
Assume that the set $\{\,\mu(A_1),\, \mu(A_2),\, \dotsc\,\}$
has an infimum, $\bw_n \mu(A_n)$.\\
Then we have $\bigcap_n A_n \in \mathcal{A}$.
Moreover,
\begin{equation*}
\mu(\, \textstyle{ \bigcap_n A_n }\,) \ =\ \bw_n \,\mu(A_n).
\end{equation*}

\item
\label{prop:measure-5}
\textit{($\Sigma$-Completeness)}\ 
Let $A_1 \subseteq A_2 \subseteq \dotsb$
from~$\mathcal{A}$ 
be such that $\bv_n \mu(A_n)$ exists.\\
Then we have $\bigcup_n A_n \in \mathcal{A}$.
Moreover,
\begin{equation*}
\mu(\, \textstyle{ \bigcup_n A_n }\,) \ =\ \bv_n \,\mu(A_n).
\end{equation*}

\item
\label{prop:measure-6}
\textit{(Convexity)}
Let $A\subseteq Z \subseteq B$ be subsets of~$\R$.\\
Assume that $A,B\in\mathcal{A}$ and $\mu(A)=\mu(B)$.\\
Then we have $Z\in \mathcal{A}$ and $\mu(A) = \mu(Z)= \mu(B)$.
\end{enumerate}
\end{dfn}

\begin{dfn}
\label{D:lebesgue-integral}
The \keyword{Lebesgue integral}
$\Lphi \colon \LF \ra \R$
is the smallest  $\varphi\colon F \ra \R$
where $F$ is a subset of~$\E^\R$
that has the following properties.
\begin{enumerate}
\item
\label{prop:integral-1}
Let $a,b,\lambda\in\R$ with $a\leq b$. 
Then $\lambda\cdot\mathbf{1}_{[a,b]}\in F$
and $\lambda\cdot\mathbf{1}_{(a,b)}\in F$,
and
\begin{equation*}
\varphi(\,\lambda\cdot \mathbf{1}_{[a,b]}\,)
\ =\ \varphi(\,\lambda\cdot \mathbf{1}_{(a,b)}\,)
\ =\ \lambda\cdot(b-a).
\end{equation*}

\item
\label{prop:integral-2}
\textit{(Monotonicity)}\ 
Let $f,g\in F$.
Then $\varphi(f)\leq \varphi(g)$
when  $f\leq g$.

\item
\label{prop:integral-3}
\textit{(Modularity)}\ 
Let $f,g\in F$.
Then $f\wedge g\in F$ and $f \vee g \in F$, and
\begin{equation*}
\varphi(\,f\wedge g\,)\,+\,\varphi(\,f\vee g\,)\ =\ \varphi(f)\,+\,\varphi(g).
\end{equation*}

\item
\label{prop:integral-4}
\textit{($\Pi$-Completeness)}\ 
Let $f_1 \geq f_2 \geq \dotsb$
from~$F$ be such that $\bw_n \varphi(f_n)$ exists.\\
Then we have $\bw_n f_n \in F$.
Moreover,
\begin{equation*}
\varphi(\, \textstyle{ \bigwedge_n f_n }\,) \ =\ \bw_n \,\varphi(f_n).
\end{equation*}
Here $\bw_n f_n$
is the infimum of $\{\,f_1,\,f_2,\,\dotsc\,\}$ in $\E^\R$;
more concretely,
it is the \emph{pointwise infimum}, i.e., $(\bw_n f_n)(x) = \bw_n f_n(x)$
for all~$x\in \R$.

\item
\label{prop:integral-5}
\textit{($\Sigma$-Completeness)}\ 
Let $f_1 \leq f_2 \leq \dotsb$
from~$F$ 
be such that $\bv_n \varphi(f_n)$ exists.\\
Then we have $\bigvee_n f_n \in F$.
Moreover,
\begin{equation*}
\varphi(\, \textstyle{ \bigvee_n f_n }\,) \ =\ \bv_n \,\varphi(f_n).
\end{equation*}

\item
\label{prop:integral-6}
\textit{(Convexity)}
Let $f\leq z \leq g$ be $\E$-valued functions on~$\R$.\\
Assume that $f,g\in F$ and $\varphi(f)=\varphi(g)$.\\
Then we have $z\in F$ and $\varphi(f) = \varphi(z)= \varphi(g)$.
\end{enumerate}
\end{dfn}

\noindent
In this thesis
we present an abstract theory based on the properties
(Monotonicity), (Modularity), ($\Pi$-Completeness),
($\Sigma$-Completeness) and (Convexity)
and we try to fit some of the results
of measure and integration theory
in this framework.

\subsection{Valuations}
We begin by considering (Monotonicity) and (Modularity).

Maps with these two properties 
are called \emph{(lattice) valuations}.
More precisely,
let~$L$ be a lattice,
and let~$E$ be an
ordered Abelian group (e.g.~$\R$, see Appendix~\ref{S:ag}).
A map $\varphi\colon L\ra E$ is a \emph{valuation}
if it is order preserving and \emph{modular}, i.e.,
\begin{equation*}
\varphi(x)\,+\,\varphi(y) \ =\ 
\varphi(x\wedge y)\,+\, \varphi(x \vee y)
\qquad(x,y\in L).
\end{equation*}

Of course,
the Lebesgue measure~$\Lmu$
and the Lebesgue integral~$\Lphi$
are valuations,
and there are many more examples.
We study valuation in Section~\ref{S:vals}.

\subsection{Valuation Systems}
Let us now look at ($\Pi$-Completeness).
For the Lebesgue measure it
involves intersections, ``$\bigcap_n A_n$'',
i.e., infima in~$\wp(\R)$.
Similarly,
for the Lebesgue integral
it
involves pointwise infima, ``$\bw_n f_n$'',
i.e., infima in~$\smash{\E^\R}$.
In order to
generalise 
the notion of 
($\Pi$-Completeness) 
to any valuation $\varphi\colon L\ra E$
we involve a `surrounding' lattice,~$V$.
That is, we will define what
it means for an object of the following shape
to be \emph{$\Pi$-complete}
(see Definition~\ref{D:system-complete}).
\begin{equation*}
\vs{V}{L}\varphi{E}
\end{equation*}
We call these objects \emph{valuation systems},
and we study them in Section~\ref{S:system}.

The Lebesgue measure and the Lebesgue integral give us valuation systems:
\begin{equation*}
\vsLA\qquad\text{and}\qquad\vsLF.
\end{equation*}
Of course
these valuation systems are $\Pi$-complete
by ($\Pi$-Completeness).\\
They are also \emph{$\Sigma$-complete},
which is a generalisation of  ($\Sigma$-Completeness).

Finally,
(Convexity)
can easily be generalised to valuation systems as well.
We will define what it means
for a valuation system to be \emph{convex}
in  Definition~\ref{D:convex}.
We study these convex valuation systems in Subsection~\ref{SS:convex}.

\vspace{1em}
\noindent
Now that we have introduced
the main objects of study,
valuations and valuation systems,
let us spend some words on the theorems that we will prove.

\subsection{Completion and Convexification}
Recall that we defined
the Lebesgue measure~$\Lmu$
as the smallest map $\mu\colon \mathcal{A}\ra \R$
that has properties~\ref{prop:measure-1}--\ref{prop:measure-6}.
It is important to note that
it is not obvious at all that such a map exists.
While it is relatively easy to see
that if there is a map $\mu\colon \mathcal{A} \ra \R$
that has properties~\ref{prop:measure-1}--\ref{prop:measure-6},
then there is also a smallest one,
it takes quite some effort
to prove that there is any map $\mu\colon \mathcal{A}\ra \R$
with properties~\ref{prop:measure-1}--\ref{prop:measure-6} to begin with.

One could call this statement the \emph{Extension Theorem
for the Lebesgue measure}.
Similarly,
to define~$\Lphi$,
we need an \emph{Extension Theorem
for the Lebesgue integral}.

We will generalise (a part) of these two theorems
to the setting of valuations.
To see how we could do this,
note that to prove
the Extension Theorem for the Lebesgue measure,
one could take the following three steps.
\begin{enumerate}
\item
\label{extension-step-1}
\emph{Find the smallest map $\Smu\colon \SA\ra \R$
that has properties~\ref{prop:measure-1}--\ref{prop:measure-3}.}

This is not too difficult.
Let $\mathcal{S}$ be the family of subsets of~$\R$ of
the form~$[a,b]$ or~$(a,b)$ where $a,b\in \R$ with $a\leq b$.
Let~$\SA$ be the set of all unions
of finite disjoint subsets of~$\mathcal{S}$,
and let $\Smu\colon \SA\ra \R$
be given by
\begin{equation*}
\Smu(\,I_1 \cup\dotsb\cup I_N\,) \ = \ |I_1| \,+\,\dotsb\,+\,|I_N|,
\end{equation*}
where $I_1,\dotsc,I_N\,\in\, \mathcal{S}$
with $I_n\cap I_m = \varnothing$ when $n\neq m$.

Of course, it requires some calculations
to see that such a map~$\Smu$ exists,
and that~$\Smu$ will have the 
properties~\ref{prop:measure-1}--\ref{prop:measure-3}
(see Example~\ref{E:smeas-val}).

\item
\label{extension-step-2}
\emph{Extend~$\Smu$
to the smallest map $\overline{\Smu}\colon\overline{\SA}\ra\R$
that has properties~\ref{prop:measure-1}--\ref{prop:measure-5}.}

This is the most interesting and the most difficult step.
To give an idea of how one could try obtain such~$\overline{\Smu}$,
consider the following `algorithm'.
\begin{equation*}
\left[\quad
\begin{minipage}{.75\columnwidth}
Let $\mu\colon\mathcal{A} \ra \R$
be a variable. To begin, set $\mu \eqdf \Smu$.\\
$\mathbf{(*)}$ For all~$A_1,A_2,\dotsc$ from~$\mathcal{A}$ 
do the following.
\begin{enumerate}
\item[$\bullet$] If $A_1 \supseteq A_2 \supseteq \dotsb$
and  $\bw_n \mu(A_n)$ exists
and  $\bigcap_n A_n \notin \mathcal{A}$,\\
then add $\bigcap_n A_n$ to~$\mathcal{A}$,
and set $\mu(\bigcap_n A_n) \eqdf \bw_n \mu(A_n)$.
\item[$\bullet$] If $A_1 \subseteq A_2 \subseteq \dotsb$
and  $\bv_n \mu(A_n)$ exists
and  $\bigcup_n A_n \notin \mathcal{A}$,\\
then add $\bigcup_n A_n$ to~$\mathcal{A}$,
and set $\mu(\bigcup_n A_n) \eqdf \bv_n \mu(A_n)$.
\end{enumerate}
If~$\mu$ was changed, repeat~$\mathbf{(*)}$.
\end{minipage}
\right.
\end{equation*}
There are many problems with this `algorithm'.
Perhaps the most serious problem
is, loosely speaking, that the same set~$A$ may be obtained
in several ways
and it is not clear that
$\mu(A)$ would be given the same value each time.

Note that the `algorithm' resembles the  definition
of the Borel sets.
In fact, $\overline{\Smu}$
will be the family of all Borel subsets of~$\R$
with finite measure.

\item
\label{extension-step-3}
\emph{Extend~$\overline{\Smu}$
to the smallest map~$\Lmu\colon \LA \ra \R$
that has properties~\ref{prop:measure-1}--\ref{prop:measure-6}.}

This is straightforward.
Simply
define $\LA$ to be the family of all subsets of~$\R$
that are `sandwiched' between elements of~$\overline{\SA}$,
that is,
all $Z\in\wp(\R)$
for which
there are $A,B\in \overline{\SA}$
such that $A\subseteq Z\subseteq B$ and 
$\overline{\Smu}(A) = \overline{\Smu}(B)$.

Now, define $\Lmu\colon \LA \ra \R$
by $\Lmu (Z) = \overline{\Smu} (A)$ for $Z$ and~$A$ as above.
\end{enumerate}

\noindent
We have sketched how 
to get
the Lebesgue measure~$\Lmu\colon \LA \ra \R$
in three steps, 
\begin{equation*}
\xymatrix{
\ar @{-->}[rr]^{\text{\ref{extension-step-1}}}
&&
\Smu
\ar @{-->}[rr]^{\text{\ref{extension-step-2}}}
&&
\ol\Smu
\ar @{-->}[rr]^{\text{\ref{extension-step-3}}}
&&
\Lmu.
}
\end{equation*}
We will generalise step~\ref{extension-step-2}
and step~\ref{extension-step-3}
to the setting of valuations.
More precisely:
\begin{enumerate}
\item
Let $\vs{V}{L}\varphi{E}$
be a valuation system.
We will give
a necessary and sufficient condition,
namely~$\vs{V}{L}\varphi{E}$
is \emph{extendible}
(see Definition~\ref{D:extendible}),
for 
the existence of a smallest
valuation~$\overline\varphi\colon\overline{L}\ra E$
which extends~$\varphi$
where~$\overline{L}$ is a sublattice of~$V$
such that the valuation system
$\vs{V}{\overline{L}}{\overline\varphi}{E}$
is both $\Pi$-complete and $\Sigma$-complete
(see Lemma~\ref{L:complete} and Proposition~\ref{P:comp-minimal}).\\
We will call~$\overline\varphi$
the \emph{completion}
of~$\varphi$ (relative to~$V$).
\item
Let $\vs{V}{L}\varphi{E}$
be a valuation sytem.
We will prove the following.\\
There is smallest valuation~$\varphi^\bullet\colon L^\bullet\ra E$
extending~$\varphi$
with~$L^\bullet$ a sublattice of~$V$
such that 
$\vs{V}{L^\bullet}{\varphi^\bullet}{E}$ is convex
(see Propisition~\ref{P:convexification}).\\
Moreover,
$\vs{V}{L^\bullet}{\varphi^\bullet}{E}$
is $\Pi$-complete and $\Sigma$-complete
provided that
$\vs{V}{L}\varphi{E}$ 
is $\Pi$-complete and $\Sigma$-complete
(see Proposition~\ref{P:convexification_versus_completion}).\\
We will call $\varphi^\bullet$
the \emph{convexification} of~$\varphi$
(relative to~$V$).
\end{enumerate}
By the discussion above
we see that 
the Lebesgue measure~$\Lmu$
is the convexification of the completion of~$\Smu$
relative to~$\wp(\R)$:
\begin{equation*}
\xymatrix{
\Smu
\ar @{-->}[rrr]^{\text{completion}}
&&&
\ol\Smu
\ar @{-->}[rrr]^{\text{convexification}}
&&&
{\Lmu}
}
\end{equation*}
Similarly,
the Lebesgue integral~$\Lphi$
is the convexification 
of the completion of~$\Sphi$
relative to~$\smash{\E^\R}$,
where~$\Sphi\colon \SF\ra \R$ is the obvious
valuation on the set of \emph{step functions}~$\SF$
(see Example~\ref{E:sint-val}).
So we get the following diagram.
\begin{equation*}
\xymatrix{
\Sphi
\ar @{-->}[rrr]^{\text{completion}}
&&&
\ol\Sphi
\ar @{-->}[rrr]^{\text{convexification}}
&&&
{\Lphi}
}
\end{equation*}
Let us note that $\ol\Sphi\colon \ol{\SF}\ra \R$
will be the restriction of the Lebesgue integral
to the set~$\ol{\SF}$
of Lebesgue integrable \emph{Baire functions}.
We will not prove this.

We belief that
the completion is the most important step,
and that the convexification
is mere decoration.
In line with this believe,
we spend most words
on the completion,
and we leave it to the reader
to think about the convexification.

\subsection{Closedness under Operations}
We have found an abstract method
to get the Lebesgue measure~$\Lmu$ and
the Lebesgue integral~$\Lphi$.
However, such a method is nothing but a curiosity
if we cannot use it to derive some
basic properties of~$\Lmu$ and~$\Lphi$.
One such property might be:
\begin{equation*}
\left[\quad
\begin{minipage}{.7\columnwidth}
If $f,g\in \R^\R$
are Lebesgue integrable,\\
then $f+g$ is Lebesgue integrable,\\
and $\Lphi(f+g) = \Lphi(f)+\Lphi(g)$.
\end{minipage}
\right.
\end{equation*}
So, roughly speaking,
$\Lphi$ is \emph{closed under the operation~``$+$''}.
Instead of this,
we will prove that~$\overline{\Sphi}$ is closed under the operation~``$+$''.
We leave it to the reader to use this to prove that the convexification
 of~$\overline{\Sphi}$, i.e.~$\Lphi$, is closed under~``$+$'' as well.

More generally,
in Section~\ref{S:closedness}
we will prove some  statements of the following shape.
If $\vs{V}{L}\varphi{E}$
is a valuation system,
and~$\varphi$ is closed under some operation in some sense,
then the completion~$\ol{\varphi}$ is closed under the same operation as well.

\subsection{Convergence Theorems}
An important part of the theory of integration
is that of the convergence theorems.
So we have studied whether
these make sense in the setting of valuations.
We will show in Subsection~\ref{SS:complete-val_convergence}
that it is possible to formulate
and prove
the Lemma of Fatou
and Lebesgue's Dominated Convergence Theorem
for complete valuation systems.
Interestingly,
the surrounding lattice~$V$ will play no role.
This leads to the study 
of \emph{complete valuations}
(as opposed to complete valuation systems),
see Section~\ref{S:complete-val} for more details.

\subsection{Fubini's Theorem}
Another important part of the theory of integration
is Fubini's Theorem.
Unfortunately,
it seems that that it not possible
to make sense of Fubini's Theorem in
the general setting of valuations.

Nevertheless,
in Section~\ref{S:fub} 
we will split the proof of Fubini's Theorem for the Lebesgue integral
into two parts.
The first part concerns step functions and is specific
to the Lebesgue integral,
while the second part is a consequence of a
general extension theorem for valuations
(see Theorem~\ref{T:fubext}).

\subsection{Extendibility}
We have remarked that a valuation system
$\vs{V}{L}\varphi{E}$
has a completion if and only if~$\varphi$ is \emph{extendible}.
As the reader will see in Subsection~\ref{SS:hierarchy-abstract}
the definition of ``$\varphi$ is extendible''
is rather involved.

Fortunately, 
the situation is simpler
for some choices of~$E$.
We say that~$E$ is \emph{benign}
if 
for every valuation system $\vs{V}{L}\varphi{E}$
we have  that $\varphi$ is extendible iff
\begin{equation*}
\left[\quad
\begin{minipage}{.7\columnwidth}
Let $a_1 \geq a_2 \geq \dotsb$ in~$L$
with $\bw_n\varphi(a_n)$ exists be given.\\
Let $b_1 \leq b_2 \leq \dotsb$ in~$L$ 
with $\bv_n \varphi(b_n)$ exists be given. \\
Then we have the following implication.
\begin{equation*}
\bw_n a_n  \,\leq\,\bv_n b_n
\quad\implies\quad
\varphi(\bw_n a_n)
\,\leq\,
\varphi(\bv_n b_n),
\end{equation*}
Here, $\bw_n a_n$
is the infimum of $a_1 \geq a_2 \geq \dotsb$ in~$V$,\\
and 
$\bv_n b_n$ is the supremum of~$b_1 \leq b_2 \leq \dotsb$ in~$V$.
\end{minipage}
\right.
\end{equation*}
We will prove that~$\R$ is benign.
More generally,
we will prove in Section~\ref{S:unif}
any ordered Abelian group~$E$
that has a suitable unformity
(see Def.~\ref{D:uniformity})
is benign.

\subsection{Attribution}
Some work of others has been included in this master's thesis.
\begin{enumerate}
\item
An early version of the theory in Section~\ref{S:unif}
has been developed together with Bas Westerbaan,
and some of his work is undoubtedly still there.
\item
The proof of the Borel Hierarchy Theorem
in Subsection~\ref{SS:borel-hierarchy}
is an adaptation of the work by Wim Veldman~\cite{Veldman08}.
\item
Countless improvements 
were suggested by dr.~A.C.M. van Rooij.
Most notably,
he strengthened Lemma~\ref{lem:main}
to its current form,
and he suggested that I should restrict the theory  to
lattice valuations.
\end{enumerate}
Aside from the things mentioned above,
and that which is common knowledge,
and unless stated otherwise,
every definition and proof in this thesis is my own.

Nevertheless,
I do not want to
claim 
that any part of my work is original as well,
because that would be mere gambling.
Indeed,
recently I discovered
an article~\cite{Alfsen63}
on the foundation of integration
in which 
 valuations
are used as well.

\subsection{Prerequisites}
We have tried to keep this text as accessible as possible.\\
We assume that
the reader is familiar
with the ordinal numbers
and is comfortable with the basic notions
of order theory (suprema, infima, lattices, etc., see~\cite{DP02}).\\
We have attached some material on \emph{ordered Abelian groups},
in Appendix~\ref{S:ag}.\\
While some knowledge
about measure, integral,
topology,
uniform spaces, 
Borel sets,
and Riesz spaces will helpful as well,
we hope this will not be necessary.

\subsection{Notation}
Let us take this opportunity to fix some notation.
\begin{enumerate}
\item
We write $\N=\{1,2,\dotsc\}$ and $\omega=\{0,1,2,\dotsc\}$.
\item
We will use the symbol ``$\bv$'' for suprema,
and the symbol ``$\bw$'' for infima,
as the symbol ``$\sum$'' is used for sums,
and the symbol ``$\prod$'' for products.
\item
Given  $x\in \R$ we say that $x$ is \emph{positive}
when $x \geq 0$.\\
Given  $x\in \R$ we say that $x$ is \emph{strictly positive}
when $x > 0$.\\
\end{enumerate}
 }
\clearpage
{ \section{Valuations}
\label{S:vals}
\noindent
Both the Lebesgue measure and Lebesgue integral
are \emph{valuations}.
In fact they are both \emph{complete valuations}
(see Section~\ref{S:complete-val}).
While the better part of this thesis involves
complete valuations,
much can already be said about valuations.

In this section we study the elementary properties of valuations.
We start with some examples in Subsection~\ref{SS:vals_intro}.
We study the distance 
induced by a valuation~$\varphi$,
\begin{equation*}
\ld\varphi(x,y) \ =\ \varphi(x\vee y) \,-\, \varphi(x\wedge y),
\end{equation*}
in Subsection~\ref{SS:vals_d}.
Finally, 
we study the equivalence induced by this distance,
\begin{equation*}
x\approx y\quad\iff\quad \ld\varphi(x,y)=0,
\end{equation*}
in Subsection~\ref{SS:vals_eq}.
The notion of distance 
is especially important.

We end the section 
some exotic
examples
(in Subsection~\ref{SS:more_examples}).
%%%%%%%%%%%%%%%%%%%%%%%%%%%%%%%%%%%%%%%%%%%%%%%%%%%%%%%%%%%%%%%%%%%%%%%%%%%%%%%
%%%%%%%%%%%%%%%%%%%%%%%%%%%%%%%%%%%%%%%%%%%%%%%%%%%%%%%%%%%%%%%%%%%%%%%%%%%%%%%
%%%%%%%%%%%%%%%%%%%%%%%%%%%%%%%%%%%%%%%%%%%%%%%%%%%%%%%%%%%%%%%%%%%%%%%%%%%%%%%
%%%%%%%%%%%%%%%%%%%%%%%%%%%%%%%%%%%%%%%%%%%%%%%%%%%%%%%%%%%%%%%%%%%%%%%%%%%%%%%
\subsection{Introduction}
\label{SS:vals_intro}
\noindent
%
%          DEFINITION OF MODULAR MAPS
%
\begin{dfn}

\label{D:val}
Let $L$ be a lattice. 
Let $E$ an ordered Abelian group
(see Section~\ref{S:ag}).\\
Let  $\varphi\colon L \ra E$ be a map.
We say that
\begin{enumerate}
\item
\label{D:val-mod}
$\varphi$ is \keyword{modular} provided that
\begin{equation*}
\varphi(a\wedge b) \,+\, \varphi(a \vee b)
\ =\ 
\varphi(a) \,+\, \varphi(b)
\qquad(a,b\in L);
\end{equation*}

\item
\label{D:val-val}
$\varphi$ is a \keyword{valuation}
provided that $\varphi$ is modular and order preserving.
\end{enumerate}
\end{dfn}

\begin{ex}
Let $\mathcal{F}$ be the set of finite subsets of~$\N$,
and for each $A\in \mathcal{F}$,
let $\#(A)$ be the number of elements of~$A$.
Then we have 
\begin{equation*}
\#(A\cap B) \,+\, \#(A\cup B) \ =\ \#(A) \,+\, \#(B)
\qquad(A,B\in\mathcal{F}),
\end{equation*}
so obviously the map $\mathcal{F} \ra \N$
given by by $A\mapsto \#(A)$ is a valuation.
\end{ex}

%
%                  LEBESGUE MEASURE IS A VALUATION
%
\begin{ex}
\label{E:lmeas-val}
Let $\LA$ be the set of Lebesgue measurable
subsets of~$\R$ with finite Lebesgue measure.
Then~$\LA$ is a lattice of subsets of~$\R$.
Given~$A\in\LA$,
let $\Lmu(A)$ denote the Lebesgue measure of~$A$.
Then $A\subseteq B \implies \Lmu(A)\leq \Lmu (B)$,
and 
\begin{equation*}
\Lmu(A\cap B) \,+\, \Lmu(A\cup B) \ =\  \Lmu(A) \,+\, \Lmu(B),
\end{equation*}
where $A,B\in \LA$.
So~$\Lmu$ is a valuation.
\end{ex}
%
%                  HISTORICAL REMARK
%
\begin{rem}
Valuations have been known for a long time,
see~\cite{Birkhoff67}.
\end{rem}

%
%                  LEBESGUE INTEGRAL IS A VALUATION
%
\begin{ex}
\label{E:int-val}
Let~$\LF$ be the set of
Lebesgue integrable functions on~$\R$.
When we write ``function on~$\R$''
we mean a map $f\colon \R \ra \E$.
Allowing the infinite values $+\infty$
and $-\infty$
might make the story a bit more complicated
in the short run,
but it will turn out to be a convenient choice later on
(see Remark~\ref{R:non-finite-functions}).

The set~$\LF$ is a lattice
under pointwise ordering,
and  $\LF\cap \R^\R$
is a lattice ordered Abelian group.
The assignment $f\mapsto \int f$ 
yields an order preserving map
\begin{equation*}
\Lphi\colon \LF\longrightarrow \R
\end{equation*}
that is group homomorphism
restricted to~$\LF\cap \R^\R$.

It takes some work see that $\Lphi$ is modular (and hence a valuation).

First,
note that for $x,y\in \R$, we have
\begin{equation*}
\min\{x,y\}\,+\,\max\{x,y\} \ =\ x\,+\,y.
\end{equation*}
So given $f,g\in \LF\cap \R^\R$,
we have $f\wedge g + f \vee g = f+ g$,
and hence
\begin{equation}
\label{eq:lint-fin-mod}
\begin{alignedat}{3}
\Lphi(f\wedge g) \,+\, \Lphi(f\vee g) 
\ &=\ \Lphi(f\wedge g + f\vee g)\\
\ &=\  \Lphi(f+g)\\
\ &=\  \Lphi (f)\,+\,\Lphi(g).
\end{alignedat}
\end{equation}
So we see that~$\Lphi$ is modular on~$\LF\cap \R^\R$.

To see that $\Lphi$ is modular on~$\LF$, we need some observations.
\begin{enumerate}
\item
\label{negl:1}
Let $f\in\LF$.
Then the set of~$x\in \R$ such that $f(x) = +\infty$ 
or $f(x)=-\infty$
is negligible.
Define $f_\R\colon \R \ra \R$ by,
for $x\in \R$,
\begin{equation*}
f_\R(x) \ =\
\begin{cases}
\ f(x) &\text{if $f(x)\in\R$,}\\
\ 0 &\text{otherwise.}
\end{cases}
\end{equation*}
Then $f(x)=f_\R(x)$
for almost all~$x\in \R$.

\item
\label{negl:2}
Let $f_1,f_2 \in \LF$
be given and assume $f_1(x) = f_2(x)$
for almost all~$x\in \R$.
(We denote this by $f_1 \approx f_2$.)
Then we have  $\Lphi(f_1) = \Lphi(f_2)$.

\item
\label{negl:3}
Let $f_1,f_2\in \LF$ with $f_1 \approx f_2$ be given,
and let $g\in \LF$.\\
Then $f_1 \wedge g \approx f_2 \wedge g$
and $f_1 \vee g \approx g_2 \vee g$.
\end{enumerate}
Now,
let $f,g\in \LF$ be given.
To prove that $\Lphi$ is modular,
we must show that 
\begin{equation*}
\Lphi(f)\,+\,\Lphi(g) \ =\  \Lphi(f\wedge g) \,+\, \Lphi(f\vee g).
\end{equation*}
Indeed, we have:
\begin{alignat*}{3}
\Lphi(f) \,+\, \Lphi(g) 
  \ &=\  \Lphi(f_\R) \,+\, \Lphi(g_\R)
    \qquad&&\text{by \ref{negl:1} and \ref{negl:2}}\\
  \ &=\ \Lphi(f_\R\wedge g_\R) \,+\, \Lphi(f_\R \vee g_\R) 
    \qquad&&\text{by Statement~\eqref{eq:lint-fin-mod}}\\
  \ &=\ \Lphi(f\wedge g_\R) \,+\, \Lphi(f \vee g_\R)
    \qquad&&\text{see \ref{negl:3}}\\
  \ &=\ \Lphi(f\wedge g) \,+\, \Lphi(f \vee g)
    \qquad&&
\end{alignat*}
Hence the Lebesgue integral $\Lphi\colon \LF\ra \R$ is a valuation.
\end{ex}

%
%                  ANY MAP ON A CHAIN IS MODULAR
%
\begin{ex}
Let $C$ be a chain,
i.e. a totally ordered set.
Then $C$ is a lattice with
\begin{equation*}
a\wedge b = \min\{a,b\},
 \qquad 
a\vee b = \max\{a,b\}.
\end{equation*}
One quickly sees that
\emph{any map $f\colon C\rightarrow E$
to an ordered Abelian group is modular.}
\end{ex}

%
%                  RING OF SUBSETS
%
\begin{ex}
\label{E:ring-val}
Let~$X$ be a set and let~$\mathcal{A}$ be a ring of subsets of~$X$.
That is,
\begin{equation*}
A\cap B,\qquad\qquad A\cup B,\qquad\qquad A\backslash B
\end{equation*}
are in~$\mathcal{A}$ for all $A,B\in\mathcal{A}$.
Then clearly~$\mathcal{A}$ is a lattice.

Let $E$ be an ordered Abelian group
and let~$\mu\colon \mathcal{A}\rightarrow E$ be a map.
Recall that~$\mu$ is \keyword{additive} if $\mu(A) + \mu(B) = \mu(A\cup B)$
for all $A,B\in\mathcal{A}$ with $A\cap B=\varnothing$.

\emph{If~$\mu$ is additive,
then $\mu$ is modular.}
Indeed,
let $A,B\in \mathcal{A}$ be given. We need to prove that
$\mu(A) + \mu(B) =\mu(A\cap B) + \mu(A\cup B)$
assuming~$\mu$ is additive.
We have
\begin{alignat*}{6}
\mu(A) + \mu(B) \,
  & =\  \mu(\,A\cap B \ \cup\ A\backslash B\,) \,+\, \mu(B) \\ 
  & =\  \mu(A\cap B) \,+\, \mu(A\backslash B)  \,+\, \mu(B)\qquad
    && \text{since } A\cap B \ \cap\ A\backslash B = \varnothing \\ 
  & =\  \mu(A\cap B) \,+\, \mu(A\backslash B \ \cup\ B ) 
    && \text{since } A\backslash B\ \cap\ B = \varnothing \\
  & =\  \mu(A\cap B) \,+\, \mu(A\cup B).
\end{alignat*}

Recall that~$\mu$ is \keyword{positive} whenever
$\mu(A)\in E^+$ for all~$A\in\mathcal{A}$.

\emph{If~$\mu$ is additive and positive,
then~$\mu$ is a valuation.}
Since~$\mu$ is additive,
$\mu$ is modular.
It remains to be shown 
(see Definition~\ref{D:val})
that~$\mu$ is order preserving.\\
Let $A\subseteq B$ from~$\mathcal{A}$ be given
in order to prove $\mu(A)\leq \mu(B)$.
We have
\begin{equation*}
(B\backslash A)\,\cup\, A\,=\,B,\qquad\qquad 
(B\backslash A)\,\cap\, A\,=\,\varnothing.
\end{equation*}
So by additivity, 
$\mu(B)=\mu(B\backslash A)+\mu(A)$.
Then $\mu(B)\geq \mu(A)$, since $\mu(B\backslash A)\geq 0$.
\end{ex}

%
%                  RING OF SIMPLE LEBESGUE SUBSETS
%
\begin{ex}
\label{E:smeas-val}
We describe a ring of subsets of~$\R$
and a positive and additive 
map~$\Smu\colon \SA\ra \R$ that 
will eventually
lead to the Lebesgue measure.

Let $\mathcal{S}$ be the set of all subsets of~$\R$ of the form,
with $a\leq b$ from~$\R$,
\begin{equation*}
(a,b)\qquad\text{or}\qquad  [a,b].
\end{equation*}
Let~$\SA$ be the ring generated by~$\mathcal{S}$.
Every element~$A$ of~$\SA$ is of the form
\begin{equation*}
I_1 \,\cup \,\dotsb \,\cup\, I_N
\end{equation*}
where $I_1,\dotsc,I_N\in \mathcal{S}$
are disjoint.
Let $\Smu(A)$ be given by
\begin{equation*}
\mu_{\mathrm L}(A) \ \eqdf\  |I_1| + \dotsb + |I_N|.
\end{equation*}
One can verify that the number $\Smu(A)$
only depends on~$A$ and not on the choice of~$I_1,\dotsc,I_N$.
Hence we obtain a map~$\Smu\colon \SA \ra \R$.
Almost by definition $\Smu$ is additive and positive.
Hence $\Smu\colon\SA\ra \R$ is a valuation 
(see Example~\ref{E:ring-val}).
\end{ex}

In Example~\ref{E:int-val},
we saw a group homomorphism that is modular,
namely
the Lebesgue integral~$\Lphi$ restricted to~$\R^\R$.
In fact any group homomorphism
on a lattice ordered Abelian group is modular
(see Corollary~\ref{C:hom-val}\ref{C:hom-val-group}).
\begin{ex}
\label{E:1-valuation}
Let~$R$ be a lattice ordered Abelian group.
Then the identity map $\mathrm{id}_R$ is a valuation.
Indeed, $\mathrm{id}_R$ is modular by Lemma~\ref{L:1-valuation},
and clearly order preserving.
\end{ex}

\begin{lem}
\label{L:mod-comp}
Suppose we have the following situation.
\begin{equation*}
\xymatrix{
L' \ar[r]^f&
L \ar[r]^\varphi&
E \ar[r]^g&
E',
}
\end{equation*}
where $L$, $L'$ are lattices,
$E$, $E'$ are ordered Abelian groups,
$f$ is a lattice homomorphism,
$\varphi$ a map,
and $g$ is a group homomorphism.
Then
\begin{enumerate}
\item
\label{L:mod-comp-mod}
$g\circ \varphi \circ f$ is modular
provided that $\varphi$ is modular;
\item
\label{L:mod-comp-val}
$g\circ \varphi \circ f$ is a valuation
provided that $\varphi$ is a valuation
and~$g$ is positive.
\end{enumerate}
\end{lem}
\begin{proof}
\noindent
\ref{L:mod-comp-mod}
\  Suppose~$\varphi$ is modular.
Let $a,b\in L'$ be given.
Writing $\varphi'= g\circ\varphi \circ f$,
we need to prove that
$\varphi'(a\wedge b)+\varphi'(a\vee b)=\varphi'(a)+\varphi'(b)$.
We have
\begin{alignat*}{3}
\varphi'(a) + \varphi'(b)
\ &=\ g(\varphi(f(a))) \,+\, g(\varphi(f(b))) \\
  &=\ g(\ \varphi( f(a)) + \varphi( f(b))\ ) \\
  &=\ g(\ \varphi(f(a)\wedge f(b)) + \varphi(f(a)\vee f(b))\ ) \\
  &=\ g(\  \varphi(f(a\wedge b)) \,+\, \varphi(f(a\vee b))\ ) \\
  &=\ g(\varphi(f(a\wedge b))) \,+\, g(\varphi(f(a\vee b))) \\ 
  &=\ \varphi'(a\wedge b) + \varphi'(a \vee b)
\end{alignat*}
\ref{L:mod-comp-val}
\  Suppose~$\varphi$ is a valuation
and~$g$ is positive.
We need to prove that~$\varphi'\eqdf g\circ \varphi\circ f$
is a valuation.
By part~\ref{L:mod-comp-mod}
we know that~$\varphi'$ is modular.
It remains to be shown that~$\varphi'$ is order preserving.
This is easy: $g$, $\varphi$, and~$f$ are all order preserving.
So $\varphi'=g\circ\varphi\circ f$ must be order preserving too.
\end{proof}

\begin{cor}
\label{C:hom-val}
Let~$R$ be a lattice ordered Abelian group.
\begin{enumerate}
\item
\label{C:hom-val-lat}
Let~$L$ be a lattice.
Any lattice homomorphism $f\colon L\ra R$ 
is a valuation.

\item
\label{C:hom-val-group}
Let $E$ be an ordered Abelian group
and $g\colon R\ra E$ a group homomorphism.
Then $g$ is modular.
Moreover,
if~$g$ is positive,
then~$g$ is a valuation.
\end{enumerate}
\end{cor}
\begin{proof}
Apply Lemma~\ref{L:mod-comp} to the following situations.
\begin{equation*}
\xymatrix{
L\ar[r]^f &
R\ar[r]^{\mathrm{id}_R} &
R\ar[r]^{\mathrm{id}_R} &
R &&
R\ar[r]^{\mathrm{id}_R} &
R\ar[r]^{\mathrm{id}_R} &
R\ar[r]^{g} &
E}
\end{equation*}
(Recall that $\mathrm{id}_R$ is a valuation,
see Example~\ref{E:1-valuation}.)
\end{proof}

%
%                  RIESZ SPACE OF FUNCTIONS
%
\begin{ex}
\label{E:val-riesz-space-of-functions}
Let $X$ be a set.
We say that $F\subseteq \R^X$
is  \emph{Riesz space of functions} if
\begin{equation*}
f\vee g,\quad\qquad 
f\wedge g,\quad\qquad
f+g,\quad\qquad 
\lambda \cdot f
\end{equation*}
are all in~$F$
where $f,g\in F$ and $\lambda \in \R$.
Then~$F$ is a lattice ordered Abelian group.

Let~$E$ be an ordered Abelian group
and let $\varphi\colon F\ra E$ be a positive linear map.
We see that~$\varphi$ is a valuation
 by Corollary~\ref{C:hom-val}\ref{C:hom-val-group}.
\end{ex}

%
%                  RIESZ SPACE OF STEP FUNCTIONS
%
\begin{ex}
\label{E:sint-val}
We  describe a Riesz space
of functions~$\SF$ on~$\R$
and a positive linear map~$\Sphi\colon \SF\ra \R$
that will eventually lead to the Lebesgue integral.

A \emph{step function} is a function~$f\colon \R\ra\R$
for which there are $s_1 < s_2 <\dotsb <s_N$ in~$\R$
such that $f$ is constant on each~$(s_n,s_{n+1})$
and $f$ is zero outside $[s_1,s_N]$.

Let $\SF$ be the set of step functions.
One can easily see that~$\SF$ is a Riesz space of functions.
Let $f\in \SF$.
Let $s_1 < s_2 <\dotsb <s_N$
be such that $f$ is constant, say $c_n\in \R$,
 on~$(s_n,s_{n+1})$
and $f$ is zero outside $[s_1,s_N]$.
One can prove that 
\begin{equation}
\label{exp:step}
\sum_{n=1}^{N-1} \, c_n\cdot(s_{n+1} - s_n)
\end{equation}
does not depend on the choice
of~$s_1 < s_2 <\dotsb <s_N$.
So Expression~\eqref{exp:step}
gives a map $\Sphi \colon \SF \ra \R$.
This map is easily seen to be linear.

Consequently, $\Sphi\colon \SF\ra \R$
is a valuation (see Example~\ref{E:val-riesz-space-of-functions}).
\end{ex}

\noindent
We end this subsection
with some tame examples of valuations we 
need later on.
\begin{ex}
\label{E:val-product}
Let $I=\{ 1,2\}$.
For each $i\in I$, 
let $L_i$ be a lattice,
$E_i$ an ordered Abelian group,
and $\varphi_i \colon L_i \ra E_i$
a valuation.
Then 
the map 
\begin{equation*}
\varphi_1 \times \varphi_2 \colon \,
L_1 \times L_2 \,\longrightarrow\, E_1 \times E_2,
\end{equation*}
given by $(\varphi_1 \times\varphi_2)(a_1,a_2) = (\varphi_1(a),\varphi_2(b))$
for all~$a_i\in L_i$, is a valuation.

We call the valuation $\varphi_1 \times \varphi_2$
the \emph{product} of $\varphi_1$ and $\varphi_2$.
Of course,
one can similarly define a product of an $I$-indexed family
of valuations for any set~$I$.
\end{ex}

\begin{ex}
\label{E-val-opposite}
Let $L$ be a lattice.
If we reverse the order on~$L$,
i.e., consider the partial order on~$L$ 
given by $\smash{a \leq_{L^\mathrm{op}} b}
\iff a\geq_L b$,
then if a subset~$S\subseteq L$
has a supremum, $\bigvee S$,
then $\bigvee S$ is the 
\emph{infimum} of~$S$
with respect to~$\leq^\mathrm{op}$.
So we see that $\leq^\mathrm{op}$
gives us a lattice, $L^\mathrm{op}$.
(The \emph{opposite} lattice.)

Let~$E$ be an ordered Abelian group.
If we reverse the order on~$E$,
we obtain an ordered Abelian group~$E^\mathrm{op}$
with the same group structure,
but whose positive elements, $\smash{(E^\mathrm{op})^+}$,
are precisely the negative elements of~$E$.

Let~$\varphi\colon L\ra E$ be a modular map
(see Definition~\ref{D:val}).
Then one quickly sees that $\varphi$ is also modular considered as a map
$L^\mathrm{op} \ra E$.
However,
$\varphi\colon L^\mathrm{op}\ra E$
is a valuation (that is, also order preserving)
if and only if $\varphi\colon L\ra E$
is order \emph{reversing},
i.e., $a\leq b\implies \varphi(a)\geq \varphi(b)$ for
all $a,b\in L$.

Of course,
if $\varphi\colon L\ra E$ is a valuation,
then $\varphi$ is a valuation $L^\mathrm{op}\ra E^\mathrm{op}$.
\end{ex}

%%%%%%%%%%%%%%%%%%%%%%%%%%%%%%%%%%%%%%%%%%%%%%%%%%%%%%%%%%%%%%%%%%%%%%%%%%%%%%%
%%%%%%%%%%%%%%%%%%%%%%%%%%%%%%%%%%%%%%%%%%%%%%%%%%%%%%%%%%%%%%%%%%%%%%%%%%%%%%%
%%%%%%%%%%%%%%%%%%%%%%%%%%%%%%%%%%%%%%%%%%%%%%%%%%%%%%%%%%%%%%%%%%%%%%%%%%%%%%%
%%%%%%%%%%%%%%%%%%%%%%%%%%%%%%%%%%%%%%%%%%%%%%%%%%%%%%%%%%%%%%%%%%%%%%%%%%%%%%%
%%%%%%%%%%%%%%%%%%%%%%%%%%%%%%%%%%%%%%%%%%%%%%%%%%%%%%%%%%%%%%%%j%
%
%
%                  DISTANCE INDUCED BY 
%                     a valuation
%
%
\subsection{Distance Induced by a Valuation}
\label{SS:vals_d}
In this subsection,
we derive some facts 
concerning the following notion
of distance induced by a valuation.
\begin{dfn}
\label{D:d}
Let $E$ be an ordered Abelian group.
Let $L$ be a lattice.\\
Let $\varphi\colon L \ra E$ be a valuation.
Define $\ld\varphi\colon L\times L \ra E$ by
\begin{equation*}
\ld\varphi(a,b)\ =\  \varphi(a\vee b) - \varphi(a \wedge b)
\qquad\quad(a,b\in L).
\end{equation*}
\end{dfn}

To give the name ``distance'' for~$\ld\varphi$ some credibility,
we will prove that $\ld\varphi$ 
is a pseudometric (see Lemma~\ref{L:d-metric}).
After that,
we turn our attention to the following fact,
which we will use often.
Given~$a\in L$, the map $x\mapsto a\wedge x$ is a \emph{contraction}, i.e.,
\begin{equation*}
\ld\varphi(\,a\wedge x,\, a\wedge y\,)
\ \ \leq\ \ \ld\varphi(x,y)
\qquad\quad(x,y\in L).
\end{equation*}
In fact,
we will prove the following, stronger, statement
(see Lemma~\ref{L:curry-wc-unif}).
\begin{equation*}
\ld\varphi(\,a\wedge x,\, a\wedge y\,) \ +\ 
\ld\varphi(\,a\vee x,\, a\vee y\,)
\ \ \leq\ \ \ld\varphi(x,y)
\qquad\quad(x,y\in L).
\end{equation*}

Before we do all this,
let us consider some examples.
%
%                  EXAMPLES OF DISTANCE
%
\begin{ex}
\label{E:d-riesz}
Let~$E$ be an ordered Abelian group.
Let $F$ be a Riesz space of functions,
and let $\varphi\colon F\ra E$ be
a positive and linear map
(see Example~\ref{E:val-riesz-space-of-functions}).

Let $f,g\in F$ be given.
The distance between $f$ and $g$ is the usual one,
\begin{equation*}
\ld\varphi(f,g) \,=\, \|f-g\|_1 \ \eqdf\ \varphi(\,|f-g|\,).
\end{equation*}
To see this,
note that since~$\varphi$ is linear,
we have 
\begin{equation*}
\ld\varphi(f,g) \,=\, 
\varphi(f\vee g) - \varphi(f\wedge g)
\,=\, \varphi(f\vee g - f\wedge g).
\end{equation*}
Further,
since we have the identity $\max\{x,y\} - \min\{x,y\} = |x-y|$
for reals $x,y$,
we have the identity $f\vee g - f\wedge g = |f-g|$
for functions.
\end{ex}
\begin{ex}
Let $E$ be an ordered Abelian group.
Let $\mathcal{A}$ be a ring of sets,
and let $\mu\colon \mathcal{A}\ra E$
be a positive additive map
(see Example~\ref{E:ring-val}).
Let $A,B\in\mathcal{A}$.
We have
\begin{equation*}
\ld\mu(A,B) \ =\ \mu(A \ominus B),
\end{equation*}
where $A\ominus B \eqdf A\backslash B \,\cup\, B\backslash A$
is the \emph{symmetric difference} of~$A$ and~$B$.
To see this, note that $A\cup B$ is the disjoint union
of $A\ominus B$ and $A\cap B$. So
since $\mu$ is additive, 
\begin{equation*}
\mu(A\cup B) \ =\ \mu(A\ominus B) + \mu(A\cap B).
\end{equation*}
\end{ex}
%
%                  LEMMA ON THE DISTANCE
%
\begin{lem}
\label{L:d-metric}
Let $E$ be an ordered Abelian group.\\
Let $L$ be a lattice,
and let $\varphi\colon L \ra E$ be a valuation.\\
Let $a,b,z\in L$ be given.
We have:
\begin{enumerate}
\item \label{d-metric_pos}
$\ld\varphi(a,b)\,\geq\, 0$
\item\label{d-metric_self} 
$\ld\varphi(a,a)\,=\,0$
\item\label{d-metric_sym}
$\ld\varphi(a,b)\,=\,\ld\varphi(b,a)$
\item\label{d-metric_triangle}
$\ld\varphi(a,b)\,\leq\,\ld\varphi(a,z)\,+\,\ld\varphi(z,b)$
\end{enumerate}
\end{lem}
\begin{proof}
Only point~\ref{d-metric_triangle} requires some work.
Let $a,b,z\in L$ be given.
We want to show that $\ld\varphi(a,b) \leq \ld\varphi(a,z)+\ld\varphi(z,b)$.
In other words:
\begin{equation}
\label{eq:d-metric}
\varphi(a\vee b) \,+\, \varphi(a\wedge z) \,+\, \varphi (z\wedge b)
\ \leq\ 
\varphi (a\vee z) \,+\, \varphi(z\vee b) \,+\, \varphi (a\wedge b)\text{.}
\end{equation}
By modularity,
the left-hand side equals
\begin{equation*}
\varphi(a\vee b) 
 \,+\, \varphi(\,(a\wedge z)\vee(b\wedge z)\,)
 \,+\, \varphi(a\wedge b\wedge z).
\end{equation*}
On the other hand,
using modularity
the right-hand side of Inequality~\eqref{eq:d-metric} becomes
\begin{equation*}
\varphi(a\vee b\vee z)
 \,+\, \varphi(\,(a\vee z)\wedge(b\vee z)\,)
 \,+\, \varphi(a\wedge b).
\end{equation*}
Note that $a\vee b \leq a\vee b\vee z$,
and $(a\wedge z)\vee (b\wedge z) \leq z \leq (a\vee z)\wedge (b\vee z)$,
and $a\wedge b\wedge z \leq a\wedge b$,
so that the monotonicity of~$\varphi$ yields Inequality~\eqref{eq:d-metric}.
\end{proof}
%
%                  DEFINITION HAUSDORFF
%
\noindent
It is possible that $\ld\varphi(a,b) = 0$ while $a\neq b$
(see Example~\ref{E:eq-int}).
So in general, $\ld\varphi$ is not a metric
(but merely a  \emph{pseudo}metric).
Those $\varphi$ for which~$\ld\varphi$ is a metric
turn out to be useful. So let us give them a name.
\begin{dfn}
\label{D:val_Hausdorff}
\label{D:hausdorff}
Let $L$ be a lattice.
Let $E$ be an ordered Abelian group.\\
Let $\varphi\colon L \ra E$ be a valuation.
We say $\varphi$ is \keyword{Hausdorff}
provided that
\begin{equation*}
\ld\varphi(a,b) = 0 
\quad \implies \quad a=b \qquad\qquad(a,b\in L).
\end{equation*}
\end{dfn}
\noindent
We return to Hausdorff valuations in Subsection~\ref{SS:vals_eq}.
%
%                  CURRY-WC-UNIF
%
\begin{lem}
\label{L:curry-wc-unif}
Let $E$ be an ordered Abelian group.\\
Let $L$ be a lattice,
and $\varphi\colon L \ra E$ a valuation.
We have,
for $a,b,z\in L$,
\begin{equation*}
\ld\varphi(a\wedge z,\,b\wedge z)\,+\,
 \ld\varphi(a\vee z,b\vee z)\ \leq\ \ld\varphi(a,b).
\end{equation*}
\end{lem}
\begin{proof}
By expanding Definition~\ref{D:d},
we see that 
we need to prove that
\begin{equation}
\label{eq:curry-wc-unif}
\begin{split}
\varphi(a\vee b\vee z) \,+\,
\varphi(\,(a\wedge z)\vee(b\wedge z)\,) & \,+\,
\varphi(a\wedge b) \\
\ \leq\ 
\varphi(a \vee b)  \,+\, &
\varphi(\,(a\vee z)\wedge(b\vee z)\,) \,+\,
\varphi(a\wedge b \wedge z)
\end{split}
\end{equation}
By modularity,
the left-hand side equals
\begin{equation*}
\varphi(a\vee b \vee z) \,+\,
\varphi(\, (a\wedge z) \vee (b\wedge z) \vee (a\wedge b)\,) \,+\,
\varphi(\,a\wedge b\wedge ( (a \wedge z) \vee (b\wedge z) )\,).
\end{equation*}
To simplify the above expression,
we prove that
$a\wedge b\wedge ( (a \wedge z) \vee (b\wedge z)=a\wedge b\wedge z$.
To this end, note that
$a\wedge z \,\leq\, (a\wedge z)\vee (b\wedge z) \leq z$
so that 
\begin{equation*}
a\wedge b\wedge z
\, =\, a\wedge b \wedge (a\wedge z)
\,\leq\, a\wedge b \wedge ((a\wedge z)\vee (b\wedge z))
\,\leq\, a\wedge b \wedge z.
\end{equation*}
Hence the left-hand side of Inequality~\eqref{eq:curry-wc-unif}
equals
\begin{equation*}
\varphi(a\vee b \vee z) \,+\,
\varphi(\, (a\wedge z) \vee (b\wedge z) \vee (a\wedge b)\,) \,+\,
\varphi(\,a\wedge b\wedge z\,).
\end{equation*}
In a similar fashion,
one can show that the right-hand side of Inequality~\eqref{eq:curry-wc-unif}
equals
\begin{equation*}
\varphi(a\vee b \vee z) \,+\,
\varphi(\, (a\vee z) \wedge (b\vee z) \wedge (a\vee b)\,) \,+\,
\varphi(\,a\wedge b\wedge z\,).
\end{equation*}
So in order to prove Inequality~\eqref{eq:curry-wc-unif},
we must show that
\begin{equation*}
\varphi(\, (a\wedge z) \vee (b\wedge z) \vee (a\wedge b)\,) 
\ \leq\ 
\varphi(\, (a\vee z) \wedge (b\vee z) \wedge (a\vee b)\,).
\end{equation*}
Since $\varphi$ is order preserving,
it suffices to show that 
\begin{equation*}
 (a\wedge z) \vee (b\wedge z) \vee (a\wedge b)
\ \leq\ 
  (a\vee z) \wedge (b\vee z) \wedge (a\vee b).
\end{equation*}
Writing $c_1 = a$, $c_2 = b$, and $c_3 = z$,
we must prove that 
\begin{equation*}
\bv_{i\neq j}\, c_i \wedge c_j \ \leq\ \bw_{k\neq \ell}\, c_k \vee c_\ell.
\end{equation*}
That is,
we must show that 
$c_i \wedge c_j \leq c_k \vee c_\ell$
for given $i\neq j$ and $k\neq \ell$.
Now, note
\begin{equation*}
\# (\{ i,j \} \cap \{ k,\ell\}) + \#\{i,j,k,\ell\}
\ = \ \# \{i,j\} + \# \{k,\ell \}
\ = \ 4.
\end{equation*}
Since  $\# \{ i,j,k,\ell \}\leq 3$,
we see that $\# \{ i,j \} \cap \{ k,\ell\} \geq 1$.
So pick $m\in  \{ i,j \} \cap \{ k,\ell\}$.
Then $c_i \wedge c_j \leq c_m \leq c_k \vee c_\ell$.
\end{proof}
%
%                  WC-UNIF
%
\begin{lem}
\label{L:wv-unif}
Let $E$ be an ordered Abelian group.\\
Let $L$ be a lattice,
and $\varphi\colon L \ra E$ a valuation.
Then we have
\begin{equation*}
\ld\varphi(a\wedge w,b\wedge z) \,+\, \ld\varphi(a \vee w,b\vee z)
\ \leq\ 
\ld\varphi(a,b) \,+\, \ld\varphi(w,z),
\end{equation*}
where $a,b,w,z\in L$.
\end{lem}
\begin{proof}
By the triangle inequality (point~\ref{d-metric_triangle}
of Lemma~\ref{L:d-metric}),
we have
\begin{equation}
\label{eq:L:d-metric-1}
\begin{alignedat}{3}
\ld\varphi(a\wedge w,\,b\wedge z)\ &\leq\ 
\ld\varphi(a\wedge w,\,b\wedge w)\,+\,
\ld\varphi(b\wedge w,\,b\wedge z), \\
\ld\varphi(a\vee w,\,b\vee z)\ &\leq\ 
\ld\varphi(a\vee w,\,b\vee w)\,+\,
\ld\varphi(b\vee w,\,b\vee z).
\end{alignedat}
\end{equation}
On the other hand, Lemma~\ref{L:curry-wc-unif} gives us
\begin{equation}
\label{eq:L:d-metric-2}
\begin{alignedat}{3}
\ld\varphi(a\wedge w,b\wedge w) + \ld\varphi(a\vee w,b\vee w)
   \ &\leq\ \ld\varphi(a,b), \\
\ld\varphi(b\wedge w,b\wedge z) + \ld\varphi(b\vee w, b\vee z)
   \ &\leq\ \ld\varphi(w,z).
\end{alignedat}
\end{equation}
The sum of the right-hand sides of Equation~\eqref{eq:L:d-metric-1}
equals the sum of the left-hand sides of Equation~\eqref{eq:L:d-metric-2}.
Hence
$\ld\varphi(a\wedge w,b\wedge z) + \ld\varphi(a \vee w,b\vee z)
\leq 
\ld\varphi(a,b) + \ld\varphi(w,z)$.
\end{proof}
%%%%%%%%%%%%%%%%%%%%%%%%%%%%%%%%%%%%%%%%%%%%%%%%%%%%%%%%%%%%%%%%%%%%%%%%%%%%%%%
%%%%%%%%%%%%%%%%%%%%%%%%%%%%%%%%%%%%%%%%%%%%%%%%%%%%%%%%%%%%%%%%%%%%%%%%%%%%%%%
%%%%%%%%%%%%%%%%%%%%%%%%%%%%%%%%%%%%%%%%%%%%%%%%%%%%%%%%%%%%%%%%%%%%%%%%%%%%%%%
%%%%%%%%%%%%%%%%%%%%%%%%%%%%%%%%%%%%%%%%%%%%%%%%%%%%%%%%%%%%%%%%%%%%%%%%%%%%%%%
%%%%%%%%%%%%%%%%%%%%%%%%%%%%%%%%%%%%%%%%%%%%%%%%%%%%%%%%%%%%%%%%%%%%%%%%%%%%%%%
%
%                  EQUIVALENCE INDUCED BY A VALUATION
%
\subsection{Equivalence Induced by a Valuation}
\label{SS:vals_eq}
In measure theory
two functions are considered equivalent
if they are equal almost everywhere.

In this subsection, we extend this notion of equivalence 
to valuations. 
%
%                  DEFINITION OF EQUIVALENCE
%
\begin{dfn}
\label{D:eq}
\label{D:approx}
Let $E$ be an ordered Abelian group.
Let $L$ be a lattice.\\
Let $\varphi\colon L\ra E$ be a valuation.
We define $\approx$ to be 
the binary relation on~$L$ given by
\begin{equation*}
a \approx b
\quad\iff\quad
\ld\varphi(a,b)=0\qquad\qquad(a,b\in L).
\end{equation*}
\end{dfn}
\begin{rem}
$\varphi$ is Hausdorff (see Definition~\ref{D:val_Hausdorff})
iff $a \approx b\iff a = b$.
\end{rem}
%
%                  EXAMPLES OF EQUIVALENCE
%
\begin{ex}
\label{E:eq-int}
We consider the Lebesgue integral $\Lphi\colon \LF\ra \R$
(see Example~\ref{E:int-val}).

Let $f,g\in \LF\cap \R^\R$ be given.
By Example~\ref{E:d-riesz} 
we know that
\begin{equation}
\label{eq:eq-int}
f\approx g\quad\iff\quad \Lphi(\,|f-g|\,)=0.
\end{equation}
In fact,
Statement~\eqref{eq:eq-int}
holds for all~$f,g\in \LF$,
as the reader can verify using the remarks
made in Example~\ref{E:int-val}.

Now, for any $h\in \LF$
with $h\geq 0$,
we have that
\begin{equation*}
\Lphi(h)\,=\,0
\quad\iff\quad \text{$h(x)=0$\quad for almost all~$x$.}
\end{equation*}
So we see that
we have, for $f,g\in \LF$,
\begin{equation*}
f\approx g\quad\iff\quad f(x) = g(x)\quad\text{for almost all~$x$}.
\end{equation*}
So ``$\approx$''
is \emph{equality almost everywhere}
when $\varphi=\Lphi$,
as was intended.
\end{ex}
%
%                  PROP ON EQUIVALENCE
%
\begin{prop}
\label{P:eq}
Let $E$ be an ordered Abelian group.
Let $L$ be a lattice.\\
Let $\varphi\colon L\ra E$ be a valuation.
Let $\approx$ be as in Definition~\ref{D:eq}.
\begin{enumerate}
\item
\label{P:eq-1}
The relation $\approx$ is an equivalence.
\item 
\label{P:eq-2}
Let~$a_1,a_2\in L$ with 
with $a_1\approx a_2$ be given.
Then $\varphi(a_1)=\varphi(a_2)$.
\item
\label{P:eq-3}
Let $a_1,a_2 \in L$ with $a_1\approx a_2$,
and let $b_1,b_2 \in L$ with $b_1 \approx b_2$ be given.
Then
\begin{equation*}
a_1 \wedge b_1 \,\approx\, a_2 \wedge b_2
\qquad\text{and}\qquad 
a_1 \vee b_1 \,\approx\, a_2 \vee b_2.
\end{equation*}
\item
\label{P:eq-4}
Let $a_1,a_2 \in L$ with $a_1\approx a_2$,
and let $b_1,b_2 \in L$ with $b_1 \approx b_2$ be given.
Then
\begin{equation*}
\ld\varphi(a_1,b_1)\ =\ \ld\varphi(a_2,b_2).
\end{equation*}
\end{enumerate}
\end{prop}
\begin{proof}
\ref{P:eq-1}\ 
The relation~$\approx$ is clearly reflexive
and symmetric. So to prove~$\approx$ is an equivalence relation,
we will only show that~$\approx$ is transitive.
Let $a,b,c\in L$ with $a\approx b\approx c$ be given.
We must show that $a\approx c$.
Or in other words, $\ld\varphi(a,c) = 0$.

By Lemma~\ref{L:d-metric}, points~\ref{d-metric_pos}
and~\ref{d-metric_triangle},
we get
\begin{equation}
\label{eq:P:eq-1}
0 \ \leq\ \ld\varphi(a,c) \ \leq\ 
\ld\varphi(a,b) + \ld\varphi(b,c).
\end{equation}
But 
$\ld\varphi(a,b)=0$ and $\ld\varphi(b,c)=0$,
since $a\approx b$ and $b\approx c$, respectively.

So we see that Statement~\eqref{eq:P:eq-1} implies $\ld\varphi(a,c)= 0$.
Hence $a\approx c$.

\vspace{.3em}
\ref{P:eq-2}\ 
Let $a_1,a_2\in L$ with $a_1\approx a_2$ be given.
We must prove $\varphi(a_1) = \varphi(a_2)$.

Let $i\in\{1,2\}$ be given. 
Note that $a_1 \wedge a_2 \leq a_i \leq a_1 \vee a_2$.
So we have 
\begin{equation}
\label{eq:P:eq-2}
\varphi(a_1 \wedge a_2) \ \leq\ \varphi(a_i) \ \leq\  \varphi(a_1 \vee a_2).
\end{equation}
Since $\ld\varphi(a_1,a_2)=0$,
we know that $\varphi(a_1 \vee a_2) = \varphi(a_1 \wedge a_2)$.
So Statement~\eqref{eq:P:eq-2}
implies that
$\varphi(a_1 \vee a_2) = \varphi(a_i) = \varphi(a_1 \wedge a_2)$.
Hence $\varphi(a_1) = \varphi(a_2)$.

\vspace{.3em}
\ref{P:eq-3}
Let $a_1,a_2 \in L$ with $a_1 \approx a_2$ be given.
Let $b_1,b_2 \in L$ with $b_1 \approx b_2$ be given.
We will only show that $a_1 \wedge b_1 \approx a_2 \wedge b_2$;
the proof of $a_1 \vee b_1 \approx a_2 \vee b_2$ is similar.

Note that we have the following inequalities by 
Lemma~\ref{L:d-metric}\ref{d-metric_pos}
and Lemma~\ref{L:wv-unif}.
\begin{equation}
\label{eq:P:eq-3}
0\ \leq\ 
\ld\varphi(a_1 \wedge b_1, a_2\wedge b_2)
\ \leq\ 
\ld\varphi(a_1,b_1) + \ld\varphi(a_2,b_2)
\end{equation}
Since $a_1 \approx a_2$ and $b_1 \approx b_2$,
we have $\ld\varphi(a_1,b_1)=0$
and $\ld\varphi(a_2,b_2)=0$,
respectively.
Hence Statement~\eqref{eq:P:eq-3} implies 
$\ld\varphi(a_1\wedge b_1, a_2 \wedge b_2)=0$.
So $a_1 \wedge b_1 \approx a_2 \wedge b_2$.

\vspace{.3em}
\ref{P:eq-4}
Let $a_1,a_2 \in L$ with $a_1 \approx a_2$ be given.
Let $b_1,b_2 \in L$ with $b_1 \approx b_2$ be given.
We must prove that $\ld\varphi(a_1,b_1)=\ld\varphi(a_2,b_2)$.
Note that
by point~\ref{P:eq-3}
we have
\begin{alignat*}{5}
a_1 \wedge b_1 &\approx a_2 \wedge b_2&
\qquad&\text{and}\qquad&
a_1 \vee b_1 &\approx a_2 \vee b_2.&
\shortintertext{%
So by point~\ref{P:eq-2} of this lemma, we get}
\varphi(a_1 \wedge b_1) &= \varphi( a_2 \wedge b_2)&
\qquad&\text{and}\qquad&
\varphi(a_1 \vee b_1)&=\varphi(a_2 \vee b_2).&
\end{alignat*}
So if we unfold Definition~\ref{D:d},
we see that
\begin{alignat*}{5}
\ld\varphi(a_1,b_1)
\ &=\ 
\varphi(a_1 \vee b_1) - \varphi(a_1 \vee b_1) \\
\ &=\ 
\varphi(a_2 \vee b_2) - \varphi(a_2 \vee b_2)
\ =\ 
\ld\varphi(a_2,b_2).\qedhere
\end{alignat*}
\end{proof}
%
%                  THE QUOTIENT LATTICE
%
\noindent
When studying the Lebesgue integrable functions, $\LF$,
it is sometimes convenient
to consider the space~$\mathrm{L}^1 = \qvL{\LF}$
of integrable functions modulo equality almost everywhere
(see Example~\ref{E:eq-int}).
Of course,
one can consider the space $\qvL{L}$ 
for any valuation~$\varphi\colon L\ra E$.
We list some of the properties of~$\qvL{L}$
in Proposition~\ref{P:quotient-lattice}.
\begin{prop}
\label{P:quotient-lattice}
Let $E$ be an ordered Abelian group.\\
Let $L$ be a lattice.
Let $\varphi\colon L \ra E$ be a valuation.
Let $\approx$ be as in Definition~\ref{D:eq}.\\
Let $\qvL{L}$ denote the quotient set,
and let $q\colon L\ra \qvL{L}$ be the
quotient map. Then:
\begin{enumerate}
\item
The set
$\qvL{L}$ is lattice if
the operations are given by
\begin{equation*}
qa \wedge qb \,=\, q(a\wedge b),
\qquad
 qa \vee qb \,=\, q(a\vee b)\qquad\quad(a,b\in L).
\end{equation*} 
Then, in particular, $q\colon L\ra \qvL{L}$ is a lattice homomorphism.

\item
There is a unique map $\qvphi\varphi\colon\qvL{L}\ra E$
such that
\begin{equation*}
(\qvphi\varphi)(q(a)) \ =\ \varphi(a)\qquad\quad(a\in L).
\end{equation*}
Moreover, the map $\qvphi\varphi$ is a valuation.

\item
We have
\ $\ld{\qvphi{\varphi}}(\,qa,\,qb\,) \,=\,\ld\varphi(a,b)$\ \ 
for all \ $a,b\in L$.

\item
\label{P:quotient-lattice-iv}
We have
\ $\ld{\qvphi{\varphi}}(\mathfrak{a},\mathfrak{b})=0
\ \implies \ 
\mathfrak{a} = \mathfrak{b}$\ \ 
for all \ $\mathfrak{a},\mathfrak{b}\in \qvL{L}$.
\end{enumerate}
\end{prop}
\begin{proof}
Follows from Proposition~\ref{P:eq}.
We leave the verification to the reader.
\end{proof}
\begin{rem}
Note that $\qvphi\varphi$ is Hausdorff 
(see Definition~\ref{D:val_Hausdorff})
by Proposition~\ref{P:quotient-lattice}\ref{P:quotient-lattice-iv}.
\end{rem}

\subsection{More Examples}
\label{SS:more_examples}
%
%                  TOTIENT FUNCTION
%
\noindent
Valuations also appear outside measure theory.\\
We begin with an example from elementary number theory.
\begin{ex}
Recall that \emph{Euler's totient function}
$\varphi$
is given by,
for $n\in \N$,
\begin{equation*}
\varphi(n) \ =\ \#\{ \  x\in \{1,\dotsc,n\}\colon \ \gcd\{x,n\}=1\ \}.
\end{equation*}
We will prove that~$\varphi$
is a valuation `with respect to the division order'.

More precisely,
we consider $\varphi$ to be a map
\begin{equation*}
\varphi\colon \N \longrightarrow  \Q^\circ,
\end{equation*}
where $\Q^\circ$
is the set of strictly positive rational numbers.
Write, for $q,r\in\Q^\circ$,
\begin{equation*}
q\dv r \quad\iff\quad \exists n\in\N\, [\ q \cdot n \ =\  r\ ].
\end{equation*}
We order the sets~$\N$ and $\Q^\circ$ by ``$\dv$''.
They are both lattices with
\begin{equation*}
a\wedge b \ =\ \gcd\{a,b\},\qquad a \vee b \ =\ \lcm\{a,b\}.
\end{equation*}
Moreover,
$\Q^\circ$
is an ordered Abelian group under the normal multiplication~``$\cdot$''.

Before we prove that $\varphi$ is a valuation,
we make a useful observation: for $n\in\N$,
\begin{equation}
\label{eq:phi-Z-inv}
\varphi(n) \ =\ \# \Zmod{n}^*.
\end{equation}
Here $\Zmod{n}$ is the set of integers modulo~$n$,
and $\Zmod{n}^*$
are the invertible elements of~$\Zmod{n}$.

To see that Equation~\eqref{eq:phi-Z-inv} holds,
note that for~$x\in \Z$, we have
\begin{alignat*}{3}
\gcd\{x,n\}\,=\,1
\quad&\iff&\quad
\exists a,b\in\Z,& \quad ax + bn \,=\,1  
\qquad&& \text{by B\'ezout's Lemma}\\
\quad&\iff&\quad
\exists a,b\in\Z,& \quad ax  \,=\,1 - bn 
\qquad&& \\
\quad&\iff&\quad
\exists a\in\Z,& \quad [a]_n\cdot[x]_n \,=\,1 \  
\qquad&& \\
\quad&\iff&\quad
& \quad [x]_n \,\in\, \Zmod{n}^*,
\qquad&& 
\end{alignat*}
where $[-]_n\colon \Z \rightarrow \Zmod{n}$
is the quotient map.

Let us now prove that~$\varphi$ is a valuation.
We first prove that~$\varphi$ is order preserving.
Let $m,n\in \N$ with $m\dv n$ be given.
We must show that $\varphi(m) \,\dv\, \varphi(n)$.

Note that there is a unique ring homomorphism
$h\colon \Zmod{n}\rightarrow \Zmod{m}$
given by,
for $x\in \Z$,
\begin{equation*}
h(\,[x]_n\,) \ =\ [x]_m.
\end{equation*}
Note that~$h$ is surjective,
and that $[x]_n$ if invertible iff $[x]_m$ is invertible
for $x\in \Z$.
So we see that if we restrict~$h$ to~$\Zmod{n}^*$,
we get a surjective group homomorphism
\begin{equation*}
\tilde{h}\colon \Zmod{n}^* \longrightarrow \Zmod{m}^*.
\end{equation*}
By Lagrange's Theorem we know that
\begin{equation*}
\# \Zmod{m}^*  \,\cdot\, \# \ker(\tilde h) \ =\ \# \Zmod{n}^*,
\end{equation*}
where $\ker(\tilde{h}) 
\eqdf\smash{\{\, a\in \Zmod{n} \colon \,\tilde{h}(a)=0\,\}}$
is the kernel of~$\tilde{h}$.
Thus
\begin{equation*}
\varphi(m) \,=\,
\# \Zmod{m}^*  \,\dv\,
 \# \Zmod{n}^*
\,=\,\varphi(n).
\end{equation*}
Hence $\varphi$ is order preserving.

It remains to be shown that~$\varphi$
is modular. That is,
for $m,n\in\N$,
\begin{equation}
\label{eq:totient-modular}
\varphi(\,\gcd\{m,n\}\,) \,\cdot\, \varphi(\,\lcm\{m,n\}\,)
\ =\ 
\varphi(m) \,\cdot\,\varphi(n).
\end{equation}

We first prove a special case,
namely,
that for $m,n\in\N$ with $\gcd\{m,n\}=1$,
\begin{equation}
\label{eq:totient-additive}
\varphi(m\cdot n) \ =\ \varphi(m)\,\cdot\,\varphi(n).
\end{equation}
By the Chinese Remainder Theorem we 
have the following isomorphism of rings.
\begin{equation*}
\Zmod{m} \times \Zmod{n} \ \cong\ \Zmod{m\cdot n}
\end{equation*}
As a consequence, we get the following
isomorphism of groups.
\begin{equation*}
\Zmod{m}^* \times \Zmod{n}^* \ \cong\ \Zmod{m\cdot n}^*
\end{equation*}
If we count the number of elements in the above
groups we see that
\begin{equation*}
\#\Zmod{m}^* \,\cdot\, \#\Zmod{n}^* \ =\ \#\Zmod{m\cdot n}^*.
\end{equation*}
Hence Equation~\eqref{eq:totient-additive}
holds (see Equation~\eqref{eq:phi-Z-inv}).

Let $m,n\in \N$ be given.
We prove that Equation~\eqref{eq:totient-modular} holds.

By the Fundamental Theorem of Arithmatic,
we have 
\begin{alignat*}{3}
m \,&=\,\prod_{p\in\Prm}\ p^{w(p)},&
\qquad
n \,&=\,\prod_{p\in\Prm}\ p^{v(p)},
\shortintertext{%
where $v,w\colon \Prm \ra \{0,1,2,\dotsc\}$
have finite support,
and $\Prm$ are the primes.
Hence}
\varphi(m) \,&=\,\prod_{p\in\Prm}\ \varphi(p^{w(p)}),&
\qquad
\varphi(n) \,&=\,\prod_{p\in\Prm}\ \varphi(p^{v(p)}),
\end{alignat*}
by Equation~\eqref{eq:totient-additive},
because $\gcd\{p_1^{k_1},p_2^{k_2}\}=1$
for all~$p_1 \neq p_2$ from~$\Prm$
and $k_1,k_2 \in \N$.

Let $p\in \Prm$ be given.
Note that either $w(p) \leq v(p)$ or $v(p)\leq w(p)$.
Hence,
\begin{alignat}{3}
\label{eq:totient-almost-done}
\varphi(p^{w(p)}) \,\cdot\, \varphi(p^{v(p)})
\ &=\ \varphi(\,p^{\min\{w(p),v(p)\}}\,)
 \,\cdot\, \varphi(\,p^{\max\{w(p),v(p)\}}\,).
\end{alignat}
This gives us the following equality.
\begin{equation*}
m\cdot n \ =\ 
\prod_{p\in\Prm} \ \varphi(\,p^{\min\{w(p),v(p)\}}\,)
 \ \cdot\  \prod_{p\in \Prm}\ \varphi(\,p^{\max\{w(p),v(p)\}}\,).
\end{equation*}
Note that\quad
$\gcd\{m,n\} = \prod_{p\in \Prm}\ p^{\min\{w(p),v(p)\}}$,\quad
so we have
\begin{equation*}
\varphi(\,\gcd\{m,n\}\,)
\  =\  \prod_{p\in \Prm}\ \varphi(\,p^{\min\{w(p),v(p)\}}\,).
\end{equation*}
Similarly, we have
\begin{equation*}
\varphi(\,\lcm\{m,n\}\,)
\  =\  \prod_{p\in \Prm}\ \varphi(\,p^{\max\{w(p),v(p)\}}\,).
\end{equation*}
If we apply the above equalities to Equation~\eqref{eq:totient-almost-done}
we get
\begin{equation*}
\varphi(m)\,\cdot\,\varphi(n) 
\ =\ \varphi(\,\gcd\{m,n\}\,)\,\cdot\, \varphi(\,\lcm\{m,n\}\,).
\end{equation*}
So $\varphi$ is modular.
Hence Euler's totient function $\varphi$ is a valuation.
\end{ex}

%
%                  EXAMPLE ON VECTOR SPACES
%
\noindent
Up to this point
we have only seen valuations on distributive lattices.
We will now give an example
of a valuation on a non-distributive lattice.
\begin{ex}
Let $W$ be a vector space.
Let $L$ be the set of finite-dimensional linear subspaces
of~$W$ ordered by inclusion.
  Then~$L$ is a lattice, and 
for all~$A,B\in L$,
\begin{equation*}
A\wedge B \ =\ A\cap B,
\qquad A\vee B \ =\ \left<A\cup B\right>,
\end{equation*}
where $\left< S \right>$ denotes the smallest
linear subspace containing $S$.
We have
\begin{equation*}
\dim (A\wedge B) \,+\, \dim(A\vee B)
\ =\ 
\dim A \,+\, \dim B 
\qquad\qquad(A,B\in L).
\end{equation*}
To see this,
apply the dimension theorem
to the map~$f\colon A\times B\ra A\vee B$ given by  $(a,b)\mapsto a+b$.
Hence the assignment $A \mapsto \dim A$
gives a valuation $\dim\colon L\ra \N$.

The lattice~$L$ might be distributive.
For instance, if~$W=\{ 0 \}$.
This occurs only seldom:
if $W$ contains two linearly independent vectors,
then~$L$ is non-distributive.

Indeed,
let $v_1,v_2\in W$ be linearly independent vectors
and consider  $w\eqdf v_1 + v_2$.
One can verify that $v_i, w$ are linearly independent too.
So $\left< v_i \right> \cap \left< w \right> = \{0\}$.
Hence
\begin{equation*}
\left< w \right> \wedge(\left< v_1 \right>\vee\left< v_2 \right>)
\ =\ 
\left< v_1, v_2 \right>
\ \neq\ 
\{0\}
\ =\ 
(\left<w\right> \wedge \left<v_1\right>)
\,\vee\, (\left<w\right> \wedge \left<v_2\right>).
\end{equation*}
\end{ex}
\noindent
It is interesting to note that
there are some `connections'
between
modular maps (see Definition~\ref{D:val})
and \emph{modular lattices}.
Recall that a lattice~$L$ is modular if
\begin{equation*}
\ell \vee (a \wedge u) \ =\ (\ell \vee a) \wedge u
\end{equation*}
for all~$\ell,u,a \in L$ with $\ell \leq u$.
One such connection is given by the following lemma.
%
%                  LEMMA ON MODULARITY FOR MODULAR MAPS
%
\begin{lem}
\label{L:modular-map-modular}
Let $E$ be an ordered Abelian group.
Let~$L$ be a lattice.\\
Let $\varphi\colon L \ra E$ 
be a modular map.
Let $\ell,u\in L$ with $\ell\leq u$ be given.
We have
\begin{equation}
\label{eq:modular-map}
\varphi(\,\ell \vee (a \wedge u)\,) 
\ =\ 
\varphi(\,(\ell\vee a)\wedge u\,)
\qquad (a\in L).
\end{equation}
\end{lem}
\begin{proof}
The trick is to consider the expression 
$\varphi(\ell) + \varphi(a) + \varphi(u)$.
On the one hand,
\begin{align*}
\varphi(\ell) + \varphi(a) + \varphi(u)
\ &=\ \varphi(\ell\wedge a) + \varphi(\ell \vee a) + \varphi(u) \\
\ &=\ \varphi(\ell \wedge a)
      + \varphi(\,(\ell\vee a)\wedge u\,)
      + \varphi(a\vee u),
\end{align*}
where we have used modularity twice.
On the other hand,
\begin{align*}
\varphi(\ell) + \varphi(a) + \varphi(u)
\ &=\ \varphi(\ell) + \varphi(a\wedge u) + \varphi(a \vee u) \\
\ &=\ \varphi(\ell \wedge a)
      + \varphi(\,\ell\vee (a\wedge u)\,)
      + \varphi(a\vee u).
\end{align*}
The difference,
$\varphi(\,(\ell\vee a)\wedge u\,)
- \varphi(\,\ell\vee (a\wedge u)\,)$,
must be zero.
\end{proof}

 }
\clearpage
{ %%%%%%%%%%%%%%%%%%%%%%%%%%%%%%%%%%%%%%%%%%%%%%%%%%%%%%%%%%%%%%%%%%%%%%%%%%
%
%
%                  C O M P L E T E      V A L U A T I O N S 
%
%
%
\section{Complete Valuations}
\label{S:complete-val}
\noindent
We now turn to the study of \emph{complete} valuations
(see Definition~\ref{D:complete-val}).
Among all valuations
the complete valuations resemble
the Lebesgue measure
and the Lebesgue integral most closely.
To support this claim,
we will
prove generalisations of
some of the classical convergence theorems of integration
in Subsection~\ref{SS:complete-val_convergence}.

But first, we give some examples of complete valuations
in Subsection~\ref{SS:complete-val_introduction}.

After that,
we  study a notion of completeness
for an ordered Abelian group~$E$,
called \emph{$R$-completeness},
in Subsection~\ref{SS:complete-val_R-completeness},
which will be useful later on.

The notion of complete valuation is not at the end of the 
road.  We will study the
slightly more sophisticated 
\emph{valuation systems}
(see Definition~\ref{D:system})
and \emph{complete valuation systems}
(see Definition~\ref{D:system-complete})
in Section~\ref{S:valuation-systems}.

%
%                  PHI-CONVERGENCE
%
\subsection{Introduction}
\label{SS:complete-val_introduction}
\begin{dfn}
\label{D:phi-conv}
Let $E$ be an ordered Abelian group.\\
Let $L$ be a lattice, and let $\varphi\colon L \ra E$ be a valuation.

Consider a sequence
$a_1 \geq a_2 \geq \dotsb$ from~$L$.
We say
\begin{equation*}
a_1 \geq a_2 \geq \dotsb \text{ is \keyword{$\varphi$-convergent}}
\qquad\text{if}\qquad \bw_n \varphi(a_n)\text{ exists.}
\end{equation*}

Similarly,
if
$b_1 \leq b_2 \leq \dotsb$ is
a sequence in~$L$, then 
\begin{equation*}
b_1 \leq b_2 \leq \dotsb \text{ is \keyword{$\varphi$-convergent}}
\qquad\text{if}\qquad \bv_n \varphi(b_n)\text{ exists.}
\end{equation*}
\end{dfn}
%
%                  COMPLETE VALUATION
%
\begin{dfn}
\label{D:complete-val}
Let $E$ be an ordered Abelian group. Let $L$ be a lattice.\\
Let $\varphi\colon L \ra E$ be a valuation.
We say $\varphi$ is \keyword{$\Pi$-complete} if
\begin{alignat*}{5}
a_1 \geq a_2 \geq \dotsb \text{$\varphi$-convergent }
  \quad &\implies \quad 
  & \bw_n a_n &\text{ exists,}\quad 
  &&\text{and}\quad
  &\varphi(\,\bw_n a_n\,) &= \bw_n \varphi(a_n).
\shortintertext{%
We say $\varphi$ is \keyword{$\Sigma$-complete} if
}
b_1 \leq b_2 \leq \dotsb \text{$\varphi$-convergent }
  \quad &\implies \quad 
  & \bv_n b_n &\text{ exists,}\quad 
  &&\text{and}\quad
  &\varphi(\,\bv_n b_n\,) &= \bv_n \varphi(b_n).
\end{alignat*}
We say $\varphi$ is \keyword{complete}
if $\varphi$ is both $\Pi$-complete and $\Sigma$-complete.
\end{dfn}
%
%                  THE LEBESGUE MEASURE IS A COMPLETE VALUATION
%
\begin{ex}
\label{E:lmeas-complete-val}
The Lebesgue measure~$\Lmu$
(see Example~\ref{E:lmeas-val})
is a complete valuation.
We must show that~$\Lmu$ is
both $\Pi$-complete and $\Sigma$-complete
(see Definition~\ref{D:complete-val}).

Let us prove~$\Lmu$ is $\Sigma$-complete.
Let $B_1 \subseteq B_2 \subseteq \dotsb$ in~$\LA$
be $\Lmu$-convergent.
We must prove that $\bv_n B_n$ exists in~$\LA$
and that 
$\Lmu(\bv_n B_n) = \bv_n \Lmu(B_n)$.

Note that $\bigcup_n B_n$
is Lebesgue measurable,
and that  $\bigcup_n B_n$ has (finite) Lebesgue measure~$\bv_n\Lmu(B_n)$.
Hence we have $\bigcup_n B_n \in \LA$, and,
\begin{equation}
\label{eq:E:lmeas-complete-val-1}
\textstyle
\Lmu(\,\bigcup_n B_n\,) \ =\  \bv_n \Lmu(B_n).
\end{equation}
So we are done if we prove that $\bigcup_n B_n = \bv_n B_n$.
Since $\bigcup_n B_n$ is the smallest subset of~$\R$
containing all~$B_n$
(i.e. $\bigcup_n B_n$ is the supremum of the~$B_n$ in~$\wp \R$),
$\bigcup_n B_n$ is also the smallest
subset \emph{of finite Lebesgue measure} containing all~$B_n$
(i.e. $\bigcup_n B_n$ is the supremum of the~$B_n$ in~$\LA$).
So
 $\bigcup_n B_n = \bv_n B_n$.
Hence $\Lmu$ is $\Sigma$-complete.

Using an easier reasoning one can prove that $\Lmu$
is $\Pi$-complete.
\end{ex}
%
%                  THE LEBESGUE INTEGRAL IS A COMPLETE VALUATION
%
\begin{ex}
\label{E:int-complete-val}
The Lebesgue integral~$\Lphi$
(see Example~\ref{E:int-val})
is a complete valuation.
We must show that~$\Lphi$ is
both $\Pi$-complete and $\Sigma$-complete
(see Definition~\ref{D:complete-val}).

Let us prove~$\Lphi$ is $\Sigma$-complete.
Let $f_1\leq f_2 \leq\dotsb$ in~$\LF$
be $\Lphi$-convergent.
We must prove that $\bv_n f_n$ exists in~$\LF$,
and that 
\begin{equation*}
\Lphi(\bv_n f_n) \ =\  \bv_n \Lphi(f_n).
\end{equation*}
Of course,
this follows immediately from
Levi's Monotone Convergence Theorem;
the supremum $\bv_n f_n$ in~$\LF$ is simply
the pointwise supremum
(which is the supremum of~$f_1\leq f_2 \leq\dotsb$ in $\E^\R$).
So we see $\Lphi$ is $\Sigma$-complete.

With a similar argument one can see that $\Lphi$ is $\Pi$-complete.
\end{ex}
%
%                  REMARK ON \E INSTEAD OF \R
%
\begin{rem}
\label{R:non-finite-functions}
Note that the restriction of~$\Lphi$ to~$\LF\cap \R^\R$
is not complete.

Indeed,
consider for instance
the following sequence.
\begin{equation*}
1\cdot\mathbf{1}_{\{0\}} 
\ \leq\   2\cdot\mathbf{1}_{\{0\}} 
\ \leq\   3\cdot\mathbf{1}_{\{0\}}
\ \leq\  \dotsb
\end{equation*}
It is~$\Lphi$-convergent
in~$\LF\cap \R^\R$,
but it has no supremum in~$\R^\R$.\\
On the other hand,
it does have a supremum in~$\LF$, 
namely~$+\infty\cdot \mathbf{1}_{\{0\}}$.

Because of the above observation,
we work with
the $\E$-valued Lebesgue integrable functions
instead of the $\R$-valued Lebesgue integrable functions.
\end{rem}
%
%                  THE SIMPLE MEASURABLE SETS ARE NOT COMPLETE
%
\begin{ex}
The valuation $\Smu$ 
(see Example~\ref{E:smeas-val}) is \emph{not} complete.

To see this,
we consider the sets~$A_1,A_2,\dotsc$ given by, for~$n\in\N$,
\begin{equation*}
A_n \ =\  \{1,\dotsc,n\}.
\end{equation*}
Then $A_n\in \SA$
and $\Smu(A_n) = 0$
for all~$n\in \N$.
So we see that
\begin{equation*}
A_1 \,\subseteq\, A_2 \,\subseteq\, \dotsb
\end{equation*}
is a $\Smu$-convergent sequence.
To prove that~$\Smu$ is not complete,
we  show that
the $\Smu$-convergent sequence~$A_1 \subseteq A_2 \subseteq \dotsb$
has no supremum in~$\SA$
(see Definition~\ref{D:complete-val}).

Suppose (towards a contradiction) that $A_1 \subseteq A_2 \subseteq\dotsb$
has a supremum~$B$ in~$\SA$.
Then in particular $A_n\subseteq B$ for all~$n\in \N$.
So we have
\begin{equation*}
\N \,=\, \textstyle{\bigcup_n A_n} \ \subseteq\  B.
\end{equation*}
Note that $B$ is the disjoint union of
elements from~$\mathcal{S}$ (see Example~\ref{E:smeas-val}).
Since all~$I\in \mathcal{S}$ are bounded, the set~$B$ is bounded.
That is, $B\subseteq [a,b]$ for some~$a,b\in \R$.

We now see that $\N\subseteq B\subseteq [a,b]$,
which is nonsense.
So $A_1 \subseteq A_2 \subseteq \dotsb$ has no supremum
in~$\SA$.
Hence $\Smu$ is not complete.
\end{ex}
%
%                  THE STEP FUNCTIONS ARE NOT COMPLETE
%
\begin{ex}
The valuation $\Sphi$ 
(see Example~\ref{E:sint-val}) is also \emph{not} complete.\\
We leave it to the reader to prove this fact.
\end{ex}

%
%                  PHI-CONVERGENCE FOR SIGMA-DEDEKIND COMPLETE VALUATIONS
%
\noindent
If $E=\R$,
or more generally,
if~$E$ is $\sigma$-Dedekind complete
(see Definition~\ref{D:sdc}),
then there is a nice description of $\varphi$-convergence,
see Proposition~\ref{P:phi-conv-dom}.
\begin{lem}
\label{L:phi-conv-dom}
Let $E$ be an ordered Abelian group.\\
Assume that~$E$ is $\sigma$-Dedekind complete
(see Definition~\ref{D:sdc}).\\
Let $L$ be a lattice,
and let $\varphi\colon L\ra E$ be a valuation.\\
Let $a_1 \geq a_2 \geq \dotsb$ be a sequence in~$L$.\\
Then $a_1 \geq a_2 \geq \dotsb$ 
is $\varphi$-convergent 
provided that~$a_1\geq a_2\geq \dotsb$ 
has a lower bound.
\end{lem}
\begin{proof}
Let $\ell\in L$ be a lower bound of~$a_1 \geq a_2 \geq \dotsb$,
that is, $\ell \leq a_n$ for all~$n\in\N$.
We must prove that~$a_1 \geq a_2 \geq \dotsb$ is
$\varphi$-convergent,
i.e., $\bw_n \varphi(a_n)$ exists.

Note that $\varphi(\ell)\leq \varphi(a_n)$ for all~$n\in \N$.
So $\varphi(a_1),\,\varphi(a_2),\, \dotsc$
has a lower bound.
But then $\bw_n \varphi(a_n)$ exists,
because~$E$ is $\sigma$-Dedekind complete
(see Remark~\ref{R:sdc}).
Hence $a_1 \geq a_2 \geq \dotsb$ is $\varphi$-convergent.
\end{proof}
\begin{prop}
\label{P:phi-conv-dom}
Let $E$ be an ordered Abelian group.\\
Assume that~$E$ is $\sigma$-Dedekind complete
(see Definition~\ref{D:sdc}).\\
Let $L$ be a lattice,
and let $\varphi\colon L\ra E$ be a complete valuation.\\
For a sequence $a_1 \geq a_2 \geq \dotsb$
in~$L$ the following are equivalent.
\begin{enumerate}
\item
\label{P:phi-conv-dom_i}
$a_1 \geq a_2 \geq \dotsb$
is $\varphi$-convergent.

\item
\label{P:phi-conv-dom_ii}
$a_1 \geq a_2 \geq \dotsb$
has a lower bound in~$L$.

\item
\label{P:phi-conv-dom_iii}
$a_1 \geq a_2 \geq \dotsb$
has an infimum, $\bw_n a_n$.
\end{enumerate}
\end{prop}
\begin{proof}
The implication 
``\ref{P:phi-conv-dom_i}$\ \Longleftarrow\ $\ref{P:phi-conv-dom_ii}''\ 
holds by Lemma~\ref{L:phi-conv-dom}.

\noindent
``\ref{P:phi-conv-dom_ii}$\ \Longleftarrow\ $\ref{P:phi-conv-dom_iii}''\ 
holds, because the infimum~$\bw_n a_n$ is a 
lower bound of $a_1 \geq a_2 \geq \dotsb$.

\noindent
``\ref{P:phi-conv-dom_iii}$\ \Longleftarrow\ $\ref{P:phi-conv-dom_i}''\ 
holds since $\varphi$ is complete
(see Definition~\ref{D:complete-val}).
\end{proof}
%
%                  EXAMPLE ON INFIMUM DOES NOT IMPLY PHI-CONVERGENCE
%
\noindent
The notion of $\varphi$-convergence
is less trivial in general
as the following example shows.
\begin{ex}
We will show that
the assumption 
that~$E$ is $\sigma$-Dedekind complete
in Proposition~\ref{P:phi-conv-dom}
is necessary for 
the implication 
``\ref{P:phi-conv-dom_iii}$\ \Longrightarrow\ $\ref{P:phi-conv-dom_i}''.\\
To this end, we extend the Lebesgue integral $\Lphi$
(see Example~\ref{E:int-complete-val})
to the set 
\begin{equation*}
\LF'\ \eqdf\  \LF\,\cup\,\{\,-\infty\cdot \mathbf{1}\,\}.
\end{equation*}
Note that~$\LF'$ is a sublattice of~$\E^\R$.
Let $\Lphi'\colon \LF'\rightarrow \Lex$
be the map,
where $\Lex$ is the \emph{lexicograpic plane} 
(see Example~\ref{E:oag}\ref{E:oag_lex}),
given by, for $f\in \LF'$,
\begin{equation*}
\Lphi'( f) \ =\ 
\begin{cases}
\quad (\,\rsub{0}{-1},\,\Lphi(f)\,)\,  \qquad
&\text{if $f\in\LF$}, \\
\quad (\,-1,\,0\,)  &\text{if $f=-\infty\cdot \mathbf{1}$}.
\end{cases}
\end{equation*}
Then $\Lphi'$ is a valuation.
In fact, 
$\Lphi'$ is a complete valuation
as the reader can verify 
using the following observation.
If $a_1 \geq a_2 \geq \dotsb$ from~$\LF'$
is $\Lphi'$-convergent, then:
\begin{enumerate}
\item
If $a_n \in \LF$ for all~$n\in\N$,
then $a_1 \geq a_2 \geq \dotsb$ is $\Lphi$-convergent.
\item
If $a_N = -\infty\cdot \mathbf{1}$
for some $N\in\N$,
then $a_n = -\infty\cdot \mathbf{1}$ for all~$n\geq N$.
\end{enumerate}
Now, 
consider the following sequence in $\LF'$.
\begin{equation*}
-1\cdot\mathbf{1}_{[0,1]} \ \geq\ 
-2 \cdot\mathbf{1}_{[0,1]} \ \geq \ 
-3 \cdot\mathbf{1}_{[0,1]}
\ \geq\ \dotsb.
\end{equation*}
This sequence has an infimum in~$\LF'$, 
namely $-\infty\cdot \mathbf{1}$.
Nevertheless,
the sequence is not $\Lphi'$-convergent(,
because $(0,-1) \,\geq\, (0,-2)\,\geq\, \dotsb$
has no infimum in~$\Lex$).
\end{ex}

%%%%%%%%%%%%%%%%%%%%%%%%%%%%%%%%%%%%%%%%%%%%%%%%%%%%%%%%%%%%%%%%%%%%%%%%%%
%%%%%%%%%%%%%%%%%%%%%%%%%%%%%%%%%%%%%%%%%%%%%%%%%%%%%%%%%%%%%%%%%%%%%%%%%%%
%
\subsection{$R$-completeness}
\label{SS:complete-val_R-completeness}
We now study a notion of completeness for ordered Abelian groups
called \emph{$R$-completeness}
that will be useful later on.

Let $\varphi\colon L \ra E$ be a valuation.
Let $a_1 \leq a_2 \leq \dotsb$
and $b_1 \leq b_2 \leq \dotsb$ be $\varphi$-convergent sequences in~$L$.
(see Definition~\ref{D:phi-conv}).

For the development of the theory,
it would be convenient if also 
\begin{equation}
\label{eq:SS-complete-val_R-completeness-1}
a_1 \vee b_1 \ \leq\  a_2 \vee b_2 \ \leq\  \dotsb
\qquad\text{is $\varphi$-convergent.}
\end{equation}
Unfortunately, 
this is not always the case (see Example~\ref{E:P:R-main}).
However,
if the space~$E$ is $\sigma$-Dedekind complete 
(see Appendix~\ref{S:ag}, Definition~\ref{D:sdc}),
for instance if~$E=\R$,
then one can prove that Statement~\eqref{eq:SS-complete-val_R-completeness-1}
holds.

In fact,
if we only assume that
$E$ is \emph{$R$-complete} (see Definition~\ref{D:R-complete}) ---
which is a weaker assumption than that~$E$ is Dedekind-complete ---
then we can still prove that 
that Statement~\eqref{eq:SS-complete-val_R-completeness-1}
holds (see Proposition~\ref{P:R-main}).
%
%                  R-completeness
%
\begin{dfn}
\label{D:R-complete}
Let $E$ be an ordered Abelian group.
Consider the following.
\begin{equation*}
\left[\quad 
\begin{minipage}{.7\columnwidth}
Let $x_1 \leq x_2 \leq \dotsb$
and $y_1 \leq y_2 \leq \dotsb$ be from~$E$
such that
\begin{equation*}
x_{n+1} - x_n \ \leq\ y_{n+1} - y_n\qquad \text{for all }n.
\end{equation*}
Then $\bv x_n $ exists whenever $\bv y_n$ exists.
\end{minipage}
\right.
\end{equation*}
If the above statement holds,
we say~$E$ is \keyword{$R$-complete}.
\end{dfn}
\begin{rem}
The name ``$R$-complete''
is due to Willem van Zuijlen~\cite{Zuijlen12}.
\end{rem}

\begin{exs}
\label{E:R-complete}
\begin{enumerate}
\item
\label{E:R-complete-R}
The ordered Abelian group $\R$ is $R$-complete.

\item
In fact, any $\sigma$-Dedekind complete 
ordered Abelian group $E$ is $R$-complete.\\
Indeed,
let $x_1 \leq x_2 \leq \dotsb$ 
and $y_1 \leq y_2 \leq \dotsb$
be from~$E$
such that
\begin{equation}
\label{eq:E:R:sdc}
x_{n+1} - x_n \ \leq\ y_{n+1}-y_n
\qquad\text{for all }n,
\end{equation}
and assume that $\bv_n y_n$ exists.
We must show that $\bv_n x_n$ exists.

Let $n\in\N$ be given.
By Statement~\eqref{eq:E:R:sdc}
we see that
\begin{equation*}
x_{n+1} - y_{n+1} \ \leq\ x_n - y_n.
\end{equation*}
So with induction on~$n$, we get $x_n - y_n \ \leq\ x_1 - y_1$.
Then
\begin{equation*}
x_n \ \leq\ 
(x_1 - y_1) \,+\, y_n  \ \leq\ 
(x_1 - y_1) \,+\, \bv_m y_m.
\end{equation*}
So we see that the sequence $x_1,x_2,\dotsc$ 
has an upper bound.\\
So $\bv_n x_n$ exists,
as~$E$ is $\sigma$-Dedekind complete.
Hence~$E$ is $R$-complete.

\item
The lexicographic plane~$\Lex$ (see Examples~\ref{E:oag}\ref{E:oag_lex})
is $R$-complete,
but~$\Lex$ is not $\sigma$-Dedekind complete
(see Examples~\ref{E:sdc}\ref{E:sdc_lex}).

\item
The ordered Abelian group~$\Q$ is \emph{not} $R$-complete.\\
To see this,
and pick $q_1 \leq q_2\leq \dotsb$ in~$\Q$
with
\begin{equation*}
q_{n+1} - q_{n} \ \leq\ 2^{-(n+1)}
\qquad\quad\text{and}\quad\qquad  \text{$\bv_n q_n = \sqrt{2}$ in~$\R$.}
\end{equation*}
Note that $q_1\leq q_2 \leq \dotsb$
has no supremum in~$\Q$.

Now, let $y_n \eqdf 1-2^{-n}$ for all~$n\in\N$.
Then $y_1 \leq y_2 \leq \dotsb$ has an supremum, namely~$1$,
and we have $y_{n+1} - y_{n} = 2^{-(n+1)}$. So we see that
\begin{equation*}
q_{n+1} - q_{n} \ \leq\ y_{n+1} - y_n
\qquad\quad(n\in\N).
\end{equation*}
If~$\Q$ were~$\R$-complete,
then the above implies 
 $q_1 \leq q_2 \leq \dotsb$
would have a supremum in~$\Q$,
which it does not.
Hence~$\Q$ is not $R$-complete.

\item
\label{E:R-complete-product}
Let $I$ be a set.
For each $i\in I$,
let $E_i$ be an $R$-complete
ordered Abelian group.
Then the product, $\prod_{i\in I} E_i$,
is $R$-complete.
\end{enumerate}
\end{exs}

\begin{rem}
\label{R:R-dual}
Let~$E$ be an ordered Abelian group.
Using the map $x\mapsto -x$,
one can easily verify
that $E$ is $R$-complete
if and only if the following statement holds.
\begin{equation*}
\left[\quad 
\begin{minipage}{.7\columnwidth}
Let $x_1 \geq x_2 \geq \dotsb$
and $y_1 \geq y_2 \geq \dotsb$ be from~$E$
such that
\begin{equation*}
x_{n} - x_{n+1}\ \leq\ y_{n} - y_{n+1}\qquad \text{for all }n.
\end{equation*}
Then $\bw x_n $ exists whenever  $\bw y_n$ exists.
\end{minipage}
\right.
\end{equation*}
\end{rem}

\begin{prop}
\label{P:R-main}
Let  $E$ be an ordered Abelian group which is $R$-complete.\\
Let~$L$ be a lattice, and let $\varphi\colon L \ra E$ be a valuation.
\begin{enumerate}
\item
\label{P:R-main-descending}
If  $a_1 \geq a_2 \geq \dotsb$,
$b_1 \geq b_2 \geq \dotsb$
are  $\varphi$-convergent
sequences from~$L$,
then
\begin{equation*}
a_1 \wedge b_1 \,\geq\, a_2 \wedge b_2 \,\geq\, \dotsb
\qquad\text{and}\qquad
a_1 \vee b_1 \,\geq\, a_2 \vee b_2 \,\geq\, \dotsb
\end{equation*}
are $\varphi$-convergent.

\item
\label{P:R-main-ascending}
If  $a_1 \leq a_2 \leq \dotsb$,
$b_1 \leq b_2 \leq \dotsb$
are  $\varphi$-convergent
sequences from~$L$,
then
\begin{equation*}
a_1 \wedge b_1 \,\leq\, a_2 \wedge b_2 \,\leq\, \dotsb
\qquad\text{and}\qquad
a_1 \vee b_1 \,\leq\, a_2 \vee b_2 \,\leq\, \dotsb
\end{equation*}
are $\varphi$-convergent.
\end{enumerate}
\end{prop}
\begin{proof}
\ref{P:R-main-descending}\ 
We prove that $a_1 \wedge b_1 \geq a_2 \wedge b_2 \geq\dotsb$
is $\varphi$-convergent.
For this we need to show that $\bw_n \varphi(a_n\wedge b_n)$ exists.
Note that since $\bw_n \varphi (a_n)$
and $\bw_n \varphi(b_n)$ exist,
we know that $\bw_n \ (\varphi(a_n) + \varphi(b_n))$
exists (by Lemma~\ref{L:addition-of-infima}).
So by $R$-completeness,
in order to show $\bw_n\varphi(a_n \wedge b_n)$ exists,
it suffices to prove that (see Remark~\ref{R:R-dual}),
\begin{equation*}
\varphi(a_{n}\wedge b_{n}) \,-\, \varphi(a_{n+1} \wedge b_{n+1}) 
\ \leq\ 
(\,\varphi(a_{n}) + \varphi(b_{n})\,) 
\,-\, (\,\varphi(a_{n+1}) + \varphi(b_{n+1})\,).
\end{equation*}
Phrased differently
using ``$d_\varphi$''
(see Definition~\ref{D:d}),
we need to prove that
\begin{equation*}
d_\varphi(a_{n}\wedge b_{n},\, a_{n+1} \wedge b_{n+1}) 
\ \leq\ 
d_\varphi(a_{n},a_{n+1}) + d_\varphi(b_{n},b_{n+1}).
\end{equation*}
This follows  from Lemma~\ref{L:wv-unif}.

The 
proof that~$a_1 \vee b_1 \geq a_2 \vee b_2 \geq \dotsb$
is $\varphi$-convergent is similar.

\ref{P:R-main-ascending}.  Again, similar.
\end{proof}
%
%                  COUNTEREXAMPLE ON PROP R-MAIn
%
\begin{ex}
\label{E:P:R-main}
We will prove that
the assumption in Proposition~\ref{P:R-main},
that~$E$ is $R$-complete, is necessary.

Let $\mathcal{A}$ be the ring of subsets (see Example~\ref{E:ring-val})
of~$\R$ generated by the non-empty closed intervals with \emph{rational}
 endpoints,
i.e., subsets of the form $[q,r]$ where $q,r\in \Q$ and $q\leq r$.
Then there is a unique positive and additive map $\mu\colon\mathcal{A} \ra\Q$
such that
\begin{equation*}
\mu(\,[q,r]\,) \ =\ r-q\qquad\quad\text{for all $q\leq r$ from $\Q$}.
\end{equation*}

Recall that~$\Q$ is not $R$-complete.
To prove that the conclusion of Proposition~\ref{P:R-main}
does not hold for~$E=\Q$, we will find $\mu$-convergent 
sequences $A_1 \subseteq A_2 \subseteq \dotsb$
and $B_1 \subseteq B_2 \subseteq\dotsb$ such that
$A_1 \cup B_1 \ \subseteq\ A_2 \cup B_2 \ \subseteq\ \dotsb$
is not $\mu$-convergent.

If we have done this,
we see that the assumption ``$E$ is $R$-complete''
is necessary.

Find rational numbers $ \dotsb \leq r_2 \leq r_1 < q_1 \leq q_2 \leq \dotsb$
such that, in~$\R$,
\begin{equation*}
\bv_n q_n \ = \ \sqrt2\qquad\text{and}\qquad
\bw_n r_n \ = \ \sqrt2-1.
\end{equation*}
Now, let us define $A_1 \subseteq A_2 \subseteq \dotsb$
and $B_1 \subseteq B_2 \subseteq \dotsb$ in~$\mathcal{A}$
by, for~$n\in\N$,
\begin{equation*}
A_n  \ = \ [0, r_1]
\qquad\text{and}\qquad
B_n \ = \ [r_n,q_n].
\end{equation*}
Then clearly $A_1 \subseteq A_2 \subseteq \dotsb$
is $\varphi$-convergent.
Note that $\mu(B_n) = q_n - r_n$. So, in~$\R$,
\begin{equation*}
\bv_n \mu(B_n) \ =\ \bv_n q_n - \bw_n r_n \ =\ 1.
\end{equation*}
Hence $\bv_n \mu(B_n)=1$ in~$\Q$.
So  $B_1 \subseteq B_2 \subseteq \dotsb$
is a $\mu$-convergent sequence.

However $A_n \cup B_n = [0,q_n]$,
and thus $\mu(A_n \cup B_n) = q_n$.
So we see that, in~$\R$,
\begin{equation*}
\bv_n \mu(\,A_n \cup B_n \,) \ =\ \bv_n q_n \ =\ \sqrt{2}.
\end{equation*}
So $\mu(A_1 \cup B_1) \ \leq\  \mu(A_2 \cup B_2) \ \leq\ \dotsb$
has no supremum in~$\Q$.

Hence $A_1 \cup B_1 \ \subseteq\ A_2 \cup B_2 \ \subseteq\ \dotsb$
is not $\mu$-convergent.
\end{ex}

%%%%%%%%%%%%%%%%%%%%%%%%%%%%%%%%%%%%%%%%%%%%%%%%%%%%%%%%%%%%%%%%%%%%%%%%%%%%%%%
%%%%%%%%%%%%%%%%%%%%%%%%%%%%%%%%%%%%%%%%%%%%%%%%%%%%%%%%%%%%%%%%%%%%%%%%%%%%%%%

\subsection{Convergence Theorems}
\label{SS:complete-val_convergence}
%
%                  CONVERGENCE OF SEQUENCES
%
The notion of a complete valuation
has been based on Levi's Monotone Convergence Theorem
(see Example~\ref{E:int-complete-val}).
In this subsection,
we prove variants of some of the other classical convergence theorems
of integration theory.
For example, 
Lebesgue's Dominated Convergence Theorem.
It states:
\begin{equation}
\label{eq:Lebesgue}
\left[\quad
\begin{minipage}{.7\columnwidth}
Let $f_1,\,f_2,\,\dotsc$ be a sequence in~$\LF$.

Assume $f_1(x),\,f_2(x),\,\dotsc$
converges for almost all~$x\in\R$.

Assume that $f_1,\,f_2,\,\dotsc$
 is dominated in the sense that
 $|f_n|\leq D$ for all~$n$
for some $D\in\LF$.

Then there is an $f\in\LF$
with $f_1(x),\,f_2(x),\,\dotsc$ converges to $f(x)$ for 
for almost all~$x\in\R$,
and $\Lphi(f) = \lim_n\Lphi(f_n)$.
\end{minipage}
\right.
\end{equation} 
The difficulty 
in the setting of  valuation systems
is not the proof of the theorem,
but its formulation.
For instance, it 
not clear how we should interpret 
\begin{equation*}
\text{``$f_1(x),\,f_2(x),\,\dotsb$
converges for almost all~$x$''}
\end{equation*}
when the objects $f_n$ are not necessarily functions,
but elements of a lattice~$V$.

%
%                  CONVERGENCE IN A LATTICE
%
Let us begin by generalising the notion of convergence in~$\R$
to any lattice~$L$.
Recall that a sequence $a_1,a_2,\dotsc$ in~$\R$
is convergent (in the usual sense) if and only if 
the \emph{limit inferior}, $\lim_{N} \inf_{n\geq N}\,a_n$,
and the \emph{limit superior}, $\lim_{N} \sup_{n\geq N}\,a_n$,
exist and are equal.
This leads us to the following definitions.
\begin{dfn}
\label{D:conv}
Let~$L$ be a lattice.
Let $a_1,\,a_2,\,\dotsc$ be a sequence in~$L$.
\begin{enumerate}
\item
We say $a_1,\,a_2,\,\dotsc$ 
is \keyword{upper convergent}
if the following exists.
\begin{equation*}
\ulim{n} a_n \ \eqdf\ \bw_N\bv_{n\geq N}\ a_N \vee \dotsb \vee a_n.
\end{equation*}
Similarly,
we say $a_1,\,a_2,\,\dotsc$ is \keyword{lower convergent}
if the following exists.
\begin{equation*}
\llim{n} a_n \ \eqdf\ \bv_N\bw_{n\geq N}\ a_N \wedge \dotsb \wedge a_n.
\end{equation*}

\item
We say $a_1,\,a_2,\,\dotsc$ is \keyword{convergent}
if it is both upper and lower convergent,
and in addition $\ulim{n}a_n = \llim{n}a_n$.
In that case,
we write $\lim_n a_n\eqdf\ulim{n} a_n$.

\item
Let $a\in L$ be given.
We say $a_1,\,a_2,\,\dotsc$ 
\keyword{converges to}~$a$
if  $a=\lim_n a_n$.
\end{enumerate}
\end{dfn}
%
%                  REMARK ON THE INEQUALITY OF LIMSUP AND LIMINF
%
\begin{rem}
\label{R:conv}
Let~$L$ be a lattice.
Let $a_1,\,a_2,\,\dotsc$ be a sequence in~$L$,
which is  upper convergent and lower convergent.
Then we have the following inequality.
\begin{equation*}
\llim{n} a_n \ \leq\ \ulim{n} a_n
\end{equation*}
Indeed,
this follows immediately from the observation
that, for every~$N\in\N$,
\begin{equation*}
\bv_{n\geq N}\ a_N \wedge \dotsb \wedge a_n
\ \leq\  a_N\ \leq\ 
\bw_{n\geq N}\ a_N \vee \dotsb \vee a_n.
\end{equation*}
\end{rem}
%
%                  EXAMPLES OF CONVERGENCE IN A LATTICE
%
\begin{exs}
\label{E:conv}
\begin{enumerate}
\item 
In~$\R$ we have:
A sequence $a_1,\,a_2,\,\dotsc$
is convergent in the sense of Definition~\ref{D:conv}
if and only if $a_1,\,a_2,\,\dotsc$
is convergent as usual.\\
Moreover, if $a_1,\,a_2.\,\dotsc$
is convergent, then $\lim_n a_n$ from Definition~\ref{D:conv}
is also the limit of $a_1,\,a_2,\,\dotsc$
in the usual sense.

\item
Similarly,
in~$\R^X$, where~$X$ is any set,
``convergent'' from Definition~\ref{D:conv}
coincides with the usual ``pointwise convergent''.
\end{enumerate}
\end{exs}

\begin{ex}
Let $X$ be a set.
Let $A_1,\,A_2,\,\dotsc$ be subsets of~$X$.
Then $A_1,\,A_2,\,\dotsc$
is upper and lower convergent in the lattice~$\wp X$, 
and we have, for~$x\in X$,
\begin{alignat*}{3}
x \in \ulim{n} A_n 
\quad&\iff\quad  \forall\,N\ \ \exists\, n\geq N \quad x\in A_n, \\
x \in \llim{n} A_n 
\quad&\iff\quad  \exists\,N\ \ \forall\, n\geq N \quad x\in A_n.
\end{alignat*}
So we see that~$A_1,\,A_2,\,\dotsc$ is \emph{not} convergent
iff there is an~$\tilde{x}\in X$ such that
\begin{equation*}
\forall\,N\ \ \exists\, n\geq N \quad \tilde{x}\in A_n
\qquad\text{and}\qquad
\forall\,N\ \ \exists\, n\geq N \quad \tilde{x}\notin A_n.
\end{equation*}
\end{ex}

\begin{ex}
\label{E:conv_leb}
For the lattice of Lebesgue integrable functions, $\LF$,
the notion of convergence from Def.~\ref{D:conv}
and the usual pointwise convergence 
do not coincide.

To see this,
consider the following sequence.
\begin{equation*}
\mathbf{1}_{[0,1]},\quad \mathbf{1}_{[1,2]},\quad \mathbf{1}_{[2,3]},\quad 
\dotsc
\end{equation*}
This sequence converges pointwise to~$\mathbf{0}$,
but it not convergent in the sense of Def.~\ref{D:conv}.
Indeed, the sequence is not even upper convergent
because 
\begin{equation*}
\mathbf{1}_{[0,1]} 
\ \leq\ \mathbf{1}_{[0,2]} 
\ \leq\ \mathbf{1}_{[0,3]}  \ \leq\ \dotsb
\end{equation*}
has no supremum in~$\LF$.

Fortunately,
the situation is better for \emph{dominated sequences}.\\
Let $f_1,f_2,\dotsc \in \LF$ and $f\in \LF$ be given.
Let $D\in \LF$ be given
such 
that $|f_n|\,\leq\,D$ for all~$n\in\N$.
(We say that $f_1,f_2,\dotsc$ is \emph{dominated} by~$D$.)\\
The reader can easily verify
the following statements
(cf. Example~\ref{E:int-complete-val}).
\begin{enumerate}
\item
\label{E:conv_leb_1}
The dominated sequence $f_1,f_2,\dotsc$
is upper convergent,
and
for all $x\in\R$,
\begin{equation*}
(\ulim{n} f_n)(x) \ =\  \ulim{n}(\,{f_n}(x)\,).
\end{equation*}

\item
\label{E:conv_leb_2}
The dominated sequence $f_1,f_2,\dotsc$
is lower convergent.
and  for all $x\in\R$,
\begin{equation*}
(\llim{n} f_n)(x) \ =\  \llim{n}(\,{f_n}(x)\,).
\end{equation*}

\item
\label{E:conv_leb_3}
The dominated sequence $f_1,f_2,\dotsc$
converges pointswise to~$f$
if and only if
$f_1,f_2,\dotsc$
converges to~$f$ 
in the sense of Definition~\ref{D:conv}.
\end{enumerate}
\end{ex}

\noindent
Dominated sequences
are useful when working with the Lebesgue integrable functions,
because~$\R$ is $\sigma$-Dedekind complete
(see Definition~\ref{D:sdc}).
However, it turns out that dominated sequences
are less useful in general.

Hence we have found a replacement 
for ``$f_1,f_2,\dotsc$ is dominated'',
namely,
\begin{equation*}
\text{``$f_1,f_2,\dotsc$ is \emph{upper and lower $\varphi$-convergent}''.}
\end{equation*}
%
%                  PHI-CONVERGENCE
%
\begin{dfn}
\label{D:seq-phi-conv}
Let $E$ be an ordered Abelian group.
Let~$L$ be a lattice.\\
Let $\varphi\colon L \ra E$ be a valuation.
Let $a_1,\, a_2\, \dotsc\in L$ be given.
\begin{enumerate}
\item
We say $a_1,\,a_2,\, \dotsc$
is \keyword{upper $\varphi$-convergent}
if the following exists.
\begin{equation*}
\pulim\varphi{n}a_n \ \eqdf\ 
\bw_N \bv_{n\geq N}\ \varphi(a_N\vee\dotsb\vee a_n)
\end{equation*}
Similarly,
we say $a_1,\,a_2,\,\dotsc$ is \keyword{lower $\varphi$-convergent} if
the following exists.
\begin{equation*}
\pllim\varphi{n}a_n \ \eqdf\ 
\bv_N \bw_{n\geq N}\ \varphi(a_N\wedge\dotsb\wedge a_n)
\end{equation*}

\item
We say $a_1,\,a_2,\,\dotsc$
is \keyword{$\varphi$-convergent}
if it is lower and upper $\varphi$-convergent,
and in addition $\pulim\varphi{n}a_n = \pllim\varphi{n}a_n$.
\end{enumerate}
\end{dfn}
%
%                  REMARK ON INEQUALITY BETWEEN UPPER AND LOWER PHI-LIM
%
\begin{rem}
\label{R:seq-phi-conv}
Let $E$ be an ordered Abelian group.
Let~$L$ be a lattice.\\
Let $\varphi\colon L \ra E$ be a valuation.
Let $a_1,\,a_2,\,\dotsc$ be a sequence in~$L$,
which is upper and lower $\varphi$-convergent.
We have the following inequality
(cf. Remark~\ref{R:conv}).
\begin{equation*}
\pllim\varphi{n}a_n \ \leq\ \pulim\varphi{n}a_n.
\end{equation*}
\end{rem}
%
%                  EXAMPLE OF PHI-CONV SEQUENCE
%
\begin{prop}
\label{P:phi-conv-dom-2}
Let $E$ be a $\sigma$-Dedekind complete
ordered Abelian group.\\
Let $L$ be a lattice,
and let $\varphi\colon L\ra E$ 
be a complete valuation.\\
Then for a sequence $a_1,a_2,\dotsc$ in~$L$
the following are equivalent.
\begin{enumerate}
\item
\label{P:phi-conv-dom-2_1}
$a_1,a_2,\dotsc$ is upper and lower $\varphi$-convergent.
\item
\label{P:phi-conv-dom-2_2}
$a_1,a_2,\dotsc$ has an upper and lower bound.
\end{enumerate}
\end{prop}
\begin{proof}
``\ref{P:phi-conv-dom-2_1}
$\Longrightarrow$
\ref{P:phi-conv-dom-2_2}''\ 
Assume that  $a_1,a_2,\dotsc$ is upper and lower $\varphi$-convergent.\\
We must find~$u,\ell\in L$ such that $\ell\leq a_n \leq u$
for all~$n\in\N$.

Since $a_1,a_2,\dotsc$ is  upper $\varphi$-convergent
(see Definition~\ref{D:seq-phi-conv}), 
we know that 
\begin{equation*}
\bv_n \ \varphi(a_1 \vee \dotsb \vee a_n)\qquad\text{exists.}
\end{equation*}
In other words, we know that
\begin{equation*}
a_1 \ \leq\ a_1 \vee a_2 \ \leq\ \dotsb \qquad\text{is $\varphi$-convergent.}
\end{equation*}
Since $\varphi$ is complete,
$u\eqdf \bv_n a_n$ exists in~$L$.
Note that $a_n \leq u$ for all~$n\in \N$.

By a similar reasoning, but using the fact that
$a_1,a_2,\dotsc$ is lower $\varphi$-convergent,
we can find an $\ell\in L$ 
such that $\ell \leq a_n$ for all~$n\in \N$.
\vspace{.3em}

\noindent
``\ref{P:phi-conv-dom-2_1}
$\Longleftarrow$
\ref{P:phi-conv-dom-2_2}''\ 
Let $\ell,u\in L$ 
be such that  $\ell \leq a_n \leq u$
for all~$n\in\N$.

We prove that $a_1,a_2,\dotsc$
is upper $\varphi$-convergent.
For this,
we must show that the following exists
(see Definition~\ref{D:seq-phi-conv}).
\begin{equation}
\label{eq:E:seq-phi-conv-1}
\bw_N \bv_{n\geq N}\ \varphi(a_N\vee\dotsb\vee a_n)
\end{equation}
Let $N\in\N$ and $n\geq N$ be given.
Note that we have 
\begin{equation*}
\ell\ \leq\ a_N \vee \dotsb\vee a_n \ \leq\ u.
\end{equation*}
Since $\varphi$ is order preserving, this gives us
\begin{equation*}
\varphi(\ell)\ \leq\ \varphi(a_N \vee \dotsb\vee a_n) \ \leq\ \varphi(u).
\end{equation*}
Since $E$ is $\sigma$-Dedekind complete
it follows that Expression~\eqref{eq:E:seq-phi-conv-1} exists.

We have proven that $a_1,a_2,\dotsc$ 
is upper $\varphi$-convergent.
With a similar reasoning one can prove that
$a_1,a_2,\dotsc$
is lower $\varphi$-convergent.
\end{proof}

\begin{ex}
\label{E:seq-phi-conv}
Let $f_1,f_2,\dotsc$ be  Lebesgue integrable functions
(see Example~\ref{E:int-val}). \\
Then by Proposition~\ref{P:phi-conv-dom-2}
the following statement holds.
\begin{equation*}
\left[\quad
\begin{minipage}{.7\columnwidth}
 The sequence $f_1,f_2,\dotsc$
is  upper and lower $\Lphi$-convergent.
\begin{center}
$\Longleftrightarrow$
\end{center}
There is a Lebesgue integrable~$D$ with $|f_n| \leq D$ for all~$n$.
\end{minipage}
\right.
\end{equation*}
\end{ex}
\vspace{.3em}

We can now prove a generalisation of the  Lemma of Fatou.
%
%                  LEMMA OF FATOU
%
\begin{lem}[Fatou]
\label{L:fatou}
Let $E$ be an ordered Abelian group.\\
Let $L$ be a lattice,
and let $\varphi\colon L \ra E$ be a 
complete valuation.\\
Let $a_1,a_2,\dotsc$ be an upper $\varphi$-convergent
sequence in~$L$
(see Definition~\ref{D:seq-phi-conv}) \\
Then $a_1,a_2,\dotsc$ is upper convergent
(see Definition~\ref{D:conv}),
and we have
\begin{equation*}
\varphi(\ulim{n}a_n) \ =\ 
\pulim\varphi{n} a_n.
\end{equation*}
Moreover,
if $E$ is a lattice, and
if $\ulim{n}\varphi(a_n)$ exists
(see Definition~\ref{D:conv}), then
\begin{equation*}
\pulim\varphi{n}a_n
\ \geq\  
\ulim{n}\varphi(a_n) .
\end{equation*}
\end{lem}
\begin{proof}
Let $a_1,\,a_2,\,\dotsc$ be an upper $\varphi$-convergent sequence.
We prove that $a_1,a_2,\dotsc$ is upper convergent
(see Definition~\ref{D:conv}),
and that $\varphi(\ulim{n}a_n) = \pulim\varphi{n}a_n$.

Let $N\in\N$ be given.
Note that $\bv_{N\geq n} \,\varphi(a_N\vee\dotsb\vee a_n)$
exists because the sequence~$a_1,\,a_2,\,\dotsc$ is upper $\varphi$-convergent
(see Definition~\ref{D:seq-phi-conv}).
So the sequence
\begin{equation*}
a_N \,\leq\ a_N \vee a_{N+1} 
    \ \leq\quad a_N\vee a_{N+1} \vee a_{N+2} 
    \quad \leq\qquad \dotsb
\end{equation*}
is $\varphi$-convergent (in the sense of
Definition~\ref{D:phi-conv}).
For brevity,
let us write
\begin{equation*}
\overline{a}_{N}^n \ \eqdf \ a_N\vee\dotsb\vee a_{N+n}.
\end{equation*}
Since~$\varphi$ is complete,
and $\overline{a}_N^0 \leq \overline{a}_N^1 \leq\dotsb$
is $\varphi$-convergent,
we get $\bv_n\,\overline{a}_N^n$ exists,
and 
\begin{equation*}
\varphi(\overline{a}_N) 
\ =\ 
 \bv_{n}\ \varphi( \overline{a}_N^n),
\end{equation*}
where 
$\overline{a}_N \eqdf \bv_n \,\overline{a}_N^n$.
Note that $\overline{a}_1 \geq \overline{a}_2 \geq\dotsb$
is $\varphi$-convergent,
because
\begin{equation*}
\pulim\varphi{n} a_n
\ =\ 
\bw_N\bv_n\, \varphi(\overline{a}^n_N)
\end{equation*}
exists since $a_1,\,a_2,\,\dotsc$ is upper $\varphi$-convergent.
Since $\varphi$ is complete,
this implies that
\begin{equation*}
\bw_n \overline{a}_N\text{ exists}
\qquad\text{and}\qquad \varphi(\bw_n \overline{a}_N)
\,=\,
\bw_n\varphi(\overline{a}_N).
\end{equation*}
Now,
note that  we have the following equality.
\begin{equation*}
\bw_N \overline{a}_N \ =\ 
\bw_N \bv_{n\geq N} \, a_n
\end{equation*}
So we see that $a_1,a_2,\dotsc$
is upper $\varphi$-convergent and that
\begin{equation*}
\varphi(\ulim{n}a_n)
\ =\ 
\bw_N\varphi(\overline{a}_N)
\ =\ 
\bw_N \bv_n \varphi(\overline{a}^n_N)
\ =\ 
\pulim\varphi{n} a_n.
\end{equation*}
We have proven the first part of the lemma.

Assume $E$ is a lattice and $\ulim{n}\varphi(a_n)$ exists
(see Definition~\ref{D:conv}).
To prove the remainder of the theorem,
we need to show that 
$\pulim\varphi{n}a_n \geq \ulim{n}\varphi(a_n)$.
That is,
\begin{equation*}
\bw_N \bv_{n\geq N} \ \varphi(a_N \vee \dotsb \vee a_n)
\ \geq \ 
\bw_N \bv_{n\geq N} \ \varphi(a_N) \vee \dotsb \vee \varphi(a_n).
\end{equation*}
This is easy.  It follows immediately
from the fact that
\begin{equation*}
\varphi(a_N\vee \dotsb\vee a_n)
\ \geq\ \varphi(a_N)\vee \dotsb \vee \varphi(a_n)
\end{equation*}
for all~$N\in\N$ and $n\geq N$.
\end{proof}
%
%                  APPROX-CONVERGENCE
%
\noindent
Let us now think about ``almost everywhere convergent''. 
\begin{dfn}
\label{D:approx-conv}
Let $E$ be an ordered Abelian group.
Let $L$ be a lattice.\\
Let $\varphi\colon L \ra E$ be a valuation.
Let $a_1,a_2,\dotsc \in L$ and $a\in L$ be given.
\begin{enumerate}
\item 
If  $a_1,a_2,\dotsc$
is upper and lower convergent, and,
with  $\approx$ as in Def.~\ref{D:approx},
\begin{equation*}
\llim{n} a_n \ \approx\ \ulim{n}a_n,
\end{equation*}
then we say that
$a_1,a_2,\dotsc$
is \keyword{$\approx$-convergent}.

\item
We say that $a_1,a_2,\dotsc$ $\approx$-\keyword{converges} to~$a$
when $\llim{n}a_n \approx a \approx \ulim{n}a_n$.
\end{enumerate}
\end{dfn}
%
%                  EXAMPLE OF APPROX CONVERGENCE
%
\begin{ex}
\label{E:approx-conv}
Unfortunately,
in the lattice of Lebesgue integrable functions, $\LF$,
the notion of $\approx$-convergence
does not coincide with convergence almost everywhere,
as can be seen using a similar argument
as before (see~Example~\ref{E:conv_leb}).

Again, the situation is better for dominated sequences.\\
Let $f_1,f_2,\dotsc \in\LF$ and $f\in\LF$.
Assume $f_1,f_2,\dotsc$ is dominated by some~$D\in\LF$,
that is, $|f_n|\leq D$ for all~$n\in\N$.
Then the following statements hold.
\begin{enumerate}
\item
\label{E:approx-conv_i}
The dominated sequence $f_1,f_2,\dotsc$  $\approx$-converges to~$f$
if and only if\\
$f_1(x),\,f_2(x),\,\dotsc$
converges to~$f(x)$ for almost all~$x\in\R$.

\item
\label{E:approx-conv_ii}
The dominated sequence
$f_1,f_2,\dotsc$ is $\approx$-convergent
if and only if\\
$f_1(x),\,f_2(x),\,\dotsc$
converges for almost all~$x\in\R$.
\end{enumerate}

We will prove implication
``$\Longleftarrow$''
of 
\ref{E:approx-conv_i},
and leave the rest to the reader.\\
We must show that 
$f_1,f_2,\dotsc$ is upper and lower convergent,
and that
\begin{equation}
\label{E:approx-conv-1}
\llim{n}f_n \ \approx \ f\ \approx\  \ulim{n}f_n.
\end{equation}
Since $f_1(x),\,f_2(x),\,\dotsc$ 
converges to~$f(x)$ for almost all~$x\in \R$,
we know that:
\begin{equation}
\label{E:approx-conv-2}
\llim{n}(\,f_n(x)\,) \ =\ f(x)\ =\ \ulim{n} (\,f_n (x)\,)
\end{equation}
for almost all~$x\in \R$.
So we see that the function given by $x\mapsto \llim{n}(\,f_n(x)\,)$
is equal almost everywhere to the Lebesgue integrable function~$f$,
and hence it is Lebesgue integrable itself.
By Example~\ref{E:conv_leb}\ref{E:conv_leb_2}
it follows that $f_1,f_2,\dotsc$ is lower convergent,
and that $(\llim{n}f_n)(x) = \llim{n}(\,f_n(x)\,)$
for all~$x\in\R$.
By a similar argument,
we see that $f_1,f_2,\dotsc$ is upper convergent as well,
and that $(\ulim{n}f_n)(x) = \ulim{n}(\,f_n(x)\,)$
for all~$x\in \R$.
Hence we get by Equation~\eqref{E:approx-conv-2},
for almost all~$x\in \R$,
\begin{equation*}
(\llim{n}f_n)(x) \ =\ f(x)\ =\ (\ulim{n} f_n) (x).
\end{equation*}
This proves Equation~\eqref{E:approx-conv-1}
(see Example~\ref{E:eq-int}).
\vspace{.3em}
\end{ex}

\noindent
Let us relate $\varphi$-convergence
and $\approx$-convergence.
%
%                  LEMMA ON APPROX-CONV VS PHI-CONV
%
\begin{lem}
\label{L:approx-phi-conv}
Let $E$ be an ordered Abelian group.\\
Let $L$ be a lattice,
and let $\varphi\colon L \ra E$ be a 
complete valuation.\\
Let $a_1,a_2,\dotsc$ be an upper and lower $\varphi$-convergent
sequence in~$L$
(see Definition~\ref{D:seq-phi-conv}). \\
Then the sequence $a_1,a_2,\dotsc$
is upper and lower convergent (see Definition~\ref{D:conv}),\\
and the following statements are equivalent.
\begin{enumerate}
\item 
\label{L:approx-phi-conv_i}
$a_1,a_2,\dotsc$ is $\approx$-convergent
(see Definition~\ref{D:approx-conv}).

\item 
\label{L:approx-phi-conv_ii}
$a_1,a_2,\dotsc$ is $\varphi$-convergent
(see Definition~\ref{D:seq-phi-conv}).
\end{enumerate}
Moreover, if 
either~\ref{L:approx-phi-conv_i}
or~\ref{L:approx-phi-conv_ii} holds, we have
\begin{equation}
\label{eq:L:approx-phi-conv}
\varphi(\llim{n}a_n) \,=\, \varphi(\ulim{n} a_n) 
\ =\ 
\plim\varphi{n} a_n.
\end{equation}
\end{lem}
\begin{proof}
Let $a_1,a_2,\dotsc$ be an upper and lower
$\varphi$-convergent sequence in~$L$.\\
By Lemma~\ref{L:fatou}
and its dual,
$a_1,a_2,\dotsc$
is both upper and lower convergent, and
\begin{equation}
\label{eq:L:approx-phi-conv_fatou}
\varphi(\llim{n} a_n) \ = \ \pllim\varphi{n}a_n,
\qquad\text{and}\qquad
\varphi(\ulim{n} a_n) \ =\ \pulim\varphi{n}a_n.
\end{equation}
By Definition~\ref{D:seq-phi-conv}
and Definition~\ref{D:approx-conv}
we see that
\begin{alignat*}{3}
\text{$a_1,a_2,\dotsc$ is $\varphi$-convergent}
\quad&\iff\quad&
\pllim\varphi{n}a_n \ &=\ \pulim\varphi{n}a_n,\\
\text{$a_1,a_2,\dotsc$ is $\approx$-convergent}
\quad&\iff\quad&
\varphi(\, \llim{n}a_n\,) \ &=\ \varphi(\,\ulim{n}a_n\,).
\end{alignat*}
Hence Equation~\eqref{eq:L:approx-phi-conv_fatou}
implies that statements~\ref{L:approx-phi-conv_i}
and~\ref{L:approx-phi-conv_ii} are equivalent.

Now, assume that~\ref{L:approx-phi-conv_i}
(or~\ref{L:approx-phi-conv_ii}) holds.
We must show  that Statement~\eqref{eq:L:approx-phi-conv}
holds.
This follows from Statement~\eqref{eq:L:approx-phi-conv_fatou}
since $\plim\varphi{n}a_n \,=\,\pllim\varphi{n}a_n \,=\, \pulim\varphi{n}a_n$.
\end{proof}

%
%                  LEBESGUE'S DOMINATED CONVERGENCE THEOREM
%
\noindent
We now prove a generalisation 
of Lebesgue's Dominated Convergence Theorem.
\begin{thm}[Lebesgue]
\label{T:lebesgue}
Let $E$ be a \emph{lattice} ordered Abelian group.\\
Let $L$ be a lattice,
and let $\varphi\colon L \ra E$ be a 
complete valuation.\\
Let $a_1,a_2,\dotsc$ be an upper and lower $\varphi$-convergent
sequence in~$L$
(see Def.~\ref{D:seq-phi-conv}). \\
Assume that $a_1,a_2,\dotsc$ is $\approx$-convergent
(see Def.~\ref{D:approx-conv}). \\
Assume that $\llim{n}\varphi(a_n)$ and $\ulim{n}\varphi(a_n)$ exist
(see Def.~\ref{D:conv}).\\
Then 
the sequence $\varphi(a_1),\,\varphi(a_2),\,\dotsc$
converges (see Def.~\ref{D:conv}), and we have
\begin{equation}
\label{eq:T:lebesgue-1}
{\lim}_n \,\varphi(a_n)
\ =\ 
\varphi(\,\llim{n} a_n\,)
\ =\ 
\varphi(\,\ulim{n} a_n\,).
\end{equation}
\end{thm}
\begin{proof}
Let us first prove that the sequence 
$\varphi(a_1),\,\varphi(a_2),\,\dotsc$
is convergent.\\
By Lemma~\eqref{L:approx-phi-conv},
we see that
 $a_1,a_2,\dotsc$
is $\varphi$-convergent (see Definition~\ref{D:seq-phi-conv}),
and that
\begin{equation}
\label{eq:T:lebesgue-2}
\varphi(\,\llim{n}a_n\,) \ =\  \varphi(\,\ulim{n} a_n\,) 
\ =\ 
\plim\varphi{n} a_n.
\end{equation}
By Lemma~\ref{L:fatou} and its dual, we see that
\begin{equation*}
\pllim\varphi{n}a_n \,\leq\,
\llim{n}\varphi(a_n) \,\leq\,
\ulim{n}\varphi(a_n) \,\leq\,
\pulim\varphi{n}a_n.
\end{equation*}
But $\pllim\varphi{n}a_n = \pulim\varphi{n}a_n$,
since $a_1,a_2,\dotsc$ is $\varphi$-convergent.
So we get 
\begin{equation}
\label{eq:T:lebesgue-3}
\pllim\varphi{n}a_n \,=\,
\llim{n}\varphi(a_n) \,=\,
\ulim{n}\varphi(a_n) \,=\,
\pulim\varphi{n}a_n.
\end{equation}
In particular,
$\varphi(a_1),\,\varphi(a_2),\,\dotsc$
is convergent
(see Definition~\ref{D:conv}).

It remains to be shown that 
Statement~\eqref{eq:T:lebesgue-1} holds.\\
To do this,
combine Statement~\eqref{eq:T:lebesgue-2}
and Statement~\eqref{eq:T:lebesgue-3}.
\end{proof}

%
%                  LEBESGUE'S THEOREM FOR SIGMA-DEDEKIND COMPLETE OAGs
%
\noindent
If we assume that~$E$ is $\sigma$-Dedekind complete 
we get a  more familiar statement.
\begin{thm}
\label{T:lebesgue_sdc}
Let $E$ be a \emph{lattice} ordered Abelian group.\\
Assume that~$E$ is $\sigma$-Dedekind complete (see Def.~\ref{D:sdc}).\\
Let $L$ be a lattice,
and let $\varphi\colon L \ra E$ be a 
complete valuation.\\
Let $a_1,a_2,\dotsc$
sequence in~$L$
which has an upper and lower bound.\\
Assume that $a_1,a_2,\dotsc$ is $\approx$-convergent
(see Def.~\ref{D:approx-conv}). \\
Then 
the sequence $\varphi(a_1),\,\varphi(a_2),\,\dotsc$
converges (see Def.~\ref{D:conv}), and we have
\begin{equation*}
{\lim}_n \,\varphi(a_n)
\ =\ 
\varphi(\,\llim{n} a_n\,)
\ =\ 
\varphi(\,\ulim{n} a_n\,).
\end{equation*}
\end{thm}
\begin{proof}
We want to apply Theorem~\ref{T:lebesgue}.
For this, we must prove that $a_1,a_2,\dotsc$
is upper and lower $\varphi$-convergent,
and that $\llim{n}\varphi(a_n)$
and $\ulim{n}\varphi(a_n)$ exist.

Note that $a_1,a_2,\dotsc$
is upper and lower $\varphi$-convergent
since $a_1,a_2,\dotsc$
has an upper and lower bound
(see Proposition~\ref{P:phi-conv-dom-2}).

Let $u,\ell\in L$
be such that $\ell\leq a_n \leq u$
for all~$n\in \N$.
Then we have, for all~$n\in \N$,
\begin{equation*}
\varphi(\ell) \ \leq\ \varphi(a_n) \ \leq\ \varphi(u).
\end{equation*}
Using this,
and the fact that~$E$ is $\sigma$-Dedekind complete,
it is not so hard to see
that $\llim{n}\varphi(a_n)$
and $\ulim{n}\varphi(a_n)$ exist
 (cf.~Proposition~\ref{P:phi-conv-dom-2}).

Now we can apply Theorem~\ref{T:lebesgue},
and we are done.
\end{proof}

%
%                  EXAMPLE ON LEBESGUE'S THEOREM
%
\begin{ex}
\label{E:lebesgue}
If we apply Theorem~\ref{T:lebesgue_sdc}
to the Lebesgue integral $\Lphi$,
we get the classical form of 
Lebesgue's Dominated Convergence Theorem
(see Statement~\eqref{eq:Lebesgue}).

Indeed,
let $f_1,\,f_2,\,\dotsc$ be a sequence of Lebesgue integrable functions.
Assume there is an Lebesgue integrable function~$D$ such that
$|f_n|\leq D$ for all~$n\in\N$,
and assume that $f_1(x),\,f_2(x),\,\dotsc$ converges
for almost all~$x\in\R$.

We must prove that there is a Lebesgue integrable~$f$
such that $f_1(x),\,f_2(x),\,\dotsc$
converges to~$f(x)$ for almost all~$x\in\R$,
and $\Lphi(f) = \lim_n \Lphi(f_n)$.

Since 
$f_1(x),\,f_2(x),\,\dotsc$ 
converges for almost all~$x\in \R$,
we know that the sequence $f_1,\,f_2,\,\dotsc$
is $\approx$-convergent
(see Example~\ref{E:approx-conv}\ref{E:approx-conv_ii}).
Now, define
\begin{equation*}
f\ \eqdf\  \llim{n}{f_n}.
\end{equation*}
It is easy to see that,
$f_1,\,f_2,\,\dotsc$
$\approx$-converges to~$f$
(see Definition~\ref{D:approx-conv}).\\
That is,
$f_1(x),\,f_2(x),\,\dotsc$
converges to $f(x)$ for almost all~$x\in\R$
(see Example~\ref{E:approx-conv}\ref{E:approx-conv_i}).\\
By Theorem~\ref{T:lebesgue_sdc}
we see that $\varphi(a_1),\,\varphi(a_2),\,\dotsc$
converges, and that
\begin{equation*}
{\lim}_n \,\Lphi(f_n)
\ =\ 
\Lphi(\,\llim{n} f_n\,)
\ =\ 
\Lphi(f).
\end{equation*}
So we are done.
\end{ex}

\noindent
There are many variants of the classical convergence theorems of
integration. \\
For instance,
a variant on Levi's Monotone Convergence Theorem
is the following.
\begin{equation*}
\left[\quad
\begin{minipage}{.7\columnwidth}
Let $f_1,\,f_2,\,\dotsc$ be
 Lebesgue integrable functions.\\
Assume that $\bv_n \Lphi(f_n)$ exists.\\
Assume that 
for every $n\in\N$,
\begin{equation*}
f_n(x) \,\leq\, f_{n+1}(x) \qquad\text{ for almost all~$x\in\R$.}
\end{equation*}
Then $\bv_n f_n$,
the pointwise supremum of $f_1,\,f_2,\,\dotsc$,\\
is Lebesgue integrable
and 
\begin{equation*}
\Lphi(\bv_n f_n) \,=\, \bv_n \Lphi(f_n).
\end{equation*}
\end{minipage}
\right.
\end{equation*}
Note that if 
we want to prove the above statement
it will be useful to know that
\begin{equation*}
\qvphi{\Lphi}\colon \ \qvL{\LF} \longrightarrow \R
\end{equation*}
from Proposition~\ref{P:quotient-lattice} 
is complete. We will prove this in Proposition~\ref{P:quot-complete}.

Of course,
if we apply
Lemma~\ref{L:fatou}
and Theorem~\ref{T:lebesgue} 
to
$\qvphi{\Lphi}$,
we obtain variants of the Lemma of Fatou
and the Dominated Convergence Theorem of Lebesgue,
respectively. We leave this to the reader.

%
%                 QUOTIENT OF COMPLETE IS COMPLETE
%
\begin{prop}
\label{P:quot-complete}
Let~$L$ be a lattice. Let $E$ be an ordered Abelian group.\\
Let $\varphi\colon L\ra E$ be a complete valuation.
Then the valuation
\begin{equation*}
\qvphi\varphi\colon\  \qvL{L} \longrightarrow  E
\end{equation*}
from Proposition~\ref{P:quotient-lattice}
is a complete valuation.
\end{prop}
\begin{proof}
We leave this to the reader.
\end{proof}

\noindent
There is a small gap that needs to filled before
we continue with another topic.\\
Let $\varphi\colon L \ra E$ be a valuation.
We have defined what it means
for a sequence $a_1,\,a_2,\,\dotsc$ in~$L$ to be
$\varphi$-convergent
(see Definition~\ref{D:seq-phi-conv}),
but we have not yet given the meaning of
``$a_1,\,a_2,\,\dotsc$ converges \emph{to} $a$''.
We will do this in Definition~\ref{D:seq-phi-conv-to}.
%
%                  DEFINITION OF PHI-CONV TO AN ELEMENT
%
\begin{dfn}
\label{D:seq-phi-conv-to}
Let $L$ be a lattice.
Let $E$ be an ordered Abelian group.\\
Let $\varphi\colon L\ra E$ be a valuation.
Let $a_1,\,a_2,\,\dotsc$ be a sequence in~$L$.\\
Let $a\in L$ be given.
We say $a_1,\,a_2,\,\dotsc$
\keyword{$\varphi$-converges} to~$a$
provided that 
\begin{equation*}
a_1,\,a,\,a_2,\,a,\,\dotsc\qquad\text{is $\varphi$-convergent.}
\end{equation*}
\end{dfn}
\begin{rem}
Let $\varphi\colon L \ra E$ be a valuation.
While Definition~\ref{D:seq-phi-conv-to} is certainly reasonable,
it is also quite silly, 
and so one wonders if there is a more direct description
of when a sequence $\varphi$-converges to an element~$a\in L$.
If we assume~$\varphi$ is complete,
then there is a slightly 
better description (see Proposition~\ref{P:seq-phi-conv-to}).

In Section~\ref{S:fub} we will study a notion of convergence
(see Definition~\ref{D:weak-phi-conv})
which was intended to be a more aesthetically pleasing
definition of $\varphi$-convergence,
but which turns out to be strictly weaker 
than $\varphi$-convergence
(see Example~\ref{E:weak-phi-conv}).
\end{rem}
%
%                  EXAMPLE PHI-CONVERGENCE
%
\begin{prop}
\label{P:seq-phi-conv-to}
Let~$E$ be an ordered Abelian group.\\
Let~$L$ be a lattice, and $\varphi\colon L \ra E$ be a complete valuation.\\
Let $a_1,\,a_2,\,\dotsc$ be a $\varphi$-convergent sequence in~$L$,
and let $a\in L$.\\
Then $a_1,\,a_2,\,\dotsc$
$\varphi$-converges to~$a$
if and only if (see Definition~\ref{D:eq})
\begin{equation*}
a\ \approx\ \ulim{n}a_n.
\end{equation*}
(Recall that
$a_1,\,a_2,\,\dotsc$
is upper convergent
(see Definition~\ref{D:conv})
by Lemma~\ref{L:fatou}.)
\end{prop}
\begin{proof}
We leave this to the reader.
\end{proof}
 }
\clearpage
{ \section{Valuation Systems}
\label{S:valuation-systems}
\label{S:system}
\noindent
Note that the Lebesgue measure $\Lmu$
is a complete valuation (see Example~\ref{E:lmeas-complete-val}),
that extends the relatively simple valuation~$\Smu$
(see Example~\ref{E:smeas-val}).

We would like to consider~$\Lmu$ to be \emph{a completion} of~$\Smu$.
What should this mean?
The following definition seems obvious
when one thinks about valuations.
\begin{equation*}
\left[ \quad
\begin{minipage}{.7\columnwidth}
Let $E$ be an ordered Abelian group.\\
Let $L$ and~$K$ be lattices.\\
Let $\psi\colon K \ra E$
and $\varphi\colon L\ra E$
be valuations.\\
We say $\psi$ is \keyword{a completion} of~$\varphi$
provided that\\
$L$ is a sublattice of~$K$,
~$\psi$ extends~$\varphi$,
and  $\psi$ is complete.
\end{minipage}
\right.
\end{equation*}
However,
in the more concrete setting 
of measure theory this broad definition
of completion is not that useful.
After all,
if we are given a completion~$\psi\colon K\ra \R$ of~$\Smu$,
then we only know that~$K$ is a sublattice of~$\SA$,
while we would prefer~$K$ to be a sublattice of sets,
or resemble it.

To mend this problem
we might try to prove 
that any completion of~$\Smu$ is 
essentially a completion on a lattice of subsets.
Of course, the meaning of the previous statement
is not clear.
We suspect that if one gives it an exact meaning,
the statement will be either false or trivial.
So we will not follow this direction.

Instead,
we consider a different notion of completion
that involves the 
the surrounding lattice, $\wp \R$.
More precisely,
we will see that $\Lmu$ is a completion of~$\Smu$
\emph{relative to}~$\wp \R$,
which means 
that
$\Lmu$ extends~$\Smu$
and that $\Lmu$ is complete \emph{relative to~$\wp\R$}
(see Example~\ref{E:complete-lmeas}).
This naturally leads to the study of the following objects.
\begin{equation*}
\vsSA \qquad\qquad \vsLA.
\end{equation*}
That is, we are interested in objects of the following shape.
\begin{equation*}
\vs{V}{L}{\varphi}{E},
\end{equation*}
where $\varphi\colon L \ra E$ is a valuation,
and where $V$ is a lattice such that~$L$ is a sublattice of~$V$.
We call such objects \emph{valuation systems}
(see Definition~\ref{D:system}).

The drawback of this approach is that it requires
quite a bit of bookkeeping,
and so this section
is filled with definitions and examples,
but there is little theory.
We hope the reader will bear with us;
we are confident the reader will be rewarded
for his/her patience
in the next sections.

Since this section
 is already  administrative in nature,
we take this chance
to put 
 some additional restraints on
the notion of valuation system
which turns out to be useful later on
(see Remark~\ref{R:use-of-sigma-and-R}).
Given a valuation system,
$\vs{V}{L}{\varphi}{E}$,
we  require that~$E$ is $R$-complete
(see Def.~\ref{D:R-complete}),
and that~$V$ is $\sigma$-distributive
(see Def.~\ref{D:sigma-distributive}).

Before we give a formal definition
of ``valuation system'' 
in Subsection~\ref{SS:valuation-systems},
and define ``complete valuation system''
in Subsection~\ref{SS:complete-valuation-systems},
we consider $\sigma$-distributive lattices
in Subsection~\ref{SS:sigma-distributive}.

We end the section
with ``convex valuation systems''
 in Subsection~\ref{SS:convex}.
%%%%%%%%%%%%%%%%%%%%%%%%%%%%%%%%%%%%%%%%%%%%%%%%%%%%%%%%%%%%%%%%%%%%%%%%%%%%%%%
%
%                  SIGMA DISTRIBUTIVITY 
%
\subsection{$\sigma$-Distributivity}
\label{SS:sigma-distributive}
%
%                  ADDITIONAL RESTRICIONS ON THE SYSTE
% 
\begin{dfn}
\label{D:sigma-distributive}
Let~$V$ be a lattice.
We say~$V$ is
\keyword{$\sigma$-distributive}
provided that
\begin{enumerate}
\item
$V$ is \keyword{$\sigma$-complete}, i.e.,
for every sequence $c_1,\,c_2,\,\dotsc$ in~$V$
we have 
\begin{equation*}
\text{ $\bw_n c_n$ exists\qquad and\qquad $\bv_n c_n$ exists, }
\end{equation*}
\item
and for every  $a\in V$ and $c_1,\,c_2,\,\dotsc\in V$,
we have,
\begin{equation*}
a \vee \bw_n c_n \,=\, \bw_n\  a\vee c_n
\qquad\text{and}\qquad
a \wedge \bv_n c_n \,=\, \bv_n\  a\wedge c_n.
\end{equation*}
\end{enumerate}
\end{dfn}
\begin{exs}
\label{E:sigma-distributive}
\begin{enumerate}
\item
\label{E:sigma-distributive-sets}
Let $X$ be a set. Then $\wp(X)$ is $\sigma$-distributive.
Indeed,  
\begin{equation*}
\textstyle{
A \cup \bigcap_n C_n \,=\, \bigcap_n A \cup C_n
\qquad
A\cap \bigcup_n C_n \,=\, \bigcup_n A \cap C_n}
\end{equation*}
for all $A,\, C_1,C_2,\dotsc \subseteq X$.
\item
Let $C$ be totally ordered
and $\sigma$-complete. Then $C$ is $\sigma$-distributive.

Indeed,
let $a,\,c_1,c_2,\dotsc \in C$ be given.
We need to prove that $a \vee \bw_n c_n$
is the supremum of~$a\vee c_1,\,a\vee c_2,\,\dotsc$.
To this end note that 
\begin{equation*}
b \leq d_1 \vee d_2 \quad\iff\quad 
b\leq d_1\quad\text{or}\quad b\leq d_2
\qquad\quad(b,d_i\in C).
\end{equation*}
(To see this, note that $d_1 \vee d_2 = \max\{d_1,d_2\}$.)
Now, for $\ell \in C$, we have
\begin{alignat*}{3}
\forall n [\ \ell \leq a \vee c_n \ ]
\quad&\iff\quad
\ell \leq a
    \quad\text{or}\quad
    \forall n[\ \ell \leq c_n\ ] \\
\quad&\iff\quad
\ell \leq a
    \quad\text{or}\quad
    \ell \leq \bw_n c_n \\
\quad&\iff\quad
\ell \leq a\vee \bw_n c_n.
\end{alignat*}
So we see 
that $a\vee\bw_n c_n$ is the greatest 
lower bound of~$a\vee c_1,\,a\vee c_2\,\dotsc$.

With the same argument,
one can prove that $a \wedge \bv_n c_n = \bv a \wedge c_n$
for all $a,\,c_1,c_2,\dotsc \in C$
such that $\bv_n c_n$ exists.
Hence $C$ is $\sigma$-distributive.

\item
The lattice of the real numbers~$\R$ is a chain
and so~$\R$ is $\sigma$-distributive
\emph{if}~$\R$ would be $\sigma$-complete.
However,
$\R$ is not $\sigma$-complete.
Indeed,
a sequence $c_1,c_2,\dotsc$ in~$\R$ has a supremum
if and only if it is bounded from above,
i.e., there is an~$a\in \R$ such that $c_n \leq a$
for all~$n$.
Similarly,
a sequence
$c_1,c_2,\dotsc\in \R$ has an infimum
if and only if it is bounded from below.

\item
Let $\E$ be the lattice of the extended real numbers.
Then~$\E$ is a chain and clearly $\sigma$-complete.
Hence $\E$ is $\sigma$-distributive.

\item
\label{E:sigma-distributive-product}
Let $I$ be a set,
and for each~$i\in I$,
let $L_i$ be a $\sigma$-distributive lattice.\\
Then the product $L\eqdf \prod_{i\in I} L_i$
is $\sigma$-distributive.

\item
\label{E:sigma-distributive-functions}
Let $X$ be a set.
Then lattice $\EX$ of functions from~$X$ to~$\E$
is $\sigma$-distributive.
\end{enumerate}
\end{exs}

\subsection{Valuation Systems}
\label{SS:valuation-systems}
%
%                  SYSTEMS
%

\begin{dfn}
\label{D:system}
We say $\vs{V}{L}{\varphi}{E}$
 is a \keyword{valuation system}
provided that
\begin{enumerate}
\item \label{D:simple-system-1}
$V$ is a $\sigma$-distributive lattce 
(see Definition~\ref{D:sigma-distributive});
\item \label{D:simple-system-2}
$L$ is a sublattice of~$V$;
\item \label{D:simple-system-3}
$E$ is an ordered Abelian group,
which is $R$-complete (see Definition~\ref{D:R-complete});
\item \label{D:simple-system-4}
$\varphi\colon L\ra E$ is a valuation.
\end{enumerate}
\end{dfn}
%
%                  RING AS SIMPLE VALUATION SYSTEM
%
\begin{ex}
\label{E:ring-system}
Let $E$ be an $R$-complete ordered Abelian group.
Let~$X$ be a set, 
$\mathcal{A}$ a ring of subsets of~$X$,
and $\mu\colon \mathcal{A}\ra E$
a positive and additive map
(see Example~\ref{E:ring-val}).

Then we have the following  valuation system.
\begin{equation*}
\vs{\wp X}{\mathcal{A}}\mu{E}
\end{equation*}
Indeed, $\wp X$ is lattice with 
$\bw_n A_n = \bigcap_n A_n$
and $\bv_n A_n = \bigcup_n A_n$ for all $A_i \in \wp X$,
the set
$\mathcal{A}$ is a sublattice of~$\wp X$
by definition,
$\wp X$ is $\sigma$-distributive
(see Examples~\ref{E:sigma-distributive}\ref{E:sigma-distributive-sets})
and we have already seen 
that $\mu\colon \mathcal{A}\ra E$ is a valuation
(in Example~\ref{E:ring-val}).

In particular,
we have the following valuation systems.
\begin{equation*}
\vsLA \qquad\qquad \vsSA
\end{equation*}
See Example~\ref{E:lmeas-val} and  Example~\ref{E:smeas-val}.
\end{ex}

%
%                  RIESZ SPACE OF FUNCTIONS AS SIMPLE VALUATION SYSTEM
%
\begin{ex}
\label{E:riesz-function-space-simple-system}
Let $E$ be an $R$-complete ordered Abelian group.
Let $F$ be a Riesz space of functions on a set~$X$
(see Example~\ref{E:val-riesz-space-of-functions}),
and let $\varphi\colon F\ra E$ be a positive and linear map.
Then we have the following  valuation system.
\begin{equation*}
\vs{[-\infty,\infty]^X}{F}\varphi{E}
\end{equation*}
Indeed, 
$[-\infty,\infty]^X$ is a $\sigma$-distributive
lattice 
(see Examples~\ref{E:sigma-distributive}\ref{E:sigma-distributive-functions}).
Further, $F$ is a sublattice of~$\R^X$
which is in turn a sublattice of $[-\infty,\infty]^X$,
and we already know that
$\varphi$ is a valuation (see Example~\ref{E:val-riesz-space-of-functions}).

In particular, since~$\R$ is $R$-complete,
we have the following valuation systems
\begin{equation*}
\vs{\E^X}{(\LF\cap\R^\R)}{\Lphi}{\R},
\qquad\qquad
\vsSF,
\end{equation*}
see Example~\ref{E:int-val} and Example~\ref{E:sint-val}.
Recall that $\LF\cap \R^\R$ is a Riesz space of functions on~$X$,
while~$\LF$ is not.
Of course,
we also have the following valuation sytem.
\begin{equation*}
\vsLF
\end{equation*}
\end{ex}

\begin{ex}
Let $I=\{1,2\}$.
For each~$i\in I$,
let $\vs{V_i}{L_i}{\varphi_i}{E_i}$
be a  valuation system.
Then we have the following  valuation system
(see Example~\ref{E:val-product}).
\begin{equation*}
\vs{V_1\times V_2}{L_1 \times L_2}{\varphi_1 \times \varphi_2}{E_1 \times E_2}
\end{equation*}
Indeed,
the lattice
$V_1\times V_2$ is $\sigma$-distributive
(see Examples~\ref{E:sigma-distributive}\ref{E:sigma-distributive-product}),
and the ordered Abelian group $E_1\times E_2$ is $R$-complete
(see Examples~\ref{E:R-complete}\ref{E:R-complete-product}).
We call this system
the \emph{product} of $\vs{V_1}{L_1}{\varphi_1}{E_1}$
and $\vs{V_2}{L_2}{\varphi_2}{E_2}$.
Of course,
one can similarly define the product of 
an $I$-indexed family of valuation systems
where~$I$ is any set.
\end{ex}

%
%                  NOTATION CONCERNING SUPREMA AND INFIMA IN L
%
\begin{nt}
\label{N:V-inf-sup}
Let $\vs{V}{L}\varphi{E}$ be a  valuation system.
Let $a_1, a_2, \dotsc$ be from~$L$.
Then $a_1, a_2,\dotsc$ has a supremum
in~$V$ and might have a supremum in~$L$.
We ignore the latter:
\textbf{with $\bv_n a_n$
we always mean the supremum of~$a_1, a_2,\dotsc $ in~$V$}.\\
Similarly, \textbf{with $\bw_n a_n$
we always mean the infimum of~$a_1, a_2,\dotsc $ in~$V$}.
\end{nt}
%%%%%%%%%%%%%%%%%%%%%%%%%%%%%%%%%%%%%%%%%%%%%%%%%%%%%%%%%%%%%%%%%%%%%%%%%%%%%%%
%
%                  COMPLETENESS
%
%
\subsection{Complete Valuation Systems}
\label{SS:complete-valuation-systems}
%
%                  COMPLETE SYSTEMS
%
\begin{dfn}
\label{D:system-complete}
Let $\vs{V}{L}\varphi{E}$ be a valuation system.
\begin{enumerate}
\item 
We say $\vs{V}{L}\varphi{E}$
is \keyword{$\Pi$-complete},
or  $\varphi$ is \keyword{$\Pi$-complete relative to} $V$,\\
or even $\varphi$ is  \keyword{$\Pi$-complete} (if no confusion should 
arise with
Definition~\ref{D:complete-val}),\\
if for every $\varphi$-convergent
$a_1\geq a_2 \geq \dotsb$ 
(see Definition~\ref{D:phi-conv}), we have,
\begin{equation*}
   \bw_n a_n \in L\quad 
  \text{and}\quad
  \varphi(\,\bw_n a_n\,) \ =\  \bw_n \varphi(a_n).
\end{equation*}
Here, $\bw_n a_n$
is the infimum of $a_1 \geq a_2 \geq \dotsb$ in~$V$
(see Notation~\ref{N:V-inf-sup}).

\item
Similarly,
we say $\vs{V}{L}\varphi{E}$
is \keyword{$\Sigma$-complete}, etc.,\\
provided that for every $\varphi$-convergent sequence
$b_1\leq b_2 \leq \dotsb$ we have
\begin{equation*}
   \bv_n b_n \in L\quad 
  \text{and}\quad
  \varphi(\,\bv_n b_n\,) \ =\  \bv_n \varphi(b_n).
\end{equation*}

\item
We say $\vs{V}{L}\varphi{E}$
is \keyword{complete}, etc.,\\
provided that 
$\vs{V}{L}\varphi{E}$
is both $\Pi$-complete and $\Sigma$-complete.
\end{enumerate}
\end{dfn}
%
%                  REMARK ON COMPLETE VALUATION SYSTEMS
%
\begin{rem}
Let $\vs{V}{L}\varphi{E}$ be a complete valuation
system
(see Definition~\ref{D:system-complete}).
Then the valuation $\varphi$
is also complete in the sense
of Definition~\ref{D:complete-val}.\\
We leave it to the reader to verify this.
\end{rem}
%
%                  THE LEBESGUE INTEGRAL IS COMPLETE
%
\begin{ex}
\label{E:complete-lint}
The Lebesgue integral
gives us the   valuation system
\begin{equation*}
\vsLF
\end{equation*}
(see Example~\ref{E:int-val} and
Example~\ref{E:riesz-function-space-simple-system});
we will prove that this system is complete.

Let $f_1 \leq f_2 \leq \dotsb$
be a $\Lphi$-convergent sequence (in $\LF$).
We must prove that 
\begin{equation*}
\bv_n f_n \in\LF
\qquad\text{and}\qquad\Lphi(\,\bv_n f_n\,) \ =\  \bv_n \Lphi (f_n).
\end{equation*}
This follows immediately from Levi's Monotone Convergence Theorem.

Similarly,
if $g_1 \geq g_2 \geq \dotsb$
is a $\Lphi$-convergent sequence,
then
\begin{equation*}
\bw_n g_n \in\LF
\qquad\text{and}\qquad\Lphi(\,\bw_n g_n\,) \ =\  \bw_n \Lphi (g_n).
\end{equation*}
So the  valuation system $\vsLF$ is complete.

Note that we have given a similar argument earlier 
(see Example~\ref{E:int-complete-val}) to prove
that the valuation $\Lphi$
is complete in the sense of Definition~\ref{D:complete-val}.
\end{ex}
%
%                  THE LEBESGUE MEASURE IS COMPLETE
%
\begin{ex}
\label{E:complete-lmeas}
The Lebesgue measure
gives us the   valuation system
\begin{equation*}
\vsLA
\end{equation*}
(see Example~\ref{E:lmeas-val} and
Example~\ref{E:ring-system});
one can prove that this system is complete.

We leave this to the reader (cf.~Example~\ref{E:lmeas-complete-val}).
\end{ex}
\begin{cnv}
The complete valuation systems play a far more important
role in the remainder of this thesis
than the complete valuations of Definition~\ref{D:complete-val}.
So:\\
\textbf{Whenever we later on write ``$\varphi$ is complete'' 
we mean that  $\varphi$ is complete relative to~$V$},
where~$V$ should be clear from the context.
\end{cnv}

\subsection{Convex Valuation Systems}
\label{SS:convex}
$\,$\\
We have remarked
that the Lebesgue measure is
 complete (see Example~\ref{E:complete-lmeas}).
It should be noted that ``the Lebesgue measure is complete''
has a different meaning in the literature,
namely
that any subset of a Lebesgue neglegible set is neglegible itself.
We call this \emph{convexity} 
and we will briefly discuss this notion
in this subsection.
\begin{dfn}
\label{D:convex}
Let $\vs{V}{L}\varphi{E}$ be a valuation system.\\
We say 
that $\vs{V}{L}\varphi{E}$ is \keyword{convex},
if the following statement holds.
\begin{equation*}
\left[\quad
\begin{minipage}{.7\columnwidth}
Let~$a\leq b$ from~$L$
with $\varphi(a) = \varphi(b)$
be given. Then
\begin{equation*}
a \,\leq\,z\,\leq\, b
\qquad\implies\qquad z\in L,
\end{equation*}
where $z\in V$.
\end{minipage}
\right.
\end{equation*}
\end{dfn}
%
%                  EXAMPLES OF CONVEX VALUATION SYSTEMS
%
\begin{exs}
\begin{enumerate}
\item
The Lebesgue measure $\vsLA$ is convex.

\item
The Lebesgue integral $\vsLF$ is convex.
\end{enumerate}
\end{exs}
%
%                  BOREL MEASURE IS NOT CONVEX
%
\begin{ex}
\label{E:not-convex}
We give a serious example of a non-convex valuation system.

Let $\mathcal{B}$ be the set  of all \emph{Borel} subsets of~$\R$,
that is, those subsets of~$\R$
that can be formed by countable intersection
and countable union starting from the open subsets of~$\R$.
Every Borel subset of~$\R$ is Lebesgue measurable.

Let $\BA$ be the set of Borel subsets of~$\R$ with finite Lebesgue measure.
Note that we have $\BA\subseteq \LA$.
Let $\Bmu\colon \BA\ra \R$
be the restriction of~$\Lmu$ to~$\BA$.

Recall that $\vsLA$ is a complete and convex valuation system.
It is not hard to see that the valuation system~$\vsBA$
is complete as well.
However, we will prove that $\vsBA$ is \emph{not} convex.

To this end,
we will find a negligible Borel set~$A$
and a subset~$\Delta$ of~$A$
such that~$\Delta$ is not a Borel set.
This is sufficient to prove that
$\vsBA$ is not convex.
Indeed,
if $\vsBA$ were convex,
then~$\Delta$ would be Borel,
because
\begin{equation*}
\varnothing\subseteq \Delta\subseteq A
\qquad\text{and}\qquad
\varnothing,A\in\BA
\qquad\text{and}\qquad
\Bmu(\varnothing) =0= \Bmu(A).
\end{equation*}
To find such~$A$ and~$\Delta$,
we need the following fact
(see Corollary~\ref{C:bhier-embed}).
\begin{equation}
\label{eq:not-convex-fact}
\left[\quad
\begin{minipage}{.7\columnwidth}
There is a set~$\mathcal{C}$,
a negligible Borel set~$A\subset\R$,
and maps 
\begin{equation*}
G\colon \mathcal{C} \longrightarrow \mathcal{B}
\qquad\text{and}\qquad
p\colon \mathcal{C} \longrightarrow A
\end{equation*}
such that~$G$ is surjective, and $p$ is injective.
\end{minipage}
\right.
\end{equation}
Now,
let $\Delta$ be the subset of~$p(\N^\N)$ given by, for all $f\in\N^\N$,
\begin{alignat*}{3}
p(f)&\in\Delta& \quad&\iff&\quad
p(f)\,&\notin\, G(f).
\shortintertext{%
We prove that~$\Delta$ is not a Borel set.
Suppose that $\Delta$ is a Borel set
in order to derive a contradiction.
Since~$G$ is surjective, we have
 $\Delta\equiv G(f_0)$ for some~$f_0\in \N^\N$.
Then}
p(f_0)\,&\in\, \Delta& \quad&\iff&\quad p(f_0)\,&\notin\,G(f_0)=\Delta.
\end{alignat*}
Which is absurd.
Thus~$\Delta$ is not Borel set.\\
So we see that~$\Delta$ is a subset
of a negligible Borel set, $A$,
while~$\Delta$ is not Borel.
Hence $\vsBA$ is not convex.
\end{ex}
%
%                  CONVEXIFICATION
%
\begin{prop}
\label{P:convex-completion}
\label{P:convexification}
Let $\vs{V}{L}\varphi{E}$ be a valuation system.
\begin{enumerate}
\item
Let $L^\bullet$ be the subset of~$V$ given by,
for all $z\in V$,
\begin{equation*}
z\in L^\bullet \quad\iff\quad
\exists\, a,b\in L\ [
\quad a\leq z\leq b \quad\wedge\quad \varphi(a) = \varphi(b)\quad].
\end{equation*}
Then $L$ is a sublattice of~$L^\bullet$, which is a sublattice of~$V$.

\item
There is a unique order preserving map $\varphi^\bullet\colon L^\bullet\ra E$
which extends~$\varphi$.\\
Moreover, $\varphi^\bullet$
is a valuation,
and  $\vs{V}{L^\bullet}{\varphi^\bullet}{E}$ is convex.
\end{enumerate}
\end{prop}
\begin{proof}
We leave this to the reader.
\end{proof}
\begin{dfn}
\label{D:convexification}
The \keyword{convexification}
of a valuation system $\vs{V}{L}\varphi{E}$ 
is the valuation system
 $\vs{V}{L^\bullet}{\varphi^\bullet}{E}$ 
described in Proposition~\ref{P:convexification}.
\end{dfn}
%
%                  CONVEXIFICATION OF A COMPLETE SYSTEM IS COMPLETE
%
\begin{prop}
\label{P:convexification_versus_pi-completion}
Let $\vs{V}{L}\varphi{E}$ be a $\Pi$-complete valuation system.\\
Then 
the valuation system 
$\vs{V}{L^\bullet}{\varphi^\bullet}{E}$ is $\Pi$-complete as well.
\end{prop}
\begin{proof}
Let $a_1 \geq a_2 \geq \dotsb$
be a $\varphi^\bullet$-convergent sequence in~$L^\bullet$.
To prove that the valuation system
$\vs{V}{L^\bullet}{\varphi^\bullet}{E}$ is $\Pi$-complete,
we must show that $\bw_n a_n \in L^\bullet$,
and
\begin{equation*}
\varphi^\bullet(\,\bw_n a_n\,) \ =\  \bw_n \varphi^\bullet(a_n).
\end{equation*}
Let $n\in \N$ be given.
Note that since $a_n \in L^\bullet$
there are $\ell_n,u_n \in L$
such that 
\begin{equation*}
\ell_n \,\leq\, a_n \,\leq\, u_n
\qquad\text{and}\qquad
\varphi(\ell_n) \,=\, \varphi(u_n).
\end{equation*}
Define $\ell\eqdf \bw_n \ell_n$
and $u\eqdf \bw_n u_n$.
Then we have $\ell \leq \bw_n a_n \leq u$.
So to prove that $\bw_n a_n \in L^\bullet$,
it suffices to show that $\ell,u\in L$,
and $\varphi(\ell)=\varphi(u)$.

The trick is
to consider
$\ell_1' \geq \ell_2 '\geq \dotsb$
and $u_1'\geq u_2' \geq\dotsb$
given by, for $n\in\N$,
\begin{equation*}
\ell_n' \,\eqdf\, \ell_1 \wedge \dotsb \wedge \ell_n
\qquad\text{and}\qquad
u_n'\eqdf u_1 \wedge \dotsb \wedge u_n.
\end{equation*}
Note that
$\ell= \bw_n \ell_n '$ and $u= \bw_n u_n '$.
Let $n\in \N$ be given.
We claim that
\begin{equation}
\label{eq:convex-completion-1}
\varphi(\ell_n') \ =\  \varphi^\bullet(a_n)\ =\ \varphi(u_n').
\end{equation}
Indeed, since $\varphi(\ell_n) = \varphi(u_n)$
and $\ell_n\leq u_n$,
we have $\ell_n \approx u_n$
(see Definition~\ref{D:eq}).
Then 
by  Proposition~\ref{P:eq}\ref{P:eq-3}
and induction, we see $\ell_n' \approx u_n'$.
Hence  $\varphi(\ell_n') = \varphi(u_n')$.
Now, as
$\ell_n'\leq a_n \leq u_n'$,
and $\varphi^\bullet$ is order preserving,
we get
 $\varphi(\ell_n') \leq \varphi^\bullet(a_n) \leq \varphi(u_n')$.
Hence we easily see that Statement~\eqref{eq:convex-completion-1} holds.

Since $a_1 \geq a_2 \geq \dotsb$
is $\varphi^\bullet$-convergent,
we know that $\bw_n \varphi^\bullet (a_n)$ exists.
Further,
\begin{equation*}
\bw_n \varphi(\ell_n') \ =\ 
\bw_n \varphi^\bullet(a_n) \ =\ 
\bw_n \varphi(u_n')
\end{equation*}
by St.~\eqref{eq:convex-completion-1}.
So we see that
$\ell_1' \geq \ell_2' \geq \dotsb$
and $u_1'\geq u_2'\geq \dotsb$
are $\varphi$-convergent.
Because $\varphi$ is $\Pi$-complete
relative to~$V$
we see that $\bw_n \ell_n'=\ell\in L$
and $u\in L$,
and 
\begin{equation}
\label{eq:convex-completion-2}
\varphi(\ell) \ = \ 
\bw_n \varphi(\ell_n')\ =\ 
\bw_n \varphi^\bullet(a_n)\ = \ 
\bw_n \varphi(u_n')\ = \ 
\varphi(u).
\end{equation}
So,
since $\ell\leq \bw_n a_n \leq u$,
Statement~\eqref{eq:convex-completion-2} implies that
$\bw_n a_n \in L^\bullet$,
and 
\begin{equation*}
\varphi^\bullet(\,\bw_n a_n\,)
\ = \ 
\bw_n \varphi^\bullet(a_n).
\end{equation*}
Hence $\varphi^\bullet$ is $\Pi$-complete relative to~$V$.
\end{proof}
\begin{prop}
\label{P:convexification_versus_completion}
Let $\vs{V}{L}\varphi{E}$ be a complete valuation system.\\
Then 
the valuation system 
$\vs{V}{L^\bullet}{\varphi^\bullet}{E}$ is complete as well.
\end{prop}
\begin{proof}
Apply Proposition~\ref{P:convexification_versus_pi-completion}
and its dual.
\end{proof}
 }
\clearpage
{ \section{The Completion}
\label{S:completion}
\subsection{Introduction}
\noindent
Valuation systems
were introduced in
 Section~\ref{S:valuation-systems}
to give meaning to the phrase 
``the Lebesgue integral $\Lphi$ is a completion of~$\Sphi$'',
namely,
\begin{equation*}
\text{$\vsLF$ is complete}\qquad\text{and}\qquad
\text{$\Lphi$ extends~$\Sphi$}.
\end{equation*}
In this section we replace~$\Sphi$
and~$\smash{\E^\R}$
by any valuation system $\vs{V}{L}\varphi{E}$
and study when, so to say, $\varphi$ has 
\emph{a completion with respect to~$V$},
i.e., when there is a sublattice~$C$ of~$V$
and a valuation~$\psi\colon C\ra E$
such that 
\begin{equation*}
\vs{V}{C}\psi{E}\text{ is complete}
\qquad\text{and}\qquad
\text{$\psi$ extends~$\varphi$}.
\end{equation*}
There is not always a completion (see Example~\ref{E:no-completion}).
However,
if~$\varphi$ has a completion with respect to~$V$,
then~$\varphi$ also has a smallest\footnote{``Smallest''
with respect to the ordering on partial functions 
given by $f\leq g$ iff $g$ extends~$f$.} %
completion
with respect to~$V$,
\begin{equation*}
\vs{V}{\ol{L}}{\ol\varphi}{E}.
\end{equation*}
We describe $\ol\varphi$
in Subsection~\ref{SS:hierarchy},
and call~$\ol\varphi$ \emph{\underline{the} completion} of~$\varphi$
with respect to~$V$.

\subsubsection*{$\Pi$-Completion}
If we replace ``complete'' by ``$\Pi$-complete''
in the above discussion
the situation is much easier
(see Definition~\ref{D:system-complete}).
Indeed,
we will see  in Lemma~\ref{L:Pi-continuity}
that~$\varphi$ has a $\Pi$-completion with respect
to~$V$ iff
\begin{equation*}
\label{eq:Pi-cont}
\left[ \quad
\begin{minipage}{.7\columnwidth}
for every $b\in L$ and $\varphi$-convergent 
 $a_1 \geq a_2 \geq \dotsb$
we have
\begin{equation*}
\bw_n a_n \leq b
\quad\implies\quad
\bw_n \varphi (a_n) \leq \varphi(b).
\end{equation*}
\end{minipage}
\right.
\end{equation*}
Similar to before, if there is a $\Pi$-completion
of~$\varphi$ with respect to~$V$,
then there is also a smallest $\Pi$-completion
with respect to~$V$,
which will be denoted by
\begin{equation*}
\vs{V}{\Pi L}{\Pi \varphi}{E}.
\end{equation*}
We will study~$\Pi\varphi$ in Subsection~\ref{S:Pi-extension}.

\subsubsection*{$\Sigma$-Completion}
We can replace ``complete'' by ``$\Sigma$-complete''
as well.\\
If it exists,
the smallest $\Sigma$-completion of~$\varphi$
with respect to~$V$ is denoted by
\begin{equation*}
\vs{V}{\Sigma L}{\Sigma \varphi}{E}.
\end{equation*}
Since~$\Pi \varphi$ and~$\Sigma\varphi$
are very similar,
we will only study~$\Pi\varphi$.
All the results that we obtain about~$\Pi\varphi$
and all the definitions  for~$\Pi\varphi$
can be easily translated
to results about~$\Sigma\varphi$
and definitions for~$\Sigma\varphi$,
respectively.
We leave this to the reader.

\subsubsection*{Hierarchy of Extensions}
Recall that~$\varphi$ is complete with respect to~$V$
if and only if~$\varphi$
is both $\Pi$-complete
and $\Sigma$-complete with respect to~$V$
(see Definition~\ref{D:system-complete}).

So if we want to find a completion of~$\varphi$
with respect to~$V$ it is natural to try and see if $\Sigma\Pi\varphi$
is complete 
(when it exists) with respect to~$V$.
Unfortunately,
while $\Pi\varphi$ is $\Pi$-complete
with respect to~$V$,
the valuation
$\Sigma\Pi\varphi$ need not be $\Pi$-complete.
Nevertheless,
we can continue to apply ``$\Sigma$'' and ``$\Pi$''
whenever possible,
and we will see in Subsection~\ref{SS:hierarchy-abstract}
that this leads to a `hierarchy' of the following shape.
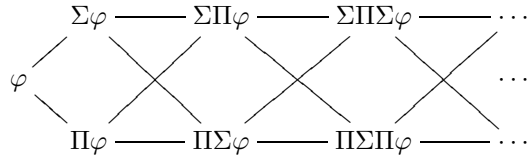
\begin{equation*}
\xymatrix @=10pt {
& \Sigma\varphi \ar @{-} [rr] \ar @{-} [rrdd]
&& \Sigma \Pi \varphi  \ar @{-} [rr] \ar @{-} [rrdd]
&& \Sigma \Pi \Sigma\varphi  \ar @{-} [rr] \ar @{-} [rrdd]
&& \dotsb
\\  
\varphi \ar @{-} [ru] \ar @{-} [rd] 
&&&&&&&\dotsb\\
& \Pi\varphi \ar @{-} [rr] \ar @{-} [rruu]
&& \Pi \Sigma \varphi \ar @{-} [rr]\ar @{-} [rruu]
&& \Pi \Sigma \Pi \varphi \ar @{-} [rr]\ar @{-} [rruu]
&& \dotsb
}
\end{equation*}
It will become clear,
that if the making of this hierarchy is hindered, e.g.,
if $\Pi\Sigma\varphi$ has no $\Sigma$-completion
with respect to~$V$
(and so $\Sigma\Pi\Sigma\varphi$
does not exist),
then $\varphi$ cannot have a completion
with respect to~$V$.
On the other hand,
we will see, loosely speaking, 
that if the making of the hierarchy
can proceed
unhindered
even if we go on endlessly using ordinal numbers,
that then we eventually obtain
the smallest completion~$\ol\varphi$
of~$\varphi$ with respect to~$V$.
This will all become clear 
in Subsection~\ref{SS:hierarchy-abstract}.

With some luck, the valuation $\Sigma\Pi\Sigma \varphi$
might already be complete with respect to~$V$.
We could then say that ``we hit $\ol\varphi$ in $3$ steps''.
We will prove in Subsection~\ref{SS:borel-hierarchy}
that in general
we need to make \emph{uncountably} many steps
before we hit~$\ol\varphi$. {}

\subsubsection*{No Completion}
Let us end the introducion
with an example
of a valuation
system $\vs{V}{L}\varphi{E}$
such that~$\varphi$ has no completion with respect to~$V$.
\begin{ex}
\label{E:no-completion}
Let $L$ be the sublattice of~$\R$ given by
\begin{equation*}
L\ \eqdf\ \{ \ n^{-1}\colon\ n\in\N\ \} \,\cup\, \{0\}.
\end{equation*}
Let $\varphi\colon L \longrightarrow \R$
be the valuation given by, for all~$n\in\N$,
\begin{alignat*}{3}
\varphi(\,n^{-1}\,) \ &=\ 1&\qquad\text{and}\qquad
\varphi(\,0\,) \ &=\ 0.
\end{alignat*}
Then we have a valuation system $\vs{\R}{L}\varphi{\R}$.

Let $C$ be a sublattice of~$\R$,
and let $\psi\colon C \longrightarrow \R$
be the valuation given by
\begin{equation*}
\vs{\R}{C}\psi\R\text{ \ is complete }
\qquad\text{and}\qquad
\psi\text{ extends }\varphi.
\end{equation*}
We will prove that this is not possible.

Consider the sequence $a_1 \geq a_2 \geq \dotsb$ in~$L$
given by $a_n = n^{-1}$.
We have 
\begin{equation*}
\psi(a_n) \ =\  \varphi(a_n) \ =\  1.
\end{equation*}
So $\bw_n \psi(a_n) = 1$.
In particular,
$a_1 \geq a_2 \geq \dotsb$ is $\psi$-convergent
(see Definition~\ref{D:phi-conv}).
Since $\psi$ is complete with respect to~$\R$
(see Definition~\ref{D:system-complete}), 
we get
\begin{equation*}
1\ = \ 
\bw_n \psi( a_n ) \ = \  \psi(\, \bw_n a_n \,)
\ =\ \psi(\,\bw_n n^{-1}\,)
\ =\ \psi(\, 0\, ) \ =\  0,
\end{equation*}
which is absurd.
Hence $\varphi$ has no completion with respect to~$\varphi$.

\end{ex}

\clearpage
%
%                     PI L
%

\subsection{The $\Pi$-Extension}
\label{S:Pi-extension}
\begin{dfn}
\label{D:PiL}
Let $\vs{V}{L}\varphi{E}$ be a valuation system.
Define 
\begin{equation*}
\Pi L \ \eqdf \ \{ \ 
\bw_n a_n \colon \ 
 a_1\geq a_2 \geq\dotsb \text{ from $L$ is }
\text{$\varphi$-convergent} \ \}.
\end{equation*}
\end{dfn}
%
%                  REMARK ON Pi(L) = L if L is Pi-complete
%
\begin{rem}
Let $\vs{V}{L}\varphi{E}$ be a valuation system.\\
Note that if~$\varphi$
is $\Pi$-complete
(see Definition~\ref{D:system-complete}),
then  $\Pi L = L$.
\end{rem}
%
%                  Pi(L) is lattice
%
\begin{lem}
\label{L:PiL-lattice}
Let $\vs{V}{L}\varphi{E}$ be a valuation system.\\
Then
$\Pi L$ is a sublattice of~$V$,
and $L$ is a sublattice of~$\Pi L$.
\end{lem}
\begin{proof}
We first prove that~$\Pi L$ is a sublattice of~$V$.
Let $a,b\in \Pi L$ be given;
we need to prove that $a\wedge b \in \Pi L$
and $a \vee b \in \Pi L$.
Choose $\varphi$-convergent 
$a_1 \geq a_2 \geq\dotsb$ and
$b_1 \geq b_2 \geq\dotsb$
with $a = \bw_n a_n$ and  $b=\bw_n b_n$.
Then $a_1 \wedge b_1 \geq a_2 \wedge b_2 \geq \dotsb$
is $\varphi$-convergent by Proposition~\ref{P:R-main},
and we have $\bw_n a_n \wedge b_n = a\wedge b$.
Hence $a\wedge b \in \Pi L$.
Similarly,
$a_1 \vee b_1 \geq a_2 \vee b_2 \geq \dotsb$
is $\varphi$-convergent by Proposition~\ref{P:R-main}
and
using $\sigma$-distributivity 
one can prove that $a \vee b  =\bw_n a_n \vee b_n$.
Hence $a\vee b \in \Pi L$.

To prove that~$L$ is a sublattice of~$\Pi L$,
we first note that~$L$ is a subset of~$\Pi L$.
Now,
since both $L$ and~$\Pi L$ are sublattices of~$V$,
and $L$ is a subset of~$\Pi L$,
we know that $L$ must be a sublattice of~$\Pi L$.
\end{proof}
\begin{rem}
\label{R:use-of-sigma-and-R}
In the proof of Lemma~\ref{L:PiL-lattice},
we have used the fact that~$V$ is $\sigma$-distributive
and the fact that~$E$ is $R$-complete
(via Proposition~\ref{P:R-main}).
\end{rem}

\begin{dfn}
\label{D:Pi-extendible}
Let $\vs{V}{L}\varphi{E}$ be a valuation system.\\
We say~$\varphi$ is \keyword{$\Pi$-extendible}
if there is a valuation $\psi \colon \Pi L \ra E$ with
\begin{equation*}
\psi ( \bw_n a_n )
\ =\ 
\bw_n \varphi(a_n) 
\qquad
\text{ for all $\varphi$-convergent }
a_1 \geq a_2 \geq \dotsb.
\end{equation*}
Clearly,
there can be at most one such map~$\psi$;
if it exists, we denote it by 
\begin{equation*}
\Pi \varphi\colon \Pi L \longrightarrow E.
\end{equation*}
Finally,
note that if $\varphi$ is $\Pi$-extendible,
then $\Pi\varphi$ extends $\varphi$ (hence the name).
\end{dfn}
%
%                  REMARK ON PI phi = phi.
%
\begin{rem}
\label{R:Pi-extendible-complete}
Let $\vs{V}{L}\varphi{E}$ be a valuation system.
\begin{enumerate}
\item
\label{R:Pi-extendible-complete-i}
Note that if~$\varphi$ is $\Pi$-complete
with respect to~$V$,\\
then~$\varphi$ is $\Pi$-extendible and~$\Pi\varphi=\varphi$.

\item
On the other hand,
if~$\varphi$ is $\Pi$-extendible,
and $\Pi\varphi=\varphi$,\\
then $\varphi$ is $\Pi$-complete
with respect to~$V$.
\end{enumerate}
\end{rem}
%
%                  EXAMPLES OF PI-EXTENDIBLE VALUATION SYSTEMS
%
\begin{ex}
The following valuation systems are $\Pi$-extendible.
\begin{equation*}
\vsLA\qquad\text{and}\qquad\vsLF
\end{equation*}
Indeed, this follows by 
Remark~\ref{R:Pi-extendible-complete}\ref{R:Pi-extendible-complete-i}
since these valuation systems are $\Pi$-complete.\\
More interestingly,
the following valuation systems are complete as well.
\begin{equation*}
\vsSA\qquad\text{and}\qquad\vsSF
\end{equation*}
This will follow from Lemma~\ref{L:Pi-minimal}.
\end{ex}

\begin{ex}
We leave it to the reader
to verify
that the  valuation system
\begin{equation*}
\vs{\R}{L}\varphi{\R}
\end{equation*}
from Example~\ref{E:no-completion}
is \emph{not}~$\Pi$-extendible.
\end{ex}
%
%                  PI EXTENDIBLE IMPLIES PI COMPLETE
%
\begin{lem}
\label{L:Pi-complete}
Let $\vs{V}{L}\varphi{E}$ be a valuation system.\\
If $\varphi$ is $\Pi$-extendible,
then $\vs{V}{\Pi L}{\Pi\varphi}{E}$ is $\Pi$-complete.
\end{lem}
\begin{proof}
Let $a^1 \geq a^2 \geq \dotsb $ from $\Pi L$
be given and suppose that $\bw_n \Pi \varphi (a^n)$ exists;
we need to prove that $\bw_n a^n \in \Pi L$
and that $\Pi\varphi(\bw_n a_n) = \bw_n \Pi \varphi(a_n )$
(see Def.~\ref{D:system-complete}).

To begin, write
$a^n = \bw_n a^n_m$
where $a^n_1 \geq a^n_2 \geq \dotsb$
is a $\varphi$-convergent sequence in~$L$
for each~$n\in\N$.
Define for each~$i\in\N$
an element $b_i \in L$ by
\begin{equation*}
b_i \ \eqdf\ \bw\{\,a^n_m\colon n,m\leq i \,\}.
\end{equation*}
Then $b_1 \geq b_2 \geq \dotsb$
and $\bw_n b_n = \bw_n a^n$.
Recall that  $\bw_n \Pi\varphi(a^n)$
exists.
\emph{%
We claim that $\bw_n \Pi\varphi(a^n)$
is the infimum of $\varphi (b_1) \geq \varphi (b_2) \geq \dotsb$.%
}
If we can prove this, we are done.
Indeed,
then $b_1 \geq b_2 \geq\dotsb$ is $\varphi$-convergent,
so $\bw_n a^n = \bw_n b_n \in \Pi L$,
and
\begin{alignat*}{3}
\Pi\varphi( \bw_n a^n ) \ &=\  \Pi\varphi( \bw_n b_n )
 \qquad&& \text{since $\bw_n a^n = \bw_n b_n$}, \\
&=\ \bw_n \varphi (b_n)  
  && \text{since $\varphi$ is $\Pi$-extendible,}\\
&=\ \bw_n \Pi\varphi(a^n) 
  && \text{by the claim.}
\end{alignat*}
Let us prove that 
 $\bw_n \Pi\varphi(a^n)$
is the infimum of $\varphi (b_1) \geq \varphi (b_2) \geq \dotsb$.

For each~$i$,
we have $b_i \geq \bw_{n\leq i} a^n = a^i$,
so $\varphi(b_i) = \Pi\varphi(b_i) \geq \Pi\varphi(a^i)$.
Hence we see that
$\bw_n \Pi\varphi(a^n)$ is a lower bound of
$\varphi(b_1)\geq \varphi(b_2) \geq \dotsb$.

On the other hand:
Let $\ell$ be a lower bound of $\varphi b_1 \geq \varphi b_2 \geq \dotsb$;
we need to prove that $\ell \leq \bw_n \Pi\varphi(a^n)$.
For all~$n$ and~$m$,
we have $a_m^n \geq b_{n\vee m}$
and so $\varphi(a_m^n) \geq \varphi(b_{n\vee m}) \geq \ell$.
Hence $\Pi\varphi(a^n) = \bw_m \varphi(a_m^n) \geq \ell$
for all~$n$.
So $\bw_n \Pi\varphi(a^n) \geq \ell$.

So $\bw_n \Pi\varphi(a^n)$ is the infimum of 
$\varphi b_1 \geq \varphi b_2 \geq \dotsb$,
and we are done.
\end{proof}
%
%                  Pi Pi = Pi
%
\begin{rem}
\label{R:PiPiisPi}
Let $\vs{V}{L}\varphi{E}$ be a $\Pi$-extendible valuation system.\\
By Lemma~\ref{L:Pi-complete}
we see that $\Pi\varphi$
is $\Pi$-complete with respect to~$V$.\\
So by
Remark~\ref{R:Pi-extendible-complete}\ref{R:Pi-extendible-complete-i}
we see that $\Pi\varphi$
is $\Pi$-extendible, and that
\begin{equation*}
\Pi(\Pi\varphi) \ =\  \Pi\varphi.
\end{equation*}
\end{rem}
%
%                  MINIMALITY OF PI L
%
\begin{lem}
\label{L:Pi-minimal}
Let $\vs{V}{L}\varphi{E}$ be a valuation system.\\
Let $C$ be a sublattice of~$V$.
Let $\psi\colon C\ra E$ be a valuation.
Assume
\begin{equation*}
\psi\text{ extends }\varphi
\qquad\text{and}\qquad
\vs{V}{C}\psi{E}\ \text{ is $\Pi$-complete.}
\end{equation*}
Then $\varphi$ is $\Pi$-extendible and
$\psi$ extends $\Pi\varphi$.
\end{lem}
\begin{proof}
We must prove that $\varphi$
is $\Pi$-extendible
and that~$\psi$ extends~$\Pi\varphi$.\\
Before we do this,
we will prove that
for every $\varphi$-convergent $a_1 \geq a_2 \geq \dotsb$,
we have
\begin{equation}
\label{eq:L:PiL-minimal}
\bw_n a_n \in C\qquad\text{and}\qquad \psi (\bw_n a_n) = \bw_n \varphi(a_n).
\end{equation}
We know that $\bw_n \varphi (a_n)$ exists
(since $a_1 \geq a_2 \geq \dotsb$ is $\varphi$-convergent),
and that $\varphi(a_n)= \psi(a_n)$
(since $\psi$ extends $\varphi$).
So $\bw_n \psi(a_n)$ exists too.
Hence $a_1 \geq a_2 \geq \dotsb$ is $\psi$-convergent.
Because $\vs{V}{C}\psi{E}$ is $\Pi$-complete
this implies that $\bw_n a_n \in C$ and 
$\psi (\bw_n a_n) = \bw_n \psi(a_n)$
(see Definition~\ref{D:system-complete}).
Hence we have proven Statement~\eqref{eq:L:PiL-minimal}.

Statement~\eqref{eq:L:PiL-minimal} implies that $\Pi L \subseteq C$.
So in order to prove that $\varphi$
is $\Pi$-extendible,
let us consider the valuation $\mu \eqdf \psi \,|\,\Pi L$.
In order to prove that~$\varphi$
is $\Pi$-extendible
we must show that
$\mu(\bw_n a_n) = \bw_n \varphi(a_n)$
for every $\varphi$-convergent sequence $a_1 \geq a_2 \geq \dotsb$
(see Def.~\ref{D:Pi-extendible}),
but this follows immediately from St.~\eqref{eq:L:PiL-minimal}.
Hence $\varphi$ is $\Pi$-extendible.

Finally,
since we know that~$\varphi$ is $\Pi$-extendible,
we can talk about~$\Pi\varphi$,
and write the second part of
St.~\eqref{eq:L:PiL-minimal}
as $\psi(\bw_n a_n) = \Pi \varphi(\bw_n a_n)$.
Hence $\psi$ extends $\Pi\varphi$.
\end{proof}
%
%                  Pi Extendible if and only if Pi completion exists
%
\begin{cor}
Let $\vs{V}{L}\varphi{E}$ be a valuation system.\\
Then $\varphi$ is $\Pi$-extendible
iff 
there is a valuation~$\psi\colon C\ra E$
such that 
\begin{equation*}
\vs{V}{C}\psi{E}\text{ \ is $\Pi$-complete}
\qquad\text{and}\qquad
\psi\text{ extends }\varphi.
\end{equation*}
(So, loosely speaking,
$\varphi$ is $\Pi$-extendible
iff $\varphi$ has a $\Pi$-complete extension.)
\end{cor}
\begin{proof}
Combine Lemma~\ref{L:Pi-minimal}
and Lemma~\ref{L:Pi-complete}.
\end{proof}
%
%                  MONOTONICITY OF PI
%
\begin{lem}
\label{L:Pi-monotonous}
Let $\vs{V}{L}\varphi{E}$ be a valuation system.\\
Let~$K$ be a sublattice of~$L$,
and let $\psi\colon K\ra E$ be a valuation
which extends $\varphi$.\\
Suppose that $\psi$ is $\Pi$-extendible.
Then $\varphi$ is $\Pi$-extendible
and  $\Pi\psi$ extends $\Pi\varphi$.
\end{lem}
\begin{proof}
Note that $\Pi\psi$ extends $\varphi$,
and $\vs{V}{\Pi K}{\Pi \psi}{E}$
is $\Pi$-complete (see Lemma~\ref{L:Pi-complete}).
So Lemma~\ref{L:Pi-minimal}
implies that
$\varphi$ is $\Pi$-extendible
and $\Pi\psi$ extends $\Pi\varphi$.
\end{proof}
%
%                  PI EXTENDIBLE IFF PI-CONTINUOUS
%
\begin{lem}
\label{L:Pi-continuity}
Let $\vs{V}{L}\varphi{E}$ be a valuation system.\\
Then $\varphi$ is $\Pi$-extendible
if and only if  $\varphi$
has the following property.
\begin{equation}
\label{eq:Pi-cont}
\left[ \quad
\begin{minipage}{.7\columnwidth}
For every $b\in L$ and $\varphi$-convergent
 $a_1 \geq a_2 \geq \dotsb$ from~$L$,
\begin{equation*}
\bw_n a_n \,\leq\, b
\quad\implies\quad
\bw_n \varphi (a_n) \,\leq\, \varphi(b).
\end{equation*}
\end{minipage}
\right.
\end{equation}
\end{lem}
\begin{proof}
``$\Longrightarrow$''\ 
Suppose $\varphi$ is $\Pi$-extendible.
Then $\varphi$ has Property~\eqref{eq:Pi-cont},
because if $b\in L$ and $\varphi$-convergent $a_1 \geq a_2 \geq \dotsb$
with $\bw_n a_n \leq b$ are given,
then we have
\begin{equation*}
\bw_n \varphi(a_n) 
\ =\ 
\Pi\varphi(\,\bw_n a_n\,)
\ \leq\ 
\Pi\varphi(b)
\ =\ 
\varphi(b).
\end{equation*}

\noindent``$\Longleftarrow$''\ 
Suppose~$\varphi$ has Property~\eqref{eq:Pi-cont};
we prove $\varphi$ is $\Pi$-extendible.
We claim that
\begin{equation}
\label{eq:Piphi-order-preserving}
\bw_n a_n \ \leq\ \bw_n b_n 
\quad\implies\quad
\bw_n \varphi(a_n) \ \leq\ \bw_n \varphi(b_n),
\end{equation}
where $a_1 \geq a_2 \geq \dotsb$ and $b_1 \geq b_2 \geq \dotsb$
are $\varphi$-convergent sequences in~$L$.

Indeed,
if $\bw_n a_n \leq \bw_n b_n$,
then $\bw_n a_n \leq b_m$ for all~$m$,
so $\bw_n \varphi(a_n) \leq \varphi(b_m)$ for all~$m$
(by Property~\eqref{eq:Pi-cont}),
and hence $\bw_n \varphi(a_n) \leq \bw_n \varphi(b_m)$.

Statement~\eqref{eq:Piphi-order-preserving} implies that
\begin{equation*}
\bw_n a_n \ =\  \bw_n b_n 
\quad\implies\quad 
\bw_n\varphi(a_n) \ =\  \bw_n \varphi(b_n),
\end{equation*}
so there is a unique map $\psi\colon \Pi L \ra E$ such that
\begin{equation}
\label{eq:def-psi}
\psi(\,\bw_n a_n\,) \ =\ \bw_n \varphi(a_n)
\qquad
\text{ for all $\varphi$-convergent }a_1 \geq a_2 \geq \dotsb.
\end{equation}
In fact, 
Statement~\eqref{eq:Piphi-order-preserving}
also implies that~$\psi$ is order preserving.

To prove that~$\varphi$ is $\Pi$-extendible
(see Definition~\ref{D:Pi-extendible}),
it suffices to show that~$\psi$ is a valuation.
For this,
it remains to be shown that~$\psi$ is modular
(see Definition~\ref{D:val}).

Let $a,b\in \Pi L$ be given.
To show that~$\psi$ is modular,
we need to prove that
\begin{equation*}
\psi(a\wedge b) \,+\, \psi(a\vee b)\ =\  \psi(a)\,+\,\psi(b).
\end{equation*}
Write  $a = \bw_n a_n$ and $b = \bw_n b_n$
where 
 $a_1 \geq a_2 \geq \dotsb$ and $b_1 \geq b_2 \geq \dotsb$
from~$L$
are $\varphi$-convergent sequences.
Then we have
\begin{alignat*}{3}
\varphi(a\wedge b) + \varphi(a \vee b) 
\ &=\ \psi(\bw_n a_n \wedge \bw_n b_n) \ +\  \psi(\bw_n a_n \vee \bw_n b_n) \\
\ &=\ \psi(\bw_n a_n \wedge b_n) \ +\  \psi(\bw_n a_n \vee b_n) \\
  &=\ \bw_n \varphi(a_n \wedge b_n) \,+\, \bw_n \varphi( a_n \vee b_n) \\
  &=\ \bw_n (\, \varphi(a_n \wedge b_n) \,+\, \varphi(a_n \vee b_n) \,)\\
  &=\ \bw_n (\, \varphi(a_n) \,+\, \varphi(b_n)\,) \\
  &\ \,\smash{\vdots}\  \\
  &=\ \psi(\bw_n a_n) \ +\  \psi(\bw_n b_n).
\end{alignat*}
Hence $\psi$ is modular, which completes the proof
that $\varphi$ is $\Pi$-extendible.
\end{proof}
\clearpage
%%%%%%%%%%%%%%%%%%%%%%%%%%%%%%%%%%%%%%%%%%%%%%%%%%%%%%%%%%%%%%%%%%%%%%%%%%%j%
%
%                  THE HIERARCHY OF EXTENSIONS
%
\subsection{The Smallest Complete Extension} $\,$\\
\label{SS:hierarchy}%
In the previous subsection,
we described the smallest $\Pi$-complete
extension
of a valuation system (when it exists).
In this subsection,
we will describe the smallest complete
extension
of a valuation system  (when it exists).
\begin{sit}
\label{Sit:smallest}
Let $\vs{V}{L}\varphi{E}$ 
and $\vs{V}{C}\psi{E}$
be valuation systems
such that 
\begin{equation*}
\text{$\psi$ extends~$\varphi$}
\qquad\text{and}\qquad
\text{$\psi$ is complete.}
\end{equation*}
\end{sit}
\noindent
In particular $\vs{V}{C}{\psi}{E}$ must be $\Pi$-complete
(see Definition~\ref{D:system-complete}).\\
Hence Lemma~\ref{L:Pi-minimal}
implies that $\varphi$ is $\Pi$-extendible 
and that $\psi$ extends  $\Pi \varphi$.\\
Thus, loosely speaking, $\Pi\varphi$ is the smallest extension of~$\varphi$
which is $\Pi$-complete.\\
In this subsection,
we identify the smallest extension~$\overline \varphi$ of~$\varphi$
which is complete.\\
We tackle this problem
in order to familiarise
the reader with the notions needed to define
``$\vs{V}{L}{\varphi}{E}$ is extendible''
(see Definition~\ref{D:extendible}).
These notions, which we introduce
 rather informally in this subsection,
will be defined rigorously and in a more general setting later on.

Let us begin. Note that
$\vs{V}{C}{\psi}{E}$ is also $\Sigma$-complete.
Hence $\varphi$ is $\Sigma$-extendible,
and $\psi$ extends $\Sigma\varphi$.
So we have the following situation.
\begin{equation*}
\text{$\psi$ extends both~$\Pi \varphi$ and~$\Sigma\varphi$}
\qquad\text{and}\qquad
\text{$\vs{V}{C}{\psi}{E}$ is complete.}
\end{equation*}
By a similar reasoning,
we see that $\Pi\varphi$ is $\Sigma$-extendible,
and that $\Sigma\varphi$ is $\Pi$-extendible
and that~$\psi$ extends both $\Sigma\Pi\varphi$ and $\Pi\Sigma\varphi$.
(Note that $\Pi(\Pi\varphi) =  \Pi\varphi$,
see Rem.~\ref{R:PiPiisPi}).
So we have the following situation.
\begin{equation*}
\text{$\psi$ extends both~$\Sigma\Pi \varphi$ and~$\Pi\Sigma\varphi$}
\qquad\text{and}\qquad
\text{$\vs{V}{C}{\psi}{E}$ is complete.}
\end{equation*}
Of course,
we can continue this proces.
More formally,
the clauses
\begin{alignat*}{5}
\Pi_{n+1} \varphi \, &=\, \Pi(\Sigma_n\varphi) &\qquad\quad
\Sigma_{n+1} \varphi \,&=\, \Sigma(\Pi_n\varphi) &\qquad\quad 
\Pi_0 \varphi \,&=\, \varphi \,&&=\, \Sigma_0 \varphi \\
\Pi_{n+1} L \, &=\, \Pi(\Sigma_n L) &\qquad\quad
\Sigma_{n+1} L \, &=\, \Sigma(\Pi_n L) &\qquad
\Pi_0 L \,&=\, L \,&&=\, \Sigma_0 L,
\end{alignat*}
give us valuation systems
$\vs{V}{\Pi_n L}{\Pi_n\varphi}{E}$ and
$\vs{V}{\Sigma_n L}{\Sigma_n\varphi}{E}$
for all $n\in\omega$.

Note that $\Pi\varphi$ extends $\varphi$.
Hence $\Sigma_2\varphi$ extends $\Sigma \varphi$ by
Lemma~\ref{L:Pi-monotonous}.
Hence $\Pi_3\varphi$ extends $\Pi_2 \varphi$.
Etcetera.
Similarly,
$\Sigma\varphi$ extends~$\varphi$,
so $\Pi_2 \varphi$ extends~$\Pi\varphi$,
and so on.

The hierarchy which we have obtained 
is shown in the following diagram.
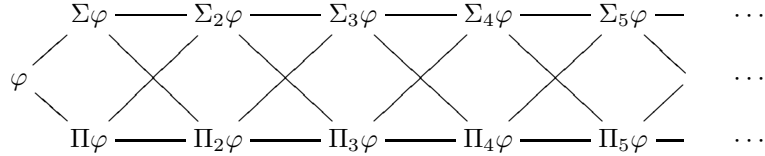
\begin{equation*}
\xymatrix @=10pt {
& \Sigma\varphi \ar @{-} [rr] \ar @{-} [rrdd]
&& \Sigma_2\varphi  \ar @{-} [rr] \ar @{-} [rrdd]
&& \Sigma_3\varphi  \ar @{-} [rr] \ar @{-} [rrdd]
&& \Sigma_4\varphi  \ar @{-} [rr] \ar @{-} [rrdd]
&& \Sigma_5 \varphi  \ar @{-} [r]\ar @{-} [rd]
&& \dotsb
\\  
\varphi \ar @{-} [ru] \ar @{-} [rd] 
&&&&&&&&&&&\dotsb\\
& \Pi\varphi \ar @{-} [rr] \ar @{-} [rruu]
&& \Pi_2\varphi \ar @{-} [rr]\ar @{-} [rruu]
&& \Pi_3\varphi \ar @{-} [rr]\ar @{-} [rruu]
&& \Pi_4\varphi \ar @{-} [rr]\ar @{-} [rruu]
&& \Pi_5 \varphi \ar @{-} [r] \ar @{-} [ru]
&& \dotsb
}
\end{equation*}
We say that the \emph{hierarchy collapsed at~$Q$}, where $Q  \in 
\{\,  L,\, \Pi_1 L, \, \Sigma_1 L,\,\Pi_2 L,\,\dotsc\,\}$, if
\begin{equation*}
\Pi( Q ) \,=\, Q \,=\, \Sigma(Q).
\end{equation*}
In that case, let~$q\colon Q\ra E$ be the associated valuation
(either $\Pi_n\varphi$ or $\Sigma_n\varphi$ for some~$n$).
Then $\vs{V}{Q}{q}{E}$ is complete,
since it is $\Pi$-complete 
and $\Sigma$-complete.

Note that the definition of $\Pi_n \varphi$
and $\Sigma_n \varphi$ does not depend
on which complete extension~$\psi$ of~$\varphi$ is given,
only on the fact that such~$\psi$ exists.
In particular,
if $\vs{V}{C'}{\psi'}{E}$ is any complete valuation system
such that~$\psi'$ extends~$\varphi$,
then $\psi'$ extends~$\Pi_n L$ and $\Sigma_n L$.
In particular,
such~$\psi'$ extends~$q$.
Hence~$q$ is the smallest complete extension of~$\varphi$ we sought.

However,
in general the hierarchy need not have collapsed at
any $\Pi_n L$ or~$\Sigma_n L$,
as we will show later on, in Subsection~\ref{SS:borel-hierarchy}.
So to find the smallest complete extension of~$\varphi$,
we need to push forwards.

To this end, consider the sets $\Pi_\omega L$
and $\Sigma_\omega L$ given by
\begin{equation*}
\Pi_\omega L \,\eqdf\,\textstyle{\bigcup_n}\, \Pi_n L
\qquad\text{and}\qquad
\Sigma_\omega L \,\eqdf\, \textstyle{\bigcup_n}\, \Sigma_n L.
\end{equation*}
Since $\Pi_{n} L \subseteq \Sigma_{n+1} L$
and $\Sigma_{n} L \subseteq \Pi_{n+1}L$ for all~$n$,
we see that $\Pi_\omega L = \Sigma_\omega L$.

Now,
since $\Pi_n\varphi$ extends $\Pi_m\varphi$
for $n\geq m$,
there is a unique map~$\Pi_\omega \varphi\colon \Pi_\omega L \ra E$
which extends all~$\Pi_n \varphi$.
One can easily see that $\vs{V}{\Pi_\omega L}{\Pi_\omega\varphi}{E}$
is a valuation system.
Similarly, there is a unique map 
$\Sigma_\omega \varphi \colon \Sigma_\omega L \ra E$
which extends all~$\Sigma_n\varphi$.
Then $\vs{V}{\Sigma_\omega L}{\Sigma_\omega\varphi}{E}$
is a valuation system.

Since $\Pi_{n+1}\varphi$ extends $\Sigma_{n}\varphi$
for all~$n$, one sees that $\Pi_\omega\varphi = \Sigma_\omega\varphi$.

Again, the hierarchy might have collapsed at~$\Pi_\omega L$,
i.e.,
\begin{equation*}
\Pi(\Pi_\omega L ) \,=\, \Pi_\omega L \,=\, \Sigma(\Pi_\omega L).
\end{equation*}
In that case $\Pi_\omega\varphi$ the minimal completion of~$\varphi$
we sought.

However,
again the hierarchy
might not have collapsed at~$\Pi_\omega L$,
so we consider the valuations
 $\Pi_{\omega+n}\varphi\eqdf\Pi_{n} (\Pi_\omega \varphi)$
and $\Sigma_{\omega+n}\varphi \eqdf \Sigma_{n}(\Pi_\omega\varphi)$.

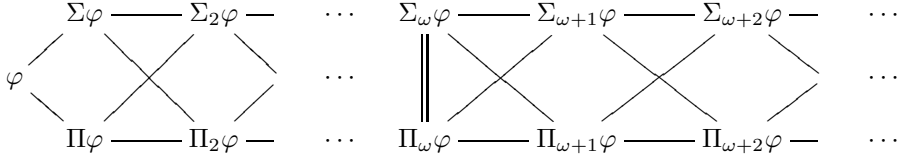
\begin{equation*}
\xymatrix @=10pt {
& \Sigma\varphi \ar @{-} [rr] \ar @{-} [rrdd]
&& \Sigma_2 \varphi  \ar @{-} [r]\ar @{-} [rd]
&& \dotsb
& \Sigma_{\omega} \varphi \ar @{-} [rr] \ar @{-} [rrdd]
                          \ar @{=} [dd]
&& \Sigma_{\omega+1}\varphi \ar @{-} [rr] \ar @{-} [rrdd]
&& \Sigma_{\omega+2} \varphi  \ar @{-} [r]\ar @{-} [rd]
&& \dotsb
\\  
\varphi \ar @{-} [ru] \ar @{-} [rd] 
&&&&&\dotsb
&&
&&&&&\dotsb\\
& \Pi\varphi \ar @{-} [rr] \ar @{-} [rruu]
&& \Pi_{2} \varphi \ar @{-} [r] \ar @{-} [ru]
&& \dotsb
& \Pi_{\omega}\varphi \ar @{-} [rr] \ar @{-} [rruu]
&& \Pi_{\omega+1}\varphi \ar @{-} [rr] \ar @{-} [rruu]
&& \Pi_{\omega+2} \varphi \ar @{-} [r] \ar @{-} [ru]
&& \dotsb
}
\end{equation*}
With induction on ordinal numbers,
we can continue this process endlessly.
However, 
the collapse of the hierarchy
can not be postponed indefinitely.

More formally,
let $\overline{L}\eqdf \{\, c\in C\colon
\exists\alpha[\,c\in \Pi_\alpha L\,]\, \}$.
Then we have $\Pi_\alpha L \subseteq \overline L$
for every (ordinaln number)~$\alpha$.
We want to prove that~$\Pi_\alpha L=\overline L$
for some~$\alpha$.
Define
\begin{equation*}
\alpha(c) \ =\ \min\, \{\  \beta\colon\   c\in \Pi_\beta L \ \}
\qquad\qquad(c\in \overline L).
\end{equation*}
Then the set of ordinal numbers $\{\,\alpha(c)\colon\,c\in \overline L\,\}$
has a supremum, say~$\xi$.
We have 
\begin{equation*}
c\in \Pi_{\alpha(c)} L \ \subseteq \ \Pi_\xi L
\qquad\qquad
(c\in \overline L).
\end{equation*}
So $\overline L \subseteq \Pi_{\xi} L$.
But we already had $\Pi_{\xi}L\subseteq \overline L$.
Hence $\Pi_{\xi} L = \overline L$.

We claim that the hierarchy has collapsed at~$\Pi_\xi L$,
i.e., 
\begin{equation*}
\Pi(\Pi_\xi L) \,=\, \Pi_\xi L \,=\, \Sigma(\Pi_\xi L).
\end{equation*}
Indeed,
we have 
\begin{equation*}
\Pi_\xi L \ \subseteq\ \Sigma(\Pi_\xi L ) \ \subseteq \ \overline L
\ = \ \Pi_\xi L.
\end{equation*}
So $\Sigma(\Pi_\xi L) = \Pi_\xi L$.
Similarly, $\Pi_\xi L = \Pi(\Pi_\xi L)$.

One can easily verify that
$\overline\varphi \eqdf \Pi_\xi \varphi$
is the smallest complete  extension of~$\varphi$.

All in all, we have proven the following.
\begin{prop}
\label{P:smallest-complete-extension}
Let $\vs{V}{L}\varphi{E}$ 
be valuation system.\\
If there is a valuation $\psi\colon C\ra E$ such that
\begin{equation*}
\text{$\psi$ extends~$\varphi$}
\qquad\text{and}\qquad
\text{$\vs{V}{C}{\psi}{E}$ is complete,}
\end{equation*}
then there is a \emph{smallest} such valuation,\\
that is, 
there is a valuation $\ol\varphi\colon \ol{L}\ra E$
such that
\begin{equation*}
\text{$\ol{\varphi}$ extends~$\varphi$}
\qquad\text{and}\qquad
\text{$\vs{V}{\ol{L}}{\ol{\varphi}}{E}$ is complete,}
\end{equation*}
and such that 
$\psi'$ extends $\ol{\varphi}$
 for every valuation $\psi'\colon C\ra E$ with
\begin{equation*}
\text{$\psi'$ extends~$\varphi$}
\qquad\text{and}\qquad
\text{$\vs{V}{C}{\psi'}{E}$ is complete.}
\end{equation*}

Moreover,
$\overline{\varphi} = \Pi_\xi \varphi$
for some ordinal number~$\xi$.
\end{prop}

\clearpage
%%%%%%%%%%%%%%%%%%%%%%%%%%%%%%%%%%%%%%%%%%%%%%%%%%%%%%%%%%%%%%%%%%%%%%%%%%%%%%%
%
%                  THE NEED FOR THE HIERARCHY
%
\subsection{The Borel Hierarchy Theorem}$\,$\\
\label{SS:borel-hierarchy}%
Before we continue our study of the hierarchy introduced in the
Subsection~\ref{SS:hierarchy}, let us take a step back and wonder: 
 is this all --- the endless hierarchy --- neccesary?

Indeed,
using the terminology of Subsection~\ref{SS:hierarchy},
it is not unthinkable
that the hierarchy is always collapsed
at, say $\Sigma_{37} L$.
In that case 
the theory would be much simpler;
we would only need to
use the symbols up to~``$\Sigma_{37}$''.
In particular,
the involvement of the (infinite) ordinal numbers would not be required.

It turns out that we \emph{do} need a large amount of symbols
to desribe the hierarchy.
In this subsection we will give an example
where the hierarchy can only be collaped at~$\Pi_\alpha L$
or at~$\Sigma_\alpha L$
if the ordinal number~$\alpha$ is \emph{uncountable},
see Proposition~\ref{P:no-countable-collapse}.

On the bright side,
it does not get worse than this:
we will see 
(in Lemma~\ref{L:aleph1})
that 
the hierarchy is always collapsed at~$\Pi_{\aleph_1} L$,
where~$\aleph_1$
is the set of all countable ordinal numbers,
i.e., the smallest uncountable ordinal number.

The material in the subsection will not be used later on.
So the reader can safely skip this subsection
and proceed to Subsection~\ref{SS:hierarchy-abstract} 
on page~\pageref{SS:hierarchy-abstract}
if so desired.

\subsubsection{Borel Subsets}$\,$\\
Our example involves Borel sets.
Recall that the \emph{Borel subsets}
of a topological space~$X$
(such as~$\R$)
are those subsets one can form using countable intersection
and countable union starting from the open subsets.

Instead of~$\R$,
we  work with the Borel subsets
of the 
 \emph{Baire space},
 $\N^\N$.
In short,
the topology on~$\N^\N$
is the product topology
when~$\N$ is given the discrete toplogy.
To understand these words,
one might look at~\cite{Willard70},
but this is not necessary as we will
give a more direct description of~$\mathcal{T}$
in Subsubsection~\ref{SSS:bhier-open}.

While we could do the following
for~$\R$ as well,
it is much easier for~$\N^\N$.
\begin{nt}
Let $\mathcal{T}$
denote the set of open subsets of~$\N^\N$, \\
and let~$\mathcal{B}$
denote the set of Borel subsets of~$\N^\N$.
\end{nt}
Note that $\mathcal{B}$
is a sublattice of~$\wp(\N^\N)$,
and~$\mathcal{B}$ is a sublattice of~$\mathcal{T}$.

\begin{dfn}
Let $\psi\colon \mathcal{B}\ra \R$
be the map given by, for all~$A\in \mathcal{B}$,
\begin{equation*}
\psi(A)\ =\ 0.
\end{equation*}
\end{dfn}
Then $\psi$ is a valuation,
and we have the following valuation system.
\begin{equation*}
\vs{\wp(\N^\N)}{\mathcal{B}}{\psi}{\R}.
\end{equation*}

\begin{lem}
The valuation $\psi$ is complete with respect to~$\wp(\N^\N)$
(see Def.~\ref{D:system-complete}).
\end{lem}
\begin{proof}
Let $A_1 \supseteq A_2 \supseteq \dotsb$
and $B_1 \subseteq B_2 \subseteq \dotsb$
be a $\psi$-convergent sequences in~$\mathcal{B}$.
To prove that~$\psi$ is complete relative to~$\wp(\N^\N)$
we must show that 
\begin{equation*}
\textstyle 
\bigcap_n A_n \,\in\, \mathcal{B}
\qquad\text{and}\qquad
\bigcup_n B_n \,\in\, \mathcal{B}.
\end{equation*}
This follows immediately by definition of
the Borel subsets.
\end{proof}

\begin{rem}
\label{R:bhier-A}
Let $\mathcal{A}$ be a sublattice of~$\mathcal{B}$
and let
\begin{equation*}
\varphi\colon\mathcal{A}\longrightarrow \R
\end{equation*}
be the restriction of~$\psi$ to~$\mathcal{A}$.
Note that we are in Situation~\ref{Sit:smallest},
\begin{equation*}
\text{$\psi$ extends $\varphi$}
\qquad\text{and}\qquad
\text{$\psi$ is complete}.
\end{equation*}
Using the notation from Subsection~\ref{SS:hierarchy},
let us see what~$\Pi\mathcal{A}$
and $\Sigma\mathcal{A}$ are.

Note that every sequence $A_1 \supseteq A_2 \supseteq \dotsb$
in~$\mathcal{A}$ is $\varphi$-convergent,
as~$\psi$ is constant.
Further, given $A_1,A_2,\dotsc \in \mathcal{A}$,
we have $\bigcap_n A_n = \bigcap_n A_n'$,
where $A_1' \supseteq A_2' \supseteq \dotsb$
are defined by $A_n' = A_1 \cap \dotsb \cap A_n$.
So we see that
\begin{alignat}{3}
\label{eq:bhier-PiA}
\Pi\mathcal{A}\ &=\ 
\{ \ \textstyle{\bigcap_n} A_n\colon\ 
A_1,A_2,\dotsc \,\in\, \mathcal{A}\ \}.
\shortintertext{%
By a similar reasoning it is not hard to see that }
\label{eq:bhier-SigmaA}
\Sigma\mathcal{A}\ &=\ 
\{ \ \textstyle{\bigcup_n} A_n\colon\ 
A_1,A_2,\dotsc \,\in\, \mathcal{A}\ \}.
\end{alignat}
\end{rem}

\subsubsection{Statement of the Borel Hierarchy Theorem}$\,$\\
Let us spend words on where we are headed.
We will define a sublattice~$\mathcal{A}$
of~$\wp(\N^\N)$
in such a way that,
using the notation of Remark~\ref{R:bhier-A},
we have that
$\Sigma \mathcal{A}$
is  precisely the family of open subsets of~$\N^\N$,
while~$\Pi \mathcal{A}$ is the family
of closed subsets of~$\N^\N$.
From this information,
the reader can deduce with induction and the principle
of the excluded middle, that for every ordinal~$\alpha>0$,
and all $A\subseteq\N^\N$, 
\begin{equation}
\label{eq:bhier-complement}
A \,\in\, \Pi_\alpha \mathcal{A}
\qquad\iff\qquad
\N^\N\backslash A \,\in\, \Sigma_\alpha \mathcal{A}.
\end{equation}
The aim of this subsection
is to prove the following statement.
\begin{equation}
\label{eq:hierarchy-theorem}
\left[\quad
\begin{minipage}{.7\columnwidth}
\textbf{%
Let~$\alpha$ be a \emph{countable} ordinal number.\\
There is a set~$S\in \Sigma(\Pi_\alpha\mathcal{A})$
with $S\notin \Pi(\Sigma_\alpha\mathcal{A})$, and \\
there is a set $P\in\Pi(\Sigma_\alpha\mathcal{A})$
with $P\notin \Sigma(\Pi_\alpha\mathcal{A})$.}
\end{minipage}
\right.
\end{equation}
In particular
this means that if the hierarchy has collapsed
at~$\Pi_\alpha\mathcal{A}$
or at~$\Sigma_\alpha\mathcal{A}$
for some ordinal number~$\alpha$
then $\alpha$ must be \emph{uncountable},
see Proposition~\ref{P:no-countable-collapse}.

Statement~\eqref{eq:hierarchy-theorem}
is known in descriptive set theory
as the \emph{Borel Hierarchy Theorem}.
We will give a proof of Statement~\eqref{eq:hierarchy-theorem}
 in this subsection that
uses the principle of excluded middle
and is based on a beautiful paper by Veldman~\cite[paragraph 5]{Veldman08}.%
\footnote{
In this paper~\cite{Veldman08},
Veldman (also) gives
a proof of a variant of the Borel Hierarchy Theorem
using Brouwer's Continuity Principle
instead of the principle of excluded middle.}

\subsubsection{Open Subsets of $\N^\N$}$\,$\label{SSS:bhier-open}\\
Before we give a definition of~$\mathcal{A}$,
and start with the proof of Statement~\eqref{eq:hierarchy-theorem}
let us describe the topology~$\mathcal{T}$ 
on the Baire space~$\N^\N$ in more detail.
\begin{dfn}
Given $m,n\in\N$, define $B^m_n$ by
\begin{equation*}
B_n^m\ \eqdf\ \{ \ f\in\N^\N\colon\  f(n) \,=\, m \ \}.
\end{equation*}
\end{dfn}

\begin{rem}
\label{R:bhier-open}
For $A\subseteq \N^\N$,
we have $A\in \mathcal{T}$
if and only if for each $f\in A$
we have,
\begin{equation*}
f\,\in\, B^{m_1}_{n_1} \cap \dotsb \cap B^{m_K}_{n_K} \ \subseteq\ A,
\end{equation*}
for some $K\in\N$ and $m_1,\dotsc,m_K\in\N$
and $n_1,\dotsc,n_K\in\N$.
\end{rem}

We can formulate~Remark~\ref{R:bhier-open} more abstractly
with some notation.
\begin{dfn}
Let $\mathcal{S}$ and $\mathcal{S}_\cap$
be families of subsets of~$\N^\N$
given by
\begin{alignat*}{3}
\mathcal{S} \ &\eqdf\  
\{\ B^m_n \colon \ m,n\in\N \ \} \\
\mathcal{S}_\cap \ &\eqdf \ 
\{ \ S_1 \cap \dotsb\cap S_K\colon\ 
K\in\N,\  S_k\in \mathcal{S}\ \}.
\end{alignat*}
\end{dfn}

\begin{rem}\label{R:bhier-basis}
By Remark~\ref{R:bhier-open},
we see that~$\mathcal{S}$
is a subbasis for the topology~$\mathcal{T}$
on~$\N^\N$,
and we see that~$\mathcal{S}_\cap$ is a basis
for~$\mathcal{T}$.
In particular, we get
\begin{equation}
\label{eq:expression-t}
\mathcal{T}\ =\ 
\{\ \textstyle{\bigcup_n} A_n\colon \ A_1,A_2,\dotsc \in \mathcal{S}_\cap\ \}.
\end{equation}
\end{rem}
\begin{rem}\label{R:bhier-clopen}
Let $m,n\in\N$. Then  $B^m_n \in \mathcal{T}$ by Remark~\ref{R:bhier-basis}.
More suprisingly,
\begin{equation*}
\N^\N \backslash B^m_n\,\in\, \mathcal{T},
\end{equation*}
that is, $B^m_n$ is closed as well.
Indeed, this follows by the following equality.
\begin{equation*}
\N^\N \backslash B^m_n \ =\ 
\textstyle{\bigcup} \{\  B^k_n\colon\ k\in\N,\ k\neq m\ \}.
\end{equation*}

\end{rem}
\subsubsection{Definition of the Sublattice~$\mathcal{A}$ of~$\N^\N$}$\,$\\
Recall that we want to define a sublattice~$\mathcal{A}$
of~$\wp(\N^\N)$
so that $\Sigma \mathcal{A}$
are the open subsets of~$\N^\N$,
while~$\Pi \mathcal{A}$  are the closed subsets
(see Remark~\ref{R:bhier-A}).

Since the elements of~$\mathcal{S}$ are both open and closed,
we let~$\mathcal{A}$ be the sub-Boolean algebra
of~$\wp(\N^\N)$
generated by~$\mathcal{S}$.
More concretely:
\begin{dfn}\label{D:bhier-A}
Let $\mathcal{S}'$, $\mathcal{S}_\cap'$
and $\mathcal{A}$ be the families of subsets of~$\N^\N$ given by
\begin{alignat*}{3}
\mathcal{S}' 
\ &\eqdf\ 
\{\ B^m_n \colon \ m,n\in\N \ \}
\ \cup\ 
\{\ \N^\N\backslash B^m_n \colon \ m,n\in\N \ \}\\
\mathcal{S}_\cap'
\ &\eqdf\ 
\{ \ S_1 \cap \dotsb\cap S_K\colon\ 
K\in\N,\  S_k\in \mathcal{S}'\ \} \\ 
\mathcal{A}
\ &\eqdf\ 
\{ \ T_1 \cup \dotsb\cup T_L\colon\ 
L\in\N,\  T_\ell\in \mathcal{S}_\cap'\ \}.
\end{alignat*}
\end{dfn}
\begin{lem}
The family
$\mathcal{A}$
is a sublattice of~$\wp(\N^\N)$,
and 
\begin{equation*}
\N^\N\backslash A\,\in\,\mathcal{A}
\qquad\iff\qquad
A\,\in\,\mathcal{A}
\end{equation*}
for every~$A\subseteq\N^\N$,
and we have the following equalities.
\begin{alignat*}{3}
\Sigma\mathcal{A} \ &=\ 
\{ \ U\subseteq \N^\N\colon\ U\text{ is open}\ \},\\
\Pi\mathcal{A} \ &=\ 
\{ \ C\subseteq \N^\N\colon \ C\text{ is closed}\ \}.
\end{alignat*}
\end{lem}
\begin{proof}
We leave this to the reader.
\end{proof}

\subsubsection{Encoding the Elements of~$\mathcal{A}$}$\,$\\
Now that we have defined~$\mathcal{A}$,
we can start the proof of
 Statement~\eqref{eq:hierarchy-theorem}.
Maybe the most important idea behind the proof presented here
is that 
we can
encode the  Borel subsets~$\mathcal{B}$ of~$\N^\N$
as elements of~$\N^\N$.

To warm up
let us see how we can encode 
a tuple $a_1 \dotsb a_n$ of natural numbers as a natural number.
We leave it to the reader to 
find a bijection
\begin{equation*}
\left< -,-\right>\colon \ \N\times \N \longrightarrow \N \backslash\{1\}.
\end{equation*}
Let $\N^*$ be the set of tuples on~$\N$.
Given a tuple $a_1\dotsb a_n \in \N^*$, define
\begin{equation*}
\left< a_1a_2  \dotsb a_n \right>\ \eqdf\ 
\left<a_1, \left<a_2,\dotsc \,\left<a_n, 1 \right>\,\dotsb\,\right>\right>.
\end{equation*}
Then the resulting map  $\left<-\right>\colon \N^* \ra \N$
is a bijection.

Let us now encode the elements of~$\mathcal{A}$
 (see Def.~\ref{D:bhier-A}).
Given $k \in \N$, let
\begin{equation*}
 \decode{k}^{\mathcal{S}'}
\ \eqdf \ 
\begin{cases}
\ \N^\N \backslash B^m_n & \text{ if $k\equiv\left<1mn\right>$}, \\
\ B^m_n & \text{ if $k\equiv\left<2mn\right>$}, \\
\ \varnothing &\text{ otherwise.}
\end{cases}
\end{equation*}
Then $\decode-^{\mathcal{S}'}\colon\N\ra \mathcal{S}'$
is a surjection.
Given~$k\in \N$ with $k\equiv \left<a_1\dotsb a_K\right>$, let
\begin{equation*}
 \decode{k}^{\mathcal{S}_\cap'}
\ \eqdf\ 
\decode{a_1}^{\mathcal{S}'} \,\cap\,\dotsb\,\cap\,\decode{a_K}^{\mathcal{S}'}.
\end{equation*}
Then $\decode-^{\mathcal{S}_\cap'}\colon \N\ra {\mathcal{S}_\cap'}$
is a surjection.
Given~$k\in \N$ with $k\equiv \left<a_1\dotsb a_K\right>$, let
\begin{equation*}
 \decode{k}^{\mathcal{A}}
\ \eqdf\ 
\decode{a_1}^{\mathcal{S}_\cap'} \,\cup\,
   \dotsb\,\cup\,\decode{a_K}^{\mathcal{S}_\cap'}.
\end{equation*}
Then $\decode-^{\mathcal{A}}\colon \N\ra {\mathcal{A}}$
is a surjection.

Let $A\in \mathcal{A}$ be given.
If
$\decode{k}^\mathcal{A}=A$
for some 
 $k\in \N$ 
we say
that~$k$ is a \emph{code} for~$A$.
Note that~$A$ might have multiple codes.
This will not be a problem.

\subsubsection{Encoding Countable Ordinal Numbers}$\,$\\
Before we can make the step
to encode all Borel subsets of~$\N^\N$
we need an encoding for
the countable ordinal numbers.
We need some notation.
\begin{nt}
Let $f\in\N^\N$
and $n\in \N$ be given.
Define $f^{[n]}\in \N^\N$ by, for $m\in\N$,
\begin{equation*}
f^{[n]} (m) \ =\ f(\,\left<n,m\right>\,).
\end{equation*}
\end{nt}
Since $\left< -,-\right>$
is a bijection from~$\N\times\N$
to $\N\backslash\{1\}$,
an element $f\in\N^\N$
is completely determined
by $f^{[1]},f^{[2]},\dotsc$
\emph{and}~$f(1)$.
More precisely,
the assignment
\begin{equation*}
f\quad\mapsto\quad
 f(1),\ f^{[1]},\  f^{[2]},\ \dotsc.
\end{equation*}
gives a bijection from
$\N^\N$ to $\N \times (\N^\N)^\N$.

To encode the countable ordinal numbers,
we use special elements of~$\N^\N$.
\begin{dfn}
\label{D:bhier-stump}
Let~$\Stump$ be the subset  of~$\N^\N$
inductively given by:
\begin{enumerate}
\item
If $f\in\N^\N$ and $f(1)\neq 1$, then $f\in \Stump$.

\item
If $f\in\N^\N$
and $f(1)=1$
and $f^{[n]}\in\Stump$ for all~$n\in\N$,
then $f\in \Stump$.
\end{enumerate}
The elements of~$\Stump$
are called \keyword{stumps}
and are used in Intuitionistic Mathematics
as a replacement for the countable ordinal numbers.
\end{dfn}
\begin{dfn}
\label{D:bhier-alpha}
Let $\alpha[-]\colon \Stump \ra \aleph_1$
be the map recursively defined by
\begin{equation*}
\alpha[f]\ =\ 
\begin{cases}
\ 0&
\text{if $f(1)\neq1$;}\\
\ \bv_{n\in \N}\ \alpha[f^{[n]}]+1&
\text{if $f(1)= 1$.}
\end{cases}
\end{equation*}
Recall that~$\aleph_1$
is the set of all countable ordinal numbers.
\end{dfn}

\begin{lem}
The map $\alpha[-]\colon \Stump\longrightarrow \aleph_1$ is surjective.
\end{lem}
\begin{proof}
To prove that $\alpha[-]$ is surjective,
we must 
show that
for each~$\alpha\in\aleph_1$
there is an~$f\in \Stump$ with $\alpha[f]=\alpha$.
We use induction on~$\alpha\in\aleph_1$.

First we must find an~$f\in\Stump$
with $\alpha[f]=0$.
Simply take the $f \in \N^\N$
with $f(n)=37$ for all~$n\in \N$.
Then $f(1)\neq 0$, so~$f\in\Stump$, and $\alpha[f]=0$.

Let $\alpha\in \aleph_1$
be given, 
and assume that $\alpha = \alpha[f]$
for some~$f\in\Stump$.
We need to find a~$g\in \Stump$
such that $\alpha[g]=\alpha+1$.
Define $g\in\N^\N$ by 
\begin{equation*}
g(1)\,=\,1,\qquad\text{and}\qquad
g^{[n]}\,=\, f\quad\text{for all~$n\in\N$}.
\end{equation*}
Then $g\in\Stump$, and $\alpha[g]=\bv_{n\in\N} \ \alpha[f]+1$.
Since $\alpha[f]=\alpha$, we have $\alpha[g]=\alpha+1$.

Let $\lambda\in\aleph_1$
be a limit ordinal,
and assume that $\alpha[-]$ is surjective on~$\lambda$.
We must find an~$f\in \Stump$ such that~$\alpha[f]=\lambda$.
Since $\lambda\in\aleph_1$
there are  $\alpha_1 ,\alpha_2 ,\dotsc \in \lambda$
such that $\lambda = \bv_{n\in\N}\ \alpha_n+1$.
Since $\alpha[-]$
is surjective on~$\lambda$,
we know that there are  $f_1,f_2,\dotsc \in \Stump$
with $\alpha[f_n] = \alpha_n$.
Define $f\in\N^\N$ by 
\begin{equation*}
f(1)\,=\,1\qquad\text{and}\qquad
f^{[n]} \,=\, f_n\quad\text{for all $n\in\N$}.
\end{equation*}
Then $f\in\Stump$,
and $\alpha[f] = \bv_{n\in\N}\  \alpha[f^{[n]}]+1 
= \bv_{n\in\N}\ \alpha_n +1 = \lambda$.
\end{proof}

\subsubsection{Encoding Borel Subsets of~$\N^\N$}$\,$\\
We are now ready to encode the Borel subsets of~$\N^\N$.
\begin{dfn}
With recursion on~$\Stump$
define for each~$f\in\Stump$
maps
\begin{equation*}
\decode{-}^\Pi_f \colon \N^\N \longrightarrow \Pi(\Sigma_{\alpha[f]}\mathcal{A})
\qquad
\text{and}
\qquad
\decode{-}^\Sigma_f \colon \N^\N 
\longrightarrow \Sigma(\Pi_{\alpha[f]}\mathcal{A})
\end{equation*}
such that the following clauses hold.
\begin{enumerate}
\item
For $f\in \Stump$ with $f(1)\neq 1$ we have
\begin{alignat*}{3}
\decode{g}^\Pi_f \ &=\  
\bigcap_{n\geq 2} \ \decode{\,g(n)\,}^\mathcal{A}&
\qquad\text{and}\qquad
\decode{g}^\Sigma_f \ &=\ 
\bigcup_{n\geq 2}\  \decode{\,g(n)\,}^\mathcal{A}.
\end{alignat*}

\item
For $f\in \Stump$ with $f(1)= 1$
we have
\begin{equation*}
\decode{g}^\Pi_f \ =\ 
\bigcap_{n\in\N}\ \decode{\,g^{[n]}\,}^\Sigma_{f^{[n]}}
\qquad\text{and}\qquad
\decode{g}^\Sigma_f \ =\ 
{\bigcup_{n\in\N}}\ \decode{\,g^{[n]}\,}^\Pi_{f^{[n]}}.
\end{equation*}
\end{enumerate}
\end{dfn}
\begin{lem}
For each~$f\in\Stump$
the maps
$\decode{-}^\Pi_f$
and 
$\decode{-}^\Sigma_f$
are surjective.
\end{lem}
\begin{proof}
We leave this to the reader.
\end{proof}
\begin{rem}
\label{R:bhier-g1}
Let $f\in \Stump$
and $g\in \N^\N$ be given.
Note that $\decode{g}^\Pi_f$
and $\decode{g}^\Sigma_f$
do not depend on~$g(1)$.
More precisely,
given $g'\in\N^\N$ 
with $g'(n)=g(n)$ for all~$n\geq 2$
--- so  possibly $g(1) \neq g'(1)$.
Then we have  $\decode{g}^\Pi_f = \decode{g'}^\Pi_f$
and $\decode{g}^\Sigma_f = \decode{g'}^\Sigma_f$.
\end{rem}
We use Remark~\ref{R:bhier-g1}
to combine the maps~$\decode{-}^\Pi_f$
and $\decode{-}^\Sigma_f$
into one map~$\decode{-}^\mathcal{B}_f$.
\begin{dfn}
Let
$\decode{-}^{\mathcal{B}}_{-}\colon 
  \N^\N\times \Stump \longrightarrow \mathcal{B}$
be given by, for $f\in\Stump$ and~$g\in\N^\N$,
\begin{equation}
\label{eq:decode-Borel}
\decode{g}^{\mathcal{B}}_f \ = \ 
\begin{cases}
\  \decode{g}^{\Pi}_f  \quad
& \text{if $g(1)= 37$},\\

\ \decode{g}^{\Sigma}_f
& \text{if $g(1)\neq 37$}.
\end{cases}
\end{equation}
\end{dfn}
We want to prove that $\decode{-}_{-}^\mathcal{B}$
is surjective. 
To do this,
we need a lemma.
\begin{lem}
\label{L:Borel-aleph1}
We have the following equality.
\begin{equation}
\label{eq:Borel-aleph1}
\mathcal{B} \ =\  \Pi_{\aleph_1} \mathcal{A}.
\end{equation}
\end{lem}
\begin{proof}
Note that we have already proven 
(at the end of
in Subsection~\ref{SS:hierarchy})
 that $\mathcal{B} = \Pi_\xi \mathcal{A}$
for \emph{some} ordinal number~$\xi$.
Statement~\eqref{eq:Borel-aleph1}
is an improvement.

Recall that $\Pi_{\aleph_1}\mathcal{A}
\equiv \bigcup_{\alpha\in\aleph_1} \Pi_\alpha \mathcal{A}$.
Since $\Pi \mathcal{A} \subseteq \Sigma_2 \mathcal{A}$,
we have
\begin{equation*}
\mathcal{T} \ \equiv\  \Sigma \mathcal{A} 
\ \subseteq\  \Pi_{\aleph_1}\mathcal{A}\ \subseteq\ \mathcal{B}.
\end{equation*}
Recall that~$\mathcal{B}$
is the family of all subsets of~$\N^\N$
that can be formed using countable unions and countable intersections
starting from~$\mathcal{T}$.
So to prove that $\Pi_{\aleph_1}\mathcal{A} = \mathcal{B}$
it suffices
to show that~$\Pi_{\aleph_1}\mathcal{A}$
is `closed' under countable unions and  intersections,
i.e.,
given $A_1,A_2,\dotsc \in \Pi_{\aleph_1}\mathcal{A}$
we must show that  $\bigcup_n A_n \in \Pi_{\aleph_1}\mathcal{A}$
and  $\bigcap_n A_n \in \Pi_{\aleph_1}\mathcal{A}$.
We will only prove that~$\bigcup_n A_n \in \Pi_{\aleph_1}\mathcal{A}$;
the proof of~$\bigcap_n A_n \in \Pi_{\aleph_1}\mathcal{A}$
is similar.

Let $A_n' \eqdf A_1 \cup \dotsb \cup A_n$
for each~$n\in\N$.
Then $\bigcup_n A_n' = \bigcup_n A_n$.
Since $\Pi_{\aleph_1}\mathcal{A}$
is a sublattice of~$\wp(\N^\N)$,
we get that 
$A_n'\in \Pi_{\aleph_1}\mathcal{A}
\equiv\bigcup_{\alpha\in\aleph_1}\Pi_{\alpha}\mathcal{A}$
for all~$n\in\N$.

Pick $\alpha_1,\alpha_2,\dotsc \in \aleph_1$
such that $A_n' \in \Pi_{\alpha_n} \mathcal{A}$
for all~$n\in\N$.
Let $\alpha\eqdf \bigvee_n \alpha_n$.
Then 
for all~$n\in\N$ we have $\alpha_n \leq \alpha$,
and  $\Pi_{\alpha_n} \mathcal{A} \subseteq \Pi_{\alpha}\mathcal{A}$,
and so $A_n'\in\Pi_{\alpha}\mathcal{A}$.
By definition of 
$ \Sigma(\Pi_{\alpha}\mathcal{A})$
we have $\bigcup_n A_n' \,\in\,
\Sigma(\Pi_{\alpha}\mathcal{A})\equiv\Sigma_{\alpha+1}\mathcal{A}$.
Now, note that since $\alpha_1,\alpha_2,\dotsc \in\aleph_1$,
also $\alpha=\bigvee_n \alpha_n \in \aleph_1$,
and so $\alpha+1 \in\aleph_1$.
Hence 
\begin{equation*}
\textstyle{\bigcup_n A_n \,=\,  \bigcup_n A_n' 
\ \in\  \Pi_{\alpha+1}\mathcal{A}
\,\subseteq\, \Pi_{\aleph_1}\mathcal{A}.}
\end{equation*}
So $\Pi_{\aleph_1} \mathcal{A}$
is closed under countable union.
Similarly, $\Pi_{\aleph_1} \mathcal{A}$
is closed under countable intersection.
It follows that  $\mathcal{B} = \Pi_{\aleph_1}\mathcal{A}$.
We have proven Statement~\eqref{eq:Borel-aleph1}.
\end{proof}

\begin{prop}
\label{P:bhier-surjective}
The map $\decode{-}^\mathcal{B}_{-}\colon \N^\N\times \Stump
\longrightarrow \mathcal{B}$
is surjective.
\end{prop}
\begin{proof}
Let $A\in\mathcal{B}$ be given.
We must find $f\in\Stump$
and $g\in\N^\N$ such that $A=\decode{g}^{\mathcal{B}}_f$.
By Lemma~\ref{L:Borel-aleph1}
we know that $A\,\in\,\Pi_{\aleph1} \mathcal{A}
\equiv \bigcup_{\alpha\in\aleph_1} \Pi_\alpha\mathcal{A}$.
Pick an~$\alpha\in \aleph_1$ with $A\in \Pi_\alpha \mathcal{A}$.
Since the map $\alpha[-]\colon \Stump \ra \aleph_1$
is a surjection,
there is an $f\in\Stump$
such that $\alpha[f] = \alpha$.
Since the map $\decode{-}^{\Pi}_{f}
\colon \N^\N \longrightarrow \Pi(\Sigma_{\alpha[f]} \mathcal{A})$
is a surjection
and $A\,\in\, \Pi_{\alpha} \mathcal{A}
\equiv \Pi_{\alpha[f]}\mathcal{A}
\subseteq \Pi(\Sigma_{\alpha[f]}\mathcal{A})$,
there is an~$h\in \N^\N$
such that $\decode{h}^{\Pi}_f = A$.

Now, let $g\in\N^\N$ be given by
$g^{[n]} = h^{[n]}$ for all~$n\in\N$
and $g(1)= 37$.
Then
\begin{equation*}
\decode{g}^{\mathcal{B}}_f
\ =\ 
\decode{g}^{\Pi}_f
\ =\ 
\decode{h}^{\Pi}_f
\ =\ 
A.
\end{equation*}
So we see that $\decode{-}^{\mathcal{B}}_{-}
\colon \N^\N\times\Stump\longrightarrow \mathcal{B}$
is surjective.
\end{proof}

\subsubsection{Outstanding Debt}$\,$\\
We take a small detour,
because in Example~\ref{E:not-convex}
we used a fact, Statement~\eqref{eq:not-convex-fact},
without proof,
and we are now in a position
to correct this situation.
\begin{nt}
Let $\mathbb{D}$ denote the \emph{Cantor set},
see~\cite[Examples~17.9c]{Willard70}.
\end{nt}
\begin{lem}
\label{L:bhier-disc-injection}
$\mathbb{D}$ is Borel negligble
subset of~$\R$,
and there is an injective map
\begin{equation*}
p\colon \N^\N \times \Stump \longrightarrow \mathbb{D}.
\end{equation*}
\end{lem}
\begin{proof}
We leave this to the reader.
\end{proof}

\begin{cor}
\label{C:bhier-embed}
We have the following situation.
\begin{equation*}
\xymatrix{
\mathcal{B} && 
\N^\N\times \Stump
  \ar[rr]^p 
  \ar[ll]_{\decode{-}^\mathcal{B}_{-}}
&& 
\mathbb{D}
}
\end{equation*}
The map~$\decode{-}^\mathcal{B}_{-}$
is surjective,
the map~$p$ is injective,
and~$\mathbb{D}$ is Borel negligible.
\end{cor}
\begin{proof}
Combine Lemma~\ref{L:bhier-disc-injection}
and Proposition~\ref{P:bhier-surjective}.
\end{proof}
\noindent
Note that Corollary~\ref{C:bhier-embed}
implies Statement~\eqref{eq:not-convex-fact}.
This ends our detour.

\subsubsection{Cataloguing Borel Subsets of~$\N^\N$}$\,$\\
We have encoded the Borel subsets
using elements of~$\N^\N$.
To prove the Borel Hierarchy Theorem
(see Statement~\eqref{eq:hierarchy-theorem})
we use the encoding to go one step further.
\begin{dfn}
For each~$f\in \Stump$,
define 
\begin{alignat*}{3}
U^\Pi_f
\ &=\ 
\{\ h\in\N^\N\colon\   h^{[1]} \in \decode{h^{[2]}}^\Pi_f\ \},&\qquad
U^\Sigma_f
\ &=\ 
\{\ h\in\N^\N\colon\ 
h^{[1]} \in \decode{h^{[2]}}^\Sigma_f \ \}.
\end{alignat*}
\end{dfn}
\noindent
One can think of the set~$\smash{U^\Pi_f}$
as a \emph{catalogue} of~$\smash{\Pi(\Sigma_{\alpha[f]}\mathcal{A}})$.\\
Indeed,
given $A\in \Pi(\Sigma_{\alpha[f]}\mathcal{A})$
with 
$A=\decode{g}^\Pi_f$
for some $g\in\N^\N$,
we have
\begin{equation*}
A \ =\ \{ \ h^{[1]}\colon \ h\in U^\Pi_f \ \text{ and } \ h^{[2]}=g\ \}.
\end{equation*}
The following lemma
might be the essential 
part of the Borel Hierarchy Theorem.
\begin{lem}
\label{L:bhier-lemma}
Let $f\in\Stump$ be given.
Then we have
\begin{equation}
\label{eq:peer-catalogue}
U_f^\Pi \,\in\, \Pi(\Sigma_{\alpha[f]}\mathcal{A})
\qquad\text{and}\qquad
U_f^\Sigma \,\in\, \Sigma(\Pi_{\alpha[f]}\mathcal{A}).
\end{equation}
\end{lem}
To give a proof of Lemma~\ref{L:bhier-lemma}
we need some notation and a lemma.
\begin{dfn}
\label{D:bhier-continuous}
Let $F\colon \N^\N\longrightarrow \N^\N$
be given.
\begin{enumerate}
\item
Given $A\subseteq \N^\N$,
let
$F^*(A)\,\eqdf\,\{\ g\in \N^\N\colon \ F(g)\in A\ \}$.
\item
We say that~$F$ is \keyword{continuous}
if~$F^*(B^m_n) \in \mathcal{T}$
for all~$m,n\in \N$.
\end{enumerate}
\end{dfn}
\begin{lem}
\label{L:bhier-continuous}
Let $F\colon \N^\N\ra \N^\N$ be a continuous function.\\
Let $\alpha$ be an ordinal number with~$\alpha>0$.
Then, for all~$A\subseteq\N^\N$,
\begin{alignat*}{3}
F^*(A) \,&\in\, \Pi_\alpha\mathcal{A}&
\qquad\text{when}\qquad A\,&\in\,\Pi_\alpha\mathcal{A}, \\
F^*(A) \,&\in\, \Sigma_\alpha\mathcal{A}&
\qquad\text{when}\qquad A\,&\in\,\Sigma_\alpha\mathcal{A}.
\end{alignat*}
\end{lem}
\begin{proof}
We leave this to the reader.
\end{proof}

\begin{proof}[Proof of Lemma~\ref{L:bhier-lemma}]
We prove Statement~\eqref{eq:peer-catalogue} using induction 
over~$f\in\Stump$.

\vspace{.6em}
\noindent
Let $f\in\Stump$ with $f(1)\neq1$ be given.
We must prove that
\begin{equation*}
U^\Pi_f \ \in\ \Pi \mathcal{A}
\qquad\text{and}\qquad
U^\Sigma_f \ \in\ \Sigma \mathcal{A}.
\end{equation*}
We will only prove that $U^\Sigma_f \in \Sigma \mathcal{A}$
since the proof of~$U^\Pi_f \in \Pi\mathcal{A}$ is similar.

Let $h\in \N^\N$ be given.
Note that the following are equivalent.
\begin{alignat*}{3}
h \ &\in\ U^\Sigma_f, \\
h^{[1]} \ &\in\ \decode{\,h^{[2]}\,}^\Sigma_f, \\
h^{[1]} \ &\in\ \textstyle{\bigcup_{n\geq 2}}\ 
           \decode{\,h^{[2]}(n)\,}^\mathcal{A}, \\
h^{[1]} \ &\in\ \decode{\,h^{[2]}(n)\,}^\mathcal{A}
\qquad &&\text{for some~$n\geq 2$},\\
h^{[1]} \ &\in\ \decode{k}^\mathcal{A}
\quad\text{and}\quad
h^{[2]}(n)\,=\,k
\qquad &&\text{for some~$n\geq 2$,\  $k\in\N$}.
\end{alignat*}
Let $F\colon \N^\N\ra\N^\N$
be given by $F(g)=g^{[1]}$.
Then $h^{[1]} \in \decode{k}^{\mathcal{A}}$
iff $h\in F^*(\decode{k}^\mathcal{A})$.
Further, note that $h^{[2]}(n)=k$ iff $h\in B^{k}_{\left<2,n\right>}$.
All in all, we get 
\begin{equation}
\label{eq:bhier-Usigmaf0}
U^\Sigma_f
\ =\ 
\bigcup_{k\in\N} \bigcup_{n\geq 2}
\ F^*(\,\decode{k}^\mathcal{A}\,) \ \cap\  B^k_{\left<2,n\right>}.
\end{equation}
By Equation~\eqref{eq:bhier-Usigmaf0}
to prove~$U^\Sigma_f \in\Sigma \mathcal{A}$
it suffices to show 
that~$F^*(\,\decode{k}^\mathcal{A}\,) \in \Sigma\mathcal{A}$.

Since it is not hard to see that~$F$ is continuous
(see Definition~\ref{D:bhier-continuous})
and $\decode{k}^\mathcal{A} \in \Sigma\mathcal{A}$
we get by Lemma~\ref{L:bhier-continuous} that
$F^*(\,\decode{k}^\mathcal{A}\,) \,\in\, \Sigma\mathcal{A}$.
Hence~$U^\Sigma_f \in \Sigma\mathcal{A}$.

\vspace{.6em}
\noindent
Recall that we are proving Statement~\eqref{eq:peer-catalogue}
using induction on~$f\in\Stump$.\\
Let $f\in\Stump$
with $f(1)=1$ be given.
Assume that for all~$n\in\N$,
\begin{equation*}
U^\Pi_{f^{[n]}} \,\in\, \Pi(\Sigma_{\alpha[f^{[n]}]}\mathcal{A})
\qquad\text{and}\qquad
U^\Sigma_{f^{[n]}} \,\in\, \Sigma(\Pi_{\alpha[f^{[n]}]}\mathcal{A}).
\end{equation*}
We must prove that 
$U^\Pi_f \,\in\, \Pi(\Sigma_{\alpha[f]} \mathcal{A})$
and $U^\Sigma_f \,\in\, \Sigma(\Pi_{\alpha[f]} \mathcal{A})$.
We will prove that
\begin{equation}
\label{eq:peer-catalogue-induction}
U^\Sigma_f \,\in\, \Sigma(\Pi_{\alpha[f]} \mathcal{A}),
\end{equation}
and we leave the proof of
$U^\Pi_f \,\in\, \Pi(\Sigma_{\alpha[f]} \mathcal{A})$
to the reader.

To proceed, we need some notation.
We will define a `pairing'
\begin{equation*}
P\colon \N^\N \times \N^\N\longrightarrow \N^\N.
\end{equation*}
Let $h_1,h_2\in\N^\N$ be given.
Define $P(h_1,h_2)\in\N^\N$
by: $P(h_1,h_2)(1) = 1$, and
\begin{equation*}
(\,P(h_1,h_2)\,)^{[1]} \,=\, h_1,
\qquad\text{and}\qquad
(\,P(h_1,h_2)\,)^{[2]} \,=\, h_2,
\end{equation*}
and $(\,P(h_1,h_2)\,)^{[n]}(m) = 1$
for all~$n,m\in\N$ with $n>2$.

Let $h\in\N^\N$ be given.
Note that the following are equivalent.
\begin{alignat*}{3}
h \ &\in\  U^\Sigma_f, \\
h^{[1]} \ &\in\ \decode{\,h^{[2]}\,}^\Sigma_f, \\
h^{[1]} \ &\in\ 
\textstyle{\bigcup_{n\in\N}}\ 
        \decode{\,h^{[2][n]}\,}^\Pi_{f^{n}}, \qquad\\
 h^{[1]} \ &\in\ 
        \decode{\,h^{[2][n]}\,}^\Pi_{f^{n}} 
\qquad&\text{for some $n\in\N$},\\
 P(h^{[1]},h^{[2][n]}) \ &\in\ 
        U^{\Pi}_{f^{[n]}} 
\qquad&\text{for some $n\in\N$}.
\end{alignat*}
Now,
for each~$n\in\N$,
let $F_n\colon \N^\N \ra\N^\N$
be given by, for $h\in\N^\N$,
\begin{equation*}
F_n(h) \ =\  P(h^{[1]},h^{[2][n]}).
\end{equation*}
Then using the notation of Definition~\ref{D:bhier-continuous}
we see that
\begin{equation}
\label{eq:expression-us}
U^\Sigma_f \ =\ 
\textstyle{\bigcup_{n\in\N}} \ F_n^*(\, U^\Pi_{f^{[n]}}\,).
\end{equation}
Recall that we must prove
that $U^\Sigma_f \in \Sigma(\Pi_{\alpha[f]}\mathcal{A})$.
It suffices to prove that
\begin{equation}
\label{eq:hierarchy-theorem-last-step}
F_n^*(\, U^\Pi_{f^{[n]}}\,) \ \in\  \Pi(\Sigma_{\alpha[f^{[n]}]}\mathcal{A}),
\end{equation}
by  Statement~\eqref{eq:expression-us},
and because we have
\begin{equation*}
\Pi(\Sigma_{\alpha[f^{[n]}]}\mathcal{A}) 
\ \subseteq \ \Sigma(\Pi_{\alpha[f]}\mathcal{A}).
\end{equation*}
By Lemma~\ref{L:bhier-continuous}
to prove
that Statement~\eqref{eq:hierarchy-theorem-last-step}
holds
it suffices to show that~$F_n$
is continuous; we leave this to the reader.
Hence we have proven Statement~\eqref{eq:peer-catalogue-induction}.

This completes the proof of Statement~\eqref{eq:peer-catalogue}.
\end{proof}

\subsubsection{Diagonalization}$\,$\\
With the catalogues $U^\Sigma_f$ and $U^\Pi_f$
at our disposal
we use Cantor's diagonal argument
to prove the Borel Hierarchy Theorem
(see Statement~\eqref{eq:hierarchy-theorem}).

\begin{dfn}
Let $f\in\Stump$ be given.
Let $D^\Pi_f$
and $D^\Sigma_f$ be give by
\begin{equation*}
D^\Pi_f \ \eqdf \ 
\{\  g\in\N^\N\colon\ g\notin \decode{g}^\Pi_f\ \},
\qquad
D^\Sigma_f \ \eqdf \ 
\{\  g\in\N^\N\colon\ g\notin \decode{g}^\Sigma_f\ \}.
\end{equation*}
\end{dfn}

\begin{thm}
\label{T:borel-hierarchy}
Let~$f\in \Stump$
be given.
Then we have
\begin{alignat*}{3}
D^\Pi_f \,&\in\, \Sigma(\Pi_{\alpha[f]} \mathcal{A})
\qquad&&\text{and}\qquad&
D^\Pi_f \,&\notin\, \Pi(\Sigma_{\alpha[f]} \mathcal{A}),\\
D^\Sigma_f \,&\in\, \Pi(\Sigma_{\alpha[f]} \mathcal{A})
\qquad&&\text{and}\qquad&
D^\Sigma_f \,&\notin\, \Sigma(\Pi_{\alpha[f]} \mathcal{A}).
\end{alignat*}
\end{thm}
\begin{proof}
Let $f\in\Stump$ be given.
We will only prove 
that
\begin{equation}
\label{eq:bhier}
D^\Pi_f \,\in\, \Sigma(\Pi_{\alpha[f]} \mathcal{A})
\qquad\text{and}\qquad
D^\Pi_f \,\notin\, \Pi(\Sigma_{\alpha[f]} \mathcal{A}),
\end{equation}
because there is a similar proof of
$D^\Sigma_f \,\in\, \Pi(\Sigma_{\alpha[f]} \mathcal{A})$
and
$D^\Sigma_f \,\notin\, \Sigma(\Pi_{\alpha[f]} \mathcal{A})$.

Let us first prove that 
$D^\Pi_f\notin  \Pi(\Sigma_{\alpha[f]}\mathcal{A})$.
So assume 
$D^\Pi_f\in \Pi(\Sigma_{\alpha[f]}\mathcal{A})$
in order to reach a contradiction.
Since the map $\decode{-}^\Pi_f\colon \N^\N \longrightarrow \Pi(\Sigma_{\alpha[f]}\mathcal{A})$
is surjective,
there is a~$g_\delta \in \N^\N$
with $\decode{g_\delta}^\Pi_f = D^\Pi_f$.
Then
by definition of~$D^\Pi_f$,
\begin{equation}
\label{eq:diagonalization}
g_\delta \,\in\, D^\Pi_f
\quad\iff\quad
g_\delta \,\notin\, \decode{g_\delta}^\Pi_f = D^\Pi_f.
\end{equation}
Statement~\eqref{eq:diagonalization}
leads to a contradiction.
Hence we conclude that $D^\Pi_f\notin \Pi(\Sigma_{\alpha[f]}\mathcal{A})$.

Let us prove that $D^\Pi_f \in \Sigma (\Pi_{\alpha[f]}\mathcal{A})$.
Let~$\Delta\colon \N^\N \ra \N^\N$
be given by, for $g\in \N^\N$,
\begin{equation*}
\Delta(g) \ =\  P(g,g).
\end{equation*}
Then $\Delta$ is continuous
(see Definition~\ref{D:bhier-continuous})
and we have, for $g\in \N^\N$,
\begin{equation*}
g\,\in\, D^\Pi_f 
\quad\iff\quad
\Delta(g) \,\notin\, U^\Pi_f
\quad\iff\quad
g\,\in\,\Delta^*(\,\N^\N \backslash U^\Pi_f \,).
\end{equation*}
So we see that $D^\Pi_f = \Delta^*(\,\N^\N\backslash U_f^\Pi\,)$.
Recall that we must prove that
\begin{equation*}
D^\Pi_f \ \in\  \Sigma (\Pi_{\alpha[f]}\mathcal{A}).
\end{equation*}
Since~$\Delta$
is continuous
it suffices to show that 
$\N^\N\backslash U_f^\Pi \ \in\  \Sigma (\Pi_{\alpha[f]}\mathcal{A})$
by Lemma~\ref{L:bhier-continuous}.
We already know that $U^\Pi_f \in \Pi(\Sigma_{\alpha[f]} \mathcal{A})$,
so 
by Statement~\eqref{eq:bhier-complement}
we get
\begin{equation*}
\N^\N\backslash U^\Pi_f \ \in\  \Sigma(\Pi_{\alpha[f]} \mathcal{A}).
\end{equation*}
So
$D^\Pi_f \in \Sigma (\Pi_{\alpha[f]}\mathcal{A})$.
Hence we have proven Statement~\eqref{eq:bhier}.
\end{proof}

\begin{cor}
\label{C:borel-hierarchy-theorem}
The Borel Hierarchy Theorem holds,
see Statement~\eqref{eq:hierarchy-theorem}.
\end{cor}

\begin{prop}
\label{P:no-countable-collapse}
Using the terminology from Subsection~\ref{SS:hierarchy},\\
if $\alpha$ is an ordinal number
such that the hierarchy collapsed at~$\Pi_\alpha \mathcal{A}$
or at $\Sigma_\alpha \mathcal{A} $,\\
then $\alpha$ must be uncountable.
\end{prop}
\begin{proof}
Let $\alpha$ be an ordinal number.
We will only prove that
$\alpha$ must be uncountable when
 the hierarchy collapsed at~$\Pi_\alpha\mathcal{A}$,
because the
proof that~$\alpha$
must be uncountable when the hierarchy collapsed
at~$\Sigma_\alpha\mathcal{A}$ is similar.

Assume that the hierarchy hierarchy collapsed at~$\Pi_\alpha\mathcal{A}$,
that is,
\begin{equation*}
\Pi(\Pi_\alpha \mathcal{A}) 
\ =\ 
\Pi_\alpha \mathcal{A}
\ =\ 
\Sigma(\Pi_\alpha \mathcal{A}).
\end{equation*}
Assume that $\alpha$ is countable in order to reach a contradiction.\\
By the Borel Hierarchy Theorem,
see Corollary~\eqref{C:borel-hierarchy-theorem},
there is an 
\begin{equation*}
S\in \Sigma(\Pi_\alpha \mathcal{A})
\qquad\text{with}\qquad
S\notin \Pi(\Sigma_\alpha \mathcal{A}).
\end{equation*}
However,
we have the following inclusion
\begin{equation*}
\Sigma(\Pi_\alpha \mathcal{A}) \ = \ 
\Pi_\alpha \mathcal{A} \ \subseteq\ \Pi(\Sigma_\alpha \mathcal{A}).
\end{equation*}
So since $S\in \Sigma(\Pi_\alpha\mathcal{A})$,
we get $S\in \Pi(\Sigma_\alpha \mathcal{A})$,
which contradicts $S\notin \Pi(\Sigma_\alpha\mathcal{A})$.

Hence $\alpha$ is not countable.
So $\alpha$ must be uncountable.
\end{proof}

\clearpage
%
%                  HIERARCHY OF EXTENSIONS
%
\subsection{The Hierarchy of Extensions}$\,$\\
\label{SS:hierarchy-abstract}%
To get the smallest complete extension
of a valuation~$\varphi$
with respect to some~$V$
(when it exists)
we can make a hierarchy 
of extensions of~$\varphi$,
see Subsection~\ref{SS:hierarchy}:
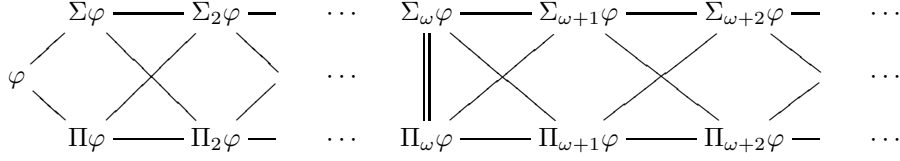
\begin{equation*}
\xymatrix @=10pt {
& \Sigma\varphi \ar @{-} [rr] \ar @{-} [rrdd]
&& \Sigma_2 \varphi  \ar @{-} [r]\ar @{-} [rd]
&& \dotsb
& \Sigma_{\omega} \varphi \ar @{-} [rr] \ar @{-} [rrdd]
                          \ar @{=} [dd]
&& \Sigma_{\omega+1}\varphi \ar @{-} [rr] \ar @{-} [rrdd]
&& \Sigma_{\omega+2} \varphi  \ar @{-} [r]\ar @{-} [rd]
&& \dotsb
\\  
\varphi \ar @{-} [ru] \ar @{-} [rd] 
&&&&&\dotsb
&&
&&&&&\dotsb\\
& \Pi\varphi \ar @{-} [rr] \ar @{-} [rruu]
&& \Pi_{2} \varphi \ar @{-} [r] \ar @{-} [ru]
&& \dotsb
& \Pi_{\omega}\varphi \ar @{-} [rr] \ar @{-} [rruu]
&& \Pi_{\omega+1}\varphi \ar @{-} [rr] \ar @{-} [rruu]
&& \Pi_{\omega+2} \varphi \ar @{-} [r] \ar @{-} [ru]
&& \dotsb
}
\end{equation*}
We have seen 
that if~$\varphi$
has a complete extension,
then~$\varphi$ also has a smallest complete extension $\overline\varphi$,
and that $\ol\varphi = \Pi_\xi \varphi$
for some ordinal number~$\xi$,
see Proposition~\ref{P:smallest-complete-extension}.

Even if we do not know whether~$\varphi$
has a complete extension,
we can still try to make the hierarchy,
and this is what we are going to do in this subsection.

It is possible that the making of the hierarchy
is hindered at some point,
e.g., if $\Sigma_2\varphi$
is not $\Pi$-extendible,
then we can not define~$\Pi_3 \varphi = \Pi(\Sigma_2 \varphi)$.

If we can define the hierarchy up to $\Pi_\alpha \varphi$ unhindered
we will say that
\begin{equation*}
\text{$\varphi$ is $\Pi_\alpha$-extendible}.
\end{equation*}
We will prove
that $\varphi$ has a complete extension
iff~$\varphi$ is $\Pi_{\aleph_1}$-extendible.
Moreover,
the valuation
$\Pi_{\aleph_1}\varphi$ will be the smallest complete extension of~$\varphi$
(see Corollary~\ref{C:aleph1}).

Let us begin by giving a formal definition
of the hierarchy
and ``$\Pi_\alpha$-extendible''.
\begin{dfn}
\label{P:hier}
Let $\vs{V}{L}{\varphi}{E}$
be a valuation system.\\
We are going to define the following
statements and valuation systems.
\begin{enumerate}
\item \label{hier:first-cond}
For each ordinal number~$\alpha$,
statements
\begin{equation*}
\text{``\keyword{$\varphi$ is $\Pi_\alpha$-extendible}''}
\qquad\text{ and }\qquad
\text{``\keyword{$\varphi$ is $\Sigma_\alpha$-extendible}''}.
\end{equation*}
\item
For each~$\alpha$ such that~$\varphi$ 
is $\Pi_\alpha$-extendible,
a valuation system 
\begin{equation*}
\vs{V}{\Pi_\alpha L}{\Pi_\alpha \varphi}{E}.
\end{equation*}

\item
For each~$\alpha$ such that~$\varphi$ is $\Sigma_\alpha$-extendible,
a valuation system 
\begin{equation*}
\vs{V}{\Sigma_\alpha L}{\Sigma_\alpha \varphi}{E}.
\end{equation*}
\end{enumerate}
We will define them in such a way
that we get a hierarchy of the following shape.
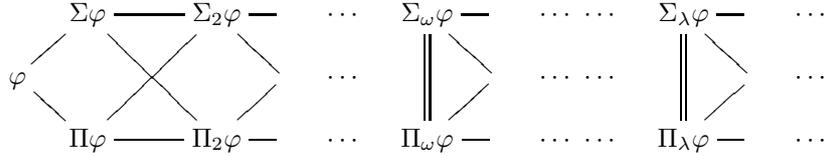
\begin{equation*}
\xymatrix @=10pt {
& \Sigma\varphi \ar @{-} [rr] \ar @{-} [rrdd]
&& \Sigma_2 \varphi  \ar @{-} [r]\ar @{-} [rd]
&& \dotsb
& \Sigma_{\omega} \varphi \ar @{-} [r] \ar @{-} [rd]
                          \ar @{=} [dd]
&& \dotsb\ \dotsb
& \Sigma_{\lambda} \varphi \ar @{-} [r] \ar @{-} [rd]
                          \ar @{=} [dd]
&& \dotsb
\\  
\varphi \ar @{-} [ru] \ar @{-} [rd] 
&&&&&\dotsb
&&&\dotsb\ \dotsb
&&&\dotsb\\
& \Pi\varphi \ar @{-} [rr] \ar @{-} [rruu]
&& \Pi_{2} \varphi \ar @{-} [r] \ar @{-} [ru]
&& \dotsb
& \Pi_{\omega}\varphi \ar @{-} [r] \ar @{-} [ru]
&& \dotsb\ \dotsb
& \Pi_{\lambda}\varphi \ar @{-} [r] \ar @{-} [ru]
&& \dotsb
}
\end{equation*}
More precisely, the following statements will be true.
\begin{enumerate}[label=(\Roman*)]
\item\label{req:hierarchy-A}
For all $\beta < \gamma$,
if  $\varphi$ is $\Pi_\gamma$-extendible,
then $\varphi$ is both $\Pi_\beta$- and $\Sigma_\beta$-extendible,
and the map $\Pi_\gamma\varphi$ extends 
both $\Pi_\beta\varphi$ and $\Sigma_\beta\varphi$.
\item\label{req:hierarchy-B}
For all $\beta < \gamma$,
if  $\varphi$ is $\Sigma_\gamma$-extendible,
then $\varphi$ is both $\Pi_\beta$- and $\Sigma_\beta$-extendible,
and the map $\Sigma_\gamma\varphi$ extends 
both $\Pi_\beta\varphi$ and $\Sigma_\beta\varphi$.
\item\label{req:hierarchy-C}
Let~$\lambda$
be zero or a a limit ordinal.
Then $\varphi$ is $\Pi_\lambda$-extendible
if and only if $\varphi$ is $\Sigma_\lambda$-extendible.
Furthermore,
if $\varphi$ is $\Pi_\lambda$-extendible,
then $\Pi_\gamma\varphi = \Sigma_\gamma\varphi$.
\end{enumerate}
Now, we define the aforementioned
statements and valuation systems
by recursion over the ordinal number
using the following clauses.
\begin{enumerate}
\item
$\varphi$ is {$\Pi_0$-extendible}
and $\varphi$ is {$\Sigma_0$-extendible}.
Moreover,
\begin{equation*}
\Pi_0 L = L\qquad 
\Sigma_0 L = L \qquad 
\Pi_0 \varphi = \varphi \qquad
\Sigma_0 \varphi = \varphi.
\end{equation*}

\item
Let $\alpha$ be an ordinal number.
Then we have
\begin{equation*}
\text{$\varphi$ is {$\Pi_{\alpha+1}$-extendible}}
\quad\iff\quad
\left[\quad
\begin{alignedat}{3}
&\varphi \text{ is $\Sigma_{\alpha}$-extendible, and}\\
&\Sigma_\alpha\varphi  \text{ is $\Pi$-extendible.}
\end{alignedat}
\right.
\end{equation*}
Moreover, if $\varphi$ is $\Pi_{\alpha+1}$-extendible,
then 
\begin{equation*}
\Pi_{\alpha+1} L \ =\  \Pi(\Sigma_\alpha L)\qquad\text{and}\qquad
\Pi_{\alpha+1} \varphi \ =\  \Pi(\Sigma_\alpha \varphi),
\end{equation*}
where 
 $\vs{V}{\Pi(\Sigma_\alpha L)}{\Pi(\Sigma_\alpha\varphi)}{E}$
is the valuation system from Definition~\ref{D:Pi-extendible}.

\item
Let $\alpha$ be an ordinal number.
Then we have
\begin{equation*}
\text{$\varphi$ is {$\Sigma_{\alpha+1}$-extendible}}
\quad\iff\quad
\left[\quad
\begin{alignedat}{3}
&\text{$\varphi$ is $\Pi_{\alpha}$-extendible, and}\\
&\text{$\Pi_\alpha \varphi$ is $\Sigma$-extendible.}
\end{alignedat}
\right.
\end{equation*}
Moreover,
if $\varphi$ is $\Sigma_{\alpha+1}$-extendible,
then 
\begin{equation*}
\Sigma_{\alpha+1} L \ =\  \Sigma(\Pi_\alpha L)\qquad\text{and}\qquad
\Sigma_{\alpha+1} \varphi \ =\  \Sigma(\Pi_\alpha \varphi).
\end{equation*}

\item
Let~$\lambda$ be a limit ordinal.
Then
we have 
\begin{equation*}
\text{$\varphi$ is {$\Pi_\lambda$-extendible}}
\quad\iff\quad
\left[\quad
\begin{alignedat}{3}
&\text{$\varphi$ is $\Pi_\alpha$-extendible,}\\
&\text{for every~$\alpha \in \lambda$}.
\end{alignedat}
\right.
\end{equation*}
Moreover,
if $\varphi$ is $\Pi_\lambda$-extendible,
then 
\begin{equation*}
\Pi_\lambda L \ =\  \textstyle{\bigcup_{\alpha \in \lambda}\, \Pi_\alpha L}
\qquad\text{and}\qquad
\Pi_\lambda \varphi (c) \ =\  \Pi_\beta \varphi(c),
\end{equation*}
where $\beta \in \lambda$ and $c\in \Pi_\beta L$.

\item
\label{hier:last-cond}
Let~$\lambda$ be a limit ordinal.
Then we have
\begin{equation*}
\text{$\varphi$ is {$\Sigma_\lambda$-extendible}}
\quad\iff\quad
\left[\quad
\begin{alignedat}{3}
&\text{$\varphi$ is $\Sigma_\alpha$-extendible,}\\
&\text{for every~$\alpha \in \lambda$.}
\end{alignedat}
\right.
\end{equation*}
Moreover,
if $\varphi$ is $\Sigma_\lambda$-extendible,
then 
\begin{equation*}
\Sigma_\lambda L \ =\  
\textstyle{\bigcup_{\alpha \in \lambda} \,\Sigma_\alpha L}
\qquad\text{and}\qquad
\Sigma_\lambda \varphi (c) \ =\  \Sigma_\beta \varphi(c),
\end{equation*}
where $\beta \in \lambda$ and $c\in \Sigma_\beta L$.
\end{enumerate}
\end{dfn}
%
%                  COLLAPSE
%
\begin{dfn}
\label{D:collapse}
Let $\vs{V}{L}\varphi{E}$ be a valuation system.
\begin{enumerate}
\item
\label{D:collapse-1}
We say the \keyword{hierarchy has collapsed at~$Q$} if
\begin{equation*}
\left[\quad
\begin{minipage}{.7\columnwidth}
$\varphi$ is $\Pi_{\alpha+1}$-extendible
 and $\Sigma_{\alpha+1}$-extendible,
and 
\begin{equation*}
\Pi(Q)\ =\ Q\ =\ \Sigma(Q),
\end{equation*}
where 
$Q=\Pi_\alpha \varphi$
or $Q=\Sigma_\alpha \varphi$.
\end{minipage}
\right.
\end{equation*}

\item
\label{D:collapse-2}
We say the \keyword{hierarchy collapses}
if the hierarchy has collapsed at some~$Q$.
\end{enumerate}
\end{dfn}
\noindent
We will prove that if the hierarchy of a valuation~$\varphi$ collapses
then~$\varphi$ has a complete extension
(see Lemma~\ref{L:complete}).
After that,
we will prove the converse,
namely,
if $\varphi$ has a complete extension,
then then the hierarchy of~$\varphi$ collapses
(see Proposition~\ref{P:comp-minimal}).

\begin{lem}
\label{L:collapse}
Let $\vs{V}{L}\varphi{E}$
be a valuation system.
\begin{enumerate}
\item
\label{L:collapse-1}
Let $\alpha<\beta$ be ordinal numbers.\\
Suppose that the hierarchy has collapsed at~$\Pi_\alpha \varphi$.\\
Then $\varphi$ is $\Pi_\beta$-extendible and $\Sigma_\beta$-extendible,
and  
\begin{equation*}
\Pi_\beta \varphi \,=\, \Pi_\alpha \varphi \,=\, \Sigma_\beta \varphi.
\end{equation*}
\item
\label{L:collapse-2}
If the hierarchy collapses
at $\Pi_\alpha \varphi$ and at $\Sigma_\alpha\varphi$
for some~$\alpha$,
then $\Pi_\alpha\varphi = \Sigma_\alpha\varphi$.

\item
\label{L:collapse-3}
Suppose the hierarchy has collapsed at~$Q_1$ and at~$Q_2$.
Then $Q_1 = Q_2$.
\item
\label{L:collapse-4}
The hierarchy collapses
if and only if 
it has collapsed at $\Pi_\alpha \varphi$
for some~$\alpha$.
\end{enumerate}
\end{lem}
\begin{proof}
We leave this to the reader.
\end{proof}
%
%                  DEFINITION OF EXTENDIBLE
%
\begin{dfn}
\label{D:extendible}
Let $\vs{V}{L}\varphi{E}$ be a valuation system.
\begin{enumerate}
\item
\label{D:extendible-1}
We say that~$\varphi$ (or  $\vs{V}{L}\varphi{E}$) is \keyword{extendible}
if the hierarchy collapses.

\item
\label{D:extendible-2}
Suppose that $\varphi$ is extendible.
Then there is precisely one valuation
at which the hierarchy has collapsed
(see Lemma~\ref{L:collapse}\ref{L:collapse-3});
we denote it by
\begin{equation*}
\vs{V}{\overline L}{\overline\varphi}{E}.
\end{equation*}
\end{enumerate}
\end{dfn}
%
%                  OVERLINE VARPHI IS COMPLETE
%
\begin{lem}
\label{L:complete}
Let $\vs{V}{L}\varphi{E}$ be an extendible valuation system.
Then
\begin{equation*}
\vs{V}{\overline L}{\overline\varphi}{E}
\text{ is complete.}
\end{equation*}
\end{lem}
\begin{proof}
To prove that~$\overline\varphi$ is complete,
it suffices to show that~$\overline\varphi$
is $\Pi$-complete and $\Sigma$-complete.
We know that $\overline\varphi = \Pi\overline\varphi$
(since the hierarchy has collapsed at~$\overline\varphi$,
see Definition~\ref{D:extendible}\ref{D:extendible-2}
and Definition~\ref{D:collapse}\ref{D:collapse-1}),
and that $\Pi\overline\varphi$ is $\Pi$-complete
(see Lemma~\ref{L:Pi-complete}).
Hence $\overline\varphi$ is $\Pi$-complete.
Similarly,
$\overline\varphi$ is $\Sigma$-complete.
So $\ol\varphi$ is complete.
\end{proof}
%
%                  REMARK ON EXTENDIBLE VERSUS COMPLETE
%
\begin{rem}
\label{R:extendible-complete}
Let $\vs{V}{L}\varphi{E}$ be a valuation system.
\begin{enumerate}
\item
Note that if~$\varphi$ is complete with respect to~$V$
(see Definition~\ref{D:system-complete}),\\
then~$\varphi$ is extendible, and $\ol\varphi=\varphi$.
\item
On the other hand,
if~$\varphi$ is extendible and $\ol\varphi=\varphi$,\\
then~$\varphi$ is complete
with respect to~$V$
(see Lemma~\ref{L:complete}).
\end{enumerate}
\end{rem}
%
%                  EXTENDIBLE IFF PI_ALEPH_1 EXTENDIBLE
%
\begin{lem}
\label{L:aleph1}
Let $\vs{V}{L}\varphi{E}$ be a valuation system.\\
If $\varphi$ is $\Pi_{\aleph_1}$-extendible,
then the hierarchy has collapsed at $\Pi_{\aleph_1}\varphi$.
\end{lem}
\begin{proof}
Suppose that $\varphi$ is $\Pi_{\aleph_1}$-extendible.
We need to prove that the hierarchy has collapsed at~$\Pi_{\aleph_1}\varphi$.
For this, we must show that
(see Definition~\ref{D:collapse}\ref{D:collapse-1}),
\begin{equation*}
\Pi(\Pi_{\aleph_1}\varphi)
= \Pi_{\aleph_1}\varphi = \Sigma(\Pi_{\aleph_1}\varphi).
\end{equation*}
Let $a_1 \geq a_2 \geq \dotsb$ 
be a $\Pi_{\aleph_1}\varphi$-convergent sequence.
In order to show that $\Pi(\Pi_{\aleph_1}\varphi)=\Pi_{\aleph_1}\varphi$,
it suffices to prove that $\bw_n a_n \in \Pi_{\aleph_1}\varphi$.

Since $\aleph_1$ is a limit ordinal,
we know that $\Pi_{\aleph_1} L = \bigcup_{\alpha < \aleph_1} \Pi_\alpha L$
(see Defintion~\ref{P:hier}).
Define for each~$n\in\N$
an ordinal number~$\alpha(n)$ by
\begin{equation*}
\alpha(n)\ \eqdf\ 
\min\,\{\, \alpha<\aleph_1\colon\, a_n \in \Pi_\alpha L\,\}.
\end{equation*}
Now,
the set $\{\,\alpha(1),\,\alpha(2),\,\dotsc\,\}$ of
ordinals has a supremum,
\begin{equation*}
\xi\,\eqdf\,\bv_n\alpha(n)\,\equiv\,\textstyle{\bigcup_n}\alpha(n).
\end{equation*}
Since~$\aleph_1$ is the smallest uncountable ordinal,
and $\alpha(n)<\aleph_1$,
we know that all~$\alpha(n)$ are countable.
Hence~$\xi$ is countable as well,
and so $\xi<\aleph_1$.

Now, we have $a_n \in \Pi_{\xi} L$
for all~$n\in \N$.  Hence
\begin{equation*}
\bw_n a_n \,\in\, \Pi(\Pi_{\xi} L) 
\ =\ \Pi_{\xi+1} L
\ \subseteq\ \Pi_{\aleph_1} L.
\end{equation*}
So we see that $\Pi(\Pi_{\aleph_1}\varphi)=\Pi_{\aleph_1}\varphi$.
Similarly, $\Sigma(\Pi_{\aleph_1}\varphi)=\Pi_{\aleph_1}\varphi$.\\
Hence the hierarchy has collapsed at~$\Pi_{\aleph_1}\varphi$.
\end{proof}
\begin{cor}
\label{C:aleph1}
Let $\vs{V}{L}\varphi{E}$ be a valuation system. Then
\begin{equation*}
\varphi\text{ is extendible }
\quad\iff\quad
\varphi\text{ is $\Pi_{\aleph_1}$-extendible.}
\end{equation*}
Moreover,
if $\varphi$ is extendible,
then $\overline\varphi = \Pi_{\aleph_1}\varphi$.
\end{cor}
\begin{proof}
Assume $\varphi$ is extendible
in order to show that $\varphi$ is $\Pi_{\aleph_1}$-extendible.
Then we know that the hierarchy collapses 
(see Definition~\ref{D:collapse}\ref{D:collapse-2}).
So it collapsed at some~$\Pi_\alpha \varphi$
(see Lemma~\ref{L:collapse}\ref{L:collapse-4}).
Pick an ordinal number~$\beta$ with $\beta > \alpha$
and $\beta > \aleph_1$.
Then $\varphi$ is $\Pi_\beta$-extendible
by Lemma~\ref{L:collapse}\ref{L:collapse-1}.
But $\aleph_1 <\beta$,
so $\varphi$ is also $\Pi_{\aleph_1}$-extendible.

Suppose $\varphi$ is $\Pi_{\aleph_1}$-extendible.
Then the hierarchy has collapsed at~$\Pi_{\aleph_1}\varphi$
by Lemma~\ref{L:aleph1}.
Hence 
$\varphi$ is extendible
and  $\overline\varphi = \Pi_{\aleph_1}\varphi$
(see Definition~\ref{D:extendible}\ref{D:extendible-2}).
\end{proof}

%
%                  Pi_\alpha is  MONOTONOUS
%
\begin{lem}
\label{L:alpha-monotonous}
Let $\vs{V}{L}\varphi{E}$
and $\vs{V}{C}\psi{E}$
 be two valuation systems.\\
Assume  that
$\psi$ extends $\varphi$.
Then for every ordinal number $\alpha$, we have
\begin{alignat*}{3}
\text{$\psi$ is $\Pi_\alpha$-extendible}
\quad&\implies\quad
\text{$\varphi$ is $\Pi_\alpha$-extendible}
\quad&&\text{and}\quad
\text{$\Pi_\alpha\psi$ }&&\text{extends $\Pi_\alpha\varphi$} \\
\text{$\psi$ is $\Sigma_\alpha$-extendible}
\quad&\implies\quad
\text{$\varphi$ is $\Sigma_\alpha$-extendible}
\quad&&\text{and}\quad
\text{$\Sigma_\alpha\psi$ }&&\text{extends $\Sigma_\alpha\varphi$.}
\end{alignat*}
\end{lem}
\begin{proof}
We prove this lemma using induction on~$\alpha$.

\noindent\emph{(Zero)} For~$\alpha=0$, the proposition is trivial.

\noindent\emph{(Successor)} 
Let~$\alpha$ be an ordinal number such that
if $\psi$ is $\Sigma_\alpha$-extendible,
then
$\varphi$ is $\Sigma_\alpha$-extendible 
and $\Sigma_\alpha\psi$ extends $\Sigma_\alpha\varphi$.
Suppose~$\psi$
is $\Pi_{\alpha+1}$-extendible.
We prove
\begin{equation}
\label{eq:P:extension-1}
\text{ $\varphi$ is $\Pi_{\alpha+1}$-extendible}
\qquad\text{and}\qquad
\text{ $\Pi_{\alpha+1}\psi$ extends $\Pi_{\alpha+1}\varphi$.}
\end{equation}
Since $\psi$ is $\Pi_{\alpha+1}$-extendible,
we know that
(see Definition~\ref{P:hier}),
\begin{equation*}
\text{ $\psi$ is $\Sigma_\alpha$-extendible}
\qquad\text{and}\qquad
\text{$\Sigma_\alpha\psi$ is~$\Pi$-extendible}.
\end{equation*}
By assumption, 
the former implies that $\varphi$ is $\Sigma_\alpha$-extendible
and $\Sigma_\alpha\psi$ extends $\Pi_\alpha\psi$;
by Lemma~\ref{L:Pi-monotonous},
the latter implies $\Sigma_\alpha\varphi$
is $\Pi$-extendible
and $\Pi(\Sigma_\alpha\psi)$ extends $\Pi(\Sigma_\alpha \varphi)$.
In other words,
we have proven Statement~\eqref{eq:P:extension-1}.

Let~$\alpha$ be an ordinal number such that
if $\psi$ is $\Pi_\alpha$-extendible,
then we have that
$\varphi$ is $\Pi_\alpha$-extendible 
and that $\Pi_\alpha\psi$ extends $\Pi_\alpha\varphi$.
Suppose that~$\psi$
is $\Sigma_{\alpha+1}$-extendible.
By a similar reasoning
as before one can prove that
\begin{equation*}
\text{$\varphi$ is $\Sigma_{\alpha+1}$-extendible}
\qquad\text{and}\qquad
\text{$\Sigma_{\alpha+1}\psi$ extends $\Sigma_{\alpha+1}\varphi$.}
\end{equation*}

\noindent\emph{(Limit)}
Let $\lambda$ be a limit ordinal such that
for all $\alpha<\lambda$, we have
\begin{equation*}
\text{$\psi$ is $\Pi_\alpha$-extendible}
\quad\implies\quad
\text{$\varphi$ is $\Pi_\alpha$-extendible }
\quad\text{and}\quad
\text{$\Pi_\alpha\psi$ extends $\Pi_\alpha\varphi$.}
\end{equation*}
Further, 
assume $\psi$ is~$\Pi_\lambda$-extendible
in order to prove that
\begin{equation}
\label{eq:P:extension-2}
\text{ $\varphi$ is $\Pi_{\lambda}$-extendible}
\qquad\text{and}\qquad
\text{ $\Pi_{\lambda}\psi$ extends $\Pi_{\lambda}\varphi$.}
\end{equation}
Let $\alpha <\lambda$ be given.
Since $\psi$ is $\Pi_\lambda$-extendible,
we know that $\psi$ is $\Pi_\alpha$-extendible.
So by assumption, 
$\varphi$ is $\Pi_\alpha$-extendible,
and $\Pi_\alpha \psi$ extends $\Pi_\alpha \varphi$.

So we see that $\varphi$ is $\Pi_\lambda$-extendible.
Further
(since $\Pi_\lambda\psi$ extends $\Pi_\alpha \psi$),
we see that $\Pi_\lambda \psi$ extends all $\Pi_\alpha\varphi$.
Hence $\Pi_\lambda \psi$ extends $\Pi_\lambda \varphi$.
So we have proven~\eqref{eq:P:extension-2}.
\end{proof}
%
%                  OVERLINE PHI IS MONOTONOUS
%
\begin{prop}
\label{P:comp-monotonous}
Let $\vs{V}{L}\varphi{E}$
and $\vs{V}{C}\psi{E}$
 be valuation systems.\\
Assume that
$\psi$ extends $\varphi$,
and that 
$\psi$ is extendible.
Then 
\begin{alignat*}{3}
\text{$\varphi$ is extendible}
\qquad&&\text{and}\qquad
\text{$\overline\psi$ }&&\text{extends $\overline\varphi$}.
\end{alignat*}
\end{prop}
\begin{proof}
By Corollary~\ref{C:aleph1},
we get the conclusion from
Lemma~\ref{L:alpha-monotonous}
with $\alpha=\aleph_1$.
\end{proof}
%
%                  OVERLINE PHI IS MINIMAL
%
\begin{prop}
\label{P:comp-minimal}
Let $\vs{V}{L}\varphi{E}$
and $\vs{V}{C}\psi{E}$
 be valuation systems.\\
Assume  that
$\psi$ extends $\varphi$
and that $\psi$ is complete
with respect to~$V$.
Then 
\begin{alignat*}{3}
\text{$\varphi$ is extendible}
\qquad&&\text{and}\qquad
\text{$\psi$ }&&\text{extends $\overline\varphi$}.
\end{alignat*}
(So, loosely speaking,
$\ol\varphi$ is the smallest complete extension of~$\varphi$.)
\end{prop}
\begin{proof}
Since~$\psi$ is complete,
$\psi$ is clearly extendible 
and $\overline\psi = \psi$
(see Remark~\ref{R:extendible-complete}).
Hence $\varphi$ is extendible
and $\psi=\overline\psi$ extends $\overline\varphi$
by Proposition~\ref{P:comp-monotonous}.
\end{proof}
%
%                  REMARK ON EXTENDIBLE
%
\begin{rem}
Let $\vs{V}{L}\varphi{E}$ be a valuation system.\\
By Lemma~\ref{L:complete} and Proposition~\ref{P:comp-minimal}
we see that
\begin{equation*}
\text{$\varphi$ is extendible}
\quad\iff\quad
\text{$\varphi$
has an complete extension}.
\end{equation*}
Hence the name ``extendible''.
\end{rem}
 }
\clearpage
{ \section{Closedness of the Completion under Operations}
\label{S:closedness}
\noindent
We have seen how we can
obtain the Lebesgue measure and the Lebesgue integral
as the (convexification of) the completion
of relatively simple valuations systems,
\begin{equation*}
\vsSA\qquad\text{and}\qquad\vsSF.
\end{equation*}
It is now time to derive some simple facts
about the completion.
In this section we will prove
statements of the following form.
\begin{enumerate}
\item
If $A,B\in\overline\SA$,
then $A\backslash B \in\ol\SA$
(see Example~\ref{E:ring-ext-2}).

\item 
If $f,g\in \ol\SF \cap \R^\R$,
then $f+g \in \ol\SF$
(see Example~\ref{E:riesz-space-ext}).
\end{enumerate}

%
%                  SIGMA-PRESERVING
%
\begin{dfn}
\label{D:sigma-preserving}
Let $P$ and $Q$ be posets.
Let $f\colon S\ra Q$ be a map, where $S\subseteq P$.
We say $f$ is \keyword{$\sigma$-preserving
with respect to~$P$}
provided that
\begin{enumerate}
\item
if 
$\bw_n{a_n}$ exists
(in~$P$)
for
$a_1 \geq a_2 \geq \dotsb$
from~$S$,
and if $\bw_n a_n\in S$,
then 
\begin{equation*}
f(\bw_n a_n) = \bw_n f(a_n);
\end{equation*}

\item
if 
$\bv_n{b_n}$ exists
for
$b_1 \leq b_2 \leq \dotsb$
from~$S$,
and if~$\bv_n b_n\in S$,
then 
\begin{equation*}
f(\bv_n b_n) = \bv_n f(b_n).
\end{equation*}
\end{enumerate}
Let $P$ and~$Q$ be posets. Let $f\colon P \ra Q$ be a map.
We say $f$ is \keyword{$\sigma$-preserving}
provided that~$f$ is $\sigma$-preserving with respect to~$P$.
\end{dfn}
\begin{rem}
If in the setting of Definition~\ref{D:sigma-preserving}
$f$ is $\sigma$-preserving
(with respect some~$S$),
then~$f$ is order preserving as well.
\end{rem}
%
%                  FIRST EXTENSION THEOREM
%
\begin{thm}
\label{T:ext1}
Let $\vs{V}{L}{\varphi}{E}$ and
 $\vs{W}{K}{\psi}{F}$ be extendible valuation systems.
Let
$A\colon V\ra W$
and 
$f\colon E\ra F$ be $\sigma$-preserving maps,
such that 
\begin{equation*}
A(L)\subseteq K
\qquad\text{and}\qquad
f \circ \varphi = \psi \circ A|L.
\end{equation*}
\begin{equation*}
\xymatrix{
V\ar[d]_{A} & 
  L\ar[r]^{\varphi} \ar[d]_{A|L} \ar @{_{(}->} [l] & 
  E\ar[d]^{f} \\
W &
  K\ar[r]^{\psi} \ar @{_{(}->} [l] & 
  F
}\end{equation*}
Then $A(\overline{L})\subseteq \overline{K}$
and $f \circ \overline\varphi = \overline\psi\circ A|\overline L$.
\end{thm}
\begin{proof}
We prove with induction that for every ordinal number~$\alpha$
we have 
\begin{equation}
\label{eq:L:ext1-1}
\begin{alignedat}{3}
A(\Pi_\alpha L)\,&\subseteq\, \overline{K}&
 \qquad f \circ \Pi_\alpha \varphi 
     \,&=\, \overline \psi \circ A | \Pi_\alpha L \\
A(\Sigma_\alpha L)\,&\subseteq\, \overline{K}&
 \qquad f \circ \Sigma_\alpha \varphi 
     \,&=\, \overline \psi \circ A | \Sigma_\alpha L. 
\end{alignedat}
\end{equation}
This is sufficient, because
$\overline{L}=\Pi_{\aleph_1} L$ and 
 $\overline \varphi = \Pi_{\aleph_1} \varphi$
(see Corollary~\ref{C:aleph1}).

\begin{enumerate}
\item We prove~\eqref{eq:L:ext1-1}
holds for $\alpha=0$.
Since $\Pi_0\varphi = \Sigma_0\varphi = \varphi$,
we need to prove that
$A(L)\subseteq \overline{K}$
and $f\circ \varphi = A| L$.
But this is valid by assumption.

\item
Let~$\alpha$ be an ordinal number
and assume~\eqref{eq:L:ext1-1}
holds for~$\alpha$;
we prove~\eqref{eq:L:ext1-1} holds for~$\alpha+1$.
We  only prove $A(\Pi_{\alpha+1} L) \subseteq \overline{K}$
and $f\circ \Pi_{\alpha+1}\varphi = \overline \psi \circ A|\Pi_{\alpha+1} L$;
the other part, $A(\Sigma_{\alpha+1} L) \subseteq \overline{K}$
and $f\circ \Sigma_{\alpha+1}\varphi 
= \overline\psi \circ A|\Sigma_{\alpha+1}L$
follows similarly.

Let $a\in \Pi_{\alpha+1} L$
be given. We need to prove that
\begin{equation}
\label{eq:L:ext1-2}
A(a)\in \overline{K}\qquad\text{and}\qquad
\overline\psi (A(a)) \,=\, f(\Pi_{\alpha+1}\varphi(a)).
\end{equation}
Recall that $\Pi_{\alpha+1}L = \Pi(\Sigma_\alpha L)$,
so write $a=\bw_n a_n$
for some 
$\Sigma_\alpha \varphi$-convergent
$a_1 \geq a_2 \geq \dotsb$
and note that $\Pi_{\alpha+1}\varphi (a) = \bw_n \Sigma_\alpha\varphi(a_n)$.
We have
\begin{alignat*}{3} 
f(\Pi_{\alpha+1}\varphi(a))
 \,&=\, f(\bw_n \Sigma_\alpha \varphi (a_n))
   && \\
 \,&=\, \bw_n f(\Sigma_\alpha\varphi(a_n))\qquad
   && \text{since $f$ is $\sigma$-preserving} \\
 \,&=\, \bw_n \overline{\psi} (A( a_n )) 
   && \text{since \eqref{eq:L:ext1-1} holds for~$\alpha$.} 
\end{alignat*}
So we see that $A(a_1)\geq A(a_2)\geq \dotsb$
is $\overline{\psi}$-convergent.
Since $\vs{W}{\overline{K}}{\overline{\psi}}{F}$ is complete,
this implies $\bw_n A(a_n) \in \overline{K}$
and $\bw_n\overline{\psi}(A(a_n))=\overline{\psi}(\bw_n A(a_n))$.
Because $A$ is $\sigma$-preserving,
we have $\bw_n A(a_n) = A(a)$.
Hence $A(a)\in \overline{K}$
and 
\begin{alignat*}{3}
f(\Pi_{\alpha+1}\varphi(a))
 \,&=\,  \bw_n \overline{\psi}(A(a_n)) \\
 \,&=\,  \overline{\psi} (\bw_n A(a_n)) \\
 \,&=\,  \overline{\psi} (A(a)).
\end{alignat*}
So we have proven Statement~\eqref{eq:L:ext1-2}.

\item
Let $\lambda$ be a limit ordinal,
and assume that \eqref{eq:L:ext1-1}
holds for all~$\alpha<\lambda$;
we prove that \eqref{eq:L:ext1-1} holds for~$\lambda$.
Since  $\Pi_\lambda \varphi = \Sigma_\lambda \varphi$,
we must prove that
\begin{equation}
A(\Pi_\lambda L)\subseteq \overline K
\qquad\text{and}\qquad
f\circ \Pi_\lambda \varphi \,=\, \overline\psi \circ A | \Pi_\lambda L.
\end{equation}
Let $a\in \Pi_\lambda L$ be given
in order to prove $A(a)\in \overline{K}$
and $\overline\psi(A(a))= f(\Pi_\lambda\varphi(a))$.
Recall that $\Pi_\lambda L = \bigcup_{\alpha < \lambda} \Pi_\alpha L$,
and $\Pi_\lambda \varphi \,|\,\Pi_\alpha L = \Pi_\alpha \varphi$
for all~$\alpha < \lambda$.
So choose $\alpha < \lambda$ such that
$a\in \Pi_\alpha L$.
Since \eqref{eq:L:ext1-1} holds for~$\alpha$,
we know that 
\begin{equation*}
A(\Pi_\alpha L) \subseteq \overline K
\qquad\text{and}\qquad
f\circ \Pi_\alpha \varphi = \overline\psi\circ A | \Pi_\alpha L.
\end{equation*}
Hence $A(a) \in A(\Pi_\alpha L) \subseteq \overline K$ and 
$f(\Pi_\lambda \varphi(a)) = f(\Pi_\alpha \varphi(a)) =
\overline\psi(A(a))$.\qedhere
\end{enumerate}
\end{proof}
%
%                  EXAMPLE ON RINGS A with X in A
%
\begin{ex}
\label{E:ring-ext-1}
Let $\mathcal{A}$
be a ring of subsets of~$X$.
Let $\mu\colon \mathcal A \ra \R$
be a positive and additive map.
Recall that~$\vs{\wp X}{\mathcal{A}}{\mu}{\R}$
is a valuation system (see Example~\ref{E:ring-system}).
Assume that~$\vs{\wp X}{\mathcal{A}}{\mu}{\R}$
is extendible.

We would like to prove that~$\overline{\mathcal{A}}$
is also a ring of subsets of~$X$ (as is~$\mathcal{A}$).
For the moment,
we will prove this under the assumption that~$X\in \mathcal{A}$,
see Example~\ref{E:ring-ext-2}.

To prove that~$\overline{\mathcal{A}}$
is a ring,
we need to show that $Z\backslash Y \in \overline{\mathcal{A}}$
for all~$Z,Y\in\overline{\mathcal{A}}$.
Note that $Z\backslash Y = (X\backslash Y)\cap Z$.
So it suffices to show that
$X\backslash Y \in\overline{\mathcal{A}}$
for all $Y\in\overline{\mathcal{A}}$.

Consider the order \emph{reversing} maps $A\colon \wp X \ra \wp X$
and $f\colon \R \ra \R$ given by
\begin{alignat*}{3}
A( Y ) &= X\backslash Y
\qquad&&(Y\subseteq X) \\
f(x) &= \mu(X) - x
\qquad&&(x\in \R).
\end{alignat*}
In order to apply Theorem~\ref{T:ext1}
to these maps,
let us rebaptise them as order preserving
maps $A\colon \wp X \ra (\wp X)^\mathrm{op}$
and $f\colon \R \ra \R^\mathrm{op}$
(see Example~\ref{E-val-opposite}).

We have the following situation.
\begin{equation*}
\xymatrix{
\wp X\ar[d]_{A} & 
  \mathcal A \ar[r]^{\mu} \ar[d]_{A|L} \ar @{_{(}->} [l] & 
  \R \ar[d]^{f} \\
(\wp X)^\mathrm{op} &
  {\mathcal{A}}^\mathrm{op} \ar[r]^{\mu} \ar @{_{(}->} [l] & 
  \R^\mathrm{op}
}\end{equation*}
We leave it to the reader to verify that
$\vs{(\wp X)^\mathrm{op}}{{\mathcal A}^\mathrm{op}}\mu{\R^\mathrm{op}}$
is again valuation system which is extendibe.
We have
 $A(\mathcal A)\subseteq \mathcal A$,
because
$X\backslash Y \in \mathcal{A}$
for all~$Y\in\mathcal{A}$
since $\mathcal{A}$ is a ring containing~$X$.
Further,
since $\mu$ is additive, we have
\begin{equation*}
\mu(A(Y))\,=\, \mu(X\backslash Y)
\,=\, \mu(X) - \mu(Y)
\,=\, f(\mu(Y)).
\end{equation*}
So $\mu\circ A$ and $f\circ \mu$ 
are identical on~$\mathcal A$.
Note that
\begin{equation*}
\textstyle{X\backslash\bigcup_n A_n \,=\, \bigcap_n X\backslash A_n}
\qquad\text{and}\qquad
\mu(X)-\bv_n x_n \,=\, \bw_n \,\mu(X) - x_n
\end{equation*}
where $Y_1 \subseteq Y_2 \subseteq \dotsb$
are from~$\wp X$
and $x_1 \leq x_2 \leq \dotsb$ 
is a bounded sequence in~$\R$.
So $A$ and $f$ are $\sigma$-preserving
(see Definition~\ref{D:sigma-preserving}).

Hence by Theorem~\ref{T:ext1},
we get $A(\overline{\mathcal{A}})\subseteq \overline{\mathcal{A}}$
and secondly
  $f\circ\overline\mu=\overline \mu \circ A$ on $\overline{\mathcal{A}}$.

From the first fact
we get that $X\backslash Y \in \overline{Y}$
for all $Y\in \mathcal{A}$,
and hence~$\overline{\mathcal{A}}$ is a ring.

From the second fact,
we get $\overline\mu (X\backslash Y)
= \overline\mu(X) - \overline\mu(Y)$
for all $Y\in \overline{\mathcal{A}}$.
From this,
one might say,
we see that~$\overline\mu$ is additive.
However, we already knew this
as $\overline\mu$ is modular
(and $\mu(\varnothing)=0$, see Definition~\ref{P:hier}).
\end{ex}
%
%                  SECOND EXTENSION THEOREM
%
\begin{thm}
\label{T:ext2}
Let $\vs{V}{L}{\varphi}{E}$ and
 $\vs{W}{K}{\psi}{F}$ be extendible valuation systems.
Let
$A\colon V\ra W$
be a $\sigma$-preserving map
such that $A(L)\subset K$.
Let $f\colon E\ra F$ be a $\sigma$-preserving group-homomorphism
such that
\begin{equation}
\label{eq:T:ext2-as}
\ld{\overline{\psi}}(A(c),A(d))\ \leq\  f(\ld{\overline{\varphi}}(c,d))
\end{equation}
for all $c,d\in V$ with $A(c),A(d)\in \overline K$.
\begin{equation*}
\xymatrix{
V\ar[d]_{A} & 
  L\ar[r]^{\varphi} \ar[d]_{A|L} \ar @{_{(}->}[l] & 
  E \ar[d]^f\\
W &
  K\ar[r]^{\psi} \ar @{_{(}->}[l] & 
  F
}\end{equation*}
Then $A(\overline{L})\subseteq \overline{K}$.
\end{thm}
\begin{proof}
With induction,
we prove that for every ordinal number~$\alpha$,
we have 
\begin{equation*}
A(\Pi_\alpha L)\subseteq \overline K
\qquad\text{and}\qquad
A(\Sigma_\alpha L)\subseteq \overline K.
\end{equation*}
This is sufficient, because $\overline L = \Pi_{\aleph_1} L$.

As one can  see,
such a proof might be quite similar to the proof of Theorem~\ref{T:ext1}.
Therefore,
we leave the details to the reader
and only prove the following statement.
\begin{equation}
\label{eq:T:ext2-1}
A(\Sigma_\alpha L) \subseteq \overline K
\quad\implies\quad
A(\Pi_{\alpha+1}L)\subseteq\overline K.
\end{equation}

Assume $A(\Sigma_\alpha L)\subseteq \overline{K}$
and let $a\in \Pi_{\alpha+1} L$ be given;
we must prove $A(a)\in \overline K$.
Write $a=\bw_n a_n$ for some $\Sigma_{\alpha}\varphi$-convergent
sequence
$a_1 \geq a_2 \geq \dotsb$
in~$\Sigma_\alpha L$.
Because we have assumed $A(\Sigma_\alpha L)\subseteq \overline{K}$,
we know that $A(a_i) \in \overline K$.
To prove $A(a)\in\overline K$,
it suffices to show that $A(a_1) \geq A(a_2) \geq \dotsb$
is $\overline{\psi}$-convergent.
Indeed, then
\begin{equation*}
A(a) \,\equiv\, A(\bw_n a_n) \,=\, \bw_n A(a_n)\,\in\, \overline{K},
\end{equation*}
because $A$ is $\sigma$-preserving and 
$\vs{W}{\overline{K}}{\overline{\psi}}{F}$
is complete.

To prove that the sequence $A(a_1)\geq A(a_2) \geq \dotsb$
is $\overline{\psi}$-convergent,
we must show that $\bw_n \overline \psi(A(a_n))$ exists.
Note that by Inequality~\eqref{eq:T:ext2-as},
we have
\begin{alignat*}{3}
\overline{\psi}(A(a_{n+1}))
\,-\, \overline{\psi}(A(a_n))
\ &=\ 
\ld{\overline\psi}(A(a_{n+1}),A(a_n))\\
\ &\leq\ 
f(\ld{\overline\varphi}(a_{n+1},a_n))
\ =\ 
f(\overline\varphi(a_{n+1},a_n)) 
\,-\, f(\overline\varphi(a_n)).
\end{alignat*}
So since~$F$ is $R$-complete (see Definition~\ref{D:R-complete}),
in order to show that $\bw_n\overline{\psi}(A(a_n))$
exists,
it suffices to prove that $\bw_n f(\overline{\varphi}(a_n))$
exists. 
For this
we need to prove that $\bw_n\overline\varphi(a_n)$ exists
(as  $f$ is $\sigma$-preserving).
That is,
we must show that $a_1 \geq a_2 \geq \dotsb$ is 
$\overline{\varphi}$-convergent.
Of course,
this follows quickly from the fact that $a_1 \geq a_2 \geq \dotsb$
is $\Sigma_{\alpha}\varphi$-convergent.
We have proven Statement~\eqref{eq:T:ext2-1}.
\end{proof}
%
%                  RESTRICTION PROPOSITION         
%
\begin{prop}
\label{P:restriction}
Let $\vs{V}{L}\varphi{E}$ be an extendible valuation system.
Note that its completion is denoted by
$\vs{V}{\overline L}{\overline \varphi}{E}$.
Given $\ell \leq u$ from~$L$,
consider
\begin{equation*}
\vs{[\ell,u]}{L\cap{[\ell, u]}}{\varphi|[\ell,u]}{E};
\end{equation*}
it is an extendible valuation system.
Note that its completion is denoted by
\begin{equation*}
\vs{[\ell,u]}{\overline{L\cap[\ell,u]}}{\overline{\varphi|[\ell,u]}}{E}.
\end{equation*}
We have $\overline{L\cap [\ell, u]} = \overline{L}\cap [\ell,u]$.
Moreover,
 $\overline{\varphi}$ and $\overline{\varphi|[\ell,u]}$
are identical on $\overline{L\cap [\ell,u]}$.
\end{prop}
\begin{proof}
One can easily see that $\overline\varphi | [\ell,u]$
extends $\varphi|[\ell,u]$
and that the valuation system
\begin{equation*}
\vs{[\ell,u]}{\overline{L}\cap [\ell,u]}{\overline\varphi | [\ell, u]}{E}
\end{equation*}
is complete.
Hence $\varphi|[\ell,u]$ is extendible,
and $\overline\varphi | [\ell,u]$ extends
$\overline{\varphi|[\ell,u]}$ (see Proposition~\ref{P:comp-monotonous}).
In particular,
$\overline{L\cap[\ell,u]}\subseteq \overline{L}\cap[\ell,u]$
and $\overline{\varphi}$ and $\overline{\varphi|[\ell,u]}$
are identical on~$\overline{L\cap[\ell,u]}$.
It remains to be shown that 
\begin{equation}
\label{eq:P:restriction}
\overline{L}\cap[\ell,u]\,\subseteq\,\overline{L\cap[\ell,u]}.
\end{equation}
To this end,
consider the map~$\varrho\colon V\ra [\ell, u]$
given by $\varrho(x)=\ell\vee (x\wedge u)$.
Note that $\varrho(x)=x$ for all~$x\in [\ell, u]$,
and thus $\varrho(\overline{L})=\overline{L}\cap[\ell,u]$.
So in order to prove~\eqref{eq:P:restriction},
we must show that
$\varrho(\overline{L})\subseteq\overline{L\cap[\ell,u]}$.
To do this, we apply Theorem~\ref{T:ext2}.
\begin{equation*}
\xymatrix{
V\ar[d]_{\varrho} &&
  L\ar[rr]^{\varphi} \ar[d]_{\varrho|L} \ar @{_{(}->}[ll] & &
  E\ar[d]^{1_E} \\
[\ell,u] &&
  L\cap[\ell,u]\ar[rr]^{\varphi|[\ell,u]} \ar @{_{(}->}[ll] &&
  E
}\end{equation*}
We must verify that $\varrho$ is
$\sigma$-preserving and that
\begin{equation}
\label{eq:P:restriction-1}
\ld{\overline{\varphi|[\ell,u]}}(\varrho(c),\varrho(d))
\ \leq\ 
\ld{\overline{\varphi}}(c,d)
\end{equation}
for all $c,d\in V$ with $\varrho(c),\varrho(d)\in \overline {L\cap [\ell,u]}$.
One can easily see that~$\varrho$ is $\sigma$-preserving,
because $V$ is $\sigma$-distributive
(see Definition~\ref{D:sigma-distributive}).
Concerning Inequality~\eqref{eq:P:restriction-1},
note that for $c,d\in V$
with $\varrho(c),\varrho(d)\in \overline {L\cap [\ell,u]}$
we have
\begin{alignat*}{3}
\ld{\overline\varphi|[\ell,u]}(\varrho(c),\varrho(d)) 
\ &=\ \ld{\overline\varphi}(\varrho(c),\varrho(d)) 
   && \text{since $\overline\varphi|\,\overline{L\cap[\ell,u]}
                           =\overline{\varphi| [\ell,u]}$}\\
\ &=\  \ld{\overline\varphi}(\,\ell\vee(c\wedge u),\,\ell\vee(d\wedge u)\,) 
   \quad&& \text{by definition of~$\varrho$} \\
\ &\leq\  \ld{\overline\varphi}(\,c\wedge u,\,d\wedge u\,) 
   && \text{by Lemma~\ref{L:wv-unif}} \\
\ &\leq\  \ld{\overline\varphi}(c,d)
   && \text{by Lemma~\ref{L:wv-unif}}.
\end{alignat*}
Hence Theorem~\ref{T:ext2} is applicable,
and we obtain Inequality~\eqref{eq:P:restriction}.
\end{proof}
%
%                  EXAMPLE ON RINGS
%
\begin{ex}
\label{E:ring-ext-2}
Let $\mathcal{A}$
be a ring of subsets of~$X$.
Let $\mu\colon \mathcal A \ra \R$
be a positive and additive map.
Recall that~$\vs{\wp X}{\mathcal{A}}{\mu}{\R}$
is a valuation system (see Example~\ref{E:ring-system}).
Assume that~$\vs{\wp X}{\mathcal{A}}{\mu}{\R}$
is extendible
(see Definition~\ref{D:extendible}).

We prove that~$\overline{\mathcal{A}}$ is a ring.
(In Example~\ref{E:ring-ext-1},
we saw that this
is the case if~$X\in \mathcal{A}$.)

Let $Y,Z\in\overline{\mathcal{A}}$ be given.
To prove that~$\overline{\mathcal{A}}$
is a ring,
we must show that
\begin{equation*}
Y\backslash Z\in\overline{\mathcal{A}}.
\end{equation*}

We restrict our attention to the interval $I\eqdf [\,\varnothing,\, Y\cup Z\,]$.
Note that $\mathcal{A}\cap I$ is a ring of subset of $Y\cup Z$
with $Y\cup Z\in\mathcal{A}\cap I$.
So by Example~\ref{E:ring-ext-1},
we know that~$\overline{\mathcal{A}\cap I}$ is a ring.
Note that 
$Y,Z\in\overline{\mathcal{A}\cap I}$
because $\overline{\mathcal A} \cap I 
= \overline{ \mathcal A \cap I}$
by Proposition~\ref{P:restriction}.
So we get $Y\backslash Z\in \overline{\mathcal{A}\cap I}$,
and thus $Y\backslash Z\in \overline{\mathcal{A}}\cap I$
by Proposition~\ref{P:restriction}.

Hence~$\overline{\mathcal{A}}$ is a ring of subsets of~$X$.
\end{ex}

%
%                  THIRD EXTENSION THEOREM
%
\begin{thm}
\label{T:ext3}
Let $\vs{V}{L}{\varphi}{E}$
and
 $\vs{W}{K}{\psi}{F}$ be extendible valuation systems.
Let $R$ be a sublattice of~$V$ with 
$L\subseteq R$.
Let
$f\colon E\ra F$ be a $\sigma$-preserving map,
and let $A\colon R\ra W$ be $\sigma$-preserving
with respect to~$V$.
Assume that $A(L)\subseteq K$
and  $f \circ \varphi = \psi \circ A|L$.
\begin{equation*}
\xymatrix{
V & 
  R\ar[d]_{A} \ar @{_{(}->}[l]& 
  L\ar[r]^{\varphi} \ar[d]_{A|L} \ar @{_{(}->}[l] & 
  E\ar[d]^{f} \\
& W &
  K\ar[r]^{\psi} \ar @{_{(}->}[l] & 
  F
}\end{equation*}
Assume that
$R$ is convex in~$V$,
and that
for every $c\in R$,
there are $\varphi$-convergent sequences $a_1 \leq a_2 \leq \dotsb$ 
and $b_1 \geq b_2 \geq \dotsb$
such that  $\bw_n a_n \leq c \leq \bv_n b_n$.

Then $A(\overline L\cap R)\subseteq \overline{K}$
and $f \circ \overline\varphi = \overline\psi\circ A$ on $\overline L \cap R$.
\end{thm}
\begin{proof}
Let us first prove the following special case.
\begin{equation}
\label{eq:T:ext3-claim-1}
\left[\quad
\begin{minipage}{.7\columnwidth}
Let $c\in \overline{L}\cap R$
with $\ell\leq c\leq u$ for some~$\ell,u\in L$.
Then 
\begin{equation*}
A(c)\in\overline{K}
\qquad\text{and}\qquad
f(\overline\varphi(c))= \overline\psi (A(c)).
\end{equation*}
\end{minipage}
\right.
\end{equation}
Let $c \in \overline{L}\cap R$
with $\ell \leq c \leq u$ for some $\ell,u\in L$ be given.
Then clearly $c\in [\ell, u]$.
Further, $[\ell,u]\subseteq R$ since~$R$ is convex
and $\ell,u\in R$ (as $L\subseteq R$).
So we have:
\begin{equation*}
\xymatrix{
[\ell,u]\ar[d]_{A|[\ell,u]} & 
  L\cap[\ell,u]\ar[r]^{\varphi|[\ell,u]} 
               \ar[d]_{A\,|\,L\cap[\ell,u]} \ar @{_{(}->}[l] & 
  E\ar[d]^{f} \\
W &
  K\ar[r]^{\psi} \ar @{_{(}->}[l] & 
  F
}\end{equation*}
Moreover, by Proposition~\ref{P:restriction},
we know that $c\in \overline{L\cap[\ell,u]}$
and $\overline{\varphi}(c) = \overline{\varphi|[\ell,u]}(c)$.
Hence Theorem~\ref{T:ext1}
yields $A(c)\in\overline{K}$
and $\overline{\psi}(A(c))
= f(\overline{\varphi|[\ell,u]}(c))$.
But then
$\overline{\psi}(A(c))
= f(\overline{\varphi}(c))$.
This proves Statement~\eqref{eq:T:ext3-claim-1}.

\vspace{.5em}
We proceed by proving another special case.
\begin{equation}
\label{eq:T:ext3-claim-2}
\left[\quad
\begin{minipage}{.7\columnwidth}
Let $c\in \overline{L}\cap R$
and suppose $c\geq \ell$ for some~$\ell\in L$.
Then 
\begin{equation*}
A(c)\in\overline{K}
\qquad\text{and}\qquad
f(\overline\varphi(c))= \overline\psi (A(c)).
\end{equation*}
\end{minipage}
\right.
\end{equation}
Let $c\in\overline{L}\cap R$ with $c\geq \ell$
for some~$\ell$ be given.
Pick $\varphi$-convergent $u_1 \leq u_2 \leq \dotsb$
such that $c \leq \bv_n u_n$.
Since $u_1 \geq u_2 \geq \dotsb$ is $\varphi$-convergent
and $c\in\overline{L}$,
we know that $c\wedge u_1 \leq c\wedge u_2 \leq \dotsb$
is $\overline{\varphi}$-convergent (see Proposition~\ref{P:R-main}).
Since $\vs{V}{\overline{L}}{\overline\varphi}{E}$ is complete,
this implies $\overline\varphi(c)
=\overline\varphi(\bv_n c\wedge u_n)
=\bv_n\overline\varphi(c\wedge u_n)$.
We get:
\begin{alignat*}{3}
f(\overline\varphi(c))
\,&=\, f(\bv_n\overline\varphi(c\wedge u_n)) &&  \\
\,&=\, \bv_n f(\overline\varphi(c\wedge u_n))
   \qquad && \text{since $f$ is $\sigma$-preserving}
\shortintertext{Note that $\ell\leq c\wedge u_n \leq u_n$.
So by \eqref{eq:T:ext3-claim-1}, we get $A(c\wedge u_n)\in\overline K$ and:}
f(\overline\varphi(c))
\,&=\, \bv_n \overline{\psi} (A (c\wedge u_n))  &&
\shortintertext{From this we 
see $A(c\wedge u_1) \leq A(c\wedge u_2)\leq \dotsb$
is $\overline\psi$-convergent.
Since $\vs{W}{\overline K}{\overline \psi}{F}$
is complete,
we get $f(\overline\varphi(c)=\bv_n A(c\wedge u_n) \in \overline K$
and}
f(\overline\varphi(c))
 \,&=\,  \overline\psi (\bv_n A(c\wedge u_n)) \\
   &=\,  \overline\psi (A(\bv_n c\wedge u_n) )
      && \text{since $A$ is $\sigma$-preserving} \\
   &=\,  \overline\psi (A(c)).
\end{alignat*}
This completes the proof of Statement~\eqref{eq:T:ext3-claim-2}.

\vspace{.5em}
We are now ready to give the proof of the general case.
Let $c\in R\cap \overline L$ be given.
We need to prove that $A(c)\in\overline K$
and $f(\overline{\varphi}(c))=\overline{\psi}(A(c))$.
Pick $\varphi$-convergent $\ell_1 \geq \ell_2 \geq \dotsb$
such that $\bw_n \ell_n \leq c$.
Since $\ell_1 \geq \ell_2 \geq \dotsb$
is $\varphi$-convergent and $c\in \overline L$,
we know that $\ell_1 \vee c \geq \ell_2 \vee c\geq \dotsb$ 
is $\overline\varphi$-convergent.
Since $\vs{V}{\overline L}{\overline\varphi}{E}$
is complete,
this implies that $\overline\varphi(c)
=\overline\varphi(\bw_n \ell_n \vee c)
=\bw_n\overline\varphi(\ell_n \vee c)$.
We get:
\begin{alignat*}{3}
f(\overline\varphi(c))
\,&=\, f(\bw_n\overline\varphi(\ell_n \vee c)) &&  \\
\,&=\, \bw_n f(\overline\varphi(\ell_n\vee c))
   \qquad && \text{since $f$ is $\sigma$-preserving}
\shortintertext{Note that $\ell_n\leq \ell_n \vee c$.
Further,
since $R$ is a sublattice of~$V$,
and $c\in R$, $\ell_n\in L\subseteq R$,
we get $\ell_n \vee c\in R$.
So by \eqref{eq:T:ext3-claim-2}, we have $A(\ell_n\vee c)\in\overline K$ and:}
f(\overline\varphi(c))
\,&=\, \bw_n \overline{\psi} (A (\ell_n \vee c))  &&
\shortintertext{From this we 
see $A(\ell_1 \vee c) \geq A(\ell_2 \vee c)\geq \dotsb$
is $\overline\psi$-convergent.
Since $\vs{W}{\overline K}{\overline \psi}{F}$
is complete,
we get $f(\overline\varphi(c)=\bw_n A(\ell_n \vee c) \in \overline K$
and}
f(\overline\varphi(c))
 \,&=\,  \overline\psi (\bw_n A(\ell_n \vee c)) \\
   &=\,  \overline\psi (A(\bw_n\ell_n \vee c) )
      && \text{since $A$ is $\sigma$-preserving} \\
   &=\,  \overline\psi (A(c)).
\end{alignat*}
We are done.
\end{proof}
%
%                  CLOSURE UNDER ADDITION
%
\begin{prop}
\label{P:subgroup}
Let $\vs{V}{L}{\varphi}{E}$ be an extendible valuation system.
Let $R$ be a sublattice of~$V$
endowed with a group structure.
Assume $L$ is a subgroup of~$R$ and that $\varphi$ is a group homomorphism
(recall that~$E$ is an ordered Abelian \emph{group}).

Further, assume 
that $R$ is convex and that for every $c\in R$,
there are $\varphi$-convergent sequences
$a_1 \geq a_2 \geq \dotsb$ and
$b_1 \leq b_2 \leq \dotsb$ 
such that $\bw_n a_n \leq c\leq \bv_n b_n$.

Then $\overline L \cap R$ is a subgroup of~$R$,
and $\overline \varphi | R$ is a group homomorphism.
\end{prop}
\begin{proof}
In order to show that $\overline L\cap R$ is a subgroup of~$R$,
we must prove the following.
\begin{enumerate}
\item \label{tp:P:subgroup-1}
If $a,b\in \overline{L}\cap R$, then $a+b\in \overline L$.

\item \label{tp:P:subgroup-2}
If $a \in \overline{L}\cap R$, then $-a\in \overline{L}$.
\end{enumerate}
We only give a proof for~\ref{tp:P:subgroup-1}.
It will then be clear how to prove~\ref{tp:P:subgroup-2}.

We aim to apply Theorem~\ref{T:ext3}.
To this end,
the reader can easily verify that
$\vs{V\times V}{L\times L}{\varphi\times\varphi}{E\times E}$
is an extendible valuation system;
that its completion is
$\vs{V\times V}{\overline L\times\overline L}
{\overline\varphi \times\overline\varphi}{E\times E}$;
that $R\times R$ is a convex sublattice of~$V\times V$;
that the assignment $(c,d)\mapsto c+d$
yields a $\sigma$-preserving map $A\colon R\times R\ra V$
 with respect to~$V\times V$;
that the map $f\colon E\times E\ra E$ given by
$f(x,y)=x+y$ is $\sigma$-preseving.

Further, note that $A(L\times L)\subseteq L$
because~$L$ is a subgroup of~$R$.
Note that for all $c_1,c_2\in R\times R$
there are $\varphi$-convergent  $\ell^i_1 \geq \ell^i_2 \geq \dotsb$
and $u^i_1 \leq u^i_2 \leq \dotsb$
such that $\bw_n \ell^i_n \leq c_i \leq \bv_n u^i_n$,
and thus $\bw_n (\ell^1_n,\ell^2_n) \leq (c_1,c_2) \leq \bv_n (u^1_n, u^2_n)$,
where $(\ell^1_1,\ell^2_1) \geq (\ell^1_2,\ell^2_2)\geq \dotsb$
and $(u^1_1,u^2_1)\leq (u^1_2,u^2_2)\leq \dotsb$
are $\varphi\times\varphi$-convergent.
Finally,
note that $f\circ (\varphi\times \varphi) = \varphi \circ A|(L\times L)$,
because $\varphi$ is a group homomorphism.
\begin{equation*}
\xymatrix{
V\times V & 
  R\times R\ar[d]_{+} \ar @{_{(}->}[l]& 
  L\times L\ar[r]^{\varphi\times\varphi} 
     \ar[d]_{+} \ar @{_{(}->}[l] & 
  E\times E\ar[d]^{+} \\
& V &
  L\ar[r]^{\varphi} \ar @{_{(}->}[l] & 
  E
}\end{equation*}

So we are in a position to apply Theorem~\ref{T:ext3}.
It gives us that $A(\overline L\times \overline L \cap R\times R)\subseteq L$ 
and 
$f\circ (\overline \varphi \times \overline \varphi)
= \overline \varphi \circ A$ on
 $\overline {L}\times \overline{L} \cap R\times R$.
In other words,
if $c,d\in \overline{L}\cap R$, then $c+d \in \overline{L}$,
and $\overline\varphi(c+d)=\overline\varphi(c)+\overline\varphi(d)$.
Hence we have proven statement~\ref{tp:P:subgroup-1},
and at the same time we have shown that~$\ol\varphi$ is a group homomorphism.
\end{proof}
%
%                  EXAMPLE ON RIESZ SPACE OF FUNCTIONS
%
\begin{ex}
\label{E:riesz-space-ext}
Let $X$ be a set.
Let $F$ be a Riesz space of functions on~$X$.
Let $\varphi\colon F\ra R$ be a positive linear map.
Recall that $\vs{[-\infty,\infty]^X}F{\varphi}\R$
is a valuation system.
Assume that $\varphi$ is extendible.

We would like to prove that $\overline F$ is a Riesz space of functions
and $\overline\varphi$ is linear.
However, since addition is only defined on $R\eqdf\R^X$,
we will instead show that
$\overline F\cap R$ is a Riesz space
of functions and that $\overline\varphi|R$ is linear.
Moreover,
we assume
that for every~$f\in\overline F \cap R$
there are $\varphi$-convergent sequences
 $\ell_1 \geq \ell_2 \geq \dotsb$
and $u_1 \leq u_2 \leq \dotsb$
such that $\bw_n \ell_n \leq f \leq \bv_n u_n$.

To prove that $\overline F \cap R$ is a Riesz space,
we must show that
\begin{enumerate}
\item
 $f+g\in \overline F$ 
for all $f,g\in \overline F \cap R$, and

\item
$\lambda \cdot f\in \overline F$
for all $\lambda\in \R$ and $f\in \overline F\cap R$.
\end{enumerate}
We only prove the first statement;
we leave it to the reader to prove the second.

Of course,
it suffices to establish that $\overline F\cap R$
is a subgroup of~$R$.
To do this,
we apply Proposition~\ref{P:subgroup}.
Indeed, one can easily
see that all the prerequisites are met.
To name a few: one sees that
$R$ is a sublattice of~$V$,
that $F$ is a subgroup of~$R$ (since $F$ is a Riesz space of functions),
that $\varphi$ is a group homomorphism (since $\varphi$ is linear),
and that $R$ is convex (since $\R$ is convex in $[-\infty,\infty]$).

Proposition~\ref{P:subgroup} not only gives us that $\overline F \cap R$
is a subgroup of~$R$, but
also that~$\varphi|R$ is a group homomorphism.
We leave it to the reader to prove that
$\overline\varphi|R$ is homogeneous,
i.e., $\overline\varphi(\lambda\cdot f) = \lambda \cdot \overline\varphi(f)$
for all~$f\in \overline F \cap R$ and $\lambda \in \R$.

Hence $\overline F\cap R$
is a Riesz space of functions, and $\overline \varphi |R$ is linear.
\end{ex}

 }
\clearpage
{ \section{Extendibility}
\label{S:benign}
\noindent
Let 
$\vs{V}{L}{\varphi}{E}$
be a valuation space,
and suppose we want to prove
that $\varphi$ can be extended to a complete valuation.
We have seen
that it suffices to prove that~$\varphi$
is $\Pi_{\aleph_1}$-extendible
(see Corollary~\ref{C:aleph1}).
However,
to prove $\varphi$ is $\Pi_{\aleph_1}$-extendible
already
seems like a monumental task
when one has only barely started
to unfold
the definition of ``$\varphi$ is $\Pi_{\aleph_1}$-extendible''
(see Definition~\ref{P:hier}):
\begin{alignat*}{3}
\varphi &\text{ is $\Pi$-extendible,}&
\quad&\text{ and }\quad& \varphi &\text{ is $\Sigma$-extendible;}\\
\Pi\varphi &\text{ is $\Pi$-extendible,}&
\quad&\text{ and }\quad& \Sigma\varphi &\text{ is $\Sigma$-extendible;}\\
\Pi_2\varphi &\text{ is $\Pi$-extendible,}&
\quad&\text{ and }\quad& \Sigma_2\varphi &\text{ is $\Sigma$-extendible;}\\
&\ \,\vdots&&\quad\vdots&&\ \,\vdots\\
\Pi_\omega\varphi &\text{ is $\Pi$-extendible,}&
\quad&\text{ and }\quad& \Sigma_\omega\varphi &\text{ is $\Sigma$-extendible;}\\
\Pi_{\omega+1}\varphi &\text{ is $\Pi$-extendible,}&
\quad&\text{ and }\quad& \Sigma_{\omega+1}\varphi
    &\text{ is $\Sigma$-extendible;}\\
&\ \,\vdots&&\quad\vdots&&\ \,\vdots
\end{alignat*}

It turns out
that for some~$E$
the situation is more tractable.
For instance,
we will see that if~$E=\R$,
then to prove that $\varphi$ is extendible
it suffices to show that~$\varphi$ is $\Pi_2$-extendible
or $\Sigma_2$-extendible.
Actually,
we have a sharper result:
it suffices to show that~$\varphi$ is \emph{continuous}
(see Definition~\ref{D:continuity}).
Those~$E$ for which we have 
\begin{equation*}
\varphi \text{ is continuous }
\quad\implies\quad
\varphi \text{ is extendible.}
\end{equation*}
will be called \emph{benign} (see Definition~\ref{D:benign}).
%%%%%%%%%%%%%%%%%%%%%%%%%%%%%%%%%%%%%%%%%%%%%%%%%%%%%%%%%%%%%%%%%%%%
%
%                  CONTINUITY
%
\subsection{Continuous Valuations}
Below we define what it means for a valuation system
to be continuous.
We will see that we have the following implications
\begin{equation*}
\xymatrix @=1em {
\Sigma_2\text{-extendible}\ar @{=>} [rd]
&&
\Sigma\text{-extendible}  \\
& \text{continuous} \ar @{=>} [rd] \ar @{=>} [ru]& \\
\Pi_2\text{-extendible}\ar @{=>} [ru]
&&
\Pi\text{-extendible}.
}
\end{equation*}
In fact,
we prove that $\varphi$ is continuous
if and only if it can be extended to~$\Pi L \cup \Sigma L$ in some sense
(see Lemma~\ref{L:continuity}),
so that we might have dubbed
it ``$\Pi\cup\Sigma$-extendible''.
\begin{dfn}
\label{D:continuity}
\label{D:continuous}
Let $\vs{V}{L}\varphi{E}$ be a valuation system.\\
We say $\varphi$ (or more precisely  $\vs{V}{L}\varphi{E}$)
is \keyword{continuous} provided that
\begin{equation*}
\bw_n a_n \,\leq\, \bv_n b_n 
\quad\implies\quad
\bw_n \varphi(a_n) \,\leq\, \bv_n \varphi(b_n)
\end{equation*}
for all $\varphi$-convergent $a_1 \geq a_2 \geq \dotsb$
and $\varphi$-convergent $b_1 \leq b_2 \leq \dotsb$.
\end{dfn}
\begin{ex}
We leave it to the reader
to verify that the valaution systems
\begin{equation*}
\vsSA\qquad\text{and}\qquad\vsSF
\end{equation*}
are continuous.
\end{ex}
\begin{lem}
\label{L:continuity}
Let $\vs{V}{L}\varphi{E}$ be a valuation system.
The following are equivalent.
\begin{enumerate}
\item
\label{L:continuity-1}
$\varphi$ is continuous.
\item
\label{L:continuity-2}
$\varphi$ is $\Pi$-extendible
and $\Sigma$-extendible,
and there is an order preserving map
\begin{equation*}
f\colon \Pi L \cup \Sigma L \ra E
\end{equation*}
 that extends both~$\Pi\varphi$ and $\Sigma\varphi$.
\end{enumerate}
\end{lem}
\begin{proof}
\noindent
\emph{\ref{L:continuity-1}
$\Longrightarrow$
\ref{L:continuity-2}}\ 
Suppose that~$\varphi$ is continuous.
By Lemma~\ref{L:Pi-continuity},
we see that $\varphi$ is $\Pi$-extendible.
Similarly, $\varphi$ must be $\Sigma$-extendible.
We need to find an order preserving map~$f\colon \Pi L \cup \Sigma L \ra E$
that extends both $\Pi\varphi$ and $\Sigma\varphi$.
We have little choice,
\begin{equation}
\label{eq:L:continuity-1}
f(c) \ \eqdf \ 
\begin{cases}
\ \Pi\varphi(c) \quad& \text{if }c\in \Pi L;\\
\ \Sigma\varphi(c) \quad& \text{if }c\in \Sigma L.
\end{cases}
\end{equation}
To see that Equation~\eqref{eq:L:continuity-1}
is a valid definition of a 
map $f\colon \Pi L \cup \Sigma L \ra E$,
we need to verify that $\Pi\varphi$ and $\Sigma\varphi$
are identical on~$\Pi L \cap \Sigma L$.
Let $c\in \Pi L \cap \Sigma L$
be given.
We must prove $\Pi\varphi(c) = \Sigma\varphi(c)$.
Choose $\varphi$-convergent
$a_1 \geq a_2 \geq \dotsb$ and
$\varphi$-convergent
$b_1 \leq b_2 \leq \dotsb$
such that $\bw_n a_n =c= \bv_n b_n$.

Then $b_n \leq a_n$ for all~$n$, so
$\varphi(b_n)\leq \varphi(a_n)$ for all~$n$.
Hence
\begin{equation*}
\Sigma\varphi(c)
=\bv_n \varphi(b_n)
\ \leq\ \bw_n \varphi(a_n)
=\Pi\varphi(c).
\end{equation*}

Conversely,
we have $\bw_n a_n \leq \bv_n b_n$,
so since $\varphi$ is continuous we get
\begin{equation*}
\Pi\varphi(c)
=\bw_n \varphi(a_n)
\ \leq\ \bv_n \varphi(b_n)
=\Sigma\varphi(c).
\end{equation*}
Hence $\Pi\varphi(c)=\Sigma\varphi(c)$.
So Equation~\eqref{eq:L:continuity-1}
is a valid definition of~$f$.

Since by defintion,
$f$ extends both $\Pi\varphi$ and $\Sigma\varphi$,
it only remains to be shown that~$f$ is order preserving.
Let $c,d\in \Pi L \cup \Sigma L$
with $c\leq d$ be given.
We prove
\begin{equation*}
f(c)\ \leq\  f(d).
\end{equation*}
Of course,
if $c,d$ are both in $\Pi L$,
then we done,
because $\Pi \varphi$
is order preserving 
and $f$ extends $\Pi\varphi$.
Similarly, if $c,d\in\Sigma L$, 
we also immediately get $f(c)\leq f(d)$.

Suppose $c\in \Pi L$ and $d\in \Sigma L$.
Choose $\varphi$-convergent
 sequences $b_1 \leq b_2 \leq \dotsb$
and  $a_1 \geq a_2 \geq \dotsb$ 
such that 
$\bv_n b_n = c$
and 
$\bw_n a_n = d$.
Then $b_m \leq \bv_n b_n \leq \bw_n a_n \leq a_m$
for all~$m$,
so $\varphi(b_m)\leq \varphi(a_m)$ for all~$m$,
and hence
\begin{equation*}
f(c) = \Sigma\varphi(c)
= \bv_n \varphi(b_n)
\ \leq\ 
\bw_n \varphi(a_n)
= \Pi \varphi(d)
= f(d).
\end{equation*}

Suppose $c\in \Pi L$ and $d\in \Sigma L$.
Choose $\varphi$-convergent sequences $a_1 \geq a_2 \geq \dotsb$ 
and $b_1 \leq b_2 \leq \dotsb$
such that $\bw_n a_n = c$ and $\bv_n b_n = d$.
Then $\bw_n a_n \leq \bv_n b_n $.
So since $f$ is continuous, 
we get 
$f(c)= \Pi\varphi(c) =\bw_n\varphi(a_n) 
\leq\bv_n\varphi(b_n) = \Sigma\varphi(d) = f(d)$.

\vspace{.3em}

\noindent\emph{\ref{L:continuity-2}
$\Longrightarrow$
\ref{L:continuity-1}}\ 
Let $f\colon \Pi L \cup \Sigma L \ra E$
be an order preserving map that extends both $\Pi\varphi$
and $\Sigma \varphi$.
We prove that~$\varphi$ is continuous.
Let $\varphi$-convergent sequences
$a_1 \geq a_2 \geq \dotsb$
and $b_1 \leq b_2 \leq \dotsb$
with $\bw_n a_n \leq \bv_n b_n$
be given.
We need to prove that  
\begin{equation*}
\bw_n \varphi(a_n)\ \leq\  \bv_n \varphi(b_n).
\end{equation*}
This is easy;
since~$f$ is order preserving
and extends $\Pi\varphi$ and $\Sigma\varphi$,
we get 
\begin{equation*}
\bw_n \varphi(a_n)
\,=\,
\Pi\varphi(\bw_n a_n)
\,=\,
 f(\bw_n a_n)
\ \leq\  
f(\bv_n b)
\,=\,
\Sigma\varphi(\bv_n b_n)
\,=\,
\bv_n \varphi(b_n).\qedhere
\end{equation*}
\end{proof}
%
%                  IMPLICATIONS
%
\begin{cor}
\label{C:cont-imp}
Let $\vs{V}{L}\varphi{E}$ be a valuation system. 
\begin{enumerate}
\item 
\label{C:cont-imp-1}
If $\varphi$ is continuous,
then $\varphi$ is $\Pi$-extendible
and $\Sigma$-extendible.

\item 
\label{C:cont-imp-2}
$\varphi$ is continuous
provided that $\varphi$ is either $\Pi_2$-extendible
or $\Sigma_2$-extendible.
\end{enumerate}
\end{cor}
\begin{proof}
Point \ref{C:cont-imp-1}
follows immediately from Lemma~\ref{L:continuity}.
Point~\ref{C:cont-imp-2} is also
a consequence of Lemma~\ref{L:continuity}.
Indeed, 
assume that $\varphi$ is $\Pi_2$-extendible.
We prove that $\varphi$ is continuous.
Note that $\Pi_2 \varphi$
is order preserving and
extends both $\Pi\varphi$ and $\Sigma\varphi$.
Hence $\varphi$ satisfies condition~\ref{L:continuity-1}
of Lemma~\ref{L:continuity}.
Thus $\varphi$ is continuous.
\end{proof}
%
%                  LEMMA ON EXT OF CONTINUITY
%
\begin{lem}
\label{L:cont-ext}
Let $\vs{V}{L}\varphi{E}$
be a valuation system.\\
Let~$K$ be a sublattice of~$L$
such that $\psi \eqdf \varphi | K$
is continuous.\\
Then $\varphi$ is continuous
under the following assumptions.
\begin{enumerate}
\item\label{L:cont-ext-1}
Given a $\varphi$-convergent sequence $a_1 \geq a_2 \geq \dotsb$ 
in~$L$, we have
\begin{equation*}
\bw_n \varphi(a_n) \ = \ 
\bv\,\bigl\{\ \Pi \psi(\ell) \colon\ 
 \ell \in S \ \bigr\},
\end{equation*}
for some $S\subseteq \Pi K$
with $\ell \leq \bw_n a_n$ for all $\ell \in S$.

\item\label{L:cont-ext-2}
Given a 
$\varphi$-convergent sequence 
$b_1 \leq b_2 \leq \dotsb$ in~$L$,
we have
\begin{equation*}
\bv_n \varphi(b_n) \ = 
\ \bv \ \bigl\{ \ \Sigma\psi(u) \colon \ u \in T\ \bigr\},
\end{equation*}
for some $T\subseteq \Sigma K$ with $\bv_n b_n \leq u$
for all~$u\in T$.
\end{enumerate}
\end{lem}
\begin{proof}
Let $\varphi$-convergent sequences $a_1 \geq a_2 \geq \dotsb$
and $b_1 \leq b_2 \leq \dotsb$ from~$L$ 
with $\bw_n a_n \leq \bv_n b_n$ be given.
To prove that~$\varphi$
is continuous
(see Definition~\ref{D:continuity}),
we must show that $\bw_n \varphi(a_n) \leq \bv_n \varphi(b_n)$.
Let $\ell\in S$ and $u\in T$
be given.
Note that 
$\ell \leq \bw_n a_n \leq \bv_n b_n \leq u$,
so
$\Pi\psi(\ell) \leq \Sigma\psi (u)$
since $\psi$ is continuous.

Hence $\bw_n \varphi(a_n) \leq \bv_n \varphi(b_n)$
by Assumptions~\ref{L:cont-ext-1}
and~\ref{L:cont-ext-2}.
\end{proof}
%%%%%%%%%%%%%%%%%%%%%%%%%%%%%%%%%%%%%%%%%%%%%%%%%%%%%%%%%%%%%%%%%%%%%%%%%%%%%
%
%                  BENIGN E
%
\subsection{Benign~$E$}
\begin{dfn}
\label{D:benign}
Let $E$ be an ordered Abelian group.
We say~$E$ is \keyword{benign} provided 
that for every valuation system
$\vs{V}{L}\varphi{E}$,
we have
\begin{equation*}
\varphi\text{ is continuous}
\quad\implies\quad
\varphi\text{ is extendible}.
\end{equation*}
\end{dfn}
\begin{ex}
We will prove that~$\R$ is benign (see Corollary~\ref{C:R-benign}).
\end{ex}
\begin{ex}
Let $I$ be a set
and let $X_i$ be a benign ordered Abelian group
for every~$i\in I$.
We leave it to the reader to verify
that $\prod_{i \in I} X_i$ is benign.
\end{ex}

 }
\clearpage
{ \section{Uniformity on $E$}
\label{S:unif}
\noindent
To prove that~$\R$
is benign (see Definition~\ref{D:benign}),
we study
ordered Abelian groups~$E$ which are endowed with
a certain uniformity (such as~$\R$)
in Subsection~\ref{SS:fitting}.
We prove that all such~$E$ are benign
(see Theorem~\ref{T:fitting-benign}),
in the following way.

Let $\vs{V}{L}\varphi{E}$ be a valuation system.
Recall that in order to prove that~$E$ is benign
 we must show that 
if $\varphi$ is continuous,
then $\varphi$ is extendible (see Definition~\ref{D:benign}).
We will first prove that
if $\varphi$ is continuous,
then both $\Pi \varphi$ and $\Sigma \varphi$ are continuous
(see Lemma~\ref{lem:cont-ext-single}).
Then, by induction,
we see that $\varphi$ is
 both $\Pi_n$-extendible and $\Sigma_n$-extendible,
and both $\Pi_n\varphi$  and $\Sigma_n \varphi$
are continuous,
for every~$n\in \N$.
Hence $\varphi$ is $\Pi_\omega$-extendible.
However,
it is not clear a priori that~$\Pi_\omega\varphi$
is continuous.

Secondly,
we prove that if 
$\varphi$ is $\Pi_\lambda$-extendible
for some ordinal number~$\lambda$,
then $\Pi_\lambda\varphi$ is continuous.
So by induction
we see that $\varphi$ is both $\Pi_\alpha$-extendible
and $\Sigma_\alpha$-extendible,
and both $\Pi_\alpha\varphi$ and $\Sigma_\alpha\varphi$
are continuous,
for every ordinal number~$\alpha$ 
(see Lemma~\ref{L:fitting-ext}).
Hence $\varphi$ is extendible
(see Corollary~\ref{C:aleph1}).

To prove the second statement
we use the fact that
elements of $\Pi_\alpha L$ (or $\Sigma_\alpha L$)
can be approximated from below
by elements of $\Pi L$, in some sense
(see Lemma~\ref{L:fitting-dense}).
We will express this by
$\Sigma L$ is \emph{lower $\Pi_\alpha \varphi$-dense}
in $\Pi_\alpha L$.
We will formally introduce this notion,
and study it, in Subsection~\ref{SS:dense}.

\subsection{Fitting Uniformity}
\label{SS:fitting}
%
%                  DEFINITION OF A GOOD UNIFORMITY
%
\begin{dfn}
\label{D:uniformity}
\label{D:fitting-uniformity}
Let~$E$ be an ordered Abelian group.
A \keyword{fitting uniformity} on~$E$
is a \emph{countable} set~$\Phi$ of binary relations on~$E$
with the following properties.
\newcounter{epropc}
\begin{enumerate}
\item 
\label{E-refl}
We have $s \se s$ for all $\ve \in \Phi$ and $s \in E$.
\item
\label{E-min}
There is a map $\wedge\colon \Phi\times\Phi \ra \Phi$
such that
\begin{equation*}
s \ \ \varepsilon\wedge\delta\ \   t
\quad\implies\quad
s \,\varepsilon\,t\ \text{ and }\ s\,\delta\,t
\qquad (\varepsilon,\delta\in\Phi,\ s,t\in E).
\end{equation*}

\item
\label{E-half}
There is a map $\dt-\colon \Phi \ra \Phi$
such that
\begin{equation*}
r\ \ \dt\varepsilon \ \ s \ \  \dt\varepsilon \ \ t
\quad\implies\quad
r \se t\qquad(\varepsilon\in\Phi,\ r,s,t\in E).
\end{equation*}

\item \label{E-ord}
Given $\varepsilon\in\Phi$ and $r,s,t\in E$ 
with $r\leq s\leq t$,
we have
\begin{equation*}
r\,\varepsilon\,t
\quad\implies\quad
r\,\varepsilon\,s
\ \text{ and }\ 
s\,\varepsilon\,t.
\end{equation*}

\item \label{E-haus}
Let $s,t\in E$ with $s\leq t$.
Then $s=t$ provided that $s\,\varepsilon\,t$ for all~$\varepsilon\in\Phi$.

\item \label{E-inf-conv}
If a sequence $s_1 \geq s_2 \geq \dotsb$ from~$E$
has an infimum~$s\in E$,
then 
\begin{equation*}
\forall\varepsilon\in \Phi
\ \ \exists N\in \N
\ \ s \, \varepsilon\, s_N.
\end{equation*}

\item  \label{E-bound-inf}
Let $s_1\geq s_2 \geq \dotsb$ be a
sequence in~$E$,
and assume that
for every $\ve\in \Phi$
there is an~$N \in \N$ such that 
$s_{n} \ \varepsilon \ s_{N_\ve}\quad (n\geq N)$.

Then $s_1 \geq s_2 \geq \dotsb$ has an infimum~$\bw_n s_n$.

\item \label{E-add}
Let $r,s,t\in E$ and $\varepsilon\in\Phi$ be given.
Then $s\,\varepsilon\,t$ implies $r+s\ \varepsilon\ r+t$.
\setcounter{epropc}{\value{enumi}}
\end{enumerate}
\end{dfn}
\begin{ex}
We define a fitting uniformity on~$\R$.
For each natural number~$n$,
let $\varepsilon_n$ be the binary relation
on~$\R$ given by 
\begin{equation*}
s \ \varepsilon_n\  t
\quad\iff\quad 
s\leq t\quad\text{and}\quad t-s \leq 2^{-n}.
\end{equation*}
Then $\Phi\eqdf \{ \varepsilon_n \colon n\in \N \}$
is a fitting uniformity on~$\R$.

(Take $\varepsilon_n \wedge \varepsilon_m \eqdf \varepsilon_{n\vee m}$
and $\dt{\varepsilon_n} \eqdf \varepsilon_{n+1}$
for all $n,m\in \N$.)
\end{ex}

\begin{rem}
The fitting uniformities
defined here
are related to
the \emph{uniform spaces}
(or more precisely,
\emph{quasi uniform spaces})
studied in topology, see~\cite{Willard70}.

However
we do not involve uniform spaces,
because  the usual way of reasoning about them
does not seem to fit well with
property~\ref{E-ord}.
Moreover,
we do not wish to assume that the reader is familiar with uniform spaces.
\end{rem}
%
%                  LEMMA ON FURTHER PROPERTIES
%
To the list of properties 
that a fitting uniformity must have (see Definition~\ref{D:uniformity}),
we add 
some easy observations
in Lemma~\ref{L:E-prop}.
When we speak of ``property~($q$)'',
where~$q$ is some Roman numeral,
we refer to this list.
\begin{lem}
\label{L:E-prop}
Let $E$ be an ordered Abelian group
with a fitting uniformity~$\Phi$.
\begin{enumerate}
\setcounter{enumi}{\value{epropc}}
\item \label{E-minus}
Let $s,t\in E$ and $\ve \in \Phi$ be given.
Then $s \se t$ implies\, $-t \ \ve\,-\!\!s$.

\item \label{E-sup-conv}
If a sequence $s_1 \leq s_2 \leq \dotsb$ from~$E$
has an supremum~$s\in E$,
then 
\begin{equation*}
\forall\varepsilon\in \Phi
\ \ \exists N\in \N
\ \ s_N \, \varepsilon\, s.
\end{equation*}

\item  \label{E-bound-sup}
Let $s_1\leq s_2 \leq \dotsb$ be a
sequence in~$E$,
and assume that
for every $\ve\in \Phi$
there is an~$N_\ve \in \N$ such that 
$s_{N_\ve} \ \varepsilon \ s_n\quad (n\geq N_\ve)$.

Then $s_1 \leq s_2 \leq \dotsb$ has a supremum~$\bv_n s_n$.
\end{enumerate}
\end{lem}
\begin{proof}
\noindent\ref{E-minus}\ 
Let $s,t\in E$ and $\varepsilon\in \Phi$ be given,
and assume $s\se t$.
We prove $-t\ \varepsilon\,-\!\!s$.
Indeed, by property~\ref{E-add}, we have
\begin{equation*}
-t\,=\,-(t+s)\,+\,s\quad\ve\quad -(t+s)\,+\,t \,=\,-s.
\end{equation*}

\noindent\ref{E-sup-conv}
Let $s_1 \leq s_2 \leq \dotsb$
be a sequence in~$E$
which has a supremum~$s$ in~$E$.
Let $\varepsilon \in \Phi$ be given.
We need to find an~$N\in\N$
such that $s_N \ \varepsilon\ s$.

Let us consider the sequence
$-s_1 \geq -s_2 \geq\dotsb$.
By Lemma~\ref{L:oag-minus-preserves}
the sequence $-s_1 \geq -s_2 \geq \dotsb$
has an infimum, $-s$.
By property~\ref{E-inf-conv}
we have 
$-s\ \varepsilon\ {-s_N}$
for some~$N$.
Then by property~\ref{E-minus},
we get $s_N\ \varepsilon\ s$,
and we are done.

\vspace{.3em}
\noindent\ref{E-bound-sup}
Similar:
apply property~\ref{E-bound-inf}
to the sequence $-s_1 \geq -s_2 \geq \dotsb$.
\end{proof}

%
%                  NOTATION
%
\begin{nt}
\label{N:unif}
Let $E$ be an ordered Abelian group
with a fitting uniformity~$\Phi$.
\begin{enumerate}
\item
\label{N:unif-leq}
Given binary relations $\varepsilon$
and $\delta$ on~$E$
(for instance, $\varepsilon,\delta\in \Phi$),
we write
\begin{equation*}
\varepsilon \ \leq\ \delta
\quad\iff\quad 
\forall s,t\in E\ 
[\ s\ \varepsilon\ t
\implies
s\ \delta\ t\ ].
\end{equation*}

\item
\label{N:unif-plus}
Given binary relations $\varepsilon$ and $\delta$ on~$E$,
let $\varepsilon + \delta$
be the relation on~$E$ given by
\begin{equation*}
s\ \ \varepsilon + \delta\ \ t
\quad\iff\quad
\exists q\in E\ 
[\ s\ \varepsilon\ q\ \delta\ t\ ].
\end{equation*}
\end{enumerate}
\end{nt}
\begin{rem}
The operation ``$+$''
defined in Notation~\ref{N:unif}\ref{N:unif-plus}
is associative,
but not in general commutative
(contrary to the expectation the symbol ``$+$'' evokes).

The chosen notation does have advantages:
property~\ref{E-half} can be written as
\begin{equation*}
\dt\varepsilon \,+\,\dt\varepsilon \ \leq\ \varepsilon
\qquad(\varepsilon \in \Phi).
\end{equation*}
\end{rem}

%
%                  WITH A FITTING UNIFORMITY,
%                  E IS R-COMPLETE
%
\begin{lem}
\label{L:E-R-complete}
Let~$E$ be an ordered Abelian group with fitting uniformity~$\Phi$.\\
Then $E$ is $R$-complete (see Definition~\ref{D:R-complete}).
\end{lem}
\begin{proof}
Let $x_1 \leq x_2 \leq \dotsb$
and $y_1 \leq y_2 \leq \dotsb$ from~$E$ be given such that
\begin{equation}
\label{eq:unif-implies-R-comp}
x_{N+1} - x_N \ \leq\ y_{N+1} - y_N
\qquad(N\in\N).
\end{equation}
Assume $\bv_n y_n$ exists.
To that $E$ is $R$-complete,
we must show that $\bv_n x_n$ exists.

Let $\ve \in \Phi$ be given.
By property~\ref{E-bound-sup},
we know that 
to prove  $\bv_n x_n$ exists,
it suffices to find $N\in \N$ such that $x_N \se x_n$
for all~$n\geq N$.

By property~\ref{E-sup-conv},
we know there is an~$N$ such that $y_N \se \bv_m y_m$.
Let $n\geq N$ be given.
We will prove that $x_N \se x_n$.
We already know  $y_N \se y_n$
by property~\ref{E-ord}
because $y_N \leq y_n \leq \bv_m y_m$
and $y_N \se \bv_n y_m$.
So 
 $0\,\se\  (y_n - y_N)$
by property~\ref{E-add}.

From Inequality~\eqref{eq:unif-implies-R-comp}
one can easily derive that
\begin{equation*}
0 \ \leq\ x_n - x_N \ \leq\ y_n - y_N.
\end{equation*}
Since $0 \,\se\,(y_n-y_N)$
we get $0 \,\se\, (x_n -x_N)$ by property~\ref{E-ord}.

Hence $x_N \se x_n$ by property~\ref{E-add}.
So we are done.
\end{proof}
%
%                  DIRECTED SET INFIMUM LEMMA
%
The following lemma will be useful.
\begin{lem}
\label{lem:conv-inf}
Let~$E$ be an ordered Abelian group with 
fitting uniformity~$\Phi$.

Let $S\subseteq E$
be non-empty and \emph{downwards directed},
i.e., for all~$s_1,s_2\in S$,
there is an~$s\in S$ such that $s\leq s_1$
and $s\leq s_2$.

Let $t\in E$ be a lower bound of~$S$
which is close to~$S$ in the sense that
\begin{equation}
\label{eq:L:conv-inf}
\forall\varepsilon\in\Phi\ \exists s\in S\quad t\,\varepsilon\, s\text{.}
\end{equation}

Then~$t$ is the infimum of~$S$.
\end{lem}
\begin{proof}
To show that $t$ is the infimum of~$S$,
we need to prove that $\ell \leq t$
for every lower bound~$\ell$ of~$S$.
To do this, we take a detour.

Let~$\varepsilon_1,\varepsilon_2,\dotsb$
be an enumeration of~$\Phi$.
Using Equation~\eqref{eq:L:conv-inf},
and the fact that~$S$ is non-empty and directed,
choose~$s_1 \geq s_2 \geq \dotsb$ in~$S$
such that 
\begin{equation}
\label{eq:L:conv-inf-1}
t\ \varepsilon_n \ s_n\qquad(n\in \N).
\end{equation}

We will prove that
$s_1 \geq s_2 \geq \dotsb$
has an infimum~$s$ and that $s=t$.

This is sufficient to prove that~$t$ is the infimum of~$S$.
Indeed,
if $\ell$ is a lower bound of~$S$,
then $\ell$ is a lower bound of~$s_1 \geq s_2 \geq\dotsb$,
and so $\ell \leq \bw_n s_n =t$.

We use property~\ref{E-inf-conv}
to show that $s_1 \geq s_2 \geq\dotsb$
has an infimum.
Given $\varepsilon\in\Phi$,
we need to find an~$N$ such that $s_n \ \ve\ s_N$
for all~$n\geq N$. Pick $k$ such that $\varepsilon=\varepsilon_k$
and take $N=k$. Let $n\geq N$ be given.
Note that $t \leq s_n \leq s_N =s_k$
and $t\ \varepsilon_k\ s_k$
by Equation~\eqref{eq:L:conv-inf-1}.
So we have $s_n\ \varepsilon_k\ s_N$ by property~\ref{E-ord}.

Hence property~\ref{E-inf-conv} implies that $s_1 \geq s_2 \geq \dotsb$
has an infimum, $s$.
It remains to be shown that $s=t$.
For this we use property~\ref{E-haus}.

Note that $t\leq s$ because $t\leq s_n$ for all~$n$.
Let $\varepsilon\in\Phi$ be given.
We need to prove that $t \se s$.
Choose $k$ such that $\ve = \ve_k$.
Then $t \leq s \leq s_k$
and $t \ \ve_k\ s_k$
by Equation~\eqref{eq:L:conv-inf-1}.
So $t\ \ve_k \  s$ by property~\ref{E-ord}.
Hence $s=t$ by property~\ref{E-haus}.
\end{proof}
%%%%%%%%%%%%%%%%%%%%%%%%%%%%%%%%%%%%%%%%%%%%%%%%%%%%%%%%%%%%%%%%%%%%%%%%%%%%%%%
%
%                  LOWER DENSENESS SUBSECTION
%
\subsection{Denseness}
\label{SS:dense}
Throughout this subsection,
$E$ will be an ordered Abelian group
endowed with a fitting uniformity~$\Phi$
(see Definition~\ref{D:uniformity}).
%
%                  LOWER DENSENESS
%
\begin{dfn}
\label{D:lower-dense}
Let $\vs{V}{L}\varphi{E}$ be a valuation system.
Let $S\subseteq T$ be subsets of~$L$.
We say $S$ is \keyword{lower $\varphi$-dense} in~$T$
provided that the following condition holds.
\begin{equation*}
\begin{minipage}{0.7\textwidth}
For every $a\in T$ and $\varepsilon\in\Phi$
there is an $\ell\in S$
such that
\begin{equation*}
\ell \leq a\quad\text{and}\quad\varphi (\ell)\,\varepsilon\, \varphi (a).
\end{equation*}
\end{minipage}
\end{equation*}
The notion of 
\keyword{upper $\varphi$-denseness} is defined similarly.
\end{dfn}
\begin{ex}
\label{E:sigma-dense}
Let $\vs{V}{L}\varphi{E}$
be a $\Sigma$-extendible valuation system
(see Def.~\ref{D:Pi-extendible}).
\emph{Then $L$ is lower $\Sigma\varphi$-dense in~$\Sigma L$.}

Indeed,
given $a\in \Sigma L$
and $\ve \in \Phi$,
we need to find an~$\ell\in L$ 
such that $\varphi(\ell) \ \ve\ \Sigma\varphi(a)$.
Write $a=\bv_n a_n$ for
some $\varphi$-convergent sequence $a_1 \leq a_2 \leq \dotsb$.
Then we have
\begin{equation*}
\Sigma\varphi(a)\ =\ \bv_n \varphi(a_n).
\end{equation*}
By property~\ref{E-sup-conv},
there is an~$N$ such that $\varphi(a_N) \ \ve\ \Sigma\varphi(a)$.
So take $\ell =a_N$.
\end{ex}
%
%                  ELEMENTARY PROPERTIES OF LOWER DENSENESS
%
\begin{lem}
\label{L:ldense-prop}
Let $\vs{V}{L}\varphi{E}$ be a valuation system.
\begin{enumerate}
\item
\label{L:ldense-prop-1}
Let $R\subseteq S\subseteq T$ be subsets of~$L$.
Suppose $R$ is lower $\varphi$-dense in~$S$,
and suppose that $S$ is lower $\varphi$-dense in~$T$.
Then $R$ is lower $\varphi$-dense in~$T$.

\item
\label{L:ldense-prop-2}
Let $R$ be a subset of~$L$,
and let $\mathcal S$ be a family of subsets of~$L$.\\
If $R$ is lower $\varphi$-dense in each~$S\in \mathcal{S}$,
then $R$ is lower $\varphi$-dense in~$\bigcup \mathcal{S}$.
\end{enumerate}
\end{lem}
\begin{proof}
\noindent\ref{L:ldense-prop-1}\ 
Let $t\in T$ and $\ve \in\Phi$ be given.
To prove $R$ is lower $\varphi$-dense in~$T$,
we need to find an $r\in R$
with $r\leq t$ and $\varphi(r) \ \ve\ \varphi(t)$.
This is easy.
Choose an $s\in S$
such that $s \leq t$ and $\varphi(s)\ \dt\ve\ \varphi(t)$
(see Definition~\ref{D:uniformity}\ref{E-half}
for the meaning of ``$\dt\ve$'').
Choose an $r\in R$
such that $r\leq s$ and $\varphi(r)\ \dt\ve\ \varphi(s)$.
Then $r\leq s$ and $\varphi(r) \ \ve\ \varphi(t)$.
\vspace{.3em}

\noindent\ref{L:ldense-prop-2}\ 
We leave this to the reader.
\end{proof}
%
%
%                  MAIN LEMMA
%
%
The proof that~$E$
 is benign
hinges on the following lemma.
\begin{lem}
\label{lem:main}
Let $\vs{V}{L}\varphi{E}$ be a valuation system.\\
Let $K$ be a lower $\varphi$-dense sublattice of~$L$.\\
Then for every $\varphi$-convergent sequence
 $a_1 \geq a_2 \geq \dotsb$ from~$L$
and $\varepsilon\in \Phi$\\
there is a $\varphi$-convergent sequence
$\tilde a_1 \geq \tilde a_2 \geq \dotsb$ from~$K$
with 
\begin{equation}
\label{eq:main}
\tilde a_n \ \leq\  a_n
\qquad\text{and}\qquad
\bw_n \varphi (\tilde a_n) \ \ \varepsilon\ \ \bw_n\varphi (a_n).
\end{equation}
\end{lem}
\begin{proof}
Let $a_1 \geq a_2 \geq \dotsb$ 
be a $\varphi$-convergent sequence in~$L$,
and let~$\varepsilon\in\Phi$ be given.
We need to find a $\varphi$-convergent sequence 
$\tilde a_1 \geq \tilde a_2 \geq \dotsb$ in~$K$
which satisfies Condition~\eqref{eq:main}.
To this end,
we seek a sequence 
$\tilde a_1 \geq \tilde a_2 \geq \dotsb$ in~$K$
such that
\begin{equation}
\label{eq:lem:main-cond}
\varphi(\tilde a_n) \ \eta\ \varphi (a_n)
\qquad\text{and}\qquad
\forall i\in\N\ \ \exists N\in\N 
\ \ \forall n\geq N
\  [\ \  \varphi \tilde a_n \ \varepsilon_i\  \varphi \tilde a_N \ \ ],
\end{equation}
where $\varepsilon_1,\,\varepsilon_2,\,\dotsc$
is
an enumeration of~$\Phi$,
and $\eta\in\Phi$ with  $2\eta \leq \varepsilon$
(see Notation~\ref{N:unif}).

Such a sequence $\tilde a_1 \geq \tilde a_2 \geq \dotsb$
is $\varphi$-convergent 
(by property~\ref{E-bound-inf}).
We prove that $\tilde a_1 \geq \tilde a_2 \geq\dotsb$
satisfies Condition~\eqref{eq:main}. 
Indeed:
We know $\bw_n \varphi(\tilde a_n)$ exists.
Hence,
there is an~$N\in\N$
with 
 $\bw_n\varphi (\tilde a_n)\ \eta\ \varphi (\tilde a_N)$
by property~\ref{E-inf-conv}.
Then
\begin{equation*}
\bw_n \varphi (\tilde a_n) \quad \eta \quad \varphi(\tilde a_N)
\quad\eta\quad \varphi(a_N).
\end{equation*}
So we have 
$\bw_n \varphi (\tilde a_n) \ \ve\  \varphi(a_N)$.
But $\bw_n \varphi(\tilde a_n) \,\leq\,
\bw_n \varphi(a_n) \,\leq\,\varphi(a_N)$.
Thus 
$\bw_n \varphi(\tilde a_n) \ \ve\ 
\bw_n \varphi(a_n)$ by property~\ref{E-ord}.
Hence $\tilde a_1  \geq \tilde a_2 \geq \dotsb$
satisfies Condition~\eqref{eq:main}.

Finding a sequence $\tilde a_1 \geq \tilde a_2 \geq \dotsb$ 
which satisfies Condition~\eqref{eq:lem:main-cond}
is a subtle affair.
Pick $\eta_1, \eta_2,\dotsc$ and $\zeta_1,\zeta_2,\dotsc$
from~$\Phi$
(using properties~\ref{E-half}
and~\ref{E-min})
such that
\begin{equation*}
2 \eta_i \leq \varepsilon_i\text{,} \qquad
\eta_i \leq \eta\text{,} \qquad
2\zeta_i \leq \eta_i\text{,} \qquad
2\zeta_{i+1} \leq \zeta_{i}\text{.}
\end{equation*}
Then we have
\begin{equation}
\label{eq:lem:main-zeta--eta}
\zeta_i + \dotsb + \zeta_j \leq \eta_i \qquad (i,j\in\N,\  i\leq j)\text{.}
\end{equation}
Pick $\ell_1,\ell_2,\dotsc$ from $K$ 
such that $\ell_n \leq a_n$ and $\varphi(\ell_n)\ \ \zeta_n\  \ \varphi(a_n)$
and define 
\begin{equation*}
\tilde{a}_{ij} \ =\ \ell_i \wedge \dotsb\wedge \ell_j\text{,}
\qquad\qquad
\tilde{a}_n \ =\ \tilde{a}_{1n} =\ell_1 \wedge\dotsb\wedge \ell_n\text{,}
\end{equation*}
where $i,j,n\in N$ with $i\leq j$.
Then $\tilde a_{ij} \in K$ and  $\tilde a_n \leq \ell_n \leq a_n$.
We will prove that
the sequence $\tilde a_1 \geq \tilde a_2 \geq \dotsb$
satisfies Condition~\eqref{eq:lem:main-cond}.

Note that for all~$i,j\in\N$ with $i\leq j$, we have,
by Lemma~\ref{L:wv-unif},
\begin{equation*}
\ld\varphi(\tilde a_{ij},a_j)\ =\ 
\ld\varphi(\ell_i \wedge\dotsb\wedge \ell_j,\  a_i \wedge \dotsb \wedge a_j)
\ \leq\ \ld\varphi(\ell_i,a_i)+\dotsb+\ld\varphi(\ell_j,a_j)\text{.}
\end{equation*}
Since $\varphi(\ell_k)\ \zeta_k\ \varphi(a_k)$
for all~$k$, 
the inequality above 
yields,
using property~\ref{E-add},
\begin{equation*}
\varphi(\tilde a_{ij}) \quad \zeta_i + \dotsb + \zeta_j\quad \varphi(a_j).
\end{equation*}
So because
$\zeta_i+\dotsb+\zeta_j \leq \eta_i$ 
(see Inequality~\eqref{eq:lem:main-zeta--eta}),
we have
\begin{equation}
\label{eq:lem:main-3}
\varphi(\tilde a_{ij} )\quad\eta_i\quad\varphi (a_j).
\end{equation}
In particular,
we get $\varphi (\tilde a_n)\equiv \varphi(\tilde a_{1n}) \ \eta_1 \ \varphi (a_n)$.
Hence $\varphi(\tilde a_n) \ \eta\ \varphi(a_n)$
as $\eta_1 \leq \eta$.

Let~$i\in \N$ be given.
To prove that $\tilde a_1 \geq \tilde a_2 \geq \dotsb$
satisfies Condition~\eqref{eq:lem:main-cond},
it remains to be shown that 
there is an~$N\in\N$ such that
\begin{equation}
\label{eq:lem:main-2}
\varphi(\tilde a_n)\ \ \varepsilon_i \ \ \varphi(\tilde a_N)
\qquad (n \geq N)
\end{equation}
Using property~\ref{E-inf-conv},
determine $N\geq i$ such that $\bw_n \varphi (a_n) \ \eta_i\ \varphi(a_N)$.
We will show that Statement~\eqref{eq:lem:main-2} holds.
Let $n\geq N$ be given.
Note that 
 by Lemma~\ref{L:curry-wc-unif},
\begin{equation*}
\ld\varphi (\tilde a_n, \tilde a_N)
\,=\,\ld\varphi(\,\tilde a_{i-1} \wedge \tilde a_{in},\,
\tilde a_{i-1} \wedge \tilde a_{iN}\,)
\ \leq\ \ld\varphi(\tilde a_{in},\tilde a_{iN}).
\end{equation*}
So to prove Statement~\eqref{eq:lem:main-2},
it suffices to show 
that $\varphi(\tilde a_{in} ) \ \ve_i\ \varphi(\tilde a_{iN})$.

Recall that $\bw_m\varphi(a_m)\ \ \eta_i\ \ \varphi(a_N)$
by choice of~$N$. Then in particular,
we get
$\varphi(a_n) \ \eta_i\ \varphi(a_N)$
by property~\ref{E-ord}.
Further, $\varphi(\tilde a_{in}) \ \eta_i\ \varphi(a_n)$
by Inequality~\eqref{eq:lem:main-3}. So
\begin{equation*}
\varphi (\tilde a_{in})\ \eta_i\ 
\varphi (a_{n})\ \eta_i\ \varphi (a_N).
\end{equation*}
Hence $\varphi (\tilde a_n) \ \ve_i\ \varphi (a_N)$,
because $2\eta_i \leq \ve_i$.
Note that $\varphi(\tilde a_{in}) \leq 
\varphi(\tilde a_{iN}) \leq \varphi(a_N)$.

So by property~\ref{E-ord},
we get
$\varphi(\tilde a_{in}) \ \ve_i\  \varphi(\tilde a_{iN})$.
\end{proof}
\begin{cor}
\label{C:main-nice}
Let $\vs{V}{L}\varphi{E}$ be a $\Pi$-extendible valuation system.\\
Let $K$ be a sublattice of~$L$. Then
\begin{equation*}
\text{ $K$ is lower dense in~$L$}
\quad\implies\quad
\text{$\Pi K$ is lower dense in~$\Pi L$.}
\end{equation*}
\end{cor}
\begin{proof}
Follows immediately from Lemma~\ref{lem:main}.
\end{proof}
%
%                  LEMMA ON DENSENESS IN EXTENSION
%
\begin{lem}
\label{L:fitting-dense}
Let $\vs{V}{L}\varphi{E}$ be a valuation system
which is both $\Sigma$-extendible
and $\Pi$-extendible.
Then for every ordinal number~$\alpha$:
\begin{enumerate}
\item 
\label{L:fitting-dense-1}
If $\varphi$ is $\Pi_\alpha$-extendible,
then 
\begin{equation*}
\text{$\Pi L$ is upper dense in $\Pi_\alpha L$,}
\quad\text{and}\qquad
\text{$\Sigma L$ is lower dense in $\Pi_\alpha L$.}
\end{equation*}

\item
\label{L:fitting-dense-2}
If $\varphi$ is $\Sigma_\alpha$-extendible,
then 
\begin{equation*}
\text{$\Pi L$ is upper dense in $\Sigma_\alpha L$},
\quad\text{and}\qquad
\text{$\Sigma L$ is lower dense in $\Sigma_\alpha L$.}
\end{equation*}
\end{enumerate}
\end{lem}
\begin{proof}
We use induction on~$\alpha$.

\vspace{.3em}
For $\alpha=0$,
Statements~\ref{L:fitting-dense-1}
and~\ref{L:fitting-dense-2} are trivial.

\vspace{.3em}
Let~$\alpha$ be an ordinal number
such that 
Statement~\ref{L:fitting-dense-1}
holds for~$\alpha$
in order to prove that
Statement~\ref{L:fitting-dense-2} holds for~$\alpha+1$.
Suppose~$\varphi$ is $\Sigma_{\alpha+1}$-extendible.
We need to prove that $\Pi L$ is upper dense in~$\Sigma_{\alpha+1} L$
and that $\Sigma L$ is lower dense in~$\Sigma_{\alpha+1} L$.

Note that 
$\varphi$ is $\Pi_\alpha$-extendible,
because  $\varphi$ is $\Sigma_{\alpha+1}$-extendible.

By Statement~\ref{L:fitting-dense-1} for~$\alpha$,
we know that $\Pi L$ is lower dense in $\Pi_\alpha L$.
Further, $\Pi_\alpha L$
is lower dense in $\Sigma(\Pi_\alpha L) = \Sigma_{\alpha+1} L$
by Example~\ref{E:sigma-dense}.
So we see that $\Pi L$ is lower dense in~$\Sigma_{\alpha+1} L$ 
by Lemma~\ref{L:ldense-prop}\ref{L:ldense-prop-1}.

By Statement~\ref{L:fitting-dense-1} for~$\alpha$,
we know that $\Sigma L$ is upper dense in $\Pi_\alpha L$.
So by the dual of Corollary~\ref{C:main-nice},
we have
$\Sigma L = \Sigma(\Sigma L)$ is upper dense in
$\Sigma(\Pi_\alpha L) = \Sigma_{\alpha+1} L$.

Hence, Statement~\ref{L:fitting-dense-2} holds for~$\alpha+1$
(if Statement~\ref{L:fitting-dense-1} holds for~$\alpha$).

Similarly,
if Statement~\ref{L:fitting-dense-2} holds for~$\alpha$,
then Statement~\ref{L:fitting-dense-1} holds for~$\alpha+1$.

\vspace{.3em}
Let
$\lambda$ be a limit ordinal
such that Statement~\ref{L:fitting-dense-1} holds
for all~$\alpha<\lambda$.
We prove that Statement~\ref{L:fitting-dense-1} holds
for~$\lambda$.
Suppose that $\varphi$ is $\Pi_\lambda$-extendible.
We need to prove that $\Pi L$ is upper dense in~$\Pi_\lambda L$
and $\Sigma L$ is lower dense in $\Pi_\lambda L$.

We know that $\varphi$ is $\Pi_\alpha$-extendible
for all~$\alpha < \lambda$.

As Statement~\ref{L:fitting-dense-1} holds for all~$\alpha <\lambda$,
we see that $\Pi L$ is upper dense in all~$\Pi_\alpha L$.
So by Lemma~\ref{L:ldense-prop}\ref{L:ldense-prop-2},
$\Pi L$ is upper dense in
 $\Pi_\lambda L = \bigcup_{\alpha <\lambda} \Pi_\alpha L$.

Similarly,
$\Sigma L$ is lower dense in
$\Sigma_\lambda L= \bigcup_{\alpha < \lambda} \Sigma_\alpha L$.
\end{proof}
%
%                  MAIN COROLLARY
%
\begin{cor}
\label{C:main}
Let $\vs{V}{L}\varphi{E}$ be a valuation system.\\
Let $K$ be a lower~$\varphi$-dense sublattice of~$L$
and assume that $\psi \eqdf \varphi | K$ is $\Pi$-extendible.\\
Let $a_1 \geq a_2 \geq \dotsb$ be a $\varphi$-convergent sequence in~$L$.
Then
\begin{equation}
\label{eq:C:main}
\bw_n \varphi(a_n) \ = \ 
\bv\ \bigl\{\ \Pi \psi(\ell) \colon \ 
 \ell \in S \ \bigr\},
\end{equation}
where 
$S \ \eqdf\ 
\bigl\{\ \bw_n \tilde a_n \colon\ 
\text{$\psi$-convergent $\tilde a_1 \geq \tilde a_2 \geq \dotsb$
with $\tilde a_n \leq a_n$ for all~$n$}\ \bigr\}$.
\end{cor}
\begin{proof}
To prove Statement~\eqref{eq:C:main},
we apply the dual of Lemma~\ref{lem:conv-inf}.
We need to verify 
that  $\Pi\psi(S)\eqdf\{\,\Pi\psi(\ell) \colon \ell \in S\,\}$
is upwards directed,
that~$\bw_n\varphi(a_n)$
is a lower bound of   $\Pi\psi(S)$,
and that
\begin{equation}
\label{eq:C:main-1}
\forall \varepsilon \in \Phi\ \ 
\exists \ell\in S\quad  \Pi\psi(\ell) \ \ \ve \ \ \bw_n \varphi(a_n).
\end{equation}
To begin,
note that Statement~\eqref{eq:C:main-1} follows immediately
from Lemma~\ref{lem:main}.

Let $\psi$-convergent 
$\tilde a_1 \geq \tilde a_2 \geq \dotsb$
with $\tilde a_n \leq a_n$ for all~$n$ be given.
Then we have  $\psi(\tilde a_n) = \varphi(\tilde a_n)\leq \varphi(a_n)$
for all~$n$, so $\Pi\psi(\bw_n\tilde a_n) = \bw_n \psi(\tilde a_n) 
\leq \bw_n \varphi(a_n)$.
Hence~$\bw_n\varphi(a_n)$
is a lower bound of $\Pi\psi(S)$.
 
To prove that~$\Pi\psi(S)$
is upwards directed,
it suffices to show that~$S$
is upwards directed 
(as $\Pi\psi$ is order preserving).
Let $\psi$-convergent sequences 
$\tilde a_1 \geq \tilde a_2 \geq \dotsb$
and
$\tilde a_1' \geq \tilde a_2' \geq \dotsb$
with $\tilde a_n \leq a_n$ and $\tilde a_n' \leq a_n$
be given.
Then 
\begin{equation*}
\tilde a_1 \vee \tilde a_1' \,\leq\, 
\tilde a_1 \vee \tilde a_2' \,\leq\,\dotsb
\end{equation*}
is again a $\psi$-convergent sequence by Proposition~\ref{P:R-main}.
Further $\tilde a_n \vee \tilde a_n' \leq a_n$
for all~$n$.
Hence $\bw_n \tilde a_n \vee \tilde a_n' \in S$.
But also $\bw_n\tilde a_n \leq \bw_n \tilde a_n \vee\tilde a_n'$
and $\bw_n\tilde a_n' \leq \bw_n \tilde a_n \vee \tilde a_n'$.
So we see that $S$ is upwards directed.
\end{proof}
%
%                  SINGLE SIDED CONTINUITY LEMMA
%
\begin{lem}
\label{lem:cont-ext-single}
Let $\vs{V}{L}\varphi{E}$ be a valuation system.\\
Assume $\varphi$ is continuous.
Then $\Pi \varphi$ is continuous.
\end{lem}
\begin{proof}
Note that~$L$ is an
upper $\Pi\varphi$-dense sublattice of~$\Pi L$
(see Example~\ref{E:sigma-dense}).
We apply Lemma~\ref{L:cont-ext}
to prove that $\varphi$ is continuous.
We must verify that Conditions~\ref{L:cont-ext-1}
and~\ref{L:cont-ext-2} of Lemma~\ref{L:cont-ext} hold.

\ref{L:cont-ext-1}\ 
Let $a_1 \geq a_2 \geq \dotsb$ be a $\Pi\varphi$-convergent
sequence in~$\Pi L$.
We need to find $S\subseteq \Pi L$
such that $\bw_n \varphi(a_n) = \bv S$
and $\ell \leq \bw_n a_n$ for all~$\ell \in S$.
By Lemma~\ref{L:Pi-complete},
we know that
$\Pi\varphi$ is $\Pi$-complete.
Hence $\bw_n a_n \in \Pi L$.
So simply take $S=\{ \bw_n a_n \}$.

\ref{L:cont-ext-2}\ 
Follows immediately from Corollary~\ref{C:main}.
\end{proof}
%
%                  DOUBLE SIDED CONTINUITY LEMMA
%
\begin{lem}
\label{lem:cont-ext-double}
Let $\vs{V}{L}\varphi{E}$ be a valuation system.

Let $K$ be a sublattice of~$L$.
Then $\varphi$ is continuous provided that:
\begin{enumerate}
\item
The restriction $\varphi|K$ of $\varphi$ to~$K$ is continuous.
\item
$K$ is lower and upper $\varphi$-dense in~$L$.
\end{enumerate}
\end{lem}
\begin{proof}
This follows from Lemma~\ref{L:cont-ext}.
Indeed, condition~\ref{L:cont-ext-1} holds
by Corollary~\ref{C:main},
and condition~\ref{L:cont-ext-2} holds
by the dual of Corollary~\ref{C:main}.
\end{proof}
%
%                  FITTING UNIFORMITY EXTENSION LEMMA
%
\begin{lem}
\label{L:fitting-ext}
Let $\vs{V}{L}\varphi{E}$ be a continuous valuation system,
and $\alpha$ an ordinal.
Then $\varphi$ is both $\Pi_\alpha$-extendible
and $\Sigma_\alpha$-extendible,
and $\Pi_\alpha\varphi$ and $\Sigma_\alpha\varphi$
are continuous.
\end{lem}
\begin{proof}
With induction on~$\alpha$.

\vspace{.3em}
For $\alpha=0$,
the statement is trivial.

\vspace{.3em}
Let~$\alpha$ be an ordinal number
and assume that $\varphi$
is $\Pi_\alpha$-extendible
and $\Pi_\alpha \varphi$ is continuous.
We prove that $\varphi$
is $\Sigma_{\alpha+1}$-extendible
and $\Sigma_{\alpha+1}\varphi$ is continuous.
Indeed,
since $\Pi_\alpha\varphi$ is continuous,
$\Pi_\alpha\varphi$ is $\Sigma$-extendible
and so $\varphi$ is $\Sigma_{\alpha+1}$-extendible.
Finally,
$\Sigma(\Pi_\alpha \varphi)=\Sigma_{\alpha+1}\varphi$
is continuous
by the dual of Lemma~\ref{lem:cont-ext-single}.

Similarly,
if $\varphi$ is $\Sigma_\alpha$-extendible
and $\Sigma_\alpha \varphi$ is continuous,
then $\varphi$
is $\Pi_{\alpha+1}$-extendible
and $\Pi_{\alpha+1}\varphi$ is continuous.

\vspace{.3em}
Let $\lambda$ be a limit ordinal
such that for each~$\alpha<\lambda$,
we have that $\varphi$ is $\Pi_\alpha$-extendible
and $\Pi_\alpha\varphi$ is continuous.
Note that $\varphi$ is $\Pi_\lambda$-extendible.
We prove that $\Pi_\lambda \varphi$ is continuous.
For this,
we use Lemma~\ref{lem:cont-ext-double}.
Consider $\psi \eqdf \Pi_2 \varphi$.
By assumption,
$\psi$ is continuous.
We know
that $\Pi_\lambda \varphi$ extends $\psi$,
and that $\psi$ extends both $\Pi\varphi$ and $\Sigma\varphi$.
Since $\Pi L$ is lower dense in $\Pi_\alpha L$,
and $\Sigma L$ is upper dense in $\Pi_\alpha L$
(by Lemma~\ref{L:fitting-dense}),
we get that $K\eqdf \Pi_2 L$ is both upper and lower dense in~$\Pi_\alpha L$.
So by Lemma~\ref{lem:cont-ext-double},
we see that~$\Pi_\lambda\varphi$ is continuous.
(Of course, 
the argument is also valid for other choices for~$\psi$,
such as $\Sigma_3\varphi$ and $\Pi_{42} \varphi$.)
\end{proof}
%
%                  THEOREM
%
\begin{thm}
\label{T:fitting-benign}
Let $E$ be an ordered Abelian group. \\
If $E$ has a fitting uniformity,
then~$E$ is benign.
\end{thm}
\begin{proof}
Let $\vs{V}{L}\varphi{E}$ be a continuous valuation system.
To prove that~$E$ is benign,
we must show that~$\varphi$ is extendible
(see Definition~\ref{D:benign}).
It suffices to prove that~$\varphi$ is $\Pi_{\aleph_1}$-extendible
by Corollary~\ref{C:aleph1}.
Now apply Lemma~\ref{L:fitting-ext}.
\end{proof}
\begin{cor}
\label{C:R-benign}
The ordered Abelian group~$\R$ is benign.
\end{cor}
 }
\clearpage
{ \section{Fubini's Theorem}
\label{S:fub}
\noindent
In this section we study Fubini's Theorem.
We have not found a satisfactory generalisation of this theorem
to the setting of valuations.
However,
we will see that it is possible to split
the proof of Fubini's Theorem
into two parts,
so that the first part (Subsection~\ref{SS:fub_part1})
is algebraic in nature
and specific to the setting of step functions,
and the second part (Subsection~\ref{SS:fub_part2}) is more analytic in nature
and a consequence of a
general extension theorem for valuations (see Theorem~\ref{T:fubext}).

\subsection{Algebraic Part}
\label{SS:fub_part1}
Let us first formulate Fubini's Theorem.
This takes time.

Let $X$ be a set,
let $\mathcal{A}_X$ be a ring of subsets of~$X$,
and let 
\begin{equation*}
\mu_X\colon \mathcal{A}_X \ra \R
\end{equation*}
be a positive and additve map
(see Example~\ref{E:ring-val}).

Similarly, let $Y$ be a set,
let $\mathcal{A}_Y$ be a ring of subsets of~$Y$,
and let 
\begin{equation*}
\mu_Y\colon \mathcal{A}_Y \ra \R
\end{equation*}
be a positive and additve map.

Now, let $\mathcal{A}_{X\times Y}$
be the ring of subsets of~$X\times Y$
generated by the subsets
\begin{equation*}
\{\ \  A \times B\colon\ \  \text{$A\in \mathcal{A}_X$, \ \ 
 $B\in\mathcal{A}_Y$} \ \ \}
\end{equation*}
Let $\mu_{X\times Y} \colon \mathcal{A}_{X\times Y} \ra \R$
be the unique positive and additive map such that
\begin{equation*}
\mu_{X\times Y} (\,A\times B\,)
\ =\ \mu_X(A)\,\cdot\,\mu_Y(B)
\end{equation*}
for all~$A\in\mathcal{A}_X$ and $B\in\mathcal{A}_Y$.
Such $\mu_{X\times Y}$ exists,
as the reader can verify.

Let $F_X$ be the set of all $\mathcal{A}_X$-stepfunctions,
i.e., functions of the form
\begin{equation*}
\textstyle{\sum_{n=1}^N \ \lambda_n \cdot \mathbf{1}_{A_n}  },
\end{equation*}
where $A_1,\dotsc,A_N \in\mathcal{A}_X$
and $\lambda_i \in \R$.
As the reader may verify,
the expression
\begin{equation*}
\varphi_X(\,\textstyle{\sum_{n=1}^N \ \lambda_n \cdot \mathbf{1}_{A_n}}\,)
\ =\ \textstyle{\sum_{n=1}^N \ \lambda_n \cdot \mu_X (A_n)}
\end{equation*}
determines a unique positive and linear
map $\varphi_X\colon F_X \ra \R$. 

Similarly, we get
a map $\varphi_Y \colon F_Y \ra \R$,
and a map $\varphi_{X\times Y} \colon F_{X\times Y} \ra \R$.

One can verify that any $f\in F_{X\times Y}$
is of the form
\begin{equation*}
\textstyle{\sum_{n=1}^N} \ \lambda_n \cdot\mathbf{1}_{A_n \times B_n},
\end{equation*}
where $A_1,\dotsc,A_N\in \mathcal{A}_X$,
and $B_1,\dotsc,B_N\in\mathcal{A}_Y$,
and $\lambda_n \in \R$.

So it is not hard to verify that the  equality
\begin{equation*}
\mathcal{F}_X \bigl(\ 
\textstyle{\sum_{n=1}^N} \ \lambda_n \cdot\mathbf{1}_{A_n \times B_n}\ \bigr)
\ \ =\ \ 
\textstyle{\sum_{n=1}^N} \ \lambda_n \cdot \mu_X(A_n) \cdot \mathbf{1}_{B_n}
\end{equation*}
gives us positive and linear map 
$\mathcal{F}_X \colon F_{X\times Y} \ra F_Y$.

Let $f\in F_{X\times Y}$ be given.
For each~$y\in Y$,
define $f^y\in F_X$ by, for all~$x\in X$,
\begin{equation*}
f^y(x) \ =\  f(x,y)
\end{equation*}
One can easily verify that we have, for all~$y\in Y$,
\begin{equation*}
\mathcal{F}_X (f) (y)\ =\ \varphi_X (f^y).
\end{equation*}
So, informally, $\mathcal{F}_{X}(f) = \int f(x,y)\,dx$.

Since $\varphi_X$ is linear
one quickly sees that $\varphi_Y \circ \mathcal{F}_X = \varphi_{X\times Y}$.
Informally,
\begin{equation*}
\int\int f(x,y)\, \,dx \,dy
\ =\ 
\int f
\qquad(f\in F_{X\times Y}).
\end{equation*}
This is Fubini's Theorem for stepfunctions, $F_{X\times Y}$.

Of course, we want to prove
Fubini's Theorem
for the extension~$\ol{F}_{X\times Y}$.

So let us assume
$\varphi_X$, $\varphi_Y$, and $\varphi_{X\times Y}$ are extendible
(see Definition~\ref{D:extendible}).

Alternatively,
we can assume that 
$\varphi_Y$ and $\varphi_Y$
are continuous (see Definition~\ref{D:continuous});
we leave it to the reader to verify that then
 $\varphi_{X\times Y}$ is continuous
(which is not too easy),
and so~$\varphi_{X\times Y}$ is extendible,
since~$\R$ is benign (see Definition~\ref{D:benign}).

Note that it is not possible to find 
an $\ol{\mathcal{F}}_{Y}\colon \ol{F}_{X\times Y} \ra \ol{F}_Y$
such that, for all~$y\in Y$,
\begin{equation*}
\ol{\mathcal{F}}_Y(f)(y) \ = \  \ol{\varphi}_X(f^y).
\end{equation*}
So to formulate Fubini's Theorem for~$\ol{F}_{X\times Y}$
we need a slightly different approach than
the one we used for the stepfunctions.

Consider the space $E_Y \eqdf \qvL{\ol{F_Y}}$
(see Proposition~\ref{P:quotient-lattice}).
We leave it to the reader to
verify that~$E_Y$ can be endowed with 
the structure of an ordered Abelian group,
and a fitting uniformity
(see Definition~\ref{D:uniformity})
such that the map $F_{X\times Y} \ra E_Y$
given by $f\mapsto \mathcal{F}_X(f) /\approx$
is a group homomorphism.

We can now formulate Fubini's Theorem as follows.
\begin{equation}
\label{eq:fub}
\left[\quad
\begin{minipage}{.7\columnwidth}
The valuation 
\begin{equation*}
\mathcal{F}_X\colon F_{X\times Y} \ra \qvL{\ol{F_Y}}
\end{equation*}
is extendible, 
and $\mathrm{dom}(\ol{\varphi}_{X\times Y})
\,\subseteq\, 
 \mathrm{dom}(\ol{\mathcal{F}}_X)$,
and 
\begin{equation*}
\ol\varphi_{X\times Y}(f) \ =\ 
(\,\qvphi{\ol\varphi} \ \circ\  \ol{\mathcal{F}_X} \,)(f)
\end{equation*}
for all $f\in\ol{F}_{X\times Y}$.
\end{minipage}
\right.
\end{equation}
Of course,
to be true to the usual formulation of Fubini's Theorem
we would need to prove that
that $\ol{\mathcal{F}}_X(f)(y) = \ol{\varphi}_X(f^y)$
for \emph{almost} all~$y\in Y$.
We will not do this.

\subsection{Extension of Operations}
\label{SS:fub_part2}
Let $\vs{V}{L}\varphi{E}$ be a valuation system.
In Section~\ref{S:closedness}
we saw that the completion~$\ol{L}$ of~$\varphi$ is \emph{closed} under
various operations.
It is also possible to \emph{extend} operations
 to~$\ol{L}$, which are (initially only) defined on~$L$.
The aim of this subsection is to prove Theorem~\ref{T:fubext}
which is an example of this principle
in case that~$E$ has a fitting uniformity~$\Phi$
(see Definition~\ref{D:uniformity}).

It should be noted that from the methods
found in the proof of Theorem~\ref{T:fubext}
one can easily obtain a stronger version of this theorem.
More interestingly,
the patterns in the proof strongly suggest that
we should make a study of the uniform structure on~$L$
given by the relations $\ol{\ve}$ (where $\ve\in\Phi$) defined by
\begin{equation*}
a \ \ol{\ve}\ b \quad\iff\quad 0\ \ve\ \ld\varphi(a,b)\qquad\quad(a,b\in L).
\end{equation*}

However, we have refrained from proving a stronger version of the theorem
and introducing yet another notion of uniform structure.
Indeed,
we have not found a clear favorite 
 among the several approaches to the strengthening
of the theorem and the axiomatisation of the uniform structure on~$L$.
Accordingly,
we introduce few new notions,
and the proofs in this subsection are sometimes ad hoc.

One new notion we do present is that of
\emph{weak $\varphi$-convergence}
(see Definition~\ref{D:weak-phi-conv}).
As the name suggests,
$\varphi$-convergence (see Definition~\ref{D:seq-phi-conv})
implies weak $\varphi$-convergence
(see Lemma~\ref{L:weak-phi-conv-implies-phi-conv}),
but the reverse implication does not hold
(see Example~\ref{E:weak-phi-conv}).
Nevertheless,
any weakly $\varphi$-convergent sequence
has a $\varphi$-convergent subsequence
(see Proposition~\ref{P:weak-phi-conv-subseq}).

Due to this all 
the notions of $\varphi$-convergent
and weakly $\varphi$-convergent
can be used somewhat interchangeably.
The main merit of ``weak $\varphi$-convergent''
is that some useful statements concerning it
(see Lemma~\ref{L:phi-conv-diag} and  Lemma~\ref{L:phi-conv-inf}) 
can be easily proven,
while it is not clear if the same statement (or a variant) holds for
``$\varphi$-convergent''.

The  main application of Theorem~\ref{T:fubext} 
is the proof of Fubini's Theorem~\ref{T:fub}.
Let us start the work towards a proof.
%
%                  WEAK PHI-CONVERGENCE
%
\begin{dfn}
\label{D:weak-phi-conv}
\label{D:phi-conv-2}
Let $E$ be an ordered Abelian group.\\
Let $\Phi$ be a fitting uniformity on~$E$.\\
Let $L$ be a lattice, and let $\varphi\colon L \ra E$ be a valuation.\\
Let $a \in L$ and let $a_1,a_2,\dotsc$ be a sequence in~$L$.\\
We say $a_1,a_2,\dotsc$ 
\keyword{weakly $\varphi$-converges} to~$a$
if
\begin{equation*}
\forall \ve \in \Phi \quad \exists N \quad \forall n \geq N \quad
[\ \ 0\quad\ve\quad \ld\varphi(a_n,a)  \ \ ].
\end{equation*}
\end{dfn}
%
%                  WEAK PHI-CONVERGENCE IMPLIES PHI-CONVERGENCE
%
\begin{lem}
\label{L:weak-phi-conv-implies-phi-conv}
Let $E$ be an ordered Abelian group.\\
Let $\Phi$ be a fitting uniformity on~$E$.\\
Let $L$ be a lattice, and let $\varphi\colon L \ra E$ be a complete valuation.\\
Let $a \in L$ and let $a_1,a_2,\dotsc$ be a sequence in~$L$. We have:
\begin{equation*}
\text{$a_1,a_2,\dotsc$ $\varphi$-converges to~$a$}
\quad\implies\quad
\text{$a_1,a_2,\dotsc$ weakly $\varphi$-converges to~$a$}.
\end{equation*}
\end{lem}
\begin{proof}
Let $\ve \in \Phi$ be given.
To prove that $a_1,a_2,\dotsc$ weakly $\varphi$-converges to~$a$
we must find an~$N\in \N$ such that
$0 \ \ve\ \ld\varphi(a_n,a)$ for all~$n\geq N$.

To find such~$N$ takes some preparation,
so bear with us.

Since $a_1,a_2,\dotsc$  $\varphi$-converges to~$a$,
i.e., $a_1,a,a_2,a,\dotsc$ is $\varphi$-convergent,
we know that $a_1,a,a_2,a,\dotsc$ is upper $\varphi$-convergent.
That is, we know the following exists.
\begin{equation}
\label{eq:L:weak-phi-conv-implies-phi-conv-1}
u\ \eqdf \ 
\bw_N \bv_{n\geq N}\  \varphi(a\vee a_N\vee \dotsb \vee a_n)
\end{equation}
In particular,
we see that for each~$N\in \N$,
the sequence 
\begin{equation*}
a\vee a_N \,\leq\, a\vee a_N\vee a_{N+1}\,\leq\,\dotsb
\end{equation*}
is $\varphi$-convergent
(in the sense of Definition~\ref{D:phi-conv}).
Since $\varphi$ is complete,
we see that $\ol{a}_N \eqdf \bv_{n\geq N} \, a \vee a_n$ exists in~$L$
and that $\varphi(\ol{a}_N) 
= \bv_{n \geq N} \,\varphi(a \vee a_N \vee\dotsb \vee a_n)$.
Now, note that we
 can phrase Statement~\eqref{eq:L:weak-phi-conv-implies-phi-conv-1}
as $u = \bw_N \varphi(\ol{a}_N)$.

Since $a_1,a_2,\dotsc$ $\varphi$-converges to~$a$,
we know that the sequence $a_1,a,a_2,a,\dotsc$ is
lower $\varphi$-convergent. That is, the following exists.
\begin{equation}
\label{eq:L:weak-phi-conv-implies-phi-conv-2}
\ell\ \eqdf \ 
\bv_N \bw_{n\geq N}\  \varphi(a\wedge a_N\wedge \dotsb \wedge a_n)
\end{equation}
In particular we see that for each~$N\in\N$ the sequence
\begin{equation*}
a\wedge a_N \,\geq\, a\wedge a_N\wedge a_{N+1}\,\geq\,\dotsb
\end{equation*}
is $\varphi$-convergent.
As before, $\ul{a}_N \eqdf \bw_{n\geq N}\, a\wedge a_n$ exists,
and $\ell = \bv_N \varphi(\ul{a}_N)$.

Now, note that for each~$N\in \N$ and $n\geq N$ we have
\begin{equation*}
\ul{a}_N \ \leq\ a\wedge a_n  \ \leq\ a\vee a_n  \ \leq\ \ol{a}_N.
\end{equation*}
In particular, we have the following inequalities.
\begin{equation}
\label{eq:L:weak-phi-conv-implies-phi-conv-3}
\varphi(\ul{a}_N) \ \leq\ \varphi(a\wedge a_n)  
\ \leq\ \varphi(a\vee a_n)  \ \leq\ \varphi(\ol{a}_N).
\end{equation}
Recall that we want to prove (for some~$N$) that
$0 \ \ve\ \ld\varphi(a,a_n)$.
That is, we must show that $\varphi(a\wedge a_n) \ \ve\ \varphi(a\vee a_n)$
(see Definition~\ref{D:uniformity}\ref{E-add}).
To prove this it suffices to show 
that $\varphi(\ul{a}_N) \ \ve\ \varphi(\ol{a}_N)$
as we can see from Statement~\eqref{eq:L:weak-phi-conv-implies-phi-conv-3}
(see  Definition~\ref{D:uniformity}\ref{E-ord}).

So to complete the proof of this lemma,
we need to find an~$N\in\N$ with
\begin{equation*}
\varphi(\ul{a}_N) \ \ \ve\ \ \varphi(\ol{a}_N).
\end{equation*}

Since $a_1,a,a_2,a,\dotsc$
is $\varphi$-convergent,
we know that $u = \ell$.
Now, recall that we have $u = \bw_N \varphi(\ol{a}_N)$
and $\ell = \bw_N \varphi(\ol{a}_N)$.
Determine
an $N$ with
\begin{equation*}
\varphi(\ul{a}_N) \ \ \dt\ve\ \ \ell
\qquad\text{and}\qquad
u \ \ \dt\ve \ \ \varphi(\ol{a}_N)
\end{equation*}
using Definition~\ref{D:uniformity}\ref{E-inf-conv}
and Lemma~\ref{L:E-prop}\ref{E-sup-conv}.
Hence we see that $\varphi(\underline{a}_N) \ \ve \ \varphi(\ol{a}_N)$.
\end{proof}
%
%                  EXAMPLE ON WEAKLY PHI-CONVERGENCE
%
\begin{ex}
\label{E:weak-phi-conv}
Let $\varphi$ be a complete valuation.
We know that $\varphi$-convergence
implies weak $\varphi$-convergence (see Lemma~\ref{D:weak-phi-conv}).
The reverse implication does not hold.
Indeed,
consider the Lebesgue integral $\Lphi\colon \LF \ra \R$
and the sequence
$f_1,f_2, \dotsc$ of functions on~$\R$
given by $f_n = \frac{1}{n}\cdot\mathbf{1}_{[n,n+1]}$.
Note that 
\begin{equation*}
\textstyle{\ld\Lphi(f_n, \mathbf{0}) \ =\  \Lphi(|f_n|) 
\ =\ \frac{1}{n}.}
\end{equation*}
So we see that $f_1,f_2,\dotsc$ weakly $\Lphi$-converges to~$\mathbf{0}$.

However, we prove that the sequence $f_1,f_2,\dotsc$ 
does not $\Lphi$-converge to~$\mathbf{0}$.
Indeed, assume
(towards a contradiction) that
 $f_1,f_2,\dotsc$ does $\Lphi$-converge to~$\mathbf{0}$.
Then $f_1,f_2,\dotsc$ is $\Lphi$-convergent.
So in particular $f_1,f_2,\dotsc$ is upper $\Lphi$-convergent
(see Definition~\ref{D:seq-phi-conv}).
So we know that the following exists.
\begin{equation*}
\pulim{\Lphi}{n} f_n 
\ = \ \bw_N \bv_{n\geq N} \,\Lphi(f_N\vee \dotsb \vee f_n)
\end{equation*}
Now, note that any $N\in\N$ and $n\geq N$ we have
$f_N\vee \dotsb \vee f_n = f_N + \dotsb + f_n$,
so 
\begin{equation*}
\textstyle{\Lphi(f_N\vee \dotsb\vee f_n) 
\ =\ \Lphi(f_N) + \dotsb + \Lphi(f_n) \ = \ 
 \frac{1}{N} + \dotsb + \frac{1}{n}}.
\end{equation*}
So we see that  $\sum_{n} \frac{1}{n} 
= \bigvee_{n} \Lphi(f_1\vee\dotsb \vee f_n)$,
 which is absurd.
Hence $f_1,f_2,\dotsc$ 
does not $\Lphi$-converge to~$\mathbf{0}$.
\end{ex}
%
%                  LEMMA: SUMMABILITY + WEAK PHI-CONV ==> PHI-CONV
%
Let $\varphi$ be a complete valuation.
If $a_1,a_2,\dotsc$ weakly
$\varphi$-converges to~$a$,
then $a_1,a_2,\dotsc$ might not $\varphi$-converge to~$a$
(as we saw in Example~\ref{E:weak-phi-conv}).
However, 
there is always a subsequence of~$a_1,a_2,\dotsc$
which does $\varphi$-converge to~$a$
(see Proposition~\ref{P:weak-phi-conv-subseq}).
To prove this, we need a lemma.
\begin{lem}
\label{L:weak-phi-conv-sum}
Let $E$ be an ordered Abelian group.\\
Let $\Phi$ be a fitting uniformity on~$E$.\\
Let $L$ be a lattice, and let $\varphi\colon L \ra E$ be a valuation.\\
Let $a \in L$ and let $a_1,a_2,\dotsc$ be a sequence in~$L$ 
that weakly $\varphi$-converges to~$a$.\\
Assume that 
$\textstyle{\sum_{n}} \,\ld\varphi(a,a_n) 
\ \eqdf \ \bv_N  \sum_{n=1}^N  \,\ld\varphi(a,a_n)$ exists.\\
Then $a_1,a_2,\dotsc$ $\varphi$-converges to~$a$.
\end{lem}
\begin{proof}
To prove that $a_1,a_2,\dotsc$ $\varphi$-converges to~$a$,
we must show that $a_1,a,a_2,a,\dotsc$ is $\varphi$-convergent
(see Definition~\ref{D:seq-phi-conv-to}).
For this, 
we must show that the following exist,
\begin{alignat*}{3}
u \ &\eqdf\  \bw_N \bv_{n\geq N} \ \varphi(a \vee a_N \vee \dotsb \vee a_n)\\
\ell \ &\eqdf\  
\bv_N \bw_{n\geq N} \ \varphi(a \wedge a_N \wedge \dotsb \wedge a_n),
\end{alignat*}
and we must prove that $\ell = u$.

Let $N\in \N$ be given.
We prove that $\bv_{n\geq N} \ \varphi(a\vee a_N \vee \dotsb \vee a_n)$ exists.
Let us write $a_n ' \eqdf a\vee a_N \vee \dotsb \vee a_n$
for brevity. 
 To prove that $\bv_{n\geq N} \ \varphi(a_n')$ exists,
we want to use the fact that~$E$ is $R$-complete
(see Proposition~\ref{L:E-R-complete}).
So the task at hand is to study given $n\geq N$ the 
value $\varphi(a_{n+1}') - \varphi(a_n')$. Note that
\begin{equation}
\label{eq:L:weak-phi-conv-sum-1}
\varphi(a_{n+1}') - \varphi(a_n')
\ \ =\  \ \ld\varphi(\,a_{n+1}',\,a_n'\,)
\ \ =\ \ \ld\varphi(\ a_{n}'\vee a_{n+1},\  a_n'\vee a\ ),
\end{equation}
since $a_{n+1}' = a_{n}'\vee a_{n+1}$
and $a_n' = a_n'\vee a $ (as $a\leq a_n'$).
By Lemma~\ref{L:curry-wc-unif} we have
\begin{equation}
\label{eq:L:weak-phi-conv-sum-2}
\ld\varphi(\ a_n' \vee a_{n+1},\ a_n' \vee a\ )
\ \ \leq\ \ 
\ld\varphi(\,a_{n+1},\,a\,).
\end{equation}
So if we combine Statement~\eqref{eq:L:weak-phi-conv-sum-1}
and Statement~\eqref{eq:L:weak-phi-conv-sum-2}
we get 
\begin{equation}
\label{eq:L:weak-phi-conv-sum-3}
\varphi(a_{n+1}') - \varphi(a_n')
\ \ \leq\ \ 
\ld\varphi(a_{n+1},a).
\end{equation}
Recall that we have assumed that $\sum_n \ld\varphi(a_n,a)$ exists.
From this,
Statement~\eqref{eq:L:weak-phi-conv-sum-3},
and the fact that~$E$ is $R$-complete
it follows that $\bv_{n\geq N} \,\varphi(a_n')$ exists.

We prove that $u\eqdf 
\bw_N \bv_{n\geq N} \ \varphi(a\vee a_N\vee \dotsb \vee a_n)$ exists.
Again we use the fact that~$E$ is $R$-complete:
it is sufficient to prove that $\xi_N - \xi_{N+1} \leq \ld\varphi(a,a_N)$
where 
\begin{equation*}
\xi_N \ \eqdf\  \bv_{n\geq N} \ \varphi(a\vee a_N \vee \dotsb \vee a_n).
\end{equation*}
Let $n\in \N$ be given.
It is useful to begin by  considering the value
$\varphi(a_N'') - \varphi(a_{N+1}'')$
where 
$a_N''\eqdf a\vee a_N \vee \dotsb \vee a_n$
for all $N\leq n$. We obtain
\begin{equation*}
\varphi(a_N'') - \varphi(a_{N+1}'')
\ \ \leq\ \ \ld\varphi(a,a_N)
\end{equation*}
using a similar reasoning as before. 
Written differently, we have
\begin{equation*}
\varphi(a \vee a_N \vee\dotsb\vee a_n) 
\ \ \leq\ \ \ld\varphi(a,a_N) \ +\ 
\varphi(a \vee a_{N+1} \vee \dotsb \vee a_n)
\end{equation*}
for all~$N\in \N$ and $n\geq N$.
This implies
\begin{alignat*}{3}
\bv_{n\geq N}\ \varphi(a \vee a_N \vee\dotsb\vee a_n) 
\ \ &\leq\ \ 
\ld\varphi(a,a_N) \ +\ 
\varphi(a \vee a_{N+1} \vee \dotsb \vee a_n) \\
\ \ &\leq\ \ 
\ld\varphi(a,a_N) \ +\ 
\bv_{n\geq N+1}\ \varphi(a \vee a_{N+1} \vee \dotsb \vee a_n).
\end{alignat*}
Or in other words, $\xi_{N} \leq \ld\varphi(a,a_N) + \xi_{N+1}$.
Hence we have proven:
 \begin{equation*}
u\ \eqdf\ 
\bw_N \bv_{n\geq N} \ \varphi(a\vee a_N\vee \dotsb \vee a_n)
\quad\text{exists.}
\end{equation*}
Of course,
the above argument can be adapted to yield:
 \begin{equation*}
\ell\ \eqdf\ 
\bv_N \bw_{n\geq N} \ \varphi(a\wedge a_N\wedge \dotsb \wedge a_n)
\quad\text{exists.}
\end{equation*}

\vspace{.3em}
\noindent It remains to be shown that $\ell = u$.
Let $\ve \in \Phi$ be given.
Since reader can easily verify that $\ell \leq u$,
to prove that $\ell = u$,
it suffices to show that $\ell \se u$
(see Definition~\ref{D:uniformity}\ref{E-haus}).
Let~$N\in\N$ be given. Note that we have the following inequalities.
\begin{equation*}
\bw_{n\geq N}\ \varphi(a\wedge a_N\wedge \dotsb \wedge a_n)
\ \ \leq\ \ 
\ell \ \leq\ u
\ \ \leq\ \ 
\bv_{n\geq N}\ \varphi(a\vee a_N\vee \dotsb \vee a_n).
\end{equation*}
So to prove $\ell \se u$,
it suffices to show that for some~$N$
(see Definition~\ref{D:uniformity}\ref{E-ord}),
\begin{equation}
\label{eq:L:weak-phi-conv-sum-4}
\bw_{n\geq N}\ \varphi(a\wedge a_N\wedge \dotsb \wedge a_n)
\quad\ve\quad
\bv_{n\geq N}\ \varphi(a\vee a_N\vee \dotsb \vee a_n).
\end{equation}
Since $\sum_n \ld\varphi(a,a_n)$ exists,
we can find an~$N\in\N$ such that
(see Lemma~\ref{L:E-prop}\ref{E-sup-conv})
\begin{equation}
\label{eq:L:weak-phi-conv-sum-4b}
0\quad \dtn\ve4 \quad \ld\varphi(a,a_N) \,+\, \dotsb \,+\, \ld\varphi(a,a_n)
\qquad\quad(n\geq N).
\end{equation}
We will prove that Statement~\eqref{eq:L:weak-phi-conv-sum-4}
holds for this~$N$.
Since $\bw_{n\geq N}\ \varphi(a\wedge a_N\wedge \dotsb \wedge a_n)$
and 
$\bv_{n\geq N}\ \varphi(a\vee a_N\vee \dotsb \vee a_n)$
exist,
we can find $n\geq N$ such that
\begin{alignat*}{3}
\bw_{n\geq N}\ \varphi(a\wedge a_N\wedge \dotsb \wedge a_n)
\quad&\dtn\ve4 \quad
\varphi(a\wedge a_N\wedge \dotsb \wedge a_n)\\
\varphi(a\vee a_N\vee \dotsb \vee a_n)
\quad&\dtn\ve4 \quad
\bv_{n\geq N}\ \varphi(a\vee a_N\vee \dotsb \vee a_n).
\end{alignat*}
So to prove that Statement~\eqref{eq:L:weak-phi-conv-sum-4}
it suffices to show that
\begin{equation}
\label{eq:L:weak-phi-conv-sum-5}
\varphi(a\wedge a_N\wedge \dotsb \wedge a_n)
\quad\dt\ve\quad
\varphi(a\vee a_N\vee \dotsb \vee a_n).
\end{equation}
Note that we have the following inequalities.
\begin{equation*}
\varphi(a\wedge a_N\wedge \dotsb \wedge a_n)
\ \ \leq\ \ \varphi(a)\ \ \leq\ \ 
\varphi(a\vee a_N\vee \dotsb \vee a_n).
\end{equation*}
So to prove that Statement~\eqref{eq:L:weak-phi-conv-sum-5}
holds, we will show that
\begin{equation}
\label{eq:L:weak-phi-conv-sum-6}
\varphi(a\wedge a_N\wedge \dotsb \wedge a_n)
\quad\dtn\ve4\quad\varphi(a)\quad\dtn\ve4\quad
\varphi(a\vee a_N\vee \dotsb \vee a_n).
\end{equation}
Now, note that vteration of Lemma~\ref{L:wv-unif} yields
\begin{alignat*}{3}
\varphi(a\wedge a_N\wedge \dotsb \wedge a_n)
\,-\, \varphi(a)
\ \ &=\ \ 
\ld\varphi(\ a\wedge a_N\wedge\dotsb \wedge a_n,\ a\ )  \\
\ \ &=\ \ 
\ld\varphi(\ a\wedge a_N\wedge\dotsb \wedge a_n,
                      \ a\wedge a\wedge \dotsb \wedge a \  ) \\
\ \ &\leq\ \ 
\ld\varphi(a_N,a) \,+\,\dotsb \,+\,\ld\varphi(a_n,a).
\end{alignat*}
So by Statement~\eqref{eq:L:weak-phi-conv-sum-4b}
and Definition~\ref{D:uniformity}\ref{E-add} we get
\begin{equation*}
\varphi(a\wedge a_N\wedge \dotsb \wedge a_n)
\quad\dtn\ve4\quad\varphi(a).
\end{equation*}
Using a similar argument,
we obtain 
\begin{equation*}
\varphi(a)\quad\dtn\ve4\quad
\varphi(a\vee a_N\vee \dotsb \vee a_n).
\end{equation*}
So we have proven Statement~\eqref{eq:L:weak-phi-conv-sum-6} and 
thereby completed the proof.
\end{proof}
%
%                  PROPOSITION ON PHI-CONV SUBSEQUENCE OF A WEAK PHI-CONV SEQ
%
\begin{prop}
\label{P:weak-phi-conv-subseq}
Let $E$ be an ordered Abelian group.\\
Let $\Phi$ be a fitting uniformity on~$E$.\\
Let $L$ be a lattice, and $\varphi\colon L \ra E$ a valuation.\\
Let $a_1,a_2,\dotsc$ be a sequence in~$L$
that weakly $\varphi$-converges to some~$a\in L$.\\
Then: there are $j_1<j_2<\dotsb$ in~$\N$
such that $a_{j_1},a_{j_2},\dotsc$
$\varphi$-converges to~$a$.
\end{prop}
\begin{proof}
It suffices to find $j_1<j_2<\dotsb$ in~$\N$ such that
$\sum_k \ \ld\varphi(a,a_{j_k})$ exists.
Indeed,
then we have $a_{j_1},a_{j_2},\dotsc$
weakly $\varphi$-converging to~$a$
since $a_1,a_2,\dotsc$ $\varphi$-converges to~$a$.
So by Lemma~\ref{L:weak-phi-conv-sum} we 
$a_{j_1},a_{j_2},\dotsc$ $\varphi$-converges to~$a$, as we must prove.

Let $\ve_1',\ve_2',\dotsc$ be an enumeration of~$\Phi$.
Pick $\ve_1,\ve_2,\dotsc$ in~$\Phi$ such that for all~$n$,
\begin{equation*}
\ve_{n} \,\leq\, \ve_n ' \qquad\text{and}\qquad  \ve_{n+1}\,\leq\,\dt{\ve_{n}}.
\end{equation*}
Note that for all $N\in \N$ and $n\geq N+1$
we have
(see Notation~\ref{N:unif}),
\begin{equation}
\label{eq:P:weak-phi-conv-subseq-1}
\ve_{N+1} + \dotsb + \ve_{n} \ \leq\  \ve_N.
\end{equation}

Pick $j_1 < j_2 < \dotsb$ in~$\N$ such that for all~$k\in\N$, 
\begin{equation*}
0\quad\ve_{k+1}\quad \ld\varphi(a, a_{j_k}).
\end{equation*}
Then by Statement~\eqref{eq:P:weak-phi-conv-subseq-1} 
and Definition~\ref{D:uniformity}\ref{E-add} for all~$N\in\N$ and $n\geq N$,
\begin{equation}
\label{eq:P:weak-phi-conv-subseq-2}
0\quad\ve_{N}\quad 
\ld\varphi(a, a_{j_N}) \,+\, \dotsb \,+\, \ld\varphi(a,a_{j_n}).
\end{equation}

Recall that we need to prove that $\sum_k\ \ld\varphi(a,a_{j_k})$ exists.
Let $\ve \in \Phi$ be given.
By Lemma~\ref{L:E-prop}\ref{E-bound-sup}
it suffices to find an $N\in\N$ such that for all~$n\geq N$,
\begin{equation}
\label{eq:P:weak-phi-conv-subseq-3}
0\quad\ve\quad \ld\varphi(a,a_{j_N}) \,+\, \dotsb \,+\, \ld\varphi(a,a_{j_n}).
\end{equation}
Since $\ve_1',\ve_2',\dotsc$ enumerates~$\Phi$
we can find an $N\in \N$ such that $\ve_N'=\ve$.
Recall that $\ve \leq \ve_N' \leq \ve_N$.
So Statement~\eqref{eq:P:weak-phi-conv-subseq-3}
follows directly from Statement~\eqref{eq:P:weak-phi-conv-subseq-2}.
\end{proof}
%
%                  LEMMA ON DIAGONALISATION OF PHI-CONV SEQUENCES
%
\begin{lem}
\label{L:phi-conv-diag}
Let $E$ be an ordered Abelian group.\\
Let $\Phi$ be a fitting uniformity on~$E$.\\
Let $L$ be a lattice, and $\varphi\colon L \ra E$ a valuation.\\
Let $a_1,a_2,\dotsc$ be a sequence in~$L$
which weakly $\varphi$-converges to some~$a\in L$\\
For each~$N\in\N$, let $b^N_1, b^N_2, \dotsc$
be a sequence in~$L$ which weakly $\varphi$-converges to~$a_N$.\\
Then there are $j_1 < j_2 < \dotsb$ in~$\N$ such that
 $b^1_{j_1},\,b^2_{j_2},\,\dotsc$
weakly $\varphi$-converges to~$a$.
\end{lem}
\begin{proof}
To find a suitable sequence $j_1 < j_2 <\dotsb$ we need some preparation.

We know that~$\Phi$ is countable
(see Definition~\ref{D:fitting-uniformity}).
Let $\ve_1',\,\ve_2',\,\dotsc$ be an enumeration of~$\Phi$.
Define a sequence $\ve_1 \geq \ve_2 \geq \dotsb$ in~$\Phi$
(see Notation~\ref{N:unif-plus})
by 
\begin{equation*}
\ve_1 \ = \ \ve_1'
\qquad\text{and}\qquad 
\ve_{n+1} \ =\ \ve_n \wedge \ve_{n+1}'.
\end{equation*}
Note that we have $\ve_n \leq \ve_n'$ for all~$n$.

Let~$N\in\N$ be given.
Since $b^N_1,b^N_2,\dotsc$ weakly $\varphi$-converges to~$a_N$,
we know by Definition~\ref{D:weak-phi-conv} that there is an~$M\in\N$
such that $\ld\varphi(b^N_n, a_N)\ \ve_N\ 0$
for all~$n\geq M$.

Now,
choose $j_1<j_2<\dotsb$ such that 
$\ld\varphi(b^N_n, a_N)\ \ve_N\ 0$
for all~$n\geq j_N$.

We will prove that $b^1_{j_1},\,b^2_{j_2},\dotsc$ weakly $\varphi$-converges
to~$a$.
Let $\ve\in\Phi$ be given.
We must find an~$\mathfrak{n}\in\N$
such that $\ld\varphi(b^N_{j_N},a_N) \ \ve\ 0$
for all~$N\geq \mathfrak{n}$
(see Definition~\ref{D:phi-conv-2}).

Find an~$k\in\N$ such that $\dt\ve = \varepsilon_k'$.
(Recall that $\ve_1',\ve_2',\dotsc$ enumerates~$\Phi$.)

Pick~$\mathfrak{n}\geq k$
such that $\ld\varphi(a_N,a) \ \dt\ve\ 0$.
We prove that $\ld\varphi(b_{j_N}^N,a) \ \ve\ 0$
for all~$N\geq \mathfrak{n}$.
Let $N\geq\mathfrak{n}$ be given.
We have $\ld\varphi(b^N_{j_N},a_N) \ \ve_N \ 0$
by choice of~$j_N$.
So in particular $\ld\varphi(b^N_{j_N},a_N) \ \dt\ve \ 0$
since $\dt\ve = \ve_k'\geq\ve_k \geq \ve_{\mathfrak{n}} \geq \ve_N$.

Further,
we have $\ld\varphi(a_N,a) \ \dt\ve\ 0$
since $N\geq \mathfrak{n}$.

So by property~\ref{E-add} 
of a fitting uniformity (see Definition~\ref{D:fitting-uniformity})
we have
\begin{equation*}
\ld\varphi(b_{j_N}^N,a_N) \,+\,  \ld\varphi(a_N,a) 
\quad\dt\ve\quad 
\ld\varphi(a_N,a)
\quad\dt\ve\quad
0.
\end{equation*}
So by property~\ref{E-half} of a fitting uniformity
we have
\begin{equation*}
\ld\varphi(b_{j_N}^N,a_N) \,+\,  \ld\varphi(a_N,a) 
\quad\ve\quad 
0.
\end{equation*}
Now,
by points~\ref{d-metric_pos} and~\ref{d-metric_triangle}
of Lemma~\ref{L:d-metric}
we get
\begin{equation*}
0 \ \leq\ 
\ld\varphi(b_{j_N}^N, a) \ \leq\ 
\ld\varphi(b_{j_N}^N,a_N) \,+\,  \ld\varphi(a_N,a).
\end{equation*}
So by property~\ref{E-ord} we get $\ld\varphi(b^N_{j_N},a) \ \ve\ 0$.
\end{proof}
%
%                  LEMMA ON WEAK PHI-CONVERGENCE AND THE BINARY INFIMUM
%
\begin{lem}
\label{L:phi-conv-inf}
Let $E$ be an ordered Abelian group.\\
Let~$\Phi$ be a fitting uniformity on~$E$.\\
Let~$L$ be a lattice, and $\varphi\colon L\ra E$ a valuation.
Let $a,b\in L$ be given.\\
Let $a_1,a_2,\dotsc$ 
be a sequence in~$L$ which weakly $\varphi$-converges to~$a$.\\
Let $b_1,b_2,\dotsc$
be a sequence in~$L$ which weakly $\varphi$-converges to~$b$.\\
Then $a_1\wedge b_1,\ a_2\wedge b_2,\ \dotsc$
weakly $\varphi$-converges to~$a\wedge b$,\\
and $a_1\vee b_1,\ a_2\vee b_2,\ \dotsc$
weakly $\varphi$-converges to~$a\vee b$.
\end{lem}
\begin{proof}
We will only prove that $a_1\wedge b_1,\ a_2\wedge b_2,\ \dotsc$
weakly $\varphi$-converges to~$a\wedge b$.

Let $\ve \in \Phi$ be given.
To prove  $a_1\wedge b_1,\ a_2\wedge b_2,\ \dotsc$
weakly $\varphi$-converges to~$a\wedge b$,
we must find an~$N\in\N$
 such that 
\begin{equation}
\label{L:phi-conv-inf-1}
0\quad \ve\quad \ld\varphi(\,a_n\wedge b_n,\ a\wedge b\,)
\qquad\qquad(n\geq N).
\end{equation}
Since $a_1,a_2, \dotsc$ weakly $\varphi$-converges to~$a$
and $b_1,b_2,\dotsc$ weakly $\varphi$-converges to~$b$
we know
that there is an~$N\in\N$ such that 
\begin{equation}
\label{L:phi-conv-inf-2}
0\ \ \dt\ve\ \ \ld\varphi(a_n, a)
\qquad\text{and}\qquad
0\ \ \dt\ve\ \ \ld\varphi(b_n,b)
\qquad\qquad(n\geq N).
\end{equation}
We will prove that Statement~\eqref{L:phi-conv-inf-1}
holds for this~$N$.

To this end, note by Lemma~\ref{L:d-metric}\ref{d-metric_pos}
and Lemma~\ref{L:wv-unif} we have, for all~$n\in \N$,
\begin{equation*}
0 \ \ \leq\ \ 
\ld\varphi(\,a_n\wedge b_n,\ a\wedge b\,)
\ \ \leq\ \ 
\ld\varphi(a_n, a) \,+\,
\ld\varphi(b_n, b)
\end{equation*}
So by property~\ref{E-ord} of~$\Phi$
(see Def.~\ref{D:uniformity})
to prove~\eqref{L:phi-conv-inf-1}
it suffices to show  that
\begin{equation}
\label{L:phi-conv-inf-3}
0\quad\ve\quad
\ld\varphi(a_n, a) \,+\,
\ld\varphi(b_n, b)
\end{equation}
for any $n\geq N$.
By Statement~\eqref{L:phi-conv-inf-2} 
and property~\ref{E-add} of~$\Phi$, we get
\begin{equation*}
0 
\ \ \dt\ve  \ \ 
\ld\varphi(a_n, a)
\quad \dt\ve \quad
\ld\varphi(a_n, a) \,+\,
\ld\varphi(b_n, b)
\end{equation*}
for all~$n\geq N$. 
So we see that Statement~\eqref{L:phi-conv-inf-3} 
holds by property~\ref{E-add} of $\Phi$.
\end{proof}
%
%                  COUNTEREXAMPLE ON THE PREVIOUS LEMMA FOR PHI-CONV
%
\begin{ex}
\label{E:phi-conv-inf}
Given Lemma~\ref{P:R-main},
one might surmise that Lemma~\ref{L:phi-conv-inf}
holds if one replaces ``weakly $\varphi$-converges''
by
``$\varphi$-converges''.
This is not the case, as we will show.

Recall that  Lebesgue integral $\Lphi\colon \LF \ra \R$
is a valuation.
For all~$n\in \N$, define
\begin{equation*}
f_n \ = \ (-1)^{n} \cdot \textstyle{\frac{1}{n}}\cdot \mathbf{1}_{[n,n+1]}.
\end{equation*}
Then $f_n\in \LF$ for all~$n$, 
and the sequence $f_1,f_2,\dotsc $ $\Lphi$-converges to~$\mathbf{0}$.
As one expects,
the sequence $-f_1,-f_2,\dotsc$ 
also
 $\Lphi$-converges to~$\mathbf{0}$.
However, the sequence 
\begin{equation*} 
f_1\vee(-f_1),\  f_2\vee(-f_2),\ \dotsc
\end{equation*}
does not $\Lphi$-converge to~$\mathbf{0}$,
because $f_n \vee (-f_n) = \frac{1}{n}\cdot \mathbf{1}_{[n,n+1]}$
(see Example~\ref{E:weak-phi-conv}).
\end{ex}
%
%                  LEMMA ON WEAK DENSENESS 
%
\begin{lem}
\label{L:fitting-phi-conv-dense}
Let $E$ be an ordered Abelian group.\\
Let $\Phi$ be a fitting uniformity on~$E$.\\
Let $\vs{V}{L}\varphi{E}$ be an extendible valuation system.
Let $a\in \ol L$ be given.\\
Then there is a sequence
$a_1,a_2,\dotsc$ in~$L$
that weakly $\ol\varphi$-converges to~$a$.
\end{lem}
\begin{proof}
By Corollary~\ref{C:aleph1}
we know that~$\overline L = \Pi_{\aleph_1} L$.
So it suffices to prove that
the following statement holds for every ordinal number~$\alpha$.
\begin{equation*}
\left[\quad
\begin{minipage}{.7\columnwidth}
Let $a\in \Pi_\alpha L\cup \Sigma_\alpha L$ be given.\\
There is a sequence $a_1,a_2,\dotsc$ in~$L$
that $\ol\varphi$-converges to~$a$.
\end{minipage}
\right.
\end{equation*}
Let us name the above statement~$P(\alpha)$.
We prove $\forall\alpha\ P(\alpha)$ with induction.

Clearly, $P(0)$ holds,
since $\Pi_0 L = L = \Sigma_0 L$.

\vspace{.3em}
Let $\alpha$ be an ordinal  such that $P(\alpha)$ holds.
We prove that $P(\alpha+1)$ holds.
Let $a\in \Pi_{\alpha+1}L \cup \Sigma_{\alpha+1} L$ be given.
We must find a sequence in~$L$ that
$\ol\varphi$-converges to~$a$.

Assume that $a\in \Pi_{\alpha+1} L$.
There is a $\Sigma_\alpha\varphi$-convergent
sequence $b_1 \geq b_2 \geq \dotsb$ in $\Sigma_\alpha L$
such that $\bw_n b_n = a$.
In particular, $a_1,a_2,\dotsc$ $\ol\varphi$-converges to~$a$.
Since $P(\alpha)$ holds,
we can find for each~$N\in \N$ a
sequence $b^N_1, b^N_2,\dotsc$ in~$L$
that $\ol\varphi$-converges to~$a_N$.
Then by Lemma~\ref{L:phi-conv-diag}
there are $j_1 < j_2 < \dotsb$ in $\N$
such that $b^1_{j_1},\, b^2_{j_2},\,\dotsc$
$\ol\varphi$-converges to~$a$.
So we see that there is a sequence in~$L$
that $\ol\varphi$-converges to~$a$.

By a similar reasoning
we see that if $a\in \Sigma_{\alpha+1}L$
then there is a sequence in~$L$ 
that $\ol\varphi$-converges to~$a$.
Hence $P(\alpha+1)$.

\vspace{.3em}
Let $\lambda$ be a limit ordinal
such that $P(\alpha)$ holds
for all~$\alpha<\lambda$.
We must prove that $P(\lambda)$ holds.
Let $a\in \Pi_\lambda L\cup \Sigma_\lambda L$ be given.
We must find a sequence in~$L$ that $\ol\varphi$-converges to~$a$.
By definition of~$\Pi_\lambda L$ and $\Sigma_\lambda L$,
there is an~$\alpha<\lambda$
such that $a \in \Pi_\alpha L \cup \Sigma_\alpha L$
(see Definition~\ref{P:hier}).
Since we know that~$P(\alpha)$ holds,
there must be a sequence in~$L$ that $\ol\varphi$-converges to~$a$.
Hence~$P(\lambda)$.
\end{proof}
%
%                  TECHNICAL COROLLARY TO THE DENSENESS LEMMA
%
\begin{cor}
\label{C:fitting-phi-conv-dense}
Let $E$ be an ordered Abelian group.\\
Let $\Phi$ be a fitting uniformity on~$E$.\\
Let $\vs{V}{L}\varphi{E}$ be an extendible valuation system.
Let $a,b\in \ol L$ with $a\leq b$ be given.\\
Then there is a sequence
$a_1,a_2,\dotsc$ in~$L$ that weakly $\ol\varphi$-converges to~$a$,\\
and there is a sequence
$b_1,b_2,\dotsc$ in~$L$ that weakly $\ol\varphi$-converges to~$b$,\\
such that $a_n \leq b_n$ for all~$n\in\N$.
\end{cor}
\begin{proof}
Let $a,b\in\ol L$ with $a\leq b$ be given.
Using Lemma~\ref{L:fitting-phi-conv-dense}
find a sequence $a_1',a_2',\dotsc$ in~$L$
that weakly $\ol\varphi$-converges to~$a$
and find a sequence $b_1',b_2',\dotsc$ in~$L$
that weakly $\ol\varphi$-converges to~$b$.
Consider the sequences
$a_1,a_2,\dotsc$ and $b_1,b_2,\dotsc$  in~$L$ 
given by  
\begin{equation*}
a_n \eqdf a_n'\wedge b_n'
\qquad\text{and}\qquad
b_n \eqdf a_n' \vee b_n'
\qquad\qquad(n\in\N).
\end{equation*}
Clearly $a_n \leq b_n$ for all~$n\in \N$.
Moreover,
by Lemma~\ref{L:phi-conv-inf} 
we know that $a_1,a_2,\dotsc$ 
weakly $\ol\varphi$-converges to~$a = a\wedge b$
and weakly $b_1,b_2,\dotsc$ $\ol\varphi$-converges to~$b= a\vee b$.
\end{proof}
%
%                  DENSENESS PROPOSITION 
%
\begin{prop}
\label{P:fitting-phi-conv-dense}
Let $E$ be an ordered Abelian group.\\
Let $\Phi$ be a fitting uniformity on~$E$.\\
Let $\vs{V}{L}\varphi{E}$ be an extendible valuation system.
\begin{enumerate}
\item
Given $a\in \ol L$,
 there is a sequence
$a_1,a_2,\dotsc$ in~$L$
which $\ol\varphi$-converges to~$a$.

\item
Let $a,b\in \ol L$ with $a\leq b$ be given.\\
There is a sequence 
$a_1,a_2,\dotsc$ in~$L$
which $\varphi$-converges to~$a$, and \\
there is a sequence
$b_1,b_2,\dotsc$ in~$L$
which $\varphi$-converges to~$b$, such that
\begin{equation*}
a_n \ \leq\ b_n \qquad\quad(n\in\N).
\end{equation*}
\end{enumerate}
\end{prop}
\begin{proof}
Combine Lemma~\ref{L:fitting-phi-conv-dense},
Corollary~\ref{C:fitting-phi-conv-dense},
and Proposition~\ref{P:weak-phi-conv-subseq}.
\end{proof}

%
%                  EXTENSION THEOREM USED BY FUBINI
%
\begin{thm}
\label{T:fubext}
Let $E$ be a lattice ordered Abelian group.\\
Let $\Phi$ be a fitting uniformity on~$E$.\\
Let $\vs{V}{L}\varphi{E}$ be an extendible valuation system.\\
Let $\psi\colon C\ra E$ be 
a complete Hausdorff valuation.\\
Let $f\colon L\ra C$ be an order preserving map such that 
$\psi\circ f = \varphi$.\\
Then there is a unique
order preserving extension~$g\colon \ol L \ra C$ of~$f$
such that $\psi\circ g = \ol\varphi$.
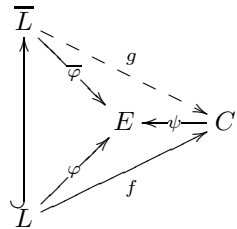
\begin{equation*}
\xymatrix{
\ol L \ar [rd] |{\ol\varphi} 
\ar @{-->} [rrd] ^{g} \\
&E&C\ar[l]|{\psi}\\
L \ar [rru] _{f}
\ar [ru] |{\varphi}
\ar @{^(->} [uu] 
}
\end{equation*}
\end{thm}
\begin{proof}
(Uniqueness)
Let $g_1,g_2\colon \ol L \ra C$
be order preserving extensions of~$f$ such that
$\psi\circ g_i = \ol\varphi$.
We prove that $g_1 = g_2$.
Let~$a\in \ol L$ be given.
By Proposition~\ref{P:fitting-phi-conv-dense}
there is a $\varphi$-convergent sequence $a_1,a_2,\dotsc$
in~$L$ which $\ol\varphi$-converges to~$a$.

We must show that~$g_1(a) = g_2(a)$.
To this end we prove that $f(a_1),\,f(a_2),\,\dotsc$
weakly $\psi$-converges to~$g_i(a)$.
Then  $g_1(a) = g_2(a)$
because~$\psi$ is Hausdorff.

Let~$i\in\{1,2\}$ be given.
Note that $g_i(a_n) = f(a_n)$ because $g_i$ extends~$f$.
So we must prove that $g_i(a_1),\,g_i(a_2),\,\dotsc$
weakly $\psi$-converges to~$g_i(a)$.
Let~$\ve\in\Phi$ be given.
We must find 
an~$N\in \N$
such that
\begin{equation}
\label{eq:T:fubext-1}
\ld\psi(\,g_i(a_n),\,g_i(a)\,) \quad\ve\quad 0
\qquad\quad(n\geq N).
\end{equation}
Since $a_1,a_2,\dotsc$ $\ol\varphi$-converges to~$a$
we know that $a_1,a_2,\dotsc$ weakly $\ol\varphi$-converges to~$a$
so we know there is an $N\in \N$ with
\begin{equation*}
\ld{\ol\varphi}(a_n,a) \quad\ve\quad 0\qquad\quad(n\geq N).
\end{equation*}
We will prove that Statement~\eqref{eq:T:fubext-1}
holds for this~$N$.

Let $n\geq N$ be given.
Since $g$ is order preserving, we have
\begin{equation*}
g_i(a_n \wedge a) \ \leq\ g_i(a_n)\wedge g_i(a)
\qquad\qquad
g_i(a_n) \vee g_i(a) \ \leq\  g_i (a_n \vee a).
\end{equation*}
In particular
(recall that $\psi\circ g_i=\ol\varphi$),
\begin{equation}
\label{eq:T-fubext-3}
\begin{alignedat}{3}
\ld\psi(\,g_i(a_n),\,g_i(a)\,) 
\ &=\ 
\psi(g_i(a_n) \vee g_i(a)) \,-\, \psi(g_i(a_n)\wedge g_i(a)) \\
\ &\leq\ 
\psi(g_i(a_n\vee a)) \,-\, \psi(g_i(a_n \wedge a)) \\
\ &=\ 
\ol\varphi(a_n \vee a) \,-\, \ol\varphi(a_n\wedge a) 
\ =\ \ld{\ol\varphi}(a_n,a).
\end{alignedat}
\end{equation}
So we
know that $0\ \ve\ \ld{\ol\varphi}(a_n,a)$
and we have the following inequalities.
\begin{equation*}
0 \ \leq\ \ld{\psi}(\,g_i(a_n),\,g_i(a)\,) \ \leq\ \ld{\ol\varphi}(a_n,a)
\end{equation*}
Hence $0 \ \ve\ \ld{\psi}(g_i(a_n),g_i(a))$
by property~\ref{E-ord} of a fitting uniformity.
So we have shown that Statement~\eqref{eq:T:fubext-1} holds.
Thus, $g_1=g_2$.

\vspace{.3em}(Existence)
We will prove the following statement.
\begin{equation}
\label{eq:T-fubext-2}
\left[\quad
\begin{minipage}{.7\columnwidth}
Let $a\in \ol L$ be given.
There is a unique $b\in C$ 
such that for every sequence $a_1,a_2,\dotsc$ in~$L$
that $\ol\varphi$-converges to~$a$,
we have $f(a_1),\,f(a_2),\,\dotsc$
$\psi$-converges to~$b$.
\end{minipage}
\right.
\end{equation}
Of course,
we will later define~$g\colon \ol L \ra C$ by  $g(a) = b$.

Let $a\in\ol L$ be given.
For each~$i\in\{1,2\}$,
let $b_i\in C$
and $a^i_1,a^i_2,\dotsc\in L$
be given,
such that 
$a^i_1,a^i_2,\dotsc$
 $\ol\varphi$-converges to~$a$,
and $f(a^i_1),\,f(a^i_2),\dotsc$
$\psi$-converges to~$b_i$.

We must prove that $b_1=b_2$.
Let $\ve\in\Phi$ be given.
Since $\psi$ is Hausdorff,
it suffices to show that $0\ \ve\ \ld\psi(b_1,b_2)$
(see~\ref{D:hausdorff}).

Note that
by
points \ref{d-metric_pos}
and \ref{d-metric_triangle}
of Lemma~\ref{L:d-metric}
we have
\begin{equation*}
0 \ \leq\ 
\ld\psi(b_1,b_2)
\ \ \leq\ \ 
\ld\psi(b_1,\,f(a^1_n)\,)
\ +\  
\ld\psi(\,f(a^1_n),\,f(a^2_n)\,)
\ +\ 
\ld\psi(\,f(a^2_n),\,b_2).
\end{equation*}
So to prove 
 $0\ \ve\ \ld\psi(b_1,b_2)$,
it sufficient to find~$N\in\N$
such that for all~$n\geq N$
the following statement holds
(see Definition~\ref{D:uniformity},
points~\ref{E-ord}, \ref{E-half} and \ref{E-add}).
\begin{alignat*}{3}
0\quad &\dtn\ve4 \quad \ld\psi(b_1,\,f(a^1_n)\,),
\quad\text{and}\\
0\quad &\dtn\ve4 \quad \ld\psi(\,f(a^1_n),\,f(a^2_n)\,),
\quad\text{and}\\
0\quad&\dtn\ve4 \quad \ld\psi(\,f(a^2_n),\,b_2).
\end{alignat*}

Recall that $f(a^i_1),\,f(a^i_2),\,\dotsc$
$\psi$-converges to~$b_i$ for all~$i$.
Hence  $f(a^i_1),\,f(a^i_2),\,\dotsc$
weakly $\psi$-converges to~$b_i$ for all~$i$.
So we know  there is an~$N\in\N$ such that
$0 \ \dtn\ve4\ \ld\psi(f(a^i_n),\,b_i)$
for all~$n\geq N$ and~$i$.

It remains to be shown that there is an~$N\in \N$
such that 
$0\  \dtn\ve4 \  \ld\psi(f(a^1_n),\,f(a^2_n))$
for all $n\geq N$.
To this end, note that~$f$ is order preserving
and that $\psi \circ f = \varphi$.
So with a similar reasoning as before (see Statement~\eqref{eq:T-fubext-3}),
we see that
\begin{equation*}
\ld\psi(\,f(a^1_n),\,f(a^2_n)\,)
\ \leq\ 
\ld\varphi(a^1_n, a^2_n).
\end{equation*}
So to complete the proof Statement~\eqref{eq:T-fubext-2}
it suffices to find an~$N\in\N$ such that 
\begin{equation}
\label{eq:T:fubext-5}
0\quad\dtn\ve4\quad\ld\varphi(a^1_n, a^2_n)
\qquad\quad(n\geq N).
\end{equation}
Note that
by
points \ref{d-metric_pos}
and \ref{d-metric_triangle}
of Lemma~\ref{L:d-metric}
we have
\begin{equation*}
0\ \leq\ 
\ld\varphi(a_n^1,a_n^2) 
\ \leq\ 
\ld\varphi(a_n^1,a)
\,+\,
\ld\varphi(a,a_n^2).
\end{equation*}
Since the sequence $a^i_1,a^i_2,\dotsc$
$\ol\varphi$-converges to~$a$
(and hence also weakly)
we can find an~$N\in\N$
such that $0\ \dtn\ve4 \ \ld\varphi(a_n^i,a)$
 for all~$n\geq N$ and~$i\in\{1,2\}$.

So by points~\ref{E-ord}, \ref{E-half} and \ref{E-add}
of Definition~\ref{D:uniformity}
we see that Statement~\eqref{eq:T:fubext-5} holds.

Hence we have proven Statement~\eqref{eq:T-fubext-2}.
So we now know there is a unique map~$g\colon \ol L \ra C$
such that 
for every $a\in\ol L$
and every sequence $a_1,a_2,\dotsc$ in~$L$ that $\ol\varphi$-converges to~$a$
we have $f(a_1),\,f(a_2),\,\dotsc$ 
$\psi$-converges to~$g(a)$.

To complete the proof of this theorem,
we show that $g$ extends~$f$,
we show that $g$ is order preserving,
and that $\psi \circ g = \ol \varphi$.

Let $a\in L$ be given. To prove that $g$ extends~$f$
we show that $g(a)=f(a)$.
Note that $a,a,\dotsc$ $\ol\varphi$-converges to~$a$.
So by definition of~$g$ we know that $f(a),\,f(a),\,\dotsc$
$\psi$-converges to~$g(a)$.
But $f(a),\,f(a),\,\dotsc$ $\psi$-converges to~$f(a)$
too, and $\psi$ is Hausdorff.
So we see that $f(a)=g(a)$.

Let $a,b\in L$ with $a\leq b$ be given.
To prove that $g$ is order preserving 
we must show that $g(a)\leq g(b)$.
By Proposition~\ref{P:fitting-phi-conv-dense}
we can find a sequence $a_1,a_2,\dotsc$ in~$L$
that $\ol\varphi$-converges to~$a$
and a sequence $b_1,b_2,\dotsc$ in~$L$
that $\ol\varphi$-converges to~$b$
such that we have $a_n \leq b_n$ for all~$n\in\N$.
Now, note that  by Lemma~\ref{L:phi-conv-inf}
we know that
\begin{equation*}
f(a_1) \wedge f(b_1),\quad f(a_2)\wedge f(b_2),\quad \dotsc
\qquad\text{\emph{weakly} $\psi$-converges to}\quad g(a)\wedge g(b).
\end{equation*}
Let $n\in \N$ be given.
Since $f$ is order preserving
and $a_n \leq b_n$
we have $f(a_n)\leq f(b_n)$
and so $f(a_n) \wedge f(b_n) = f(a_n)$.
Hence  $f(a_1),\,f(a_2),\,\dotsc$
weakly $\psi$-converges to both $g(a)$ and $g(a)\wedge g(b)$.
So we see that $g(a) = g(a)\wedge g(b)$ 
and thus $g(a)\leq g(b)$.

Let $a\in\ol L$ be given.
We show that $\psi(g(a)) = \ol\varphi(a)$.
Find a sequence $a_1,a_2,\dotsc$ in~$L$
that $\ol\varphi$-converges to~$a$
(see Lemma~\ref{L:fitting-phi-conv-dense}).

Recall that $E$ is a \emph{lattice} ordered Abelian group.
By Theorem~\ref{T:lebesgue}
we see that
\begin{equation}
\label{eq:T:fubext-7}
\ol\varphi(a) \ = \ \textstyle{\lim_n} \varphi(a_n).
\end{equation}
By definition of~$g$
we have $f(a_1),\,f(a_2),\,\dotsc$
$\psi$-converges to~$g(a)$.
So by Theorem~\ref{T:lebesgue}
we have  $\psi(g(a)) = \lim_n \psi(f(a_n)) $.
But $\psi(f(a_n)) = \varphi(a_n)$ for all~$n\in\N$,
so we have
\begin{equation}
\label{eq:T:fubext-8}
\psi(f(a))
\ =\  \textstyle{\lim_n} \varphi(a_n).
\end{equation}
If we combine Equalities~\eqref{eq:T:fubext-7}
and~\eqref{eq:T:fubext-8} we get
$\ol\varphi(a) = \psi(f(a))$.
\end{proof}
%%%%%%%%%%%%%%%%%%%%%%%%%%%%%%%%%%%%%%%%%%%%%%%%%%%%%%%%%%%%%%%%%%%%%%%%%%%%%%%
%
 %                 FUBINI'S THEOREM
\begin{thm}
\label{T:fub}
Statement~\eqref{eq:fub} holds.
\end{thm}
\begin{proof}
We only give hints and leave the details to the reader.
With the notation of Subsection~\ref{SS:fub_part1}
apply Theorem~\eqref{T:fubext}
to the following situation.
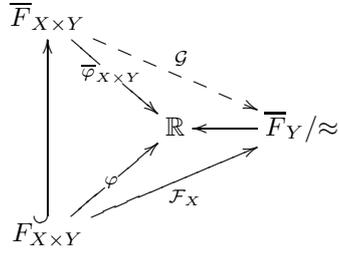
\begin{equation*}
\xymatrix{
\ol{F}_{X\times Y} \ar [rd] |{\ol\varphi_{X\times Y}} 
\ar @{-->} [rrd] ^{\mathcal{G}} \\
&\R& \qvL{\ol{F}_Y} \ar [l]\\
F_{X\times Y} \ar [rru] _{\mathcal{F}_X}
\ar [ru] |{\varphi}
\ar @{^(->} [uu] 
}
\end{equation*}
Now, note that $\mathcal{G}$
is a complete valuation
which extends $\mathcal{F}_X$.
\end{proof}

 }
\clearpage
{ \section{Epilogue}
\noindent
Starting from the similarity
between
 the Lebesgue measure
and the Lebesgue integral
as shown  on page~\pageref{S:intro}
I have tried to 
rebuild
a small part of the theory of
measure and integration 
in a more general setting.
When I look back at the result
I am most pleased that it was possible
to introduce the Lebesgue measure and the Lebesgue integral
with such natural and old primitives.
Indeed,
completeness and convexity
together
is nothing more than
the method of exhaustion
used by the ancient Greeks 
to determine the area of the disk (see title page).

The price for simple primitives
seems to be that much more effort is
required to prove even the simplest statements,
as attested by the size of this text.
Of course,
the number of pages could be
greatly reduced if we worked with~$\R$
instead of any~$E$,
but even then
I doubt that the approach taken
in this thesis would
be suitable
for a first course on the Lebesgue measure
and the Lebesgue integral.

Whether the theory in this thesis
will bear any fruit
I cannot tell,
but nevertheless I am content,
because I have enjoyed writing it,
and I hope that you have enjoyed reading it as well.
\label{S:conclusion}
 }

\appendix
\clearpage
{ \section{Ordered Abelian Groups}
\label{S:ag}
\noindent
In this thesis
we do not only consider $\R$-valued measures and integrals,
but also~$E$-valued ones, 
where $E$ is an \emph{ordered Abelian group}.
Since we do not expect reader to be familiar with this
particular generalisation of~$\R$,
we have collected the relevant definitions
and basic results in this appendix.
\begin{dfn}
\label{D:oag}
An \keyword{ordered Abelian group} $E$
is a set that is endowed with an Abelian group operation, $+$,
and a partial order, $\leq$, 
such that 
\begin{equation*}
x \ \leq\ y \quad\implies\quad w+x \ \leq\ w+y
\qquad\quad(w,x,y\in E).
\end{equation*}
\end{dfn}
%
%                  EXAMPLES OF AOGs
%
\begin{exs}
\label{E:oag}
\begin{enumerate}
\item
\label{E:aog-1}
The integers, $\Z$, the rationals, $\Q$, and the reals, $\R$,
endowed with  the usual addition and order
are ordered Abelian groups.

\item
\label{E:oag_div}
Let
$\Q^{\circ}$
be the set of rational numbers $q$ with $q>0$.
Order $\Q^\circ$ by
\begin{equation*}
q \ \preccurlyeq \ r
\quad\iff\quad 
\exists n\in \N \ [\ q \cdot n = r \ ].
\end{equation*}
Then $\Q^\circ$ with 
the usual multiplication is an  ordered Abelian group.

\item
\label{E:aog-2}
Let $E_1$ and $E_2$ be ordered Abelian groups.
Then $E_1\times E_2$ with pointwise order
and pointwise group operation is an ordered Abelian group.

\item
\label{E:oag_lex}
Consider~$\R^2$ with the pointwise addition.
By point~\ref{E:aog-2},
$\R^2$ with the usual order
is an  ordered Abelian group.
Further,
$\R^2$ with the
 \emph{lexicograpgic order},
\begin{equation*}
(x_1,x_2)\ \leq\ (y_1,y_2)
\quad\iff\quad
\left[ \ \ 
\begin{alignedat}{3}
&x_1 < y_1 \quad\text{or}\\
&x_1 = y_1  \quad\text{and}\quad  x_2 \leq y_2 
\end{alignedat}
\right.,
\end{equation*}
is also an ordered Abelian group,
called the \keyword{lexicographic plane},
 $\Lex$.
\end{enumerate}
\end{exs}

\noindent
Let us prove some simple statements
concerning ordered Abelian groups.
%
%                  + IS AN ORDER ISOMORPHISM
%
\begin{lem}
\label{L:oag-plus-iso}
Let $E$ be an ordered Abelian group.
Then, for $x,y,w\in E$,
\begin{equation*}
x \,\leq\, y \quad\iff\quad w+x \,\leq\, w+y.
\end{equation*}
\end{lem}
\begin{proof}
``$\Longrightarrow$''\ 
By the definition of ordered Abelian group.

\noindent ``$\Longleftarrow$''\ 
If $w+x\leq w+y$, then $x \,=\,-w\,+\, (w+x) \ \leq\ 
-w\,+\,(w+y) \,=\, y$.
\end{proof}
%
%                  PRESERVATION /\ BY +
%
\begin{lem}
\label{L:oag-plus-preserves}
Let $E$ be an ordered Abelian group.
Let $A\subseteq E$ and $x\in E$ be given.
\begin{enumerate}
\item \label{L:oag-plus-preserves-meet}
If $A$ has an infimum,
then so has $x+A \eqdf \{ x+a\colon a\in A\}$,
and 
\begin{equation*}
\bw \,\,x+A \ =\  x+\bw A.
\end{equation*}
\item \label{L:oag-plus-preserves-join}
If $A$ has a supremum,
then so has $x+A$,
and 
\begin{equation*}
\bv \,\,x+A \ =\  x+\bv A.
\end{equation*}
\end{enumerate}
\end{lem}
\begin{proof}
It suffices to prove that the map $E\ra E$
given by $u\mapsto x + u$ is an order isomorphism.
This follows easily using Lemma~\ref{L:oag-plus-iso}.
\end{proof}
%
%                  - IS AN ORDER ISOMORPHISM
%
\begin{lem}
\label{L:oag-minus-iso}
Let $E$ be an ordered Abelian group.
Then, for $x,y\in E$,
\begin{equation*}
x \,\leq\, y \quad\iff\quad -x \,\geq\, -y.
\end{equation*}
\end{lem}
\begin{proof}
``$\Longrightarrow$''\ 
If $x\leq y$, then
$-y \,=\, (-x-y)\,+\,x
\,\leq\, (-x-y)\,+\, y
\,=\, -x$.

\noindent ``$\Longleftarrow$''\ 
If $-y \leq -x$,
then $x\,=\,-(-x) \,\leq\, -(-y) \,=\, y$
by ``$\Longrightarrow$''.
\end{proof}
%
%                  PRESERVATION OF /\ BY -
%
\begin{lem}
\label{L:oag-minus-preserves}
Let $E$ be an ordered Abelian group,
and $A\subseteq E$.
\begin{enumerate}
\item
If $A$ has an infimum,
then $-A\eqdf \{-a\colon a\in A\}$
has a supremum,
and 
\begin{equation*}
-\bw A \,=\, \bv -\hspace{-3pt}A.
\end{equation*}
\item
If $A$ has an supremum,
then $-A$
has an infimum, and
\begin{equation*}
-\bv A \,=\, \bw -\hspace{-3pt}A.
\end{equation*}
\end{enumerate}
\end{lem}
\begin{proof}
The map $E\rightarrow E$ given by
$u\mapsto -u$ is an order reversing isomorphism.
\end{proof}

\noindent
We use the following
lemma
regularly.
%
%                  ADDITION OF SEQUENCES
%
\begin{lem}
\label{L:addition-of-suprema}
Let~$E$ be an ordered Abelian group.\\
Let $x_1 \leq x_2 \leq \dotsb$
be from~$E$ such that $\bv_n x_n$ exists.\\
Let $y_1 \leq y_2 \leq \dotsb$
be from~$E$ such that $\bv_n y_n$ exists.
Then
\begin{equation}
\label{eq:addition-of-suprema}
(\bv_n x_n) \,+\, (\bv_n y_n)
\ =\ 
\bv_k\, x_k + y_k.
\end{equation}
\end{lem}
\begin{proof}
By Lemma~\ref{L:oag-plus-preserves}
we know that 
\begin{equation*}
(\bv_n x_n) \,+\, (\bv_m y_m) \ =\  \bv_{n,m}\ x_n + y_m.
\end{equation*}
So to prove Equation~\eqref{eq:addition-of-suprema} holds,
it suffices to show that 
\begin{equation*}
\bv_{n,m}\ x_n + y_m \ =\ \bv_k \ x_k + y_k.
\end{equation*}
That is,
writing $z\eqdf \bv_{n,m}\ x_n + y_m$,
we must show that~$z$ is the supremum of 
\begin{equation*}
S \ \eqdf\ \{\ x_1+y_1,\ x_2 + y_2,\  \dotsc \ \}.
\end{equation*}
That is,
we must show that~$z$ is the smallest upper bound of~$S$.

Given $s\in S$, we have $s\equiv x_k + y_k$ for some  $k\in \N$,
and 
\begin{equation*}
x_k + y_k \ \leq\  \bv_{n,m}\ x_n + y_m\,\equiv\, z.
\end{equation*}
So we see that~$z$ is an upper bound of~$S$.

Let $u\in E$ be an upper bound of~$S$.
To prove that~$z$ is the smallest upper bound
of~$S$, we must show that $z\leq u$.
It suffices to prove that,
for all $n,m\in \N$,
\begin{equation}
\label{L:addition-of-suprema_final-bit}
x_n + y_m \ \leq \ u.
\end{equation}
Let $n,m\in\N$ be given,
and define  $k\eqdf \max\{n,m\}$.
Then we see that
\begin{equation*}
x_n + y_m \ \leq\  x_k + y_k \ \leq\  u.
\end{equation*}
Hence Statement~\eqref{L:addition-of-suprema_final-bit}
holds,
and we are done.
\end{proof}
\noindent
Of course,
we have a similar statement concerning infima.
\begin{lem}
\label{L:addition-of-infima}
Let~$E$ be an ordered Abelian group.\\
Let $x_1 \geq x_2 \geq \dotsb$
be from~$E$ such that $\bw_n x_n$ exists.\\
Let $y_1 \geq y_2 \geq \dotsb$
be from~$E$ such that $\bw_n y_n$ exists.
Then
\begin{equation*}
(\bw_n x_n) \,+\, (\bw_n y_n)
\ =\ 
\bw_k\, x_k + y_k.
\end{equation*}
\end{lem}
\begin{proof}
Similar to the proof of Lemma~\ref{L:addition-of-suprema}.
\end{proof}
%
%                  E^+ and E^-
%
\noindent
We will occasionally use the following notation.
\begin{dfn}
Let $E$ be an ordered Abelian group.
We write
\begin{align*}
E^+ \ &\eqdf\ \{\  a\in E\colon\  a\geq 0\ \},&
E^- \ &\eqdf\ \{\  a\in E\colon\  a\leq 0\ \}.
\end{align*}
\end{dfn}
\noindent
Let us now turn to a special class of ordered Abelian groups.
%
%                  Lattice ordered Abelian groups
%
\begin{dfn}
\label{D:loag}
A \keyword{lattice ordered Abelian group}
is an ordered Abelian group~$E$,
such that the order~$\leq$ 
makes~$E$ a lattice,
i.e., each pair $x,y\in E$
has an infimum, $x\wedge y$,
and a supremum, $x\vee y$.
\end{dfn}
\begin{exs}
\label{E:loag}
\begin{enumerate}
\item
The sets
$\Z$, $\Q$, and $\R$ are lattices
under the usual order.
The supremum of two elements is their maximum,
the infimum is the minimum.

\item
More generally,
any partially ordered set~$E$
that is \emph{totally ordered},
i.e.,
\begin{equation*}
\text{either } \quad x\leq y\quad\text{ or } 
\quad y\leq x\quad \text{ for all }x,y\in E,
\end{equation*}
is a lattice.
The supremum of $x,y\in E$ is simply the maximum of~$x$ and $y$,
the infimum $x$ and $y$ is the minimum of~$x$ and~$y$.

\item
The space $\Lex$
(see Example~\ref{E:oag}\ref{E:oag_lex})
is totally ordered and hence a lattice.

\item
\label{E:loag_div}
The set~$\Q^\circ$ ordered by~$\preccurlyeq$
(see Examples~\ref{E:oag}\ref{E:oag_div})
is a lattice.

Let $m,n\in \Q^\circ$ be given.
If $m,n\in \Z$, then the supremum of $m$ and~$n$
is the least common multiple of~$m$ and~$n$,
and the infimum of~$m$ and~$n$ is the greatest common divisor
of~$m$ and~$n$.

\end{enumerate}
\end{exs}

\noindent
The following result is quite suprising.
%
%                  MODULARITY EQUATION FOR LATTICE ORDERED ABELIAN GROUPS
%
\begin{lem}
\label{L:1-valuation}
Let $E$ be a lattice ordered Abelian group.
Then we have 
\begin{equation*}
a\wedge b  \,+\, a\vee b \,=\, a\,+\,b \qquad(a,b\in E).
\end{equation*}
\end{lem}
\begin{proof}
$a\vee b - a - b
= (a - a - b) \vee (b - a - b)
= (-b)\vee(-a) = -(a\wedge b)$.
\end{proof}
\begin{exs}
\begin{enumerate}
\item
Let $x,y\in \R$ be given.
Then Lemma~\ref{L:1-valuation}
gives us
\begin{equation*}
x+y\ =\ \min\{x,y\} \,+\, \max\{x,y\}.
\end{equation*}
Of course, this is trivial.

\item
Let $m,n\in \Z$ with $m,n\geq 0$ be given.
Then Lemma~\ref{L:1-valuation}
gives us
\begin{equation*}
m\cdot n \ =\ \gcd\{m,n\} \,\cdot\,\lcm\{m,n\}.
\end{equation*}
The above equality is more difficult to derive directly.
\end{enumerate}
\end{exs}
%
%                  SIGMA-DEDEKIND COMPLETE ORDERED ABELIAN GROUPS
%
\noindent
We now turn to `complete' ordered Abelian groups.
\begin{dfn}
\label{D:sdc}
Let~$E$ be an ordered Abelian group.\\
We say~$E$ is \keyword{$\sigma$-Dedekind complete}
if the following statement holds.
\begin{equation*}
\left[\quad
\begin{minipage}{.7\columnwidth}
Let $x_1,x_2,\dotsc$ be a sequence in~$E$.\\
Assume $x_1,x_2,\dotsc$ has an upper bound.\\
Then $\bv_n x_n$ exists.
\end{minipage}
\right.
\end{equation*}
\end{dfn}
%
%                  EXAMPLES CONCERNING SIGMA-DEDEKIND COMPLETE OAGs
%
\begin{exs}
\label{E:sdc}
\begin{enumerate}
\item 
The ordered Abelian group $\R$ is $\sigma$-Dedekind complete.
\item 
The ordered Abelian $\Q$ is \emph{not} $\sigma$-Dedekind complete.
\item
\label{E:sdc_lex}
The lexicographic plane~$\Lex$
(see Examples~\ref{E:oag}\ref{E:oag_lex})
is \emph{not} $\sigma$-Dedekind complete.\\
Indeed,
consider the following elements of~$\Lex$.
\begin{equation*}
(0,0)\ \leq\ (0,1)\ \leq\ (0,2)\ \leq\ \dotsb
\ \leq \ (1,0)
\end{equation*}
If $\Lex$ were $\sigma$-Dedekind complete,
then $S\eqdf \{ \, (0,n)\colon\, n\in\N\,\}$
would have a supremum;
we will prove that~$S$ does not have a supremum.

Suppose (towards a contradiction) that~$S$ has a supremum, $(x,y)$.\\
Then we have $(0,n)\leq (x,y)$ for all~$n\in \N$.
In other words, for all $n\in\N$,
\begin{equation*}
0<x \qquad\text{or}\qquad 
(\ 0=x\quad\text{and}\quad n\leq y\ ).
\end{equation*}
Hence $0<x$,
because there is no $y\in \R$ such that $n\leq y$
for all~$n\in\N$.\\
But then $(x,y-1)$ is an upper bound of~$S$ as well.\\
Since $(x,y)$ is the smallest upper bound of~$S$,
we have $(x,y)\leq(x,y-1)$. 
So $y\leq y-1$,
which is absurd.
Hence $S$ has no supremum.
\end{enumerate}
\end{exs}
%
%                  REMARK ON BOUNDED
%
\begin{rem}
The requirement
in Definition~\ref{D:sdc}
 that $x_1,x_2,\dotsc$
has an upper bound 
is essential
to make the notion
of $\sigma$-Dedekind completeness non-trivial.

Indeed,
if~$E$ is an ordered Abelian group
in which \emph{every} sequence $x_1,x_2,\dotsc$ has a supremum~$\bv_n x_n$,
then we have~$E=\{0\}$~!

Let $a\in E^+$ be given. We prove that $a=0$.
Note that the sequence
\begin{equation*}
1\cdot a  \ \leq\ 2\cdot a \ \leq\  3\cdot a \ \leq\ \dotsb
\end{equation*}
has a supremum, $\bv_n \ n\cdot a$.
Note that by Lemma~\ref{L:oag-plus-preserves}, we have
\begin{equation*}
(\bv_n\ n\cdot a)\,-\,a\ =\ 
\bv_n \ (n-1) \cdot a \ = \ 
\bv_n\ n\cdot a.
\end{equation*}
So we see that $b-a = b$,
where $b\eqdf \bv_n \ n\cdot a$.
Hence $a=0$.

Let $a\in E$ be given.
We must prove that~$a=0$.
Note that
by Lemma~\ref{L:1-valuation},
\begin{equation}
\label{eq:R:sdc}
a\ = \ 0\wedge a \,+\, 0\vee a.
\end{equation}
We have $0\wedge a =0$,
since $0\wedge a \in E^+$.
We also have $0\vee a=0$,
because $0\vee a\in E^-$,
so $-(0\vee a)\in E^+$, so $-(0\vee a)=0$,
and thus $0\vee a=0$.

So we see that $a=0$ by Equation~\eqref{eq:R:sdc}.
Hence $E=\{0\}$.
\end{rem}
%
%                  REMARK ON EQUIVALENT DEFINITION
%
\begin{rem}
\label{R:sdc}
Let~$E$ be an ordered Abelian group.\\
Using the order reversing isomorphism $x\mapsto -x$,
the reader can easily verify that
$E$ is $\sigma$-Dedekind complete
if and only if the following statement holds.
\begin{equation*}
\left[\quad
\begin{minipage}{.7\columnwidth}
Let $x_1,x_2,\dotsc$ be a sequence in~$E$.\\
Assume $x_1,x_2,\dotsc$ has a lower bound.\\
Then $\bw_n x_n$ exists.
\end{minipage}
\right.
\end{equation*}
\end{rem}
 }

\clearpage
\bibliographystyle{amsplain}
\bibliography{main}

\end{document}